\documentstyle[12pt]{amsart}

\newtheorem{theorem}{Theorem}[section] 
\newtheorem{corollary}[theorem]{Corollary} 
 
\newtheorem{lemma}[theorem]{Lemma} 
\newtheorem{proposition}[theorem]{Proposition}

\theoremstyle{definition}
\newtheorem{definition}{Definition}[section]
\newtheorem{example}{Example}[section]

\theoremstyle{remark}
\newtheorem{remark}{Remark}[section]

\numberwithin{equation}{section} 

\def\ord{{\mathop{\rm ord}}}

\def\supp{{\mathop{\rm supp}}}
\def\div{{\mathop{\rm div}}}

\def\det{{\mathop{\rm det}}}

\def\deg{{\mathop{\rm deg}}}

\def\Spec{{\mathop{\rm Spec}}}
\def\CPA{{\mathop{\rm CPA}}}

\def\Res{{\mathop{\rm Res}}}

\def\CH{{\mathop{\it CH}}}
\def\Max{{\mathop{\rm Max}}}
\def\Ker{{\mathop{\rm Ker}}}
\def\Ann{{\mathop{\rm Ann}}}
\def\radius{{\mathop{\rm radius}}}
\def\diam{{\mathop{\rm diam}}}
\def\Zh{{\mathop{\rm Zh}}}
\def\BDV{{\mathop{\rm BDV}}}
\def\SH{{\mathop{\rm SH}}}

\def\Berk{{\rm Berk}}

\def\<{\langle }
\def\>{\rangle }

\def\({(\!(}
\def\){)\!)}

\def\[{[\![}
\def\]{]\!]}

\def\v1{{\vec{1}}}

\def\BF1{{\mathbf 1}}

\def\AA{{\Bbb A}}

\def\CC{{\Bbb C}}
\def\FF{{\Bbb F}}

\def\QQ{{\Bbb Q}}

\def\NN{{\Bbb N}}
\def\RR{{\Bbb R}}
\def\ZZ{{\Bbb Z}}
\def\PP{{\Bbb P}}

\def\cA{{\cal A}}
\def\cB{{\cal B}}
\def\cC{{\cal C}}

\def\cF{{\cal F}}

\def\cI{{\cal I}}

\def\cK{{\cal K}}

\def\cM{{\cal M}}

\def\cO{{\cal O}}

\def\cV{{\cal V}}

\def\xbar{{\overline{x}}}
\def\ybar{{\overline{y}}}
\def\wbar{{\overline{w}}}
\def\zbar{{\overline{z}}}
\def\Ubar{{\overline{U}}}
\def\Vbar{{\overline{V}}}
\def\Wbar{{\overline{W}}}
\def\Zbar{{\overline{Z}}}

\def\inftybar{{\widetilde{\scriptstyle{\infty}}}}

\def\tK{{\tilde{K}}}

\def\tx{{\tilde{x}}}

\def\tM{{\tilde{M}}}

\def\tV{{\tilde{V}}}
\def\tW{{\tilde{W}}}
\def\tGamma{{\tilde{\Gamma}}}

\def\blacksquare{\Box}

\newenvironment{proof}{\noindent{\textbf{Proof: }}}{\hfill $\blacksquare$}

\DeclareMathSymbol{\varnothing} {\mathord}{AMSb}{"3F}

\begin{document}

\title[Berkovich Line]{Analysis and Dynamics on the Berkovich Projective Line}

\author[Robert Rumely]{Robert Rumely}
\author[Matthew Baker]{Matthew Baker}
\email{rr@@math.uga.edu \\
mbaker@@math.uga.edu}
       
\address{Department of Mathematics,
         University of Georgia, Athens, GA 30602-7403, USA}
\address{Department of Mathematics,
         University of Georgia, Athens, GA 30602-7403, USA}

\date{July 19, 2004}



\thanks{Work supported in part by NSF grant DMS-0300784.
We thank Robert Varley for several useful suggestions.  The idea for
using the Hsia kernel as the fundamental kernel for potential theory
on the Berkovich line was inspired by a manuscript of L. C. Hsia (\cite{Hsia}).
}

\begin{abstract}  \ \

This is a set of expanded lecture notes from the Berkovich Space seminar
held at the University of Georgia during Spring, 2004.  
The purpose of the notes is to provide a non-technical introduction to 
Berkovich spaces, and to develop the foundations for analysis on the 
Berkovich projective line, with a view toward applications in dynamics.  
After describing the underlying topological space and the
sheaf of functions on the Berkovich line, we introduce the Hsia kernel, 
the fundamental kernel for potential theory.  We develop a theory of 
capacities, define a Laplacian operator, and construct a theory of harmonic functions.
We then develop the theory of subharmonic functions and give applications
to dynamics, including a construction of the Lyubich measure attached
to a rational function.

These notes are still in a preliminary form; we plan in the future to revise and
expand them into a research monograph.  We are making them accessible
now because they provide proofs of some results needed in \cite{B-R2}. 
The paper \cite{B-R2} establishes an adelic equidistribution theorem for 
points of small dynamical height; the 
nonarchimedean part of this theorem uses in an essential way 
the theory developed in these notes.

We have been informed by A.~Chambert-Loir that his student A.~Thuillier
has independently proved a number of results concerning potential
theory on Berkovich curves, including the construction of a Laplacian
operator and a theory of harmonic functions.  
Chambert-Loir and Thuillier, and independently
C.~Favre and J.~Rivera-Letelier, have also recently given
constructions of a measure on $\PP^1_{\Berk}$
attached to a rational function which presumably coincides with our 
Lyubich measure.

\end{abstract}

\maketitle

\newpage

\tableofcontents

\section{Topological description of the Berkovich Unit Disc.}
\label{Section A}

Let $K$ be a global field and let $v$ be a nonarchimedean place of $K$.
Let $\CC_v$ be the completion of the algebraic closure of $K$.  
In this section we recall Berkovich's theorem that the unit Berkovich Disc 
over $\CC_v$  can be identified with the collection 
of all equivalence classes of sequences of nested discs  
$\{B(a_i,r_i)\}_{i=1,2,\ldots}$ contained in $B(0,1)$.  
Here $B(a,r) = \{z \in \CC_v : |z-a|_v \le r \}$.  
This leads to an explicit description of
the Berkovich Disc in terms of an infinitely-branched metrized tree.   

\vskip .1 in
Let  $\cA = \CC_v\{\{T\}\}$ be the ring of all formal power series 
with coefficients in $\CC_v$, converging on the unit  $B(0,1)$.  
That is, $\cA$ is the ring of 
of all power series $f(T) = \sum_{i=0}^\infty a_i T^i \in \CC_v[[T]]$ such that
$\lim_{i \rightarrow \infty} |a_i|_v = 0$.  Equip $\cA$ with the Gauss norm
\begin{equation*}
\|f\|_v \ = \ \max_i (|a_i|_v)  \ .
\end{equation*}
With this norm, it is a Banach algebra over  $\CC_v$. 

A multiplicative seminorm on $\cA$ is a function 
$[ \cdot ]_x : \cA \rightarrow \RR_{\ge 0}$ such that $[0]_x = 0$, 
$[1]_x = 1$, $[f \cdot g]_x = [f]_x \cdot [g]_x$ 
and $[f-g]_x \le [f]_x + [g]_x$ for all $f, g \in \cA$.  It is called 
{\it bounded} if there is a constant $C_x$ such that $[f]_x \le C_x \|f\|$ 
for all $f \in \cA$.  Boundedness is equivalent to continuity.  

       It can be deduced from these properties that a bounded multiplicative
seminorm $[ \cdot ]_x$ on $\cA$ behaves like a nonarchimedean absolute value,
except that its kernel may be nontrivial.           
The multiplicativity implies that $C_x$ can be
taken to be $1$, that is, $[f]_x \le \|f\|$ for all $f$.  Indeed  
for each $n$, $([f]_x)^n = [f^n]_x \le C_x \|f^n\| = C_x \|f\|^n$,
so $[f]_x \le C_x^{1/n} \|f\|$.  For a constant $c \in \CC_v$, necessarily
$[c]_x = |c|_v$.  For this, note that by the definition of the Gauss norm,
$\|c\| = |c|_v$.  If $c = 0$ then trivially $[c]_x = 0$;  otherwise,
$[c]_x \le \|c\| = |c|_v$ and $[c^{-1}|_v \le \|c^{-1}\| = |c^{-1}|_v$,
while multiplicativity gives $[c]_x \cdot [c^{-1}]_x = [c \cdot c^{-1}]_x = 1$.
Combining these gives  $[c]_x = |c|_v$.   The usual proof
of the ultrametric inequality now carries over to show that 
$[f+g]_x \le \max([f]_x,[g]_x)$ for all $f, g \in \cA$, 
with equality if $[f]_x \ne [g]_x$.  Indeed, the binomial theorem shows that
for each $n$
\begin{eqnarray*} 
([f+g]_x)^n & = & [(f+g)^n]_x 
            \ = \ [ \sum_{k=0}^n {{n} \choose {k}} f^k g^{n-k}]_x \\
            & \le & \ \sum_{k=0}^n 
                    |{{n} \choose {k}}|_v \cdot [f]_x^k [g]_x^{n-k} 
            \ \le \ \sum_{k=0}^n [f]_x^k [g]_x^{n-k} \\
            & \le & (n+1) \cdot \max([f]_x,[g]_x)^n \ .
\end{eqnarray*}
Taking $n^{th}$ roots and passing to a limit gives the inequality.  
If in addition $[f]_x < [g]_x$ then
$[g]_x \le \max([f+g]_x, [-f]_x)$ implies that 
$[f+g]_x = \max([f]_x,[g]_x)$.  

\vskip .1 in
By definition, the Berkovich Disc $\cB(0,1)$
is the functional-analytic {\it spectrum} of $\cA$,  the 
set of all bounded multiplicative seminorms $[ \cdot ]_x$ 
on $\cA$.  By abuse of notation, we will sometimes write $x$ for 
$[ \cdot ]_x  \in \cB(0,1)$.  The Berkovich Disc is equipped with the 
Gel'fond topology, the weakest topology such that for all $f \in A$ and all
$\alpha \in \RR$, the sets
\begin{eqnarray*}
U(f,\alpha) & = & \{ x \in \cB(0,1) : [f]_x < \alpha \} \ , \\
V(f,\alpha) & = & \{ x \in \cB(0,1) : [f]_x > \alpha \}    
\end{eqnarray*}
are open.  This makes it a nonempty, compact Haudorff space 
(\cite{Berk}, Theorem 1.2.1, p.13).  It is also connected, even path-connected 
(\cite{Berk}, Corollary 3.2.3, p.52).
 
\vskip .1 in
Before proceeding further, we should note some 
elements of $\cB(0,1)$.  For each $a \in B(0,1)$ we have the
evaluation seminorm,
\begin{equation*}
[f]_a \ = \ |f(a)|_v \ .
\end{equation*}  
The boundedness of $[f]_a$ follows from the maximum modulus principle
in nonarchimedean analysis, 
which says that $\sup_{z \in B(0,1)} |f(z)|_v = \|f\|$.
Also, for each subdisc $B(a,r) \subset B(0,1)$, we have the sup norm
\begin{equation*}
[f]_{B(a,r)} \ = \ \sup_{z \in B(a,r)} (|f(z)|_v) \ .
\end{equation*}
Note that by the maximum modulus principle,
the sup norm over the unit disc coincides with the Gauss norm:
\begin{equation*}
[f]_{B(0,1)} = \|f\| \ .
\end{equation*}
For any decreasing sequence of discs 
$x = \{B(a_i,r_i)\}_{i \ge 1}$
we have the limit seminorm
\begin{equation*}
[f]_x \ = \ \lim_{i \rightarrow \infty} [f]_{B(a_i,r_i)} \ .
\end{equation*}
Berkovich showed that in fact each $x \in \cB(0,1)$ is of this type. 

\vskip .1 in
Fix $x \in \cB(0,1)$.  By the Weierstrass Preparation Theorem, 
each $f \in \cA$ can be uniquely written as 
\begin{equation*}
f \ = \ c \cdot \prod_{j=1}^m (T-a_j) \cdot u(T)
\end{equation*}
where $c \in \CC_v$, $a_j \in B(0,1)$ for each $j$, and $u(T)$ is a unit power
series, that is $u(T) = 1 + \sum_{i=1}^{\infty} a_i T^i \in \cA$ 
where $|a_i|_v < 1$ for
all $i$ and $\lim_{i \rightarrow \infty} |a_i|_v = 0$.  It is easy to see 
that $[u]_x = 1$.  
Indeed, $u(T)$ has a multiplicative inverse $u^{-1}(T)$
of the same form, and by the definition of the Gauss norm 
$\|u\| = \|u^{-1}\| = 1$. 
Since $[u]_x \le \|u\| = 1$, $[u^{-1}]_x \le \|u^{-1}\| \le 1$, and
$[u]_x \cdot [u^{-1}]_x = [u \cdot u^{-1}]_x = 1$, necessarily
$[u]_x = 1$.  
 
It follows that
\begin{equation*}
[f]_x \ = \ |c|_v \cdot \prod_{j=1}^m [T-a_j]_x \ .
\end{equation*}
Thus, each $x \in \cB(0,1)$ is determined by its 
values on the linear polynomials  $T-a$, $a \in B(0,1)$.  

\begin{proposition} \label{BerkProp2} {\rm (Berkovich \cite{Berk}, p.18)} 
Each $x \in \cB(0,1)$ can be realized in the form 
\begin{equation} \label{FF0}
[f]_x \ = \ \lim_{i \rightarrow \infty} [f]_{B(a_i,r_i)}
\end{equation}
for some sequence of nested discs 
$B(a_1,r_1) \supseteq B(a_2,r_2) \supseteq \cdots$.  If this sequence  
has a nonempty intersection, then either

$1)$ the intersection is a single point $a$, in which case
 $[f]_x = |f(a)|_v$, or 
 
$2)$ it is a  disc $B(a,r)$ $($where $r$ may or may not belong to the value 
group of $\CC_v$$)$, in which case $[f]_x = [f]_{B(a,r)}$.    
\end{proposition}

\begin{proof}  Fix $x \in \cB(0,1)$, and  
consider the family of discs 
\begin{equation*}
\cF \ = \ \{ B(a,[T-a]_x) : a \in B(0,1) \} \ .
\end{equation*}
If $a, b \in B(0,1)$ and $[T-a]_x \ge [T-b]_x$, then 
\begin{eqnarray}
|a-b|_v & = & [a-b]_x \ = \ [(T-b)-(T-a)]_x \notag \\
            & \le & \max([T-a]_x,[T-b]_x) \ = \ [T-a]_x \ , \label{FF1}
\end{eqnarray}
with equality if $[T-a]_x > [T-b]_x$.  In particular
$b \in B(a,[T-a]_x)$, and  
\begin{equation*}
B(a,[T-a]_x) \ \supseteq \ B(b,[T-b]_x) \ .
\end{equation*}
This means that the family $\cF$ is totally ordered by containment.  
Put $r = \inf_{a \in B(0,1)} ([T-a]_x)$  and choose a sequence of points
$a_i \in B(0,1)$ such that the numbers $r_i = [T-a_i]_x$ satisfy 
$\lim_{i \rightarrow \infty} r_i = r$.  

We claim that for each polynomial $T-a$ with $a \in B(0,1)$, 
we have  
\begin{equation} \label{FF2}
[T-a]_x \ = \ \lim_{i \rightarrow \infty} [T-a]_{B(a_i,r_i)} \ .
\end{equation}
If $[T-a]_x = r$ then for each $a_i$ we have $r_i = [T-a_i]_x \ge |a_i -a|_v$
by (\ref{FF1}), so $a \in B(a_i,r_i)$.  Hence   
\begin{equation*}
[T-a]_{B(a_i,r_i)} \ = \ \sup_{z \in B(a_i,r_i)} |z-a|_v  
                   \ = \ r_i  \ .
\end{equation*}   
Since $\lim_{i \rightarrow \infty} r_i = r$,   
(\ref{FF2}) holds in this case.  

If $[T-a]_x > r$ then for each
$a_i$ with $[T-a]_x > [T-a_i]_x$, we have $[T-a]_x = |a-a_i|_v$ by 
the strict case in (\ref{FF1}), 
which means that $|a-a_i|_v > [T-a_i]_x = r_i$.  Hence
\begin{equation*}
[T-a]_{B(a_i,r_i)} \ = \ \sup_{z \in B(a_i,r_i)} |z-a|_v  
                   \ = \ |a-a_i|_v \ = \ [T-a]_x \ .  
\end{equation*}
Thus the limit on the right side of (\ref{FF2}) stabilizes at 
$[T-a]_x$, and (\ref{FF2}) holds in this case as well.
As noted previously, 
$[ \cdot ]_x$ is determined by its values on the polynomials $T-a$, so 
for all $f \in \cA$, 
\begin{equation} \label{FF3}
[f]_x \ = \ \lim_{i \rightarrow \infty} [f]_{B(a_i,r_i)} \ .
\end{equation}

Now suppose the intersection of the family $\cF$ is non-empty, and contains
the point $a$.  Formula (\ref{FF2}) gives 
\begin{equation*}
[T-a]_x \ = \ \lim_{i \rightarrow \infty} [T-a]_{B(a_i,r_i)} 
\ \le \ \lim_{i \rightarrow \infty} r_i  \ = \ r \ ,
\end{equation*}
while the definition of $r$ shows that $[T-a]_x \ge r$.  Thus
$[T-a]_x = r$. Hence the  $B(a,r)$ (which may consist of a single point,
if $r = 0$) is a minimal element of $\cF$. 
The arguments above show that formula (\ref{FF3}) holds  
for {\it any} sequence of discs $B(a_i,r_i)$ such that  $r_i = [T-a_i]_x$
satisfies $\lim r_i = r$.  If we take $a_i = a$ for each $i$, 
then $r_i = [T-a]_x = r$, and (\ref{FF3}) gives $[f]_x = [f]_{B(a,r)}$. 
If $r = 0$, it gives $[f]_x = |f(a)|_v$.  
\end{proof}

\vskip .1 in

It is important to note that there {\it are}
sequences of nested discs $\{B(a_i,r_i)\}$ with empty intersection.  Such
sequences necessarily satisfy $r = \lim r_i > 0$, since if $r = 0$ the 
completeness of $\CC_v$ shows the intersection is a point $a \in \CC_v$.  

To construct one, fix $0 < r < 1$, and choose a sequence $\{r_i\}$ which decreases 
monotonically to $r$, with $r < r_i \le 1$ for each $i$.  The algebraic
closure $\tK$ is countable, and is dense in $\CC_v$.  Enumerate the elements
of $\tK \cap B(0,1)$ as $\{\alpha_j\}_{j \ge 1}$.  Define a sequence of
s $B(a_i,r_i)$ as follows.  Take $B(a_1,r_1) = B(\alpha_1,r_1)$.  Suppose
$B(a_i,r_i)$ has been constructed.  Let $j_i$ be the least index of an element
with $\alpha_j \in B(a_i,r_i)$.  Since $r_i > r_{i+1}$, 
$B(a_i,r_i) \backslash B(\alpha_{j_i},r_{i+1})$ is nonempty;  let $a_{i+1}$
be any element of it.   Then $\alpha_{j_i} \notin B(a_{i+1},r_{i+1})$.  

Clearly the sequence $j_i$ increases 
to $\infty$, indeed one sees inductively that $j_i \ge i$.  For each $i$,
the construction has arranged that the $\alpha_j$ with  $j \le j_i$ do 
not belong to $B(a_{i+1},r_{i+1})$.  
It follows that $\cap_{i=1}^{\infty} B(a_i,r_i)$ 
contains no elements of $\tK$.  If it were nonempty, it would be a 
disc $B(a,r)$.  However, every such disc contains elements of $\tK$.  Hence 
$\cap_{i=1}^{\infty} B(a_i,r_i)$  must be empty.   

\vskip .1 in
This brings us to Berkovich's classification of elements of $\cB(0,1)$:  

Points corresponding to sequences $\{B(a_i,r_i)\}$ with $\lim r_i = 0$ 
are said to be of Type I.  As noted above, the completeness of $\CC_v$ assures
that the intersection of such a sequence is a point $a \in \CC_v$, and the 
corresponding seminorm is $[ \cdot ]_a$.  We will call these points the
``classical points''

Points corresponding to sequences $\{B(a_i,r_i)\}$ with nonempty intersection,
for which $r = \lim r_i > 0$ belongs to the value group of $\CC_v$, are said
to be of type II.  These correspond to a sup norm  $[ \cdot ]_{B(a,r)}$. 
Rivera-Letelier calls these ``rational points''. 

Points corresponding to sequences $\{B(a_i,r_i)\}$ with nonempty intersection,
but for which $r = \lim r_i > 0$ does not belong to the value group of $\CC_v$, 
are said to be of type III.  These also correspond to a sup norm  
$[ \cdot ]_{B(a,r)}$;  however the  $B(a,r)$ is not a 
``rational domain'' in the sense of classical rigid analysis.  
Rivera-Letelier calls these ``irrational points''.

Points corresponding to sequences $\{B(a_i,r_i)\}$ with empty intersection 
are said to be of type IV.  As noted before, necessarily $\lim r_i > 0$.  
These are the truly ``new'' points in the Berkovich space, ones that aren't 
seen classically.  Rivera-Letelier calls these ``singular points''.

\vskip .05 in
There is a distinguished point $\zeta_0$ in $\cB(0,1)$,  
namely, the point corresponding to the Gauss norm $\|f\|_v = [ \ ]_{B(0,1)}$.  
Chambert-Loir calls this the ``Gauss point'', 
and we will adopt his terminology. 

\vskip .1 in

We will call two sequences of nested discs {\it equivalent} if they define the 
same point in $\cB(0,1)$.  Since the limit (\ref{FF0}) is a decreasing one, it 
is clear that two sequences $\{B(a_i,r_i)\}$ and 
$\{B(a_i^{\prime},r_i^{\prime})\}$ are equivalent if
 
\vskip .05 in 
 a)  each has a nonempty intersection, and their intersections are the same;  \ or
            
 b)  both have empty intersection, and each sequence is cofinal in the other.
         
\vskip .05 in
\noindent{(Cofinal} means that for each $i$, there is a $j$ such that 
$B(a_i,r_i) \supseteq B(a_j^{\prime},r_j^{\prime})$, 
and for each $j$ there is an $i$ such that 
$B(a_j^{\prime},r_j^{\prime}) \supseteq B(a_i,r_i)$.) 
          
These conditions are necessary for equivalence, as well as sufficient.  
From Proposition \ref{BerkProp2} it is clear that two sequences with nonempty 
intersection are equivalent if and only if they have the same intersection.  
A sequence $x$ with nonempty intersection $B(a,r)$ (possibly $r = 0$) 
cannot be equivalent to any sequence $y = \{B(a_i,r_i)\}$ 
with empty intersection, since for any $i$ with 
$a \notin B(a_i,r_i)$ 
\begin{equation*}
[T-a_i]_x = [T-a_i]_{B(a,r)} = |a-a_i|_v > r_i,  \quad
[T-a_i]_y \le [T-a_i]_{B(a_i,r_i)} = r_i 
\end{equation*}
Finally, two sequences $x = \{B(a_i,r_i)\}$ and 
$y = \{B(a_i^{\prime},r_i^{\prime})\}$ with empty intersection 
which are not cofinal cannot be equivalent.  Since $x$ is not cofinal in 
$y$, after removing some initial terms of $x$ we can assume that   
$B(a_1,r_1) \not\supseteq B(a_i^{\prime},r_i^{\prime})$
for any $i$.   Since $\cap_{i=1}^{\infty} B(a_i^{\prime},r_i^{\prime}) = \phi$,
after deleting some initial terms of $y$ we can assume that
$a_1 \notin B(a_1^{\prime},r_1^{\prime})$.  As any two discs are either
disjoint, or one contains the other, it must be that $B(a_1,r_1)$ and 
$B(a_1^{\prime},r_1^{\prime})$ are disjoint.  
Since $B(a_i^{\prime},r_i^{\prime}) \subseteq B(a_1^{\prime},r_1^{\prime})$
for all $i$, we have 
$[T-a_1]_{B(a_i^{\prime},r_i^{\prime})}  = |a_1-a_1^{\prime}|_v > r_1$ 
for each $i$.  Hence 
\begin{eqnarray*}
& & \qquad \qquad [T-a_1]_x \ \le \ [T-a_1]_{B(a_1,r_1)} \ = \ r_1 \ , \\
& & [T-a_1]_{y} \ = \ \lim_{i \rightarrow \infty} 
[T-a_1]_{B(a_i^{\prime},r_i^{\prime})} \ = \ |a_1 -a_1^{\prime}|_v \ > \ r_1 \ .
\end{eqnarray*}

\vskip .1 in
We will use this description of seminorms on $\cA$ to construct 
a topological model of $\cB(0,1)$.  
Let $x = B(a,r) \subset B(0,1)$ be any .  (Here $r$ may or may not belong 
to the value group of $\CC_v$.  We also allow the possibility that $r = 0$, so
that $B(a,r)$  degenerates to the single point $a$.)  
By the {\it line of discs}
$[r,1]_x$ we mean the set of discs $\{B(a,t) : r \le t \le 1\}$.  
We view this set as having the structure of a line segment, with
the discs $B(a,t)$ as points of the segment.  

Now let $S = \{B(a_1,r_1), \ldots, B(a_n,r_n)\}$ be any finite set of discs
contained in $B(0,1)$.  Define the {\it graph of discs} $\Lambda_S$
to be the union of the associated lines of discs $[r_i,1]_{B(a_i,r_i)}$, 
\begin{equation*}
\Lambda_S \ = \ \bigcup_{i=1}^n [r_i,1]_{B(a_i,r_i)} \ .
\end{equation*}
In forming the union, we identify points of segments whose associated
discs coincide.  The graph $\Lambda_S$ is a tree rooted at the point $B(0,1)$,
and has a natural metric on its edges gotten from the distance function on
its component segments, making it a {\it metrized graph}.   

An alternate description of $\Lambda_S$, which may make its structure clearer,
is as follows.  Define the {\it saturation } of $S$
to be the set $\hat{S}$ gotten by adjoining to $S$   
all discs $B(a_i,|a_i-a_j|_v)$ with $B(a_i,r_i), B(a_j,r_j) \in S$, and 
also the  $B(0,1)$.  
Note that the radius of each  in $\hat{S} \backslash S$ belongs to 
the value group of $\CC_v$.  
Then $\Lambda_S$ is the finite metrized graph 
whose nodes are the discs in $\hat{S}$, 
and which has an edge of length $|r_i-r_j|$
between each pair of nodes $B(a_i,r_i), B(a_j,r_j) \in \hat{S}$ for which
$B(a_i,r_i) \supseteq B(a_j,r_j)$ or $B(a_j,r_j) \supseteq B(a_i,r_i)$.  
Clearly only discs $B(a_i,r_i)$ with radii belonging to the value group of $\CC_v$ 
branch points of the graph.  

If $S_1$ and $S_2$ are any two finite sets of discs, then $\Lambda_{S_1}$
and $\Lambda_{S_2}$ are metrized subgraphs of $\Lambda_{S_1 \cup S_2}$.  Thus,
the collection of graphs $\Lambda_S$ is a directed set.  Let
\begin{equation*}
\Lambda \ = \ \bigcup_S \Lambda_S
\end{equation*}
be their union.  Each  each  
$B(a,r) \subset B(0,1)$ corresponds to a unique point in $\Lambda$.  
That is, each seminorm $[ \cdot ]_x$ of type I, II, or III 
corresponds in a natural way to a unique point of $\Lambda$.  

To incorporate the points of type IV we must enlarge $\Lambda$ by 
adding ``ends''.  Consider a decreasing sequence of nested discs 
$x = \{B(a_i,r_i)\}$, and put $r = \lim r_i$.  
The union of the lines of discs $[r_i,1]_{B(a_i,r_i)}$ is a ``half-open'' 
line of discs, which we will write as $(r,1]_x$.  If $\cap_i B(a_i,r_i)$ 
is nonempty, then the intersection is a  $B(a,r)$ and $(r,1]_x$ 
extends in a natural way to the closed line $[r,1]_x = [r,1]_{B(a,r)}$.   
However, if $\cap_i B(a_i,r_i) = \phi$, 
we must adjoin a new ``end'' in order to close up $(r,1]_x$.  
By abuse of notation
we denote this point $x$, and write $[r,1]_x = (r,1]_x \cup \{x\}$.  
It is clear that cofinal sequences define the same the half-line $(r,1]_x$, 
so the point closing up the half-line depends only on the seminorm $[ \cdot ]_x$ 
and not on the sequence defining it.    

Our model of $\cB(0,1)$ is the space $\overline{\Lambda}$ 
 gotten by adjoining to $\Lambda$ 
all the ends $x$ corresponding to points of type IV.  Like $\Lambda$, 
it is a metrized tree rooted at $B(0,1)$.  
It has countably many branches emanating
from the root, each branch corresponding to discs contained in an open
 $B(a,1)^- = \{z \in \CC_v: |z-a|_v < 1\}$.  Each branch splits into 
into countably many branches at each point $B(a,r)$ of type II 
(for which $r$ belongs to the value group of $\CC_v$), 
and each new branch behaves in the same way.  
This incredible collection of splitting branches forms a sort of 
``witch's broom''.  However, the witch's broom has some structure:  it 
splits {\it only} at the points $B(a,r)$ of type II, not those of type III; 
and there are {\it only} countably many branches at each point of type II, 
corresponding to the open discs $B(p,r)^-$ with $p \in B(a,r)$.  
Some of the branches extend all the way to the bottom 
(terminating in points of type I), while others are ``cauterized off'' earlier
and terminate at points of type IV, but every branch terminates either at 
a point of type I or type IV.

\vskip .1 in
Now consider the topology on $\cB(0,1)$.  By definition,  
the Gel'fond topology is generated by the open sets 
$U(f,\alpha)  =  \{ x \in \cB(0,1) : [f]_x < \alpha \}$ and  
$V(f,\alpha)  =  \{ x \in \cB(0,1) : [f]_x > \alpha \}$,     
for $f \in \cA$ and $\alpha \in \RR$.  
Since the value group of $\CC_v$ is dense in $\RR_{> 0}$,
it suffices to consider $\alpha$ belonging to the value group of $\CC_v$.  
Also, by the Weierstrass preparation theorem and the fact that any
unit power series $u(T)$ satisfies $[u]_x \equiv 1$, we can restrict 
to polynomials $f(T) \in \CC_v[T]$ with roots in $B(0,1)$.  

Given a nonconstant polynomial 
$f(T) = c \prod_{i=1}^n (T-a_i)^{m_i} \in \CC_v[T]$, and $\alpha > 0$ 
belonging to the value group of $\CC_v$, there is a well-known 
description of the set $\{z \in \CC_v : |f(z)|_v \le \alpha\}$   
as a finite union of closed discs 
$\bigcup_{i=1}^N B(a_i,r_i)$ 
(c.f. \cite{Cantor}, Theorem 3.1.2, p.180).  
Here the centers can be taken to be roots of $f(T)$, 
and each $r_i$ belongs to the value group of $\CC_v$.  
If desired, one can assume that the discs
in the decomposition are pairwise disjoint.  However, for us it will be more
useful to assume that all the roots occur as centers, so $N = n$,  
and that if $a_i$ and $a_j$ are roots with $a_j \in B(a_i,r_i)$, 
then $B(a_i,r_i) = B(a_j,r_j)$.  

Taking the union over an increasing 
sequence of $\alpha$, we can lift the requirement that 
$\alpha$ belongs to the value group of $\CC_v$ (note that if $\alpha$
is not in the value group, then $\{z \in \CC_v : |f(z)|_v = \alpha\}$ 
is empty.)  Thus, any $\alpha > 0$ determines a collection of numbers $r_i > 0$,
which belong to the value group of $\CC_v$ if $\alpha$ does, such that 
\begin{equation*}
\{z \in \CC_v : |f(z)|_v \le \alpha\} \ = \ \bigcup_{i=1}^n B(a_i,r_i) \ .   
\end{equation*}  
Here we assume as before that $\alpha_1, \ldots, \alpha_n$ are the roots of 
$f(z)$, and if $\alpha_j \in B(a_i,r_i)$ then $r_j = r_i$.  
Using this and the factorization of $f(T)$, it is easy to 
see that  
\begin{equation*}
\{z \in \CC_v : |f(z)|_v < \alpha\} \ = \ \bigcup_{i=1}^n B(a_i,r_i)^- \ ,    
\end{equation*}       
where $B(a_i,r_i)^- = \{z \in \CC_v : |z-a_i|_v < r_i \}$.

For any closed disc $B(b,t) \subset B(a_i,r_i)^-$ with $t$ in the value 
group of $\CC_v$, one sees readily that 
\begin{equation*}
\sup_{z \in B(b,t)} |f(z)|_v < \alpha \ .  
\end{equation*}
Likewise, on any  
$B(b,t) \subset B(a_i,r_i) \backslash (\cup_{a_j \in B(a_i,r_i)} B(a_j,r_i)^-)$ 
one has $|f(z)|_v \equiv \alpha$, and on any  $B(b,t)$ disjoint 
from $\cup_{i=1}^n B(a_i,r_i)$ one has $|f(z)|_v \equiv \beta$
for some $\beta > \alpha$.  

Suppose $x \in \cB(0,1)$ corresponds to a sequence of nested discs
$\{B(b_j,t_j)\}$. Without loss, we can assume that each $t_j$ 
belongs to the value group of $\CC_v$.
We will say that $x$ is associated to an open disc
$B(a,r)^-$ if there is some $j$ such that  $B(b_j,t_j) \subset B(a,r)^-$. 
We say that $x$ is associated to 
a closed disc $B(a,r)$ if 
there is some $j$ such that 
$B(b_j,t_j) \subset B(a,r)$, or if $\cap_{j=1}^{\infty} B(b_j,t_j) = B(a,r)$.
From the assertions in the previous paragraph, 
it follows that $[f]_x < \alpha$
if and only if $x$ is associated to some $B(a_i,r_i)^-$.  
Likewise, $[f]_x > \alpha$ if and only $x$ is not associated to any of the
$B(a_i,r_i)$.  

It is not hard to see that $x$ is associated to an open disc 
$B(a,r)^-$ if and only if
$[T-a]_x < r$.  Indeed, $[T-a]_{B(b_j,t_j)} = \max(t_j,|b_j-a|_v)$.
Thus, if $[T-a]_x < r$, 
then there is some $j$ for which $\max(t_j,|b_j-a|_v) < r$ and this implies
$B(b_j,t_j) \subset B(a,r)^-$.  Conversely, if $B(b_j,t_j) \subset B(a,r)^-$
then $t_j < r$ and $|b_j-a|_v < r$ so $[T-a]_x \le [T-a]_{B(b_j,t_j)} < r$. 

Similarly, $x$ is associated to a closed disc
$B(a,r)$ if and only if $[T-a]_x \le r$.  
First suppose $x$ is associated to $B(a,r)$.
If some $B(b_j,t_j) \subset B(a,r)$ then clearly 
$[T-a]_x \le [T-a]_{B(b_j,t_j)} \le r$.  By Proposition \ref{BerkProp2}
if $\cap_{j=1}^{\infty} B(b_j,t_j) = B(a,r)$ then 
$[T-a]_x = [T-a]_{B(a,r)} = r$.  Conversely, suppose $[T-a]_x \le r$.
If $B(b_j,t_j) \subset B(a,r)$ for some $j$ then $x$ is certainly
associated to $B(a,r)$.  Otherwise, for each $j$ either $B(b_j,t_j)$ is
disjoint from $B(a,r)$ or $B(b_j,t_j) \supset B(a,r)$.  If some 
$B(b_{j_0},t_{j_0})$ is disjoint from $B(a,r)$ then by the ultrametric inequality 
$|z-a|_v = |b_{j_0}-a|_v > r$ for all $z \in B(b_{j_0},t_{j_0})$.  Since the
discs $B(b_j,t_j)$ are nested, this means that 
$[T-a]_{B(b_j,t_j)} = |b_{j_0}-a|_v$ for all $j \ge j_0$, and hence
that $[T-a]_x = |b_{j_0}-a|_v$.  This contradicts $[T-a]_x \le r$, 
so it must be that each $B(b_j,t_j)$ contains $B(a,r)$.  Thus,
$\cap_{j=1}^{\infty} B(b_j,t_j) = B(a,t)$ for some $t \ge r$.  
By Proposition \ref{BerkProp2}, $[T-a]_x = t$.  Since we have assumed
$[T-a]_x \le r$ this gives $t = r$, so $x$ is associated to $B(a,r)$.

This leads us to define 
open and closed ``Berkovich discs'', as follows.  
For $a \in B(0,1)$ and $r > 0$, write
\begin{eqnarray*}
\cB(a,r)^- & = & \{x \in \cB(0,1) : [T-a]_x < r \} \ , \\
\cB(a,r)  & = & \{x \in \cB(0,1) : [T-a]_x \le r \} \ . 
\end{eqnarray*}
With this notation, our discussion above shows that  
\begin{equation*}
U(f,\alpha)  =  \bigcup_{i=1}^N \cB(a_i,r_i)^- , \quad
V(f,\alpha)  =  \cB(0,1) \backslash \bigcup_{i=1}^N \cB(a_i,r_i) \ .
\end{equation*}

In terms of our model $\overline{\Lambda}$, a closed Berkovich disc
$\cB(a,r)$ consists of all points in a branch on or below 
$B(a,r)$.  If $r$ belongs to the value group of $\CC_v$, 
then the open Berkovich disc $\cB(a,r)^-$ is one of the 
countably many open branches emanating from $B(a,r)$ 
(more precisely, the one containing $a$), 
while if $r$ does {\it not} belong to the value group of $\CC_v$, then
$\cB(a,r)^- = \cB(a,r) \backslash \{B(a,r)\}$ consists of all points 
in the open branch below $B(a,r)$.  

Taking finite intersections of sets of the form $U(f,\alpha)$ and $V(f,\alpha)$
gives a basis for the open sets in the Gel'fond topology.  Thus, 

\begin{proposition} \label{PP3}
A basis for the open sets of $\cB(0,1)$ is given by the sets 
\begin{equation*}
\cB(a,r)^-, \ \   
\cB(a,r)^- \backslash \bigcup_{i=1}^N \cB(a_i,r_i), \ \ \text{and} \ \
\cB(0,1) \backslash \bigcup_{i=1}^N \cB(a_i,r_i), 
\end{equation*}
where $a$ and the $a_i$ range over  $B(0,1)$,
and where each $r$, $r_i > 0$.     
\end{proposition}  
\vskip .1 in
  
Clearly this basis has a countable sub-basis, gotten by restricting
to discs $\cB(a,r)^-$ and $\cB(a_i,r_i)$ whose centers belong to 
$\tK \cap B(0,1)$ and whose radii belong to the value group of $\CC_v$. 

\vskip .1 in
\begin{corollary} \label{CC3}
$\cB(0,1)$ is a metric space.
\end{corollary}

\begin{proof}
A compact Hausdorff space is ``T3'', 
that is, each point is closed, 
and for each point $x$ and each closed set $A$ with $x \notin A$, 
there are disjoint open neighborhoods $U$ of $x$ and $V$ of $A$.  
Urysohn's Metrization theorem (\cite{Kelley}, p.125) 
says that any T3 space with a countable basis is metrizable.  
\end{proof}

\vskip .1 in
It is important to note that the path distance function $\rho(x,y)$,
the length of the shortest path from $x$ to $y$ in $\overline{\Lambda}$,
is {\it not} a metric defining the Gel'fond topology.  
For example, if $p = B(0,1)$ is the root of $\overline{\Lambda}$, then
$\{x \in \overline{\Lambda} : \rho(x,p) < 1/2\}$ does not contain any points
$a$ of type I, while every neighborhood of $p$ in the 
Gel'fond topology contains infinitely many such points.  This same
example shows $\rho(x,y)$ is not even continuous for the Gel'fond topology.
(In fact, the topology defined by $\rho(x,y)$ is strictly finer than the 
Gel'fond topology.)

Tracing through the proof of Urysohn's theorem, one can construct
a metric defining the Gel'fond topology as follows.  
First, define a `separation kernel' for discs $B(a,r)$ 
and points $a^{\prime} \in B(0,1)$ by 
\begin{equation*}
\Delta(B(a,r),a^{\prime}) \ = \ 
\sup_{z \in B(a,r)} (|z-a^{\prime}|_v) \ = \ \max(r,|a-a^{\prime}|_v) \ .
\end{equation*}
Extend it to pairs of discs $B(a,r)$, $B(a^{\prime},r^{\prime})$ by 
\begin{equation*}
\Delta(B(a,r),B(a^{\prime},r^{\prime})) \ = \ 
\sup \begin{Sb} z \in B(a,r) \\ w \in B(a^{\prime},r^{\prime}) \end{Sb} 
              (|z-w|_v) 
\ = \ \max(r,r^{\prime},|a-a^{\prime}|_v) \ ,
\end{equation*}
and then to arbitrary points $x, y \in \cB(0,1)$ by  
\begin{equation*}
\Delta(x,y) \ = \ \lim_{i \rightarrow \infty}
\max(r_i,r_i^{\prime},|a_i-a_i^{\prime}|_v) \ .
\end{equation*}
if $x$, $y$ correspond to  
sequences of nested discs $\{B(a_i,r_i)\}$, $\{B(a_i^{\prime},r_i^{\prime})\}$.

Let $\{\alpha_j\}$ be an enumeration of $\tK \cap B(0,1)$ (or more generally,
take any countable dense subset of $\cB(0,1)$), and define a map 
$\varphi$ from $\cB(0,1)$ to the infinite-dimensional unit cube $[0,1]^{\NN}$ 
by putting
\begin{equation*}
\varphi(x) \ = \ (\Delta(x,\alpha_j))_{j \in \NN}
\end{equation*}
It can be checked that $\varphi$ is a topological isomorphism from $\cB(0,1)$ 
onto its image, equipped with the induced topology.  
Pulling back the metric on $[0,1]^{\NN}$, 
we obtain a metric defining the topology on $\cB(0,1)$:  
\begin{equation*}
d(x,y) \ = \ \sum_{j=1}^{\infty} 
     \frac{1}{2^j} |\Delta(x,\alpha_j)-\Delta(y,\alpha_j)| \ .
\end{equation*}

However, this formula seems nearly useless for understanding the Gel'fond 
topology.
It is much better to visualize the open sets in $\cB(0,1)$ in one of the
following ways:  
\vskip .05 in

A)  In terms of $\overline{\Lambda}$, the basic open sets are the
      sets of the following three types:  open branches of $\overline{\Lambda}$, 
      open branches with a finite number of closed branches removed, 
   and $\overline{\Lambda}$ with a finite number of closed branches removed.
\vskip .1 in
   
B)  In terms of the classical disc $B(0,1)$, the basic open sets are the 
subsets of $\cB(0,1)$ 
consisting of all points associated to an open disc $B(a,r)^-$, 
to a `punctured open disc'
$B(a,r)^- \backslash \cup_{i=1}^N B(a_i,r_i)$, 
or to a `punctured whole disc'  
$B(0,1) \backslash \cup_{i=1}^N B(a_i,r_i)$.  
(We say $x$ is associated to  
$B(a,r)^-\backslash \cup_{i=1}^N B(a_i,r_i)$ if it is associated to $B(a,r)^-$
but it is not associated to any $B(a_i,r_i)$).

\section{The Berkovich Projective Line.}
\label{Section B}

\vskip .1 in

In this section we will describe the sheaf of functions on the
Berkovich disc, making it a locally ringed space.
Then we will discuss the gluing process used
to assemble the Berkovich affine and projective lines.

We follow Berkovich's original gluing procedure, described in
(\cite{Berk}, Chapter 3), which requires gluing on quasi-affinoid open sets.  
This suffices for constructing the Berkovich
analytic spaces corresponding to algebraic varieties.  Later (\cite{Berk2})
Berkovich gave a more sophisticated gluing procedure using {\it nets},
which allows the construction of Berkovich analytic spaces
corresponding to arbitrary rigid analytic spaces.

We simplify Berkovich's exposition by restricting the gluing process to  
affinoids corresponding to classical Tate algebras;  actually
(\cite{Berk}) permits gluing along a larger class of affinoid-like spaces.

\vskip .1 in
Let us begin by reinterpreting the points of the Berkovich disc.
By the results in Section~\ref{Section A}, each $x \in \cB(0,1)$ 
corresponds to an equivalence class of nested discs $\{B(a_i,r_i)\}$ 
in $B(0,1)$, and for $f \in \CC_v\{\{T\}\}$ the corresponding seminorm $[f]_x$ 
is a limit of the sup norms $[f]_{B(a_i,r_i)}$.   
However, this association of $[f]_x$ with sup norms is misleading.   
A more accurate assertion is that $[f]_x$ 
is the {\it generic value} of $|f(z)|_v$ at $x$.

To see this, suppose $x$ is a point of type II,
so that $x$ corresponds to a disc $B(a,r)$
with $r$ in the value group of $\CC_v^{\times}$.
Proposition~\ref{BerkProp2} asserts that for $f \in \CC_v\{\{T\}\}$,
\begin{equation*}
[f]_x \ = \ \max_{z \in B(a,r)} |f(z)|_v \ .
\end{equation*}
If the zeros of $f(z)$ in $B(a,r)$ are $a_1, \ldots, a_m$,
then by the Weierstrass Preparation Theorem
$|f(z)|_v$ is takes on its maximum value on $B(a,r)$
at each point of $B(a,r) \backslash \cup_{i=1}^m B(a_i,r)^-$.
In other words, $[f]_x$ is the constant value which $|f(z)|_v$ assumes
`almost everywhere' on $B(a,r)$.  The multiplicative seminorm
$[ \cdot ]_x$ extends in a unique way to the quotient
field of $\CC_v\{\{T\}\}$, with
\begin{equation*}
[f/g]_x \ = \ \frac{[f]_x}{[g]_x} \ .
\end{equation*}
However, this extended seminorm is definitely {\it not} the sup norm:  if
$(f/g)(z)$ has poles in $B(a,r)$ then 
$\sup_{z \in B(a,r)} |(f/g)(z)|_v = \infty$.
Rather, if $f(z)$ has zeros $a_1, \ldots, a_m$ and $g(z)$ has zeros
$b_1, \ldots, b_n$, then $[f/g]_x$ is the constant value
which $|f(z)/g(z)|_v$ assumes everywhere on the `punctured disc'
$B(a,r) \backslash (\cup_{i=1}^m B(a_i,r)^-) \cup \cup_{j=1}^n B(b_j,r)^-)$.
This is best understood as the `generic value' of $|(f/g)(z)|_v$ on $B(a,r)$.

For points $x$ of type I, III or IV, the notion of a generic value
of $|f(z)|_v$ at $x$ has to be interpreted in a slightly
broader way.  Let $x \in \cB(0,1)$ be arbitrary.
By continuity, for each $\varepsilon > 0$ there
is a neighborhood $U$ of $x$ in $\cB(0,1)$ such that for each $t \in U$,
$|[f]_t - [f]_x| < \varepsilon$.  In particular, for each type I point
$z \in U$,
\begin{equation*}
|\ |f(z)|_v - [f]_x \ | \ < \ \varepsilon \ \ .
\end{equation*}
By the description of the topology of $\cB(0,1)$ in
Section~\ref{Section A},
sets of the form $\cB(a,r) \backslash \cup_{i=1}^m \cB(a_i,r_i)^-$
are cofinal in the set of closed neighborhoods of $x$.
Thus, $[f]_x$ is the unique number such that for each
$\varepsilon > 0$, there is a punctured disc
$B(a,r) \backslash \cup_{i=1}^m B(a_i,r_i)^-$ corresponding to a closed
neighborhood of $x$ such that $|f(z)|_v$ is within $\varepsilon$ of $[f]_x$
on that punctured disc.  In this sense $[f]_x$ is the generic
value of $|f(z)|_v$ at $x$.

\vskip .1 in
In rigid analysis, each punctured disc
$V = B(a,r) \backslash (\cup_{i=1}^m B(a_i,r_i)^-)$ with $r \le 1$ and 
$r, r_1, \ldots, r_m$ in the value group of $\CC_v^{\times}$
corresponds to an {\it affinoid subdomain} of $B(0,1)$.
More precisely, if $b, b_1, \ldots b_m \in \CC_v$ are such that $|b|_v = r$
and $|b_i|_v = r_i$ for each $i$, then $V$ is a {\it Laurent domain}, 
isomorphic to the set of maximal ideals $\Max(\cA_V)$ of the Tate algebra
\begin{equation}
\cA_V \ = \ \CC_v\{\{T,T_1, \ldots, T_m\}\} [X]/\cI_V \ .  \label{AFX1}
\end{equation}
Here 
\begin{equation*}
\CC_v\{\{T,T_1,\ldots,T_m\}\} \ = \ 
\{ \sum_{\vec{i} \ge 0} c_{\vec{i}} T^{i_0} T_1^{i_1} \cdots T_m^{i_m} 
\in \CC_v[[\vec{T}]] : \lim_{|\vec{i}| \rightarrow \infty} 
|c_{\vec{i}}|_v = 0 \} 
\end{equation*}
is the ring of power series converging on the unit polydisc 
$\{(z_0,z_1,\ldots,z_m) \in \CC_v^{m+1} : \max(|z_i|_v) \le 1\}$, 
and 
\begin{equation*}
\cI_V \ = \ (bT -(X-a), (X-a_1)T_1 - b_1, \ldots, (X-a_m)T_m -b_m) \ .
\end{equation*}
The relations generating $\cI_V$ mean that 
$\Max(\cA_V)$ is isomorphic to the set of $x \in \CC_v$ for which
$(\frac{x-a}{b},\frac{b_1}{x-a_1},\ldots,\frac{b_m}{x-a_m})$ belongs
to the unit polydisc, or equivalently, that $|x-a|_v \le r$
and $|x-a_i|_v \ge r_i$ for $i = 1, \ldots, m$.  This is precisely
the punctured disc $V$.
  
The localization of $\CC_v\{\{T\}\}$ at 
$(T-a_1)/b_1, \ldots, (T-a_m)/b_m$ is dense in $\cA_V$.
Hence, writing $X$ for the element $bT+a \in \CC_v\{\{T\}\}$, 
each bounded multiplicative seminorm $[ \cdot ]_x$ on $\CC_v\{\{T\}\}$ 
for which 
\begin{equation} \label{EXB2}
[X]_x \le 1 \ \ \text{and} \ \ 
[b_1/(X-a_1)]_x \le 1, \ \ldots, \ [b_m/(X-a_m)]_x \le 1 
\end{equation} 
extends to a bounded multiplicative seminorm on $\cA_V$.  Conversely,
each bounded multiplicative seminorm on $\cA_V$ restricts to a bounded
multiplicative seminorm on $\CC_v\{\{T\}\}$ for which (\ref{EXB2}) holds.  
Thus, the function analytic spectrum $\cM(\cA_V)$ can be identified 
with the subset $\cB(a,r) \backslash \cup_{i=1}^m \cB(a_i,r_i)^-$ 
of $\cB(0,1)$.  

We have emphasized the reinterpretation of seminorms $[ \cdot ]_x$ 
as `generic values' partly to explain why such an inclusion is reasonable.  
Since $V$ is a more complicated domain than $\cB(0,1)$, 
one should not expect to identify seminorms $[ \cdot ]_x$ 
on $\cA_V$ with sup norms on nested sequences of discs;  
however, one might hope to identify them with sup norms on nested sequences of
punctured discs contained in $B(a,r) \backslash \cup_{i=1}^m B(a_i,r_i)^-$.  
Such an identification can in fact be made,
though we will not prove it here.  However, we do note that 
each punctured disc $x = B(b,t) \backslash \cup_{i=1}^n B(b_i,t)^-$, in which
the deleted discs have the same radius as the outer disc, 
determines a bounded multiplicative seminorm 
\begin{equation*}
[f]_x \ = \ \sup_{z \in B(b,t) \backslash \cup_{i=1}^m B(b_i,t)^-} |f(z)|_v 
\end{equation*}
on any ring $\cA_V$ for which 
$V = B(a,r) \backslash \cup_{i=1}^m B(a_i,r_i)^-$ contains 
$B(b,t) \backslash \cup_{i=1}^n B(b_i,t)^-$.

\vskip .1 in
We will now discuss the sheaf of functions on $\cB(0,1)$ making it into
a locally ringed space.  We give the general definition of 
the sheaf of functions on an affinoid Berkovich space, 
then illustrate it for $\cB(0,1)$.

Note that $\CC_v\{\{T_1, \ldots, T_m\}\}$, given the Gauss norm 
$\|f\| = \max_{\vec{i}}(|c_{\vec{i}}|_v)$ (which coincides with the sup norm
over the unit polydisc), is a Banach algebra over $\CC_v$.  
Let $\cA$ be a Tate algebra, that is, a quotient  
$\CC_v\{\{T_1, \ldots, T_m\}\}/\cI$ for some $m$.  Equipped with the 
quotient norm $\|f\|_{\cI} = \inf_{g_0 \in \cI} \|f+g_0\|$, it too 
is a Banach algebra over $\CC_v$.  The maximal spectrum $\Max(\cA)$ is a 
(Tate) affinoid rigid analytic space.  By definition, 
the affinoid Berkovich space associated to $\cA$ is the functional analytic 
spectrum $X = \cM(\cA)$, the set of all bounded multiplicative seminorms
$[ \cdot]_x : \cA \rightarrow \RR_{\ge 0}$, equipped with the Gel'fand topology.

A closed set $V \subseteq X$ will be called a (strict)  
{\it affinoid Berkovich subdomain} 
of $X$ if there exist a Tate algebra $\cA_V$ and a continuous homomorphism 
$\varphi: \cA \rightarrow \cA_V$ satisfying the following universal property:  
for each continuous homomorphism of Tate algebras $\tau : \cA \rightarrow \cB$
such that the pullback $\tau^{*}(\cM(\cB))$ is contained in $V$, there 
is a unique continuous homomorphism $\tilde{\tau} : \cA_V \rightarrow \cB$ 
with $\tau = \tilde{\tau} \circ \varphi$.  In that case, it can be shown that
$V \cong \cM(\cA_V)$ (\cite{Berk}, Proposition 2.2.4, p.29). 

The closure of any affinoid subdomain of $\Max(\cA)$ in $X$ 
is an affinoid Berkovich subdomain of $X$.  
Such domains form a basis for the  
closed neighborhoods of any point $x \in X$ 
(\cite{Berk}, Proposition 2.2.3, p.28;  recall that a closed neighborhood of
$x$ means the closure of an open set containing $x$).    
The intersection of two affinoid Berkovich subdomains is
again an affinoid Berkovich subdomain, but the union of two 
affinoid Berkovich subdomains need not be one.  Nonetheless, 
if $V_1, \ldots, V_m$ are affinoid Berkovich subdomains of $X$, 
and if $\cV = V_1 \cup \ldots \cup V_m$, there is a natural $\CC_v$-algebra
$\cA_{\cV}$ associated to $\cV$, 
\begin{equation*}
\cA_{\cV} \ = \ \Ker(\prod_{i=1}^m \cA_{V_i} \rightarrow 
                  \prod_{i,j} \cA_{V_i \cap V_j}) \ .  
\end{equation*}
Such a $\cV$ will be called a {\it special} subset of $X$.
Given an open set $U \subset X$, put
\begin{equation} \label{FTB1}
\Gamma(U,\cO_X) \ = \ \lim \begin{Sb} \longleftarrow \\ \cV \subset U \\
\text{$\cV$ special} \end{Sb} \cA_{\cV}\ .
\end{equation}  
It is a consequence of Tate's acyclicity theorem that this construction defines
a sheaf of rings $\cO_X$, the {\it Berkovich structure sheaf}, on $X$.  

For each $x \in X$, the stalk of $\cO_X$ at $x$ is defined by
\begin{equation*}
\cO_{X,x} \ = \ 
   \lim \begin{Sb} \longrightarrow \\ x \in U \\ \text{$U$ open} \end{Sb}
              \Gamma(U,\cO_X) \ ;
\end{equation*}              
since strict affinoid Berkovich subdomains are cofinal in the 
neighborhoods of $x$, we also have
\begin{equation}  \label{AFX2}
\cO_{X,x} \ = \  
 \lim \begin{Sb} \longrightarrow \\ x \in V \\ 
                 \text{$V$ strict affinoid} \end{Sb} \cA_V \ .
\end{equation}
The seminorms $[ \cdot ]_x$ on the $\cA_V$ with $x \in V$ are compatible, 
so they induce a seminorm $[ \cdot ]_x$ on $\cO_{X,x}$.  
If $\wp_x$ is the kernel of $[ \cdot ]_x$ on $\cO_{X,x}$,  
then $[ \cdot ]_x$ descends to a norm on the residue field $\cO_{X,x}/\wp_x$. 
Let $\cK(x)$ be the completion of $\cO_{X,x}/\wp_x$ under $[ \cdot ]_x$.  
It is a valuation field, with norm again denoted $[ \cdot ]_x$.  

There is a natural homomorphism from $\cO_{X,x}$ to $\cK(x)$.
For any function $f \in \cA$, its {\it value} $f(x)$ at $x$ 
is defined to be the image of $f$ in $\cK(x)$.  The norm of $f(x)$ under
the absolute value on $\cK(x)$ is simply $[f(x)]_x = [f]_x$.  
\vskip .1 in

We will now consider how these definitions play out for the Berkovich disc.  

Let $U = \cB(a,r)^- \backslash \cup_{j=1}^m \cB(a_j,r_j)$
be a basic open set in $\cB(0,1)$.
Take a sequence of numbers $s_1,s_2, \ldots$ in the value group of
$\CC_v^{\times}$ which increase monotonically to $r$, and for each
$j = 1, \ldots, m$ take a sequence of numbers $s_{1j}, s_{2j}, \ldots$ 
in the value group of $\CC_v$ which decrease monotonically to $r_i$.
Then the affinoid Berkovich domains
$V_i = B(a,s_i) \backslash \cup_{i=1}^m B(a_j,s_{ij})^-$ form an
exhaustion of $U$.  Let $A_{V_i}$ be the Tate algebra associated to $V_i$
as in (\ref{AFX1}).  If $k > i$, the natural restriction map
$\rho_{k,i} : A_{V_k} \rightarrow A_{V_i}$ is injective, so we can view
$A_{V_k}$ as a subring of $A_{V_i}$.  Thus the inverse limit in (\ref{AFX2})
becomes an intersection:   
\begin{equation*}
\Gamma(U,\cO_X) \ = \ \bigcap_{i=1}^{\infty} A_{V_i} \ .
\end{equation*}
For example, if $U = \cB(0,1)^-$, let $s_1, s_2, \ldots$ increase 
monotonically to $1$, and put $V_i = \cB(0,s_i)$.  Then 
\begin{equation*}
A_{V_i} \ = \ \{ \sum_{k=0}^{\infty} c_k T^k \in \CC_v[[T]] :
       \lim_{k \rightarrow \infty}  s_i^k |c_k|_v = 0 \} 
\end{equation*}
and
\begin{equation}  \label{AFX3}
\Gamma(U,\cO_X) \ = \
\{ \sum_{k=0}^{\infty} c_k T^k  \in \CC_v[[T]] :
       \text{$\lim_{k \rightarrow \infty} s^k |c_k|_v = 0$
              for each $0 < s < 1$} \} \ .
\end{equation}
Note that $\Gamma(U,\cO_X)$ is strictly larger than the ring of power
series with bounded coefficients, the ring usually associated
with the open unit ball $B(0,1)^-$.  In particular, $\Gamma(U,\cO_X)$ contains
all power series $\sum_{k=0}^{\infty} c_k T^k$ with $|c_k|_v \le \log(k)$ for
each $k$.

\vskip .1 in
Next take $f \in \cA = \CC_v\{\{T\}\}$, and consider the values 
$f(x) \in \cK(x)$ for the various types of $x \in \cB(0,1)$.  

If $x$ is of type I, so it corresponds to a classical point $z \in B(0,1)$, 
then $\cO_{X,x}$ is the ring of germs of power series converging in a 
neighborhood of $z$, and $\wp_x$ is the ideal of germs of power 
series vanishing at $z$.  Thus $\cK(x)$ is canonically isomorphic to $\CC_v$,
$f(x)$ is the usual value $f(z) \in \CC_v$, and $[f(x)]_x = [f]_x = |f(z)|_v$. 
The subring of $\CC_v(T)$ consisting of rational functions with 
no poles at $z$ is dense in $\cO_{X,x}$. 
 
If $x$ is of type II, III, or IV, then the seminorm $[ \cdot ]_x$ is 
actually a norm.  The stalk $\cO_{X,x}$ contains the field of 
rational functions $\CC_v(T)$ as a dense subset, and $\wp_x = (0)$.  
The ring $\cA$ injects into $\cO_{X,x}$, 
which injects into its completion $\cK(x)$.   
The `value' $f(x)$ is simply $f$ itself, viewed as an element of $\cK(x)$.    
However, the structure of the valued field $\cK(x)$ varies with $x$, 
as does the norm of $f(x)$, which is $[f(x)]_x = [f]_x$.  We will now
attempt to make the fields $\cK(x)$ more explicit.  

If $x$ is of type II, it corresponds to a disc $B(a,r)$ with $r = |b|_v$
in the value group of $\CC_v^{\times}$, and the seminorm $[f]_x$ is
is given by the `generic value' of $|f(z)|_v$ on $B(a,r)$.  For polynomials 
$f(T) \in \CC_v[T]$, 
\begin{equation*}
[f]_x \ = \ \sup_{t \in B(a,r)} |f(t)|_v 
\ = \ \sup_{z \in B(0,1)} |f((z-a)/b)|_v \ , 
\end{equation*}
which coincides with the Gauss norm of $f((T-a)/b)$.  The field of 
rational functions $\CC_v(T)$ is dense in $\cO_{X,x}$, so $\cK(x)$ is 
isomorphic to the completion of $\CC_v(T)$ with respect to a  
Gauss norm.  This is a field whose value group coincides with that of 
$\CC_v^{\times}$, and whose residue field is the field of rational 
functions $\overline{\FF}_v(T)$, where $\overline{\FF}_v$ 
is the residue field of $\CC_v$.  (Note that although all the fields $\cK(x)$
for points of type II are isomorphic to each other, they are completions
of $\CC_v(T)$ with respect to different Gauss norms.)  

If $x$ is of type III, it corresponds to a disc $B(a,r)$ with $r$ not in the
value group of $\CC_v^{\times}$, and
for linear polynomials, $[T-b]_x = r$ if $b \in B(a,r)$,
and $[T-b]_x = |b-a|_v$ if $b \notin B(a,r)$.  The value group of $\cK(x)$ 
is generated by $r$ and the value group of 
$\CC_v^{\times}$;  it is strictly larger than the value group of 
$\CC_v^{\times}$. If $f \in \CC_v(T)$
and $[f]_x = 1$, then $f$ must have the same number of zeros and poles in 
$B(a,r)$.  Given $b,c \in B(a,r)$, the image of $(T-b)/(T-c)$ in the residue
field of $\cK(x)$ is easily seen to be $1$.
It follows that the residue field of $\cK(x)$ is $\overline{\FF}_v$,
corresponding to the reduction of constants.  It is shown in 
(\cite{Berk}, p.21) that $\cK(x)$ is isomorphic to the field of 
convergent Laurent series 
\begin{equation*}
K_r \ = \ \{ f = \sum_{i=-\infty}^{\infty} a_i T^i \in \CC_v[[T,1/T]] : 
                \lim_{|i| \rightarrow \infty} r^i |a_i|_v = 0 \}  
\end{equation*}                 
equipped with the norm $\|f\| = \max_i (r^i |a_i|_v)$.  Two fields $K_{r_1}$
and $K_{r_2}$ are isomorphic if and only if $r_2 = |b|_v r_1$ for some 
$b \in \CC_v^{\times}$;  the isomorphism takes  
$\sum a_i T^i \rightarrow \sum a_i b^i T^i$.    
Thus, $\cK(x)$ is determined up to isomorphism by its value group.

If $x$ of type IV, it corresponds to a nested sequence of discs $\{B(a_i,r_i)\}$
with empty intersection and $\lim r_i = r > 0$, 
and for any $f \in \cA$ one has $[f]_x = [f]_{B(a_i,r_i)}$ 
for sufficiently large $i$.  Hence the value group 
of $\cK(x)$ coincides with the value group of $\CC_v^{\times}$. 
For each linear polynomial $T-a$, if $i$ is large enough that 
$a \notin B(a_i,r_i)$, then $[(T-a)/(a-a_i)]_x < 1$, so 
$1 - (T-a)/(a-a_i)$ specializes to $1$ in the residue field of $\cK(x)$.  
Since 
\begin{equation*}
T-a \ = \ (a_i-a) \cdot (1 - \frac{T-a_i}{a-a_i}) \ , 
\end{equation*}
the residue field of $\cK(x)$ coincides with $\overline{\FF}_v$,
the reduction of the field of constants $\CC_v$.  This means $\cK(x)$ is
an {\it immediate extension} of $\CC_v$, a complete valued field containing
$\CC_v$, having the same value group and residue field as $\CC_v$.  
Tautologically, this field contains an element $\alpha$ (the image of $T$ in 
the completion of $\CC_v(T)$ under $[ \cdot ]_x$) which satisfies
$|a_i - \alpha|_v = r_i$ for each $i$, and hence belongs to the
intersection of the discs $B(a_i,r_i)$ in $\cK(x)$.

For our purposes, it will not be necessary to deal with the `exotic' fields 
$\cK(x)$ for points of types II, III, or IV.  We will mainly  
be interested in the value of $[f]_x$, which is the generic
value of $|f(z)|_v$ at type I points near $x$.   

\vskip .1 in

In gluing Berkovich spaces,
one faces a problem analogous to constructing open manifolds
in a category where the basic objects are closed manifolds with boundary.
Eliminating the effects of the boundary requires care.  For example,
consider the Berkovich disc $X = \cB(0,1)$ and its open subset $U = \cB(0,1)^-$.
The affinoid algebra $\CC_v\{\{T\}\} = \Gamma(X,\cO_X)$ is very different
from the ring $\Gamma(U,\cO_X)$ described in (\ref{AFX3}):  
the first is a Banach algebra;  the second is not.  We will 
give the general definitions for gluing Berkovich spaces, 
then construct the Berkovich projective line.  

\vskip .1 in

A morphism of affinoid Berkovich domains
$\varphi : (X_1,\cO_{X_1}) \rightarrow (X_2,\cO_{X_2})$  is 
a morphism of locally ringed spaces induced by a bounded homomorphism
$\varphi^* : \cA_{X_2} \rightarrow \cA_{X_1}$ between the
Tate algebras $\cA_{X_2} = \Gamma(X_2,\cO_{X_2})$ and
$\cA_{X_1} = \Gamma(X_1,\cO_{X_1})$.

A {\it quasi-affinoid Berkovich domain} $(U,\cO_U)$
is an open subset $U$ of an affinoid Berkovich domain $X$, equipped 
with the structure sheaf $\cO_U = \cO_X|_U$.  A morphism of quasi-affinoid
Berkovich domains $\varphi : (U_1,\cO_{U_1}) \rightarrow (U_2,\cO_{U_2})$
is a morphism of locally ringed spaces such that on each pair of affinoid
Berkovich subdomains $V_1 \subset U_1$, $V_2 \subset U_2$ with
$\varphi(V_1) \subset V_2$, the induced morphism of affinoid algebras
$\varphi^* : \cA_{V_2} \rightarrow \cA_{V_1}$ is bounded.  This condition
is automatically satisfied if $\varphi$ is induced by a morphism
of affinoid Berkovich domains
$\hat{\varphi} : (X_1, \cO_{X_1}) \rightarrow (X_2, \cO_{X_2})$,
where $U_1$ is an open subset of $X_1$, and $U_2$ is an open subset of $X_2$.

A {\it Berkovich analytic space} is a locally
ringed space $(X,\cO_X)$ such that each $x \in X$ has a neighborhood
isomorphic to a quasi-affinoid Berkovich domain.  More precisely,
a locally ringed space $(X,\cO_X)$ is a Berkovich analytic space 
if there is an open covering $\{U_i\}$ of $X$  such that

\vskip .05 in
$1)$ for each $i$ there is an isomorphism $\varphi_i$ of $(U_i,\cO_X|_{U_i})$
with a quasi-affinoid Berkovich domain $(\tilde{U}_i,\cO_{\tilde{U_i}})$ 
(a {\it chart});

$2)$ the charts are compatible in the following sense.  
For each pair $(i,j)$ with $U_i \cap U_j \neq \phi$, if
$\tilde{U}_{ij} = \varphi_i(U_i \cap U_j) \subset \tilde{U}_i$ and 
$\tilde{U}_{ji} = \varphi_j(U_i \cap U_j) \subset \tilde{U}_j$, then 
$(\tilde{U}_{ij},\cO_{\tilde{U}_i}|_{\tilde{U}_{ij}})$ and
$(\tilde{U}_{ji},\cO_{\tilde{U}_j}|_{\tilde{U}_{ji}})$ are quasi-affinoid
Berkovich subdomains of $(\tilde{U}_i,\cO_{\tilde{U_i}})$, 
$(\tilde{U}_j,\cO_{\tilde{U_j}})$ respectively, 
and the morphism $\tilde{\varphi}_{ij} = \varphi_j \circ \varphi_i^{-1}$ 
between them is an isomorphism of quasi-affinoid Berkovich domains.

\vskip .05 in
Such collection of charts is called an {\it atlas}.  
There is a obvious notion of compatibility of atlases.  
An equivalence class of atlases defines an analytic structure.

Conversely, given a collection of quasi-affinoid Berkovich domains 
$(\tilde{U}_i,\cO_{\tilde{U_i}})$ with distinguished quasi-affinoid
subdomains $(\tilde{U}_{ij},\cO_{\tilde{U}_i}|_{\tilde{U}_{ij}})$ and
isomorphisms
$\tilde{\varphi}_{ij} : (\tilde{U}_{ij},\cO_{\tilde{U}_i}|_{\tilde{U}_{ij}})
\rightarrow (\tilde{U}_{ji},\cO_{\tilde{U}_j}|_{\tilde{U}_{ji}})$
satisfying evident compatibility conditions, 
one can glue the $(\tilde{U}_i,\cO_{\tilde{U_i}})$ to make a Berkovich 
analytic space $(X,\cO_X)$.  

A morphism of Berkovich analytic spaces is a morphism of locally ringed
spaces $\varphi : (X_1,\cO_{X_1}) \rightarrow (X_2,\cO_{X_2})$, 
such that for each $x \in X_1$ there are affinoid (closed) neighborhoods 
$U$ of $x$ and $V$ of $\varphi(x)$ such that $\varphi(U) \subset V$ and 
$\varphi|_U : U \rightarrow V$ is a morphism of affinoid Berkovich domains.

\vskip .1 in
We will now construct the Berkovich affine line and projective line.

There are two constructions for the Berkovich affine line.  
The first views $\AA^1_{\Berk}$ as a union of discs.
Given $r$ in the value group of $\CC_v^{\times}$, let
$\CC_v\{\{r^{-1} T \}\} =
\{ \sum c_k T^k \in \CC_v[[T]] : \lim_{k \rightarrow \infty} r^k |c_k|_v \}$
be the ring of power series converging on $B(0,r)$. The Berkovich
disc of radius $r$ is $\cB(0,r) = \cM(\CC_v\{\{r^{-1}T\}\})$.
If $r_1 < r_2$, then $\cB(0,r_1)$ is an affinoid subdomain of $\cB(0,r_2)$,
and $\cB(0,r_1)^-$ is a quasi-affinoid subdomain of $\cB(0,r_2)^-$.  Put
\begin{equation*}
\AA^1_{\Berk} \ = \ \cup_{r > 0} \cB(0,r)^- 
\end{equation*}
and give $\AA^1_{\Berk}$ the structure sheaf $\cO_X$ 
defined in (\ref{FTB1}).
As a topological space, $\AA^1_{\Berk}$ 
is also given by $\cup_{r > 0} \cB(0,r)$.  Each $\cB(0,r)$ can be
identified with the set of equivalence classes of sequences of nested discs
contained in $B(0,r)$, giving it a tree 
structure like that of $\cB(0,1)$.  These tree structures combine
to give a tree structure for $\AA^1_{\Berk}$.

The second construction is global:
it defines $\AA^1_{\Berk} = \cM(\CC_v[T])$, the functional analytic
spectrum of the polynomial ring $\CC_v[T]$.
The same argument which identifies points of $\cB(0,r)$ with equivalence classes
of sequences of nested discs $\{B(a_i,r_i)\}$ in $B(0,r)$ 
shows that $\cM(\CC_v[T])$ can be identified with
equivalence classes of arbitrary sequences of nested discs $\{\cB(a_i,r_i)\}$.  
Thus the local and global constructions lead to the same space.

\vskip .1 in
The most direct construction of the Berkovich projective line 
is to glue two open discs $\cB(0,r)^-$ with $r > 1$.
Formally, this is done as follows.  If $r = |b|_v > 1$,
then the closed annulus $\Ann(r,1/r) = \{ x \in \cB(0,r) : [T]_x \ge 1/r \}$
is the affinoid Berkovich domain $\cM(\cA_{r,1/r})$ where
\begin{equation*}
\cA_{r,1/r} \ = \ \CC_v\{\{T_0,T_1\}\}[X](X-bT_0,bXT_1-1) \ .
\end{equation*}
This ring has an involution $\varphi$ interchanging $X$ and $1/X$.
The corresponding involution 
$\tilde{\varphi} : \Ann(r,1/r) \rightarrow \Ann(r,1/r)$
takes the open annulus 
\begin{equation*}
\Ann(r,1/r)^- = \{x \in \Ann(r,1/r) : 1/r < [X]_x < r\}
\end{equation*}
to itself.  At the level of discs,
$\tilde{\varphi}$ takes each punctured disc $B(0,t) \backslash B(0,t)^-$
to the punctured disc $B(0,1/t) \backslash B(0,1/t)^-$, and takes each disc
$B(a,t) \subset A(r,1/r)$ with $0 \notin B(a,t)$
to $B(1/a,t/|a|_v^2)$.  At the level of points, it takes
$[ \cdot ]_{B(0,t)}$ to $[ \cdot ]_{B(0,1/t)}$,
and $[ \cdot ]_{B(a,t)}$ to $[ \cdot ]_{B(1/a,t/|a|_v^2)}$.
To construct $\PP^1_{\Berk}$, glue two copies of $\cB(0,r)^-$
along the quasi-affinoid subdomain $\Ann(r,1/r)^-$,
using $\tilde{\varphi} : \Ann(r,1/r)^- \rightarrow \Ann(r,1/r)^-$
as the gluing isomorphism.

It is easy to see that different choices of $r$ lead to isomorphic
copies of $\PP^1_{\Berk}$.  It follows that 
$\AA^1_{\Berk} = \cup_{r > 0} \cB(0,r)^-$ 
is an open analytic subspace of $\PP^1_{\Berk}$,
whose complement consists of a single point ${\infty}$.  
One can also construct $\PP^1_{\Berk}$ by gluing two copies
of $\AA^1_{\Berk}$.

\vskip .1 in
Each rational map $h \in \CC_v(T)$ induces a morphism of Berkovich analytic
spaces  $h : \PP^1_{\Berk} \rightarrow \PP^1_{\Berk}$.  At the level of
points, it is defined by 
\begin{equation*}
[f]_{h(x)} \ = \ [f \circ h]_x
\end{equation*}
for $f \in \CC_v(T)$.

\vskip .1 in
If one is only concerned with the underlying topological space, 
$\PP^1_{\Berk}$ can be obtained by gluing two copies of the closed disc
$\cB(0,1)$ along the annulus $\Ann(1,1) = B(0,1) \backslash B(0,1)^-$,
using the gluing morphism $X \leftrightarrow 1/X$.
This preserves path distances in the tree
underlying $\Ann(1,1)$, since it takes each disc
$B(a,t) \subset \Ann(1,1)$ to $B(1/a,t)$.
The resulting model of $\PP^1_{\Berk}$ is relatively easy to visualize:  
it consists of a copy of $\cB(0,1)$ with an additional branch 
(containing the point $\infty$) leading off 
`vertically' from the root point $[ \cdot ]_{B(0,1)}$.
We will call this model, with its natural path length structure,
the {\it small model} of $\PP^1_{\Berk}$.
In it, each point has distance at most $1$ from the root $[ \cdot ]_{B(0,1)}$.

\vskip .1 in
There is another model for $\PP^1_{\Berk}$ which will be more useful
to us.  It is obtained from the small model by reparametrizing the edges
in such a way that if $B(a,r_1) \subset B(a,r_2)$, then 
the edge from $[ \cdot ]_{B(a,r_1)}$ to $[ \cdot ]_{B(a,r_2)}$ has length 
$|\log_v(r_2) - \log_v(r_1)|$, where $\log_v(t)$ is the logarithm to the
base $q_v$ for which $|x|_v = q_v^{-\ord_v(x)}$.  

Recall that the model of the Berkovich disc $\cB(0,1)$ 
constructed in Section~\ref{Section A} was obtained by adjoining `ends' to
the union of all graphs of discs $\Gamma_S$, 
where $S = \{B(a_1,t_1), \ldots, B(a_n,t_n)\}$ 
runs over all finite sets of discs contained in $B(0,1)$.  
A similar construction applies for any $\cB(0,r)$, taking  
\begin{equation*}
\Gamma_S \ = \ \bigcup_{i=1}^n [t_i,r]_{B(a_i,r_i)}
\end{equation*}
where $[t_i,r]_{B(a_i,t_i)} = \{ B(a_i,t) : t_i \le t \le r \}$ is a 
`line of discs'.  Under this parametrization, the edge between $B(a,r_1)$
and $B(a,r_2)$ has length $|r_1-r_2|$.  However, we can logarithmically
reparametrize each line of discs $[r_1,r_2]_{B(a,r_1)}$, putting 
\begin{equation*}
\<\log(r_1),\log(r_2)\>_{B(a,r_1)} \ = \ 
\{ B(a_i,q_v^t) : \log_v(r_1) \le t \le \log_v(r_2) \} 
\end{equation*}
and giving it length $|\log_v(r_2)-\log_v(r_1)|$.  Carrying out
this reparametrization for each $\cB(0,r)$, we obtain a new `path
distance' $\rho(x,y)$ between each pair of points in $\AA^1_{\Berk}$.
Any two points of type II, III, IV are at finite distance from each other,
while points of type I are at infinite distance from
each other and from the points of type II, III, IV.  Extend $\rho(x,y)$
to  $\PP^1_{\Berk}$ by setting $\rho(x,\infty) = \rho(\infty,x) = \infty$
for each $x \in \AA^1_{\Berk}$.  We will call the model of $\PP^1_{\Berk}$
obtained in this way {\it big model}.  

Note that the big and small models of $\PP^1_{\Berk}$ have the same underlying
point set, and the same topology.  Only our way of visualizing 
the distance structure on the edges has changed.  

The path distance in the big model is canonical, in the following sense.  

\begin{proposition} \label{PropB1}
  The path distance $\rho(x,y)$ on $\PP^1_{\Berk}$ 
  is independent of the choice of homogeneous coordinates on $\PP^1/\CC_v$.
\end{proposition}

\begin{proof}  Let  $h(z) = (az+b)/(cz+d) \in \CC_v(z)$ be a linear fractional
transformation.  It induces a continuous automorphism  $h$ of
of $\PP^1_{\Berk}$.  We claim that $\rho(h(x),h(y)) = \rho(x,y)$ for all $x, y$.  

Since points of type II are dense on the edges,
it suffices to show that $\rho(x,y)$ is preserved for such points.
By the remarks earlier in this section,
if $x$ corresponds to a disc $B(a,r)$, then $h(x)$ corresponds to
a disc $B(b,R)$ if there is a punctured subdisc
$B(b,R) \backslash \cup_{i=1}^m B(b_i,R)^-$
which is the image of a punctured $B(a,r) \backslash \cup_{i=1}^m B(a_i,r)^-$
for an appropriate choice of $a_1, \ldots, a_m$.
It also suffices to show that the distance is preserved when
$x$, $y$ correspond to discs $B(p,r_1) \subseteq B(p,r_2)$,
since $h$ preserves `betweenness' in the path structure on $\PP^1_{\Berk}$.

An arbitrary linear fractional transformation 
can be written as a composition of affine maps and inversions.
Hence it suffices to show that the path distance is preserved by
these types of maps.

First suppose $h(z) = az+b$ is affine.
The image of $B(p,r)$ under $h(z)$ is $B(ap+b,|a|_v r)$.
If $x, y \in \PP^1_{\Berk}$
correspond to $B(p,r_1) \subseteq B(p,r_2)$, it follows that
\begin{eqnarray*}
\rho(h(x),h(y)) & = & |\log_v(|a|_v r_1)-\log_v(|a|_v r_2)| \\
                & = & |\log_v(r_1) - \log_v(r_2)| \ = \ \rho(x,y) \ .
\end{eqnarray*}

Next suppose $h(z) = 1/z$.  Let $x \in \PP^1_{\Berk}$ correspond to a disc
$B(p,r)$.  If $B(p,r)$ does not contain $0$,
then $h(B(p,r)) = B(1/p,r/|p|_v^2)$.
If $B(p,r)$ does contain $0$, put $\dot{B}(0,r) = B(0,r) \backslash B(0,r)^-$;
then $h(\dot{B}(0,r)) = \dot{B}(0,1/r)$.  In either case, this 
determines $h(x)$.  Now let $x, y \in \PP^1_{\Berk}$
correspond to discs $B(p,r_1) \subseteq B(p,r_2)$.

If neither disc contains $0$,
then $B(1/p,r_1/|p|_v^2) \subset B(1/p,r_2/|p|_v^2)$ and so
\begin{eqnarray*}
\rho(h(x),h(y)) & = & |\log_v(r_1/|p|_v^2)-\log_v(r_2/|p|_v^2)| \\
                & = & |\log_v(r_1) - \log_v(r_2)| \ = \ \rho(x,y) \ .
\end{eqnarray*}
If both $B(p,r_1)$ and $B(p,r_2)$ contain $0$,
then $B(0,1/r_2) \subseteq B(0,1/r_1)$, so
\begin{eqnarray*}
\rho(h(x),h(y)) & = & |\log_v(1/r_1)-\log_v(1/r_2)| \\
                & = & |\log_v(r_1) - \log_v(r_2)| \ = \ \rho(x,y) \ .
\end{eqnarray*}
Finally, if $B(p,r_2)$ contains $0$ but $B(p,r_1)$ does not,
let $z \in \PP^1_{\Berk}$ correspond to $B(p,|p|_v)$.
Then $B(p,r_1) \subseteq B(p,|p|_v) \subseteq B(p,r_2)$,
so $\rho(x,y) = \rho(x,z) + \rho(z,y)$.
Since $B(p,|p|_v) = B(0,|p|_v)$ and $B(p,r_2) = B(0,r_2)$,
and since the smallest disc containing $B(1/p,r_1/|p|^2)$ and
$B(0,1/r_2)$ is $B(0,1/|p|_v)$,
\begin{eqnarray*}
\rho(h(x),h(y)) & = & |\log_v(r_1/|p|_v^2)-\log_v(1/|p|_v)|
                        + |\log_v(1/|p|_v) - \log_v(1/r_2)|_v \\
                & = & | \log_v(|p|_v) - \log_v(r_1)|
                        + | \log_v(r_2) - \log_v(|p|_v) | \\
                & = & \rho(x,y) \ .
\end{eqnarray*}
\end{proof}

\section{The Hsia kernel.}
\label{Section C}

\vskip .05 in
In this section we introduce the Hsia kernel, 
the basic kernel for potential theory on $\AA^1_{\Berk}$.  It was  
initially defined in (Hsia, \cite{Hsia}) as a kernel on trees of discs.
We will relate it to other objects in potential theory:  
the potential kernel $j_z(x,y)$ on graphs,  
the spherical $v$-adic metric $\|x,y\|_v$ on $\PP^1(\CC_v)$,
and the canonical distance $[x,y]_{\zeta}$.

\subsection{Definition of the Hsia kernel.}

If $x \in \AA^1_{\Berk}$ corresponds to a sequence of nested discs 
$\{B(a_i,r_i)\}$, we call $r = \lim_{i \rightarrow \infty} r_i$
the {\it radius} of $x$, and write $r = \radius(x)$.  
Note that the radius depends on the choice of coordinates on $\AA^1$.  

The Hsia kernel $\delta(x,y)_{\infty}$ is defined for
$x, y \in \AA^1_{\Berk}$ as follows.
If $x$ corresponds to a sequence of nested discs $\{B(a_i,r_i)\}$
and $y$ to $\{B(b_i,s_i)\}$, then
\begin{equation}  \label{FZC1} 
\delta(x,y)_{\infty} 
\ = \ \lim_{i \rightarrow \infty} \ \max(r_i,s_i,|a_i-b_i|_v) \ .
\end{equation}
Clearly the Hsia kernel is symmetric, 
and  $\delta(x,x)_{\infty} = \radius(x)$ for each $x$.
If $x, y \in \AA^1(\CC_v)$ are points of type I, 
then $\delta(x,y)_{\infty} = |x-y|_v$.  
If $x, y \in \AA^1_{\Berk}$ are points of types I, II or III, 
with $x$ corresponding to $B(a,r)$ and $y$ corresponding to $B(b,s)$, then
$\delta(x,y)_{\infty} = \max(r,s,|a-b|_v)$.
  
The Hsia kernel has the following geometric interpretation.  Consider the paths 
from $x$ to $\infty$ and $y$ to $\infty$, and let $w$ be the point where
they first meet.  Then 
\begin{equation*}
\delta(x,y)_{\infty} \ = \ \delta(w,w)_{\infty} \ = \ \radius(w) \ .
\end{equation*}
To see this, note that if $x$ lies on the path from $y$ to $\infty$ then 
$w = x$ and  
$\delta(x,y)_{\infty} = \lim_{i \rightarrow \infty}  r_i  = \radius(x)$.
Similarly if $y$ lies on the path from $x$ to $\infty$, then $w = y$ and
$\delta(x,y)_{\infty} = \radius(y)$.  Otherwise, $w \ne x, y$, so there 
there are disjoint discs $B(a_i,r_i)$ and $B(b_i,s_i)$ in the sequences 
defining $x$, $y$.   
Then $|a_i-b_i|_v = |a_j-b_j|_v > \max(r_j,s_j)$ for all $j \ge i$, 
and the point $w$ where the paths from $x$ and $y$ to $\infty$ meet 
corresponds to the disc $B(a_i,|a_i-b_i|_v) = B(b_i,|a_i-b_i|_v)$,
whose radius is  $|a_i-b_i|_v$.

\vskip .1 in
Recall that a real-valued function $f(x,y)$ is upper semicontinuous if 
for each $x_0, y_0$ 
\begin{equation*}
\limsup_{(x,y) \rightarrow (x_0,y_0)} f(x,y) \ \le \ f(x_0,y_0) \ .
\end{equation*}
This is equivalent to requiring that $f^{-1}((-\infty,b))$
be open for each $b \in \RR$.  
We will say that $f(x,y)$ is {\it strongly upper semicontinuous} if
for each $x_0, y_0$  
\begin{equation*}
\limsup_{(x,y) \rightarrow (x_0,y_0)} f(x,y) \ = \ f(x_0,y_0) \ .
\end{equation*}

\begin{proposition} \label{PropC1}
\ \ 

$A)$  The Hsia kernel is nonnegative, symmetric, 
and continuous in each variable separately.
As a function of two variables, it is strongly upper semicontinuous.
It is continuous off the diagonal, and continuous   
at $(x_0,x_0)$ for each point $x_0$ of type I, 
but is discontinuous at $(x_0,x_0)$ for each point of type II, III, or IV. 

$B)$ The Hsia kernel is the unique extension of $|x-y|_v$ to $\AA^1_{\Berk}$ 
such that 
\begin{equation} \label{FZD0A}
\delta(x,y)_{\infty} \ = \ 
\limsup \begin{Sb} (a,b) \rightarrow  (x,y) \\ 
        a, b \in \AA^1(\CC_v) \end{Sb} |a-b|_v   
\end{equation}
for each $x, y \in \AA^1_{\Berk}$.     

$C)$  For all $x, y, z \in \AA^1_{\Berk}$
\begin{equation*}
\delta(x,y)_{\infty} 
 \ \le \ \max( \delta(x,z)_{\infty}, \delta(y,z)_{\infty} )  
\end{equation*}  
with equality if $\delta(x,z)_{\infty} \ne \delta(y,z)_{\infty}$.

$D)$ For each $a \in \AA^1_{\Berk}$, $r \in \RR$, the `open disc' 
$\cB(a,r)_{\infty}^- := \{ x \in \AA^1_{\Berk} : \delta(x,a)_{\infty} < r \}$
is open.  It is empty if $r \le \radius(a)$, and coincides with an open
disc $\cB(b,r)^-$ for some $b \in \AA^1(\CC_v)$ if $r > \radius(a)$.

For each $a \in \AA^1_{\Berk}$, $r \in \RR$, the `closed disc' 
$\cB(a,r)_{\infty} := \{ x \in \AA^1_{\Berk} : \delta(x,a)_{\infty} \le r \}$
is closed.  It is empty if $r < \radius(a)$, and coincides with a closed 
disc $\cB(b,r)$ for some $b \in \AA^1(\CC_v)$ if $r > \radius(a)$ or if
$r = \radius(a)$ and $a$ is of type  II or III.  If $r = \radius(a)$
and $a$ is of type I or IV, then $\cB(a,r)_{\infty} = \{a\}$.  
\end{proposition}

\begin{proof} {\ \ }

(A) Consider $\delta(x,y)_{\infty}$ as a function of two variables. 
We first show it is continuous off the diagonal.

Take $x_0, y_0 \in \AA^1_{\Berk}$ with $x_0 \ne y_0$.  
If there are disjoint open discs with  
$x_0 \in \cB(a,r)^-$ and $y_0 \in \cB(b,s)^-$, then $|a-b|_v > \max(r,s)$, 
so for each $x \in \cB(a,r)^-$, $y \in \cB(b,s)^-$, 
\begin{equation*}
\delta(x,y)_{\infty} \ = \ |a-b|_v \ = \ \delta(x_0,y_0)_{\infty} \ ,
\end{equation*}
and $\delta(x,y)_{\infty}$ is continuous at $(x_0,y_0)$.  

Otherwise, either each open disc $\cB(a,r)^-$ containing $x_0$ contains $y_0$, 
or each open disc $\cB(b,s)^-$ containing $y_0$ contains $x_0$.  
Suppose the former holds.  Then $x_0$ must be is a point of type II or III
corresponding to a disc $B(a,r)$ with $y_0 \in \cB(a,r)$.  
After relabeling the center, we can assume that $y_0 \in \cB(a,r)^-$, 
and that $y_0 \in B(a,s)^-$ for some $s < r$.
Take $0 < \varepsilon < r-s$.  
Then $U := \cB(a,r+\varepsilon)^- \backslash \cB(a,r-\varepsilon)$ 
and $V := \cB(a,s)^-$ are disjoint open neighborhoods 
of $x_0$, $y_0$ respectively.  
Given $x \in U$, put $r_1 = \inf \{ t : x \in \cB(a,t) \}$.
Then  $r - \varepsilon < r_1 < r+\varepsilon$, and the geometric description
of the Hsia kernel implies that $\delta(x,y)_{\infty} = r_1$ 
for each $y \in V$, so 
\begin{equation*} 
|\delta(x,y)_{\infty} - \delta(x_0,y_0)_{\infty}| 
               \ = \ |r_1-r| \ < \ \varepsilon \ .
\end{equation*}               
Thus, again $\delta(x,y)_{\infty}$ is continuous at $(x_0,y_0)$.  

Now consider points on the diagonal.  Suppose $x_0 \in \AA^1_{\Berk}$
has radius $r$.  For each $\varepsilon > 0$, there is an open disc 
$\cB(a,r+\varepsilon)^-$ containing $x_0$, and the geometric description
of the Hsia kernel shows that  $\delta(x,y)_{\infty} < r+\varepsilon$
for each $x, y \in \cB(a,r+\varepsilon)^-$.  
If $r = 0$, this implies $\delta(x,y)_{\infty}$ is continuous at $(x_0,x_0)$.  

If $r > 0$, then $\delta(x,y)_{\infty}$ is not continuous at $(x_0,x_0)$, since 
every neighborhood of $(x_0,x_0)$ contains points $(a,a)$ with $a$ of type I, 
and for such points $\delta(a,a)_{\infty} = 0$, 
while $\delta(x_0,x_0)_{\infty} = r$.  
On the other hand, the discussion in the previous paragraph shows that
\begin{equation} \label{FZD2} 
\limsup_{(x,y) \rightarrow (x_0,x_0)} \delta(x,y)_{\infty} 
     \ = \ r \ = \ \delta(x_0,x_0)_{\infty} \ .
\end{equation}
Thus $\delta(x,y)_{\infty}$ is stongly upper semicontinuous as
a function of two variables.

\vskip .05 in
Now fix $x$, and consider $\delta(x,y)_{\infty}$ as a function of $y$.  
By what has been shown above, $\delta(x,y)_{\infty}$ 
is continuous for  $y \ne x$.  
For continuity at $x$, put $r = \radius(x)$ and consider a neighborhood
$\cB(a,r+\varepsilon)^-$ of $x$ as discussed above.  For each 
$y \in \cB(a,r+\varepsilon)^-$, the paths from $x$ and $y$ to $\infty$
meet at a point $w \in \cB(a,r+\varepsilon)^-$, so
\begin{equation*}
r \ = \ \radius(x) \ = \ \delta(x,x)_{\infty}
  \ \le \ \delta(x,y)_{\infty} \ = \ \radius(w) \le r + \varepsilon \ .
\end{equation*}
Thus $|\delta(x,y)_{\infty}-\delta(x,x)_{\infty}| < \varepsilon$.

\vskip .05 in  
(B) Note that each neighborhood $\cB(a,r+\varepsilon)^-$ of $x_0$ contains 
type I points $a \ne b$ with $r < |a-b|_v < r+\varepsilon$, 
so (\ref{FZD2}) remains true even if the $\limsup$ is restricted to  
$x,y \in \AA^1(\CC_v)$.  
The characterization of $\delta(x,y)_{\infty}$ follows from this, 
together with continuity off the diagonal.     

\vskip .05 in
(C)  For discs $B(a,r)$, $B(b,s)$ it is easy to see that
\begin{equation*} 
\max(r,s,|a-b|_v) \ = \ 
\sup \begin{Sb} p \in B(a,r), \  q \in B(b,s) \end{Sb} (|p-q|_v) \ .                            
\end{equation*}   
The ultrametric inequality for $\delta(x,y)$ follows from this, (\ref{FZC1}),
and the classical ultrametric inequality.  

\vskip .05 in
(D)  Given $a \in \AA^1_{\Berk}$ and $r \in \RR$, consider the `open disc'  
$\cB(a,r)_{\infty}^- := \{ x \in \AA^1_{\Berk} : \delta(x,a)_{\infty} < r \}$.
If $r \le \radius(a)$, it is clearly empty.  If $r > \radius(a)$ it contains
$a$, and is open by the upper semicontinuity of $\delta(x,a)_{\infty}$. 
Each nonempty open set contains points of type I.  
If $b \in \cB(a,r)^- \cap \AA^1(\CC_v)$, the ultametric inequality 
for $\delta(x,y)_{\infty}$ shows that $\cB(a,r)_{\infty}^- = \cB(b,r)^-$.

Likewise, consider the `closed disc' 
$\cB(a,r)_{\infty} := \{ x \in \AA^1_{\Berk} : \delta(x,a)_{\infty} \le r \}$.
If $r < \radius(a)$, it is empty.  
If $r = \radius(a)$ and $a$ is of type I or IV it is easy to see that
$\cB(a,r)_{\infty} = \{a\}$.
If $r > \radius(a)$ then it contains $\cB(a,r)_{\infty}^-$ and hence contains
points $b$ of type I;  this also holds if $r = \radius(a)$ 
and $a$ is of type II or III.  In either case the ultametric inequality 
for $\delta(x,y)_{\infty}$ shows that $\cB(a,r)_{\infty} = \cB(b,r)$,
and in particular that $\cB(a,r)_{\infty}$ is closed.    
\end{proof}

\vskip .1 in
The function-theoretic meaning of the Hsia Kernel is as follows:

\begin{corollary} \label{CorC2}
For each $a \in \CC_v$,  the function $T-a \in \CC_v(T)$
satisfies 
\begin{equation*} [T-a]_x  \ = \ \delta(x,a)_{\infty} \ .
\end{equation*}
\end{corollary}

\begin{proof}
By Proposition \ref{PropC1} (B),
if $x \in \AA^1_{\Berk}$ corresponds to the nested
sequence of discs $\{B(a_i,r_i)\}$ then
\begin{equation*}
[T-a]_x \ = \ \lim_{i \rightarrow \infty} \|T-a\|_{B(a_i,r_i)} 
        \ = \ \limsup \begin{Sb} b \rightarrow x \\ b \in \CC_v \end{Sb}
                    |b-a|_v 
        \ = \ \delta(x,a)_{\infty} \ .
\end{equation*}
\end{proof}

\vskip .1 in

The Hsia kernel $\delta(x,y)_{\infty}$ has a pole at $\infty$.
We now set out to generalize it to a kernel $\delta(x,y)_{\zeta}$
with a pole at an arbitrary point $\zeta \in \PP^1_{\Berk}$.  
To do so we will need some facts from the theory of metrized graphs.

\vskip .1 in

\subsection{Metrized Graphs.}

A metrized graph $\Gamma$ is a finite, connected graph with
a distinguished parametrization of each of its edges.
The natural path-length function defines a metric on $\Gamma$,
making it a compact metric space.  
By a {\it vertex set\ } for $\Gamma$ we mean a
finite set of points $S$ such that $\Gamma \backslash S$ 
is a union of open intervals whose closures have distinct endpoints.
A vertex set necessarily contains all endpoints and branch points of $\Gamma$.
If $\Gamma$ has loops, it also contains at least one interior point from each
loop.

We will primarily be interested in metrized graphs which are subsets of the
big model of $\PP^1_{\Berk}$, with the induced path length  $\rho(x,y)$.

For each $p \in \Gamma$, there are a finite number of edges $e_i$
emanating from $p$.  For each such edge, let $\gamma_i(t)$ be the arclength
parametrization of $e_i$ with $\gamma_i(0) = p$.  It will be useful to
introduce a formal `unit vector $\vec{v}_i$ emanating from $p$
in the direction of $e_i$', and write $p + t \vec{v_i}$ for $\gamma_i(t)$.
Given $f : \Gamma \rightarrow \RR$, let
\begin{equation*}
d_{\vec{v}_i}(f)(p)
\ = \ \lim_{t \rightarrow 0^+} \frac{f(p +t \vec{v}_i) - f(p)}{t}
\end{equation*}
be the one-sided `directional derivative' of $f$ at $p$ along $e_i$,
provided the limit exists.  

Let $\CPA(\Gamma)$ be the space of continuous, piecewise-affine,
real-valued functions on $\Gamma$.  (By a piecewise affine function $f$,
we mean that there is a vertex set $S_f$ for $\Gamma$ such that $f$
is affine on each edge in $\Gamma \backslash S_f$, with respect to an 
arclength parametrization of that edge.)
If $f \in \CPA(\Gamma)$, clearly the directional derivatives $d_{\vec{v}}(p)$
are defined for all $p$ and all $\vec{v}$ at $p$.  

Chinburg and Rumely (\cite{C-R}) introduced a Laplacian on $\CPA(\Gamma)$.
It is a map from $\CPA(\Gamma)$ to the space of discrete, signed measures
on $\Gamma$.  We will take the Laplacian to be the negative of theirs,
and put 
\begin{equation*}
\Delta(f) \ = \ \sum_{p \in \Gamma}
    (- \sum_{ \text{$\vec{v}$ at $p$} } d_{\vec{v}}(p) ) \ \delta_p(x)
\end{equation*}
where $\delta_p(x)$ is the Dirac measure at $p$.  By abuse of notation,
we will write $\Delta(f)(p)$ for
$\Delta(f)(\{p\}) = - \sum_{ \text{$\vec{v}$ at $p$} } d_{\vec{v}}(p)$.
Here are some elementary properties of $\Delta$, which show it
behaves as a Laplacian should.

\begin{proposition} \label{PropC3}  Let $f,g \in \CPA(\Gamma)$.  Then

$A)$ $\Delta(f) \equiv 0$ if and only if $f = C$ for some constant $C$.  

$B)$ $\Delta(f) = \Delta(g)$ if and only if $f = g + C$ for some constant $C$.

$C)$  If $f$ is nonconstant, then $f(x)$ achieves its maximum at a point $p$ 
where $\Delta(f)(p) > 0$, and its minimum at a point $q$
where $\Delta(f)(q) < 0$ .

$D)$ $\int_{\Gamma} f(x) \, \Delta(g)(x) = \int_{\Gamma} g(x) \, \Delta(f)(x)$.  

$E)$  The total mass $\Delta(f)(\Gamma) = 0$.

\end{proposition}

\begin{proof}

For (A), if $f = C$, clearly $\Delta(f) \equiv 0$.
Conversely, suppose $\Delta(f) \equiv 0$.
Put $M = \max_{x \in \Gamma} f(x)$, and let
$G = \{ x \in \Gamma : f(x) = M \}$.
Since $f$ is continuous,
$G$ is closed.  To see that $G$ is also open,
note that if $p \in G$, then $d_{\vec{v}}f(p) \le 0$
for each $\vec{v}$ at $p$, since $M$ is the maximum value of $f$.
Since $f$ is piecewise affine, if $f$ were not locally constant at $p$
then $d_{\vec{v}}(f)(p) < 0$ for some $\vec{v}$, and so
$\Delta(f)(p) > 0$ contrary to our assumption.  Thus, $G$ is both open
and closed, and since $\Gamma$ is connected, $G = \Gamma$.

Trivially (A) implies (B).

For (C), suppose $f$ is nonconstant, and let $M$ and $G$ be as in (A).
If $p \in \partial G$, then the argument above shows that $\Delta(f)(p) > 0$.
Similarly, $f$ achieves its minimum at a point where $\Delta(f)(q) < 0$.

For (D), let $S$ be a vertex set for $\Gamma$ such that 
$f$ and $g$ are both affine on each edge of $\Gamma \backslash S$.  
Given an edge $e_i$ of $\Gamma \backslash S$, fix an orientation of
it, and let $\gamma_i : [0,L_i] \rightarrow e_i$ be an arclength 
parametrization.  By abuse of notation, write 
$f(x) = f(\gamma_i(x))$ and $g(x) = g(\gamma_i(x))$ on $e_i$. 
Integration by parts shows that
\begin{equation*} 
\int_{\Gamma} f(x) \Delta(g)(x) 
= \sum_{i} \int_{e_i} f^{\prime}(x) g^{\prime}(x) \, dx 
= \int_{\Gamma} g(x) \Delta(f)(x) \ .
\end{equation*}
  
Part (E) follows from (D), taking $g(x) \equiv 1$.  
\end{proof}

\subsection{The potential kernel on a metrized graph.}

There is a potential kernel $j_z(x,y)$ on $\Gamma$
which inverts the Laplacian.
It is defined by the properties that for each fixed $y$, $z$,
\begin{eqnarray} \label{FZD4}
\left\{ \begin{array}{l}
           \text{$\Delta_x(j_z(x,y)) \ = \ \delta_y(x) - \delta_z(x)$ \ ,} \\
           \text{$j_z(z,y) \ = \ 0$ \ .} \end{array} \right.
\end{eqnarray}
Its uniqueness follows from Proposition \ref{PropC3} (B).   
Its existence can be shown using circuit theory (see \cite{C-R}), or 
using linear algebra (see \cite{Zh} or \cite{B-F}).  
When $\Gamma$ is a tree, its existence is trivial:  given $x, y, z \in \Gamma$,
let $w$ be the point where the path from $x$ to $z$ and the path
from $y$ to $z$ first meet, and put $j_z(x,y) = \rho(z,w)$, the path length
from $z$ to $w$.  Then along the path from $z$ to $y$, $j_z(x,y) = \rho(z,x)$;
on branches off that path, $j_z(x,y)$ is constant.
One easily checks that the function thus defined
satisfies the equations (\ref{FZD4}).

The potential kernel has the following physical interpretation.
View $\Gamma$ as an electric circuit
with terminals at $y$ and $z$, and with the resistance of each edge given
by its length.  Then $j_z(x,y)$ is the voltage at $x$ when current $1$ 
enters at $y$ and exits at $z$, with reference voltage $0$ at $z$.  
By its construction, $j_z(x,y)$ belongs to $\CPA(\Gamma)$, 
and its slope along any edge of $\Gamma$ has absolute value at most $1$.

\begin{proposition} \label{PropC4}
\ \

$A)$  $j_z(x,y)$ is non-negative, bounded, symmetric in $x$ and $y$, 
and jointly continuous in $x$, $y$, $z$.

$B)$  If $f \in \CPA(\Gamma)$ satisfies  
$\Delta(f) = \sum_{i=1}^n c_i \delta_{p_i}(x)$
then there is a constant such that
\begin{equation*}
f(x) \ = \ \sum_{i=1}^n c_i \, j_z(x,p_i) + C \ .
\end{equation*}

$C)$  For each $z, \zeta \in \Gamma$, \quad
\begin{equation} \label{FFD1} 
j_{\zeta}(x,y) \ = \ j_z(x,y) - j_z(x,{\zeta}) - j_z({\zeta},y)
                      + j_z(\zeta,\zeta) \ .
\end{equation}  
\end{proposition}

\begin{proof}

(A)  The non-negativity of $j_z(x,y)$ follows from its defining properties
(\ref{FZD4}) and Proposition \ref{PropC3} (C).   

For symmetry, fix $a, b, z$ and apply
Proposition \ref{PropC3} (D) to $f(x) = j_z(x,b)$ and $g(x) = j_z(x,a)$.
Then $\Delta(f) = \delta_b(x) - \delta_z(x)$ and 
$\Delta(g) = \delta_a(x)-\delta_z(x)$.  Since $f(z) = g(z) = 0$, 
\begin{equation*}
j_z(a,b) \ = \ \int_\Gamma f(x) \Delta(g)(x) 
         \ = \ \int_\Gamma g(x) \Delta(f)(x) \ = \ j_z(b,a) \ .
\end{equation*}

Clearly $j_z(x,y)$ is continuous in $x$ for each fixed $y$, $z$, with
$|j_z(x,y)-j_z(p,y)| \le \rho(x,p)$ since the absolute value of the  
slope of $j_z(x,y)$ along each edge is at most $1$.  From this, one
deduces that it is jointly continuous in $x, y$ for each fixed $z$, 
and indeed that for given $p, q$ 
\begin{eqnarray} \label{FZD6}
|j_z(x,y)-j_z(p,q)| & \le & |j_z(x,q)-j(p,q)| + |j_z(x,q)-j_z(x,y)| \notag \\
                    & \le & \rho(x,p) + \rho(y,q)  
\end{eqnarray}                    
using the symmetry shown above.  Finally, for any fixed $z_0$, we claim that
\begin{equation} \label{FZD5}
j_z(x,y) \ = \ j_{z_0}(x,y) - j_{z_0}(x,z) - j_{z_0}(z,y) + j_{z_0}(z,z) \ .      
\end{equation}
The joint continuity of $j_z(x,y)$ in $x$, $y$, and $z$ follows from this
and (\ref{FZD6}).  Its boundedness follows from continuity and the compactness
of $\Gamma$.  

To prove (\ref{FZD5}), note that for fixed $y$, $z$ 
\begin{equation*}
\Delta_x(j_z(x,y)) \ = \ \delta_y(x)-\delta_z(x) 
                  \ = \ \Delta_x(j_{z_0}(x,y) - j_{z_0}(x,z))  \ , 
\end{equation*} 
so by Proposition \ref{PropC3} (B) there is a constant $C_{y,z}$
such that $j_z(x,y) = j_{z_0}(x,y) - j_{z_0}(x,y) + C_{y,z}$.  
Taking $x = z_0$ shows $C_{y,z} = j_z(z_0,y)$.  
Similarly, applying $\Delta_y$ to $j_z(z_0,y)$ and $-j_{z_0}(z,y)$ 
(and using symmetry), 
one finds there is a constant $C_z$ such that 
$j_z(z_0,y) = -j_{z_0}(z,y) + C_z$.  
Taking $y = z$ shows $C_z = j_{z_0}(z,z)$.  
Combining these gives (\ref{FZD5}).  This also proves (C).   

\vskip .1 in
For part (B), first note that if $\Delta(f) = \sum_{i=1}^n c_i \delta_{p_i}(x)$ 
then $\sum_{i=1}^n c_i = 0$ by Proposition \ref{PropC3} (E).  Using this, 
one sees that 
\begin{equation*}
\Delta(\sum_{i=1}^n c_i j_{z_0}(x,p_i)) 
\ = \ \sum_{i=1}^n c_i \delta_{p_i}(x) \ , 
\end{equation*}
and the result follows from Proposition \ref{PropC3} (B).
\end{proof}   

\vskip .1 in
By a {\it subgraph} $\Gamma$ of $\PP^1_{\Berk}$ we will mean a connected
closed subgraph with a finite number of vertices and edges, which has finite 
total length under the logarithmic path distance $\rho(x,y)$.
Such a graph is necessarily a tree.  We view $\Gamma$
as a metrized graph, with the metric $\rho(x,y)$.
Since there is a unique path between any two points of $\PP^1_{\Berk}$, there
is a natural retraction map $r_{\Gamma} : \PP^1_{\Berk} \rightarrow \Gamma$.  

\vskip .1 in
The relation between the Hsia kernel and the
potential kernel $j_z(x,y)$ is as follows.
Fix a coordinate system on $\PP^1$,
so that $\PP^1_{\Berk} = \AA^1_{\Berk} \cup \{ \infty \}$.
Recall that $\zeta_0$ denotes the ``Gauss point'', 
the point in $\PP^1_{\Berk}$ corresponding to the $\sup$ norm over $B(0,1)$.  

\begin{proposition}  \label{PropC5}
Let $z \in \PP^1_{\Berk} \backslash \PP^1(\CC_v)$, and let $\Gamma$
be any subgraph of $\PP^1_{\Berk}$  containing $z$.  
Put $\inftybar = r_{\Gamma}(\infty)$.
Then for all $x, y \in \Gamma$,
\begin{equation} \label{FPV0}
-\log_v(\delta(x,y)_{\infty}) 
\ = \ j_z(x,y) - j_z(x,\inftybar) - j_z(y,\inftybar) - \log_v(\radius(z))
\end{equation}
In particular, if $z = {\zeta_0} \in \PP^1_{\Berk}$ is the Gauss point, 
then $\radius({\zeta_0}) = 1$, so
\begin{equation} \label{FPQ1}
-\log_v(\delta(x,y)_{\infty}) 
\ = \ j_{{\zeta_0}}(x,y) - j_{{\zeta_0}}(x,\inftybar)
                   - j_{{\zeta_0}}(y,\inftybar) \ .
\end{equation}
\end{proposition} 

\begin{proof}  Fix $y$, and set $s = \radius(y)$.
The intersection of the path from $y$ to $\infty$ with
$\Gamma$ is a closed segment $[s,S]$,
whose endpoints are $y$ and $\inftybar = r_{\Gamma}(\infty)$.
Put $\Gamma_y = [s,S]$;  it is a subgraph of $\Gamma$.

Consider the function
$f_y(x) = -\log_v(\delta(x,y)_{\infty})$ on $\Gamma$.
The geometric description of $\delta(x,y)_{\infty}$ shows that for
each $x \in \Gamma$
\begin{equation*}
\delta(x,y)_{\infty} \ = \ \delta(r_{\Gamma_y}(x),y)_{\infty} \ .
\end{equation*}
Thus, $f_y(x)$ is constant on branches off $\Gamma_y$.
For $x \in \Gamma_y$,
\begin{equation*}
-\log_v(\delta(x,y)_{\infty}) \ = \ -\log_v(t)
\qquad \text{where $t = \radius(x)$.}
\end{equation*}
The arclength parameter along $\Gamma_y$ is $\log_v(t)$ for $s \le t \le S$,
so the restriction of $f_y(x)$ to $\Gamma_y$ has constant slope $-1$.
Hence $f_y(x) \in \CPA(\Gamma)$, and 
\begin{equation*}
\Delta_x(f_y(x)) \ = \ \delta_y(x) - \delta_{\inftybar}(x) \ .
\end{equation*}
By (\ref{FZD4}),
\begin{equation*}
\Delta_x(j_z(x,y)-j_z(x,\inftybar)) \ = \
\delta_y(x) - \delta_{\inftybar}(x) \ .
\end{equation*}
It follows from Proposition \ref{PropC3} (B) that there is a constant $C_z(y)$
such that
\begin{equation} \label{FZZ2}
-\log_v(\delta(x,y)_{\infty}) \ = \ f_y(x)
\ = \ j_z(x,y)-j_z(x,\inftybar) - C_z(y) \ .
\end{equation}

Fixing $x$, and letting $y$ vary, we see from (\ref{FZZ2}) that the function
$h(y) = C_z(y)$ belongs to $\CPA(\Gamma)$.  Let $\Delta_y$ be the Laplacian
with respect to the variable $y$.  Applying $\Delta_y$ to both sides
of (\ref{FZZ2}) gives
\begin{equation*}
\delta_x(y) - \delta_{\inftybar}(y)
\ = \ (\delta_x(y) - \delta_z(y)) - 0 - \Delta_y(h(y)) \ ,
\end{equation*}
so $\Delta_y(h(y)) = -(\delta_{\inftybar}(y) - \delta_z(y))$.
By Proposition \ref{PropC3} (B),
$h(y) = -j_{{\zeta_0}}(y,\inftybar) + C_z$ for some constant $C_z$.
Hence
\begin{equation*}
-\log_v(\delta(x,y)_{\infty}) \ = \ f_y(x)
\ = \ j_z(x,y)- j_z(x,\inftybar)
                  - j_z(y,\inftybar) + C_z
\end{equation*}
Taking $x = y = z$,
we see that $C_z = -\log_v(\delta(z,z)_{\infty}) = -\log_v(\radius(z))$.
\end{proof}

\vskip .1 in
We will now extend $j_z(x,y)$ to all of $\PP^1_{\Berk}$.

Fix $z \in \PP^1_{\Berk} \backslash \PP^1(\CC_v)$, 
and let $\Gamma$ vary over subgraphs of $\PP^1_{\Berk}$ containing $z$; 
temporarily write $j_z(x,y)_{\Gamma}$ for the potential kernel
on $\Gamma$.  Suppose $\Gamma_1 \subset \Gamma_2$.
Since $\Gamma_1$ and $\Gamma_2$ are both trees,
their potential kernels are compatible:  
by the description of the potential kernel on a tree,
given $x, y \in \Gamma_1$,
\begin{equation*}
j_z(x,y)_{\Gamma_1} \ = \ \rho(z,w) \ = \ j_z(x,y)_{\Gamma_2} 
\end{equation*}
where $w = w_z(x,y)$ is the point where the paths from $x$ and $y$ to 
$z$ meet. Thus the functions $j_z(x,y)_{\Gamma}$ cohere to give
a well-defined function $j_z(x,y)$ on $\PP^1_{\Berk} \backslash \PP^1(\CC_v)$.   

We can extend $j_z(x,y)$ to all $x, y \in \PP^1_{\Berk}$ by setting
\begin{equation*}
j_z(x,y) \ = \
\left\{ \begin{array}{ll} \rho(z,w_z(x,y)) & \text{if $x \ne y$, or } \\
                          \infty & \text{if $x = y$}
        \end{array} \right.
\end{equation*}
for $x,y \in \PP^1(\CC_v)$, where as before  $w = w_z(x,y)$ 
is the point where the paths from $x$ and $y$ to $z$ meet. 
If $x \ne y$, and if $\Gamma$ is any subgraph containing $z$ and $w$, then
\begin{equation}  \label{FZW1}
j_z(x,y) \ = \ j_z(r_{\Gamma}(x),r_{\Gamma}(y))_{\Gamma} \ .
\end{equation}
By Proposition \ref{PropC5} and the continuity of $\delta(x,y)_{\infty}$
off the diagonal, for all $x, y \in  \AA^1_{\Berk}$
\begin{equation} \label{FPQ2A}
-\log_v(\delta(x,y)_{\infty})
\ = \ j_z(x,y) - j_z(x,\infty) - j_z(y,\infty) - \radius(z) \ .
\end{equation}
In particular, if $z = {\zeta_0}$ is the Gauss point, then
\begin{equation} \label{FPQ2}
-\log_v(\delta(x,y)_{\infty})
\ = \ j_{{\zeta_0}}(x,y) - j_{{\zeta_0}}(x,\infty) - j_{{\zeta_0}}(y,\infty) \ .
\end{equation}
Similarly, by Proposition \ref{PropC4} (C),
for each $z, \zeta \in \PP^1_{\Berk} \backslash \PP^1(\CC_v)$ 
\begin{equation} \label{FPQ3}
j_{\zeta}(x,y) \ = \ 
j_z(x,y) - j_z(x,{\zeta}) -j_z(\zeta,y) + j_z(\zeta,\zeta) \ .
\end{equation}

\begin{corollary} \label{CorC6}
Let $0 \ne f \in \CC_v(T)$ be a rational function 
with divisor $\div(f) = \sum_{i=1}^m n_i (a_i)$.

$A)$ Fix $z \in \PP^1_{\Berk} \backslash \PP^1(\CC_v)$.
Then for all $x \in \PP^1_{\Berk}$,
\begin{equation} \label{FMR1}
-\log_v([f]_x)
 \ = \ -\log_v([f]_z) + \sum_{i=1}^m n_i \,j_z(x,a_i) \ .
\end{equation}

$B)$ Let $\Gamma$ be a subgraph of $\PP^1_{\Berk}$, take $z \in \Gamma$,
and put $\tilde{a}_i = r_{\Gamma}(a_i)$.
Then for all \ $x \in \Gamma$,
\begin{equation} \label{FMR2}
-\log_v([f]_x)
 \ = \ -\log_v([f]_z) + \sum_{i=1}^m n_i \, j_z(x,\tilde{a}_i) \ .
\end{equation}
\end{corollary}

\begin{proof}
(A)  There is a $B \in \CC_v$ such that
$f(T) = B \cdot \prod_{a_i \ne \infty} (T-a_i)^{n_i}$.
By Corollary \ref{CorC2}, for each $x \in \Gamma$
\begin{equation*}
-\log_v([f]_x)
\ = \ -\log_v(|B|_v) + \sum_{a_i \ne \infty} -\log_v(\delta(x,a_i)_{\infty}) \ .
\end{equation*}
Inserting formula (\ref{FPQ2A}) and using $\sum_{i=1}^m n_i = 0$
gives
\begin{equation*}
-\log_v([f]_x) \ = \ -\sum_{i=1} n_i \log_v(j_z(x,a_i)) \ + \ C
\end{equation*}
for some constant $C$.  Taking $x = z$ gives $C = -\log_v([f]_z)$.

Part (B) follows from (\ref{FMR1}), using (\ref{FZW1}).
\end{proof}

\vskip .1 in
For future reference, we note one situation where (\ref{FZW1}) holds
automatically:

\begin{proposition} \label{PropC7} {\rm (Retraction Formula)}
Let $\Gamma$ be a subgraph of $\PP^1_{\Berk}$, and suppose $x, z \in \Gamma$.
Then for any $y \in \PP^1_{\Berk}$
\begin{equation*}
j_z(x,y) \ = \ j_z(x,r_{\Gamma}(y)) \ .
\end{equation*}
\end{proposition}

\begin{proof}  Since $x$ and $z$ belong to $\Gamma$, so does the path from 
$x$ to $z$.  Hence the point $w$  where the paths from $x$ and $y$ to $z$ meet,
which lies on the path from $x$ to $z$, belongs to $\Gamma$.
\end{proof} 

\subsection{The spherical meet.}

Write $\hat{\cO}_v$ for the ring of integers of $\CC_v$.
The $v$-adic spherical metric $\|x,y\|_v$ is
the unique $GL_2(\hat{\cO}_v)$-invariant metric on $\PP^1(\CC_v)$ such that
$\|x,y\|_v = |x-y|_v$ for $x, y \in B(0,1)$.    
In homogeneous coordinates,
if $x = (x_0:x_1)$ and $y = (y_0:y_1)$, then
\begin{equation} \label{FBN1}
\|x,y\|_v \ = \ \frac{|x_0 y_1 - x_1 y_0|_v}
                    {\max(|x_0|_v, |x_1|_v) \max(|y_0|_v, |y_1|_v)} \ .
\end{equation}                    
In affine coordinates, for $x, y \in \AA^1(\CC_v)$ this becomes 
\begin{equation} \label{FBN2}
\|x,y\|_v \ = \ \frac{|x-y|_v}{\max(1,|x|_v) \max(1,|y|_v)} \ .
\end{equation}
Thus if we identify $\PP^1(\CC_v)$ with $\AA^1(\CC_v) \cup \{\infty\}$,
and put $1/\infty = 0$,  
\begin{equation} \label{FCM1}
\|x,y\|_v \ = \ \left\{ \begin{array}{ll}
         |x-y|_v  &  \text{if $x, y \in B(0,1)$,} \\
         |1/x - 1/y| & \text{if $x, y \in \PP^1(\CC_v) \backslash B(0,1)$,} \\
         1  &  \text{if exactly one of $x, y \in B(0,1)$.}
                        \end{array} \right.
\end{equation}
It is well known, and easy to check using (\ref{FCM1}), 
that the $v$-adic spherical metric satisfies the strict ultrametric inequality.    

\vskip .1 in
We will now extend $\|x,y\|_v$ to $\PP^1_{\Berk}$.

Let $q_v$ be the base of the logarithm $\log_v(t)$,
and let ${\zeta_0} \in \PP^1_{\Berk}$ be the Gauss point,
as before.  For $x, y \in \PP^1(\CC_v)$, we claim that 
\begin{equation*}
\|x,y\|_v \ = \ q_v^{-j_{{\zeta_0}}(x,y)} \ .
\end{equation*}  
To see this, consider the various cases in (\ref{FCM1}).
If $x, y \in B(0,1)$, the point $w$ where the paths from $x$ and $y$
to ${\zeta_0}$ meet is the ball $B(x,r) = B(y,r)$ with $r = |x-y|_v$.  The
path distance $\rho({\zeta_0},w)$ is $\log_v(1/r)$, 
so $j_{{\zeta_0}}(x,y) = -\log_v(r)$
and $q_v^{-j_{{\zeta_0}}(x,y)} = |x-y|_v = \|x,y\|_v$.
If $x, y \in \PP^1(\CC_v) \backslash B(0,1)$, a similar argument applies,
using the local parameter $1/T$ at $\infty$.  Finally, if one of $x, y$
belongs to $B(0,1)$ and the other to $\PP^1(\CC_v) \backslash B(0,1)$, then
the paths from $x$ and $y$ to ${\zeta_0}$ meet at $w = {\zeta_0}$, 
so $j_{{\zeta_0}}(x,y) = 0$
and $q_v^{-j_{{\zeta_0}}(x,y)} = 1 = \|x,y\|_v$.

This motivates us to extend $\|x,y\|_v$ to $\PP^1_{\Berk}$ by putting
\begin{equation} \label{FRT1}
\|x,y\|_v \ = \ q_v^{-j_{{\zeta_0}}(x,y)} \ .
\end{equation}
for all $x, y \in \PP^1_{\Berk}$.  We will call this extended function the
{\it spherical meet}.  For  $x \in \PP^1_{\Berk}$,
write \ $\diam(x) = \|x,x\|_v 
= q_v^{-j_{{\zeta_0}}(x,x)} = q_v^{-\rho({\zeta_0},x)}$.

The spherical meet
has the following geometric interpretation, which motivates its name.
Consider the paths from $x$ and $y$ to ${\zeta_0}$, and let $w$ be the 
first point where they meet.  Then by (\ref{FRT1})
\begin{equation*}
\|x,y\|_v \ = \ \|w,w\|_v \ = \ \diam(w) \ = \ q_v^{-\rho({\zeta_0},w)} \ .
\end{equation*} 
Note that although $\|x,y\|_v$ is a metric on $\PP^1(\CC_v)$,
the spherical meet $\|x,y\|_v$ is not a metric on $\PP^1_{\Berk}$
because $\|x,x\|_v > 0$ if $x$ is not of type I.
Nonetheless, for each $a \in \PP^1(\CC_v)$, the balls 
\begin{eqnarray*}
\cB(a,r)_0^- & = & \{ x \in \PP^1_{\Berk} : \|x,a\|_v < r \} \ , \\
\cB(a,r)_0 & = & \{ x \in \PP^1_{\Berk} : \|x,a\|_v \le r \} \ .
\end{eqnarray*}
are indeed open (resp. closed) in the Berkovich topology.  
If $r > 1$ then $\cB(a,r)_0^- = \PP^1(\CC_v)$.
If $r \le 1$ and $a \in B(0,1)$ then $\cB(a,r)_0^- = \cB(a,r)^-$, 
while if $a \notin B(0,1)$ then 
$\cB(a,r)_0^- = \{x \in \PP^1_{\Berk} : [1/T-1/a]_x < r \}$.  Similar
formulas hold for the closed balls.     

\begin{proposition} \label{PropC8}
\ \ 

$A)$ The spherical meet $\|x,y\|_v$ is nonnegative, symmetric, continuous
in each variable separately, and bounded above by $1$. 
As a function of two variables,
it is strongly upper semicontinuous. 
It is continuous off the diagonal, and continuous
at $(x_0,x_0)$ for each point $x_0$ of type I, 
but is discontinuous at $(x_0,x_0)$ for each point of type II, III, or IV. 

$B)$ For all $x, y \in \PP^1_{\Berk}$
\begin{equation} \label{FTD1B}
\|x,y\|_v \ = \
\limsup \begin{Sb} (a,b) \rightarrow  (x,y) \\ 
        a, b \in \PP^1(\CC_v) \end{Sb} \|a,b\|_v \ .
\end{equation}

$C)$  For all $x, y, z \in \PP^1_{\Berk}$
\begin{equation*}
\|x,y\|_v \ \le \ \max( \|x,z\|_v, \|y,z\|_v )
\end{equation*}  
with equality if $\|x,z\|_v \ne \|y,z\|_v$.

$D)$  For each $a \in \PP^1_{\Berk}$ and $r \in \RR$, the `open ball'
$\cB(a,r)_0^- := \{ x \in \PP^1_{\Berk} : \|x,a\|_v < r \}$
is open.  It is empty if $r \le \diam(a)$, and coincides with an open
ball $\cB(b,r)^-$ for some $b \in \PP^1(\CC_v)$ if $r > \diam(a)$.

Likewise, the `closed ball'
$\cB(a,r)_0 := \{ x \in \PP^1_{\Berk} : \|x,a\|_v \le r \}$
is closed.  It is empty if $r < \|a,a\|_v$, and coincides with
$\cB(b,r)_0$ for some $b \in \PP^1(\CC_v)$ if $r > \diam(a)$ or if
$r = \diam(a)$ and $a$ is of type II or III.  If $r = \diam(a)$
and $a$ is of type I or IV, then $\cB(a,r)_0 = \{a\}$.
\end{proposition}

\begin{proof}  Similar to Proposition \ref{PropC1}.
\end{proof} 

\vskip .1 in
The Hsia kernel and the spherical meet can be obtained from each other.

\begin{proposition} \label{PropC9} \ \

$A)$  For \ $x, y \in \AA^1_{\Berk}$, \quad
$\displaystyle{\delta(x,y)_{\infty} 
\ = \ \frac{\|x,y\|_v}{\|x,\infty\|_v \|y,\infty\|_v}}$ \ .

$B)$ For \ $x, y \in \PP^1_{\Berk}$,
\begin{equation*} 
\|x,y\|_v \ =  \
\frac{\delta(x,y)_{\infty}}{\max(1,\radius(x)) \, \max(1,\radius(y))} \qquad
    \text{if $x, y \ne \infty$,} 
\end{equation*}
with $\|x,\infty\|_v = 1/\max(1,\radius(x))$ and  
$\|\infty,y\|_v = 1/\max(1,\radius(y))$. 
\end{proposition}

\begin{proof} 
Part (A) follows from (\ref{FPQ2}).  Part (B) follows from (\ref{FBN2})
using Proposition \ref{PropC1} (B) and Proposition \ref{PropC8} (B).  
\end{proof} 

\subsection{The generalized Hsia kernel.} 

Proposition \ref{PropC9} motivates us to define the Hsia kernel 
for an arbitrary point $\zeta \in \PP^1_{\Berk}$.

For $x,y \in \PP^1_{\Berk} \backslash \{\zeta\}$, we define the generalized
Hsia kernel by 
\begin{equation} \label{FCJ1}
\delta(x,y)_{\zeta} \ = \ \frac{\|x,y\|_v}{\|x,\zeta\|_v \,\|y,\zeta\|_v} \ .
\end{equation}
If $\zeta \notin \PP^1(\CC_v)$,
this makes sense for all $x, y \in \PP^1_{\Berk}$,
since $\|x,\zeta\|_v, \|\zeta,y\|_v \ge \diam(\zeta) > 0$, giving 
$\delta(x,\zeta)_{\zeta} = \delta(\zeta,y)_{\zeta} = 1/\diam(\zeta)$.  
If $\zeta \in \PP^1(\CC_v)$,
we put $\delta(x,\zeta)_{\zeta} = \delta(\zeta,y)_{\zeta} = \infty$.
In this way, we can regard $\delta(x,y)_{\zeta}$ as being
defined (as an extended real) for all $x, y \in \PP^1_{\Berk}$.

Actually, it is better to regard the Hsia kernel as only defined up to scaling.
By the definition of the spherical meet, 
\begin{equation*}
\delta(x,y)_{\zeta}
\ = \ q_v^{-j_{{\zeta_0}}(x,y) + j_{{\zeta_0}}(x,\zeta) 
+ j_{{\zeta_0}}(y,\zeta)} \ .
\end{equation*} 
However, for any $z \in \PP^1_{\Berk} \backslash \PP^1(\CC_v)$, formula
(\ref{FPQ3}) shows that there is a constant $C_{\zeta}$ such that
\begin{equation} \label{FXX5}
C_{\zeta} \cdot \delta(x,y)_{\zeta}  \ = \
q_v^{-j_z(x,y) + j_z(x,\zeta) + j_z(y,\zeta)}.
\end{equation}
For each $\zeta$, and each $C_{\zeta} > 0$,
we will also call $C_{\zeta} \cdot \delta(x,y)_{\zeta}$
a Hsia kernel.

\vskip .1 in
When $\zeta = \infty$, Proposition \ref{PropC9} shows that (\ref{FCJ1}) is
consistent with our earlier definition of $\delta(x,y)_{\infty}$.
When $\zeta = {\zeta_0}$ is the Gauss point,
then $\|x,{\zeta_0}\|_v = \|{\zeta_0},y\|_v = 1$, so
\begin{equation*}
\delta(x,y)_{{\zeta_0}} \ = \ \|x,y\|_v \ .
\end{equation*}
For an arbitrary $\zeta \in \PP^1_{\Berk} \backslash \PP^1(\CC_v)$,
Proposition \ref{PropC4} (C) shows that  
\begin{equation*}
\delta(x,y)_{\zeta} \ = \ C_{\zeta} \cdot q_v^{-j_{\zeta}(x,y)}  
\end{equation*}
where $C_{\zeta} = q_v^{j_{{\zeta_0}}(\zeta,\zeta)}$.  
Thus for $\zeta \notin \PP^1(\CC_v)$, \ 
$\delta(x,y)_{\zeta}$ is a generalized spherical meet.

For $\zeta \in \PP^1(\CC_v)$, the reader familiar with (\cite{R1})
will recognize (\ref{FCJ1}) as the `canonical distance' $[x,y]_{\zeta}$
for $x,y \in \PP^1(\CC_v)$.  Thus, the Hsia kernel is the natural
extension of the canonical distance to the Berkovich line. 

\vskip .1 in
For each $x \in \PP^1_{\Berk}$, put \ 
$\diam_{\zeta}(x) = \delta(x,x)_{\zeta}$.
The generalized Hsia kernel has the usual geometric interpretation:  
given $x, y \in \PP^1_{\Berk}$, let $w$ be the point where the paths
from $x$ and $y$ to $\zeta$ meet.  Then 
\begin{equation*}
\delta(x,y)_{\zeta} \ = \ \delta(w,w)_{\zeta} \ = \ \diam_{\zeta}(w) \ .
\end{equation*}
  
\begin{proposition} \label{PropC10}
\ \ 

$A)$  For each $\zeta$, the generalized Hsia kernel is nonnegative, 
symmetric and continuous in each variable separately.
If $\zeta \in \PP^1_{\Berk} \backslash \PP^1(\CC_v)$ it is bounded.
If $\zeta \in \PP^1(\CC_v)$ it is unbounded, and extends the canonical
distance $[x,y]_{\zeta}$.
  
As a function of two variables, it is strongly upper semicontinuous.
It is continuous off the diagonal, and is continuous
at $(x_0,x_0)$ for each point $x_0$ of type I, 
but is discontinuous at $(x_0,x_0)$ for each point of type II, III, or IV.

$B)$ For each $x, y \in \PP^1_{\Berk}$
\begin{equation} \label{FTD0B}
\delta(x,y)_{\zeta} \ = \
\limsup \begin{Sb} (a,b) \rightarrow  (x,y) \\ 
        a, b \in \PP^1(\CC_v) \end{Sb} \delta(a,b)_{\zeta} \ .
\end{equation}

$C)$  For all $x, y, z \in \PP^1_{\Berk}$,  
\begin{equation*}
\delta(x,y)_{\zeta} \ \le \ \max( \delta(x,z)_{\zeta}, \delta(y,z)_{\zeta} ) \ , 
\end{equation*}
with equality if \ $\delta(x,z)_{\zeta} \ne \delta(y,z)_{\zeta}$.

$D)$ For each $a \in \PP^1_{\Berk}$ and $r \in \RR$, the `open ball'
$\cB(a,r)_{\zeta}^- := \{ x \in \PP^1_{\Berk} : \delta(x,a)_{\zeta} < r \}$
is open.  It is empty if $r \le \diam_{\zeta}(a)$, and coincides with an open
ball $\cB(b,r)_{zeta}^-$
for some $b \in \PP^1(\CC_v)$ if $r > \diam_{\zeta}(a)$.

Likewise, the `closed ball'
$\cB(a,r)_{\zeta} := \{ x \in \PP^1_{\Berk} : \delta(x,a)_{\zeta} \le r \}$
is closed.  It is empty if $r < \diam_{\zeta}(a)$, and coincides with
$\cB(b,r)_{\zeta}$ for some $b \in \PP^1(\CC_v)$
if $r > \diam_{\zeta}(a)$ or if
$r = \diam_{\zeta}(a)$ and $a$ is of type II or III.
If $r = \diam_{\zeta}(a)$ and $a$ is of type I or IV,
then $\cB(a,r)_{\zeta} = \{a\}$.
\end{proposition}

\begin{proof}
Parts (A), (B), and (D) follow by arguments similar to those in the proof
of Proposition \ref{PropC1}.  
Part (C) follows by consideration of the geometric 
interpretation, or can be shown using Proposition \ref{PropC8} (C)
by a case by case analysis similar to the proof of
(\cite{R1}, Theorem 2.5.1, p.125).
\end{proof}

\vskip .1 in
If $r \le \diam_{\zeta}(a)$ the balls $\cB(a,r)_{\zeta}$
and $\cB(a,r)_{\zeta}^-$ have been described in Proposition \ref{PropC10} (D).
If $r > \diam_{\zeta}(\zeta)$ then
$\cB(a,r)_{\zeta} = \cB(a,r)_{\zeta}^- = \PP^1_{\Berk}$.
Suppose $\diam_{\zeta}(a) \le r \le \diam_{\zeta}(\zeta)$.
Then the balls have the following geometric interpretation.  Consider
the function $\diam_{\zeta}(x)$ for $x$ in the path from $a$ to $\zeta$.
By Proposition \ref{PropC10} (A), it is continuous on that path,
since $\diam_{\zeta}(x) = \delta(x,x)_{\zeta} = \delta(x,a)_{\zeta}$.
By the geometric interpretation it is monotone increasing.  
Hence there is a unique $x$ in the path with $\diam_{\zeta}(x) = r$.
The closed ball $\cB(a,r)_{\zeta}$ is the set of all $z \in \PP^1_{\Berk}$
such that the path from $z$ to $\zeta$ passes through $x$, and the
open ball $\cB(a,r)_{\zeta}^-$ is the connected component of
$\cB(a,r)_{\zeta} \backslash \{x\}$ which contains $a$.

\vskip .1 in
Finally, we note that the generalized Hsia kernel can be used to decompose 
absolute values of rational functions on $\PP^1_{\Berk}$, 
just as the canonical distance does on $\PP^1(\CC_v)$:  

\begin{corollary} \label{CorC11}
Let $0 \ne f \in \CC_v(\PP^1)$ have divisor 
$\div(f) = \sum_{i=1}^m n_i (a_i)$.  
For each $\zeta \in \PP^1_{\Berk}$, there is a constant $C = C(f,\zeta)$
such that for all $x \in \PP^1_{\Berk}$
\begin{equation*}
[f]_x \ = \ C \cdot \prod_{a_i \ne \zeta} \delta(x,a_i)_{\zeta}^{n_i} \ .
\end{equation*}
\end{corollary}

\begin{proof}  Similar to the proof of Corollary \ref{CorC6}.
\end{proof}

\section{Capacities.}
\label{Section D}

Fix $\zeta \in \PP^1_{\Berk}$, and let $E \subset \PP^1_{\Berk}$ be a
set not containing $\zeta$.  In this section we will develop 
the theory of the {\it logarithmic capacity} of $E$ with
respect to $\zeta$ (or more correctly, with respect to a choice of the
Hsia kernel $\delta(x,y)_{\zeta}$).  The exposition below is 
adapted from (\cite{R1}, Chapter 4.1).  

\subsection{Logarithmic capacities.}

Recall that a probability measure is a non-negative Borel measure
of total mass $1$.  Given a probability measure $\nu$ with support
contained in $E$, define the energy integral
\begin{equation*}
I_{\zeta}(\nu) \ = \ \iint_{E \times E}
              -\log_v(\delta(x,y)_{\zeta}) \, d\nu(x) d\nu(y) \ .
\end{equation*}
Here, the integral is a Lebesgue integral.
The kernel $-\log_v(\delta(x,y)_{\zeta})$ is lower semicontinuous, and hence 
Borel measurable, since $\delta(x,y)_{\zeta}$ is upper semicontinuous.

Let $\nu$ vary over probability measures with support contained in $E$,
and define the {\it Robin constant}
\begin{equation*}
V_{\zeta}(E) \ = \ \inf_{\nu} \ I_{\zeta}(\nu) \ .
\end{equation*}
Define the {\it logarithmic capacity}
\begin{equation*}
\gamma_{\zeta}(E) \ = \ q_v^{-V_{\zeta}(E)} \ .
\end{equation*}

By its definition, the logarithmic capacity is monotonic in $E$:
if $E_1 \subset E_2$, then $\gamma_{\zeta}(E_1) \le \gamma_{\zeta}(E_2)$.
It also follows from the definition that for each $E$
\begin{equation}  \label{FVM1}
\gamma_{\zeta}(E) \ = \
\sup \begin{Sb} K \subset E \\ \text{$K$ compact} \end{Sb}
       \gamma_{\zeta}(K)
\end{equation}
Indeed, the support of each probability measure $\nu$ is compact,
so tautologically for each $\nu$ supported on $E$
there is a compact set $K \subset E$ with
$I_{\zeta}(\nu) \le \gamma_{\zeta}(K) \le \gamma_{\zeta}(E)$.

An important distinction is between sets of capacity $0$ and sets of positive
capacity.  If $E$ contains a point $a$ of types II, III, or IV,
then $\gamma_{\zeta}(E) > 0$, since the point mass $\delta_a(x)$ satisfies
$I_{\zeta}(\delta_a) = -\log(\diam_{\zeta}(a)) < \infty$.
Thus, every set of capacity $0$ is contained in $\PP^1(\CC_v)$. 
The converse is not true;  there are many sets in $\PP^1(\CC_v)$
with positive capacity.

The property that a set has capacity $0$ or positive capacity
is independent of the point $\zeta$.

\begin{proposition} \label{PropD1}
Suppose $E \subset \PP^1_{\Berk}$.  Then $\gamma_{\zeta}(E) = 0$ for some 
$\zeta \notin E$ if and only if $\gamma_{\xi}(E) = 0$ for each
$\xi \in \PP^1_{\Berk} \backslash E$.
\end{proposition}

\begin{proof} By (\ref{FVM1}) it suffices to consider the case where $E$
is compact.  Take $\xi \notin E$.  Since $\|x,\xi\|_v$ is continuous on $E$,
there is a constant $K_{\xi} > 0$ such that
\begin{equation*}
1/K_{\xi} \ \le \ \|x,\xi\|_v \ \le \ K_{\xi}
\end{equation*}
for all $x \in E$.  Since
\begin{equation*}
\delta(x,y)_{\xi} = \frac{\|x,y\|_v}{\|x,\xi\|_v \|y,\xi\|_v} \ , \quad
\delta(x,y)_{\zeta} = \frac{\|x,y\|_v}{\|x,\zeta\|_v \|y,\zeta\|_v} \ .
\end{equation*}
it follows that
\begin{equation*}
\frac{1}{(K_{\zeta} K_{\xi})^2} \delta(x,y)_{\zeta}\ \le \
\delta(x,y)_{\xi} \ \le \ (K_{\zeta} K_{\xi})^2 \delta(x,y)_{\zeta} \ .
\end{equation*}
Hence for each probability measure $\nu$ on $E$,
either $I_{\zeta}(\nu) = I_{\xi}(\nu) = \infty$, or $I_{\zeta}(\nu)$ and
$I_{\xi}(\nu)$ are both finite.
\end{proof}

\vskip .1 in
Here are some examples of capacities.

\begin{example}
If $E \subset \PP^1(\CC_v)$ is a countable set,
and $\zeta \notin E$, then $\gamma_{\zeta}(E) = 0$.  
To see this, not first that if 
$\nu$ is a probability measure supported on $E$, then necessarily $\nu$ 
has point masses
(if $\nu(\{x\}) = 0$ for each $x \in E$, then by countable additivity
$\nu(E) = 0$, which contradicts $\nu(E) = 1$).  If $p \in E$ is a point 
with $\nu(\{p\}) > 0$, then
\begin{eqnarray*}
I_{\zeta}(\nu) 
& = & \iint_{E \times E} -\log_v(\delta(x,y)_{\zeta}) \, d\nu(x) d\nu(y) \\
& \ge & -\log_v(\delta(p,p)_{\zeta}) \cdot \nu(\{p\})^2 \ = \infty 
\end{eqnarray*}
so $V_{\zeta}(E) = \infty$.
\end{example}

\begin{example}
If $E = \{a\}$ is a point not of type I, 
and $\zeta \ne a$, then 
\begin{equation*}
\gamma_{\zeta}(E) \ = \ \diam_{\zeta}(a) \ .
\end{equation*} 

The only probability measure supported on $E$ is the point mass 
$\nu = \delta_{a}(x)$, for which  
\begin{eqnarray*}
I_{\zeta}(\nu) 
& = & \iint_{E \times E} -\log_v(\delta(x,y)_{\zeta}) \, d\nu(x) d\nu(y) \\
& = & -\log_v(\delta(a,a)_{\zeta}) \ = -\log_v(\diam_{\zeta}(a)) \ . 
\end{eqnarray*}     
Hence $V_{\zeta}(E) = -\log_v(\diam_{\zeta}(a))$ 
and $\gamma_{\zeta}(a) = \diam_{\zeta}(a)$.
\end{example}

\begin{example}
If $\zeta = \infty$ and 
$E = \ZZ_p \subset \AA^1(\CC_v)$, then 
$\gamma_{\infty}(E) = p^{-1/(p-1)}$.  More generally, if $E = \cO_v$ is 
the ring of integers of a finite extension $K_v/Q_p$ with ramification
index $e$ and residue degree $f$, 
then $\gamma_{\infty}(E) = (p^f)^{-1/e(p^f-1)}$.

For $x, y \in \AA^1(\CC_v)$, $\delta(x,y)_{\infty} = |x-y|_v$, 
so the capacity is given by the same computation as in the classical 
case;  see (\cite{R1}, Example 5.2.13, p. 347).  
\end{example}

\subsection{The equilibrium distribution.}

If $E$ is compact, and if $\gamma_{\zeta}(E) > 0$,
there is a probability measure $\mu = \mu_{\zeta}$ on $E$ for which
$I_{\zeta}(\mu) = V_{\zeta}(E)$.  In Section~\ref{Section E} 
of these notes, we will show
it is unique;  it will be called the {\it equilibrium distribution},
or the {\it equilibrium measure}, of $E$ with respect to $\zeta$.  
To prove its existence, we need some preliminary lemmas.

\begin{lemma} \label{LemD2}
\text{\rm (Baire)}  Let $X$ be a metric space, and let $A \subset X$.
Let $f : A \rightarrow \RR \cup \{\infty\}$ be a lower semicontinuous
function which is bounded below by a number $m \in \RR$.
Then there is a sequence of continuous functions
$f_k : X \rightarrow \RR$ such that

$a)$  For each $x \in X$, \quad $m \le f_1(x) \le f_2(x) \le \cdots$

$b)$  For each $a \in A$, \quad $\lim_{k \rightarrow \infty} f_k(a) = f(a)$.
\end{lemma}

\begin{proof}  We adapt the proof of (\cite{Ts}, Theorem II.5, p. 36).
If $f(a) \equiv \infty$ on $A$, we can take $f_k(x) \equiv m+k$.
Hence we can assume without loss that $f(a) \not\equiv \infty$.
Let $d(x,y)$ be a metric for $X$.

Define $f_k(x)$ on $X$ by
\begin{equation} \label{FGL1}
f_k(x) \ = \ \inf_{a \in A} (f(a) + k \cdot d(a,x))
\end{equation}
Clearly $m \le f_k(x) < \infty$ for all $x$, and $f_k(x) \le f_{k+1}(x)$.

We will now show each $f_k$ is continuous.  Take $x_1 \ne x_2 \in X$, and let
$d(x_1,x_2) = \delta$.  By (\ref{FGL1}) there is an $a_0 \in A$
such that
\begin{equation} \label{FGL2}
f_k(x_1) > f(a_0) + k \cdot d(a_0,x_1) - \delta \ .
\end{equation}
Again by (\ref{FGL1})
\begin{equation} \label{FGL3}
f_k(x_2) \le f(a_0) + k \cdot d(a_0,x_2) \ .
\end{equation}
Since $d(a_0,x_2) \le d(a_0,x_1) + d(x_1,x_2) = d(a_0,x_1) + \delta$,
\begin{equation*}
f_k(x_2) \ \le \ f_k(a_0) + k \cdot d(a_0,x_1) + k \delta
         \ < \ f_k(x_1) + (k+1) \delta \ .
\end{equation*}
Similarly $f_k(x_1) \le f_k(x_2) + (k+1) \delta$, so
\begin{equation*}
|f_k(x_1) - f_k(x_2)| \ < \ (k+1) d(x_1,x_2) \ .
\end{equation*}
Thus, $f_k(x)$ is continuous on $X$.

Next we will show that $\lim_{k \rightarrow \infty} f_k(a_0) = f(a_0)$
for each $a_0 \in A$.  By (\ref{FGL1})
\begin{equation*}
f_k(a_0) \ = \ \inf_{a \in A} (f(a) + k \cdot d(a,a_0)) \
         \ \le \ f(a_0) + k d(a_0,a_0) \ = \ f(a_0) \ .
\end{equation*}
First suppose $f(a_0) < \infty$.  By lower semicontinuity, for each
$\varepsilon > 0$ there is a $\rho > 0$ such that
\begin{equation} \label{FGL4}
f(a) \ > \ f(a_0) - \varepsilon,
     \qquad \text{if $a \in A$, $d(a,a_0) < \rho$ \ .}
\end{equation}
On the other hand, if $a \in A$ and $d(a,a_0) \ge \rho$ then
\begin{equation} \label{FGL5}
f(a) + k \cdot d(a,a_0) \ \ge \ m + k \rho
\end{equation}
>From (\ref{FGL4}), (\ref{FGL5}) it follows that for sufficiently large $k$,
$f(a) + k d(a,a_0) \ge f(a_0) - \varepsilon$ for all $a \in A$, and hence
$f_k(a_0) \ge f(a_0) - \varepsilon$.

Next suppose $f(a_0) = \infty$.  Take $M$ arbitrarily large.
By lower semicontinuity, there is a $\rho > 0$ such that
\begin{equation} \label{FGL6}
f(a) \ > \ M \qquad \text{if $a \in A$, $d(a,a_0) < \rho$ \ .}
\end{equation}
If $a \in A$ and $d(a,a_0) \ge \rho$ then (\ref{FGL5}) holds, so for
sufficiently large $k$, $f(a) + k d(a,a_0) > M$ for all $a \in A$,
and hence $f_k(a_0) \ge M$.
\end{proof}

\vskip .1 in

\begin{lemma} \label{LemD3}
Let $E$ be a compact Hausdorff space, and suppose $\nu_1, \nu_2, \ldots$
are probability measures on $E$ which converge weakly to a measure $\mu$.
Then $\nu_1 \times \nu_1, \nu_2 \times \nu_2, \ldots$ converges weakly to
$\mu \times \mu$ on $E \times E$.
\end{lemma}

\begin{proof}  We fill in the details of the argument given in 
(\cite{Hille2}, p.283).  Since $E$ is a compact Hausdorff space,
the Stone-Weierstrass theorem asserts that linear combinations of functions
of the form $f(x) g(y)$,  with $f, g \in \cC(E)$,
are dense in $\cC(E \times E)$ under the $\sup$ norm $\| \cdot \|_{E \times E}$.
For such products,
\begin{eqnarray*}
\lim_{n \rightarrow \infty} \iint_{E \times E} f(x) g(y) \, d\nu_n(x) d\nu_n(y)
& = & \int_E f(x) \, d\mu(x) \cdot \int_E g(y) \, d\mu(y) \\
& = & \iint_{E \times E} f(x) g(y) \, d\mu(x) d\mu(y) 
\end{eqnarray*}
so the same holds for their linear combinations.

Now let $F(x,y) \in \cC(E \times E)$ be arbitrary.  
Given $\varepsilon > 0$, take $ f(x,y) = \sum_i c_i f_i(x) g_i(y)$ 
so that $h(x,y) = F(x,y) - f(x,y)$ 
satisfies $\|h\|_{E \times E} < \varepsilon$.
For any probability measure $\nu$ on $E$, clearly 
\begin{equation*}
| \iint_{E \times E} h(x,y) \, d\nu(x) d\nu(y)|
 \ \le \ \|h\|_{E \times E} \ .
\end{equation*}
Let $N$ be large enough that 
\begin{equation*}
|\iint_{E \times E} f(x,y) \, d\nu_n(x) d\nu_n(y) -   
\iint_{E \times E} f(x,y) \, d\mu(x) d\mu(y)| \ < \ \varepsilon
\end{equation*}
for $n \ge N$.  By a three-epsilons argument, for such $n$
\begin{equation*}
| \iint_{E \times E} F(x,y) \, d\nu_n(x) d\nu_n(y)
 - \iint_{E \times E} F(x,y) \, d\mu(x) d\mu(y) |
 \ \le \ 3 \varepsilon \ .
\end{equation*}
Since $\varepsilon$ is arbitrary,
\begin{equation*}
\lim_{n \rightarrow \infty} \iint_{E \times E} F(x,y) \, d\nu_n(x) d\nu_n(y)
 \ = \ \iint_{E \times E} F(x,y) \, d\mu(x) d\mu(y) \ .
\end{equation*}
Thus, the sequence $\{\nu_n \times \nu_n\}$ converges weakly to $\mu \times \mu$.
\end{proof}
\vskip .1 in

We can now show the existence of the equilibrium measure.

\begin{proposition} \label{PropD4}
Let $E \subset \PP^1_{\Berk} \backslash \{\zeta\}$ be a compact set with 
positive capacity.  Then there is a probability measure $\mu$ supported 
on $E$ such that $I_{\zeta}(\mu) = V_{\zeta}(E)$.
\end{proposition}

\begin{proof}
Take a sequence of probability measures $\nu_n$ on $E$ for which
$\lim_{n \rightarrow \infty} I_{\zeta}(\nu) = V_{\zeta}(E)$.  After
passing to a subsequence, if necessary, we can assume that $\{\nu_n\}$
converges weakly to a measure $\mu$ .  Clearly $\mu$ is a probability measure
supported on $E$.

Applying Lemma \ref{LemD2} to $-\log(\delta(x,y)_{\zeta})$,
we obtain a sequence of continuous functions $f_k(x,y)$ on $(\PP^1_{\Berk})^2$
converging monotonically to $-\log(\delta(x,y)_{\zeta})$ on $E \times E$.
By Lemma \ref{LemD3}, for each $k$
\begin{equation*}
\iint_{E \times E} f_k(x,y) \, d\mu(x) d\mu(y) \ = \
\lim_{n \rightarrow \infty}
\iint_{E \times E} f_k(x,y) \, d\nu_n(x) d\nu_n(y) \ .
\end{equation*}
On the other hand $f_k(x,y) \le -\log(\delta(x,y)_{\zeta})$ 
so for each $n$ and $k$
\begin{eqnarray*}
\iint_{E \times E} f_k(x,y) \, d\nu_n(x) d\nu_n(y) & \le &
\iint_{E \times E} -\log(\delta(x,y)_{\zeta}) \, d\nu_n(x) d\nu_n(y) \\
& = & I_{\zeta}(\nu_n) \ ,
\end{eqnarray*}
and for each $k$
\begin{equation*}
\iint_{E \times E} f_k(x,y) \, d\mu(x) d\mu(y) \ \le \
\lim_{n \rightarrow \infty} I_{\zeta}(\nu_n) \ = \ V_{\zeta}(E) \ .
\end{equation*}
Using this and the monotone convergence theorem,
\begin{eqnarray*}
I_{\zeta}(\mu)
& = & \iint_{E \times E} -\log(\delta(x,y)_{\zeta}) \, d\mu(x) d\mu(y) \\
& = & \lim_{k \rightarrow \infty} \iint_{E \times E} f_k(x,y) \, d\mu(x) d\mu(y) \\
& \le & V_{\zeta}(E) \ .
\end{eqnarray*}
The opposite inequality is trivial, so $I_{\zeta}(\mu) = V_{\zeta}(E)$.
\end{proof}

\begin{remark}
In the classical proof over $\CC$, one considers the truncated logarithm
$-\log^{(t)}(|x-y|) = \min(t,-\log(|x-y|)$ which is continuous, 
and which approaches $-\log(|x-y|)$ pointwise as $t \rightarrow \infty$.
Here $-\log_v^{(t)}(\delta(x,y)_{\zeta})$ is lower semicontinuous,
but not continuous, so it was necessary to introduce the functions $f_k(x,y)$.
\end{remark}

\vskip .1 in
Let $U_{\zeta}$ be the connected component of
$\PP^1_{\Berk} \backslash E$ containing $\zeta$.
Write $\partial E_{\zeta} = \partial U_{\zeta}$
for the boundary of $U_{\zeta}$, 
the part of $\partial E$ in common with $\Ubar$.

\begin{proposition} \label{PropD5}
Let $E \subset \PP^1_{\Berk} \backslash \{\zeta\}$ be
compact with positive capacity.  Then the equilibrium distribution 
$\mu_{\zeta}$ is supported on $\partial E_{\zeta}$.
\end{proposition}

\begin{proof}
Suppose $\mu_{\zeta}$ is not supported on $\partial E_{\zeta}$, 
and fix $x_0 \in \supp(\mu_{\zeta}) \backslash \partial E_{\zeta}$.
Let $r_{\zeta} : E \rightarrow \partial E_{\zeta}$ be the
retraction map which takes each $x \in E$ to the last point in $E$
on the path from $x$ to $\zeta$.  Then $r_{\zeta}$ is continuous. 
Put $\mu_0 = (r_{\zeta})_*(\mu_\zeta)$.
We claim that $I_{\zeta}(\mu_0) < I_{\zeta}(\mu_{\zeta})$.

For each $x, y \in E$,
if $\xbar = r_{\zeta}(x)$, $\ybar = r_{\zeta}(y)$,
then $\delta(x,y)_{\zeta} \le \delta(\xbar,\ybar)_{\zeta}$.
This is follows from the geometric interpretation of $\delta(x,y)_{\zeta}$:
if $w$ is the point where the paths from $x$ to $\zeta$
and $y$ to $\zeta$ meet,
and $\wbar$ is the point where the paths from $\xbar$ to $\zeta$
and $\ybar$ to $\zeta$ meet,
then $\wbar$ lies on the path from $w$ to $\zeta$.

Now consider the point $x_0$, and put $x_1 = r_{\zeta}(x_0)$.
Then $\diam_{\zeta}(x_0) < \diam_{\zeta}(x_1)$.
Fix $r$ with $\diam_{\zeta}(x_0) < r < \diam_{\zeta}(x_1)$, 
and put $U = \cB(x_0,r)_{\zeta}^-$.
For each $x \in U$, $r_{\zeta}(x) = x_1$.  If $x, y \in U$,
\begin{equation*}
\delta(x,y)_{\zeta} \ < \ r
\ < \ \diam_{\zeta}(x_1) \ = \ \delta(x_1,x_1)_{\zeta}
\end{equation*}
Since $x_0 \in \supp(\nu)$, necessarily $\mu_{\zeta}(U) > 0$.  Hence
\begin{eqnarray*}
I_{\zeta}(\mu_{\zeta}) & = &
\int -\log_v(\delta(x,y)_{\zeta}) \, d\mu_{\zeta}(x) d\mu_{\zeta}(y) \\
& > & \int -\log_v(\delta(r_{\zeta}(x),r_{\zeta}(y))_{\zeta})
                                 \, d\mu_{\zeta}(x) d\mu_{\zeta}(y) \\
& = & \int -\log_v(\delta(z,w)_{\zeta})
       \, d\mu_0(z) d\mu_0(w) \ = \ I_{\zeta}(\mu_0) \ .
\end{eqnarray*}
This contradicts the minimality of $I_{\zeta}(\mu_{\zeta})$,
so $\mu_{\zeta}$ is supported on $\partial E_{\zeta}$.
\end{proof}

\vskip .1 in
The existence of the equilibrium measure has several consequences.

\begin{corollary} \label{CorD6}
Let $E \subset \PP^1_{\Berk} \backslash \{\zeta\}$ be compact.
Then for any $\varepsilon > 0$, there is a closed neighborhood $W$ of
$E$ such that
$\gamma_{\zeta}(E^{\prime}) \le \gamma_{\zeta}(E) + \varepsilon$
for each $E^{\prime} \subset W$.

In fact, there is such a neighborhood which is a finite union of discs 
$\bigcup_{k=1}^m \cB(a_k,r_k)_{\zeta}$, where each $a_k$ is of type II 
or type III, and where $r_k = \diam_{\zeta}(a_k)$. 
\end{corollary}

\begin{proof}
Take a cofinal sequence of closed neighborhoods $\{W_n\}$ of $E$.
Without loss, we can assume $\zeta \notin W_n$ for each $n$.  Since $W_n$
necessarily contains points not of type I, $\gamma_{\zeta}(W_n) > 0$.
Since $W_n$ is compact, it has an equilibrium measure $\mu_n$;  thus
$I_{\zeta}(\mu_n) = V_{\zeta}(W_n)$.  After passing to a subsequence
if necessary, we can assume the $\mu_n$
(which are all supported on $W_1$) converge weakly to a probability
measure $\mu^*$.  Since $\cap W_n = E$, clearly $\mu^*$ is supported on $E$.
By an argument similar to the one above, 
\begin{equation*}
\lim_{n \rightarrow \infty} V_{\zeta}(W_n) \ = \ I_\zeta(\mu^*)
\ \ge \ V_{\zeta}(E) \ .
\end{equation*}
Since $V_{\zeta}(W_n) \le V_{\zeta}(E)$ for each $n$, equality must hold
throughout.  It follows that
$\lim_{n \rightarrow \infty} \gamma_{\zeta}(W_n) = \gamma_{\zeta}(E)$,
and there is an $n$ with
$\gamma_{\zeta}(W_n) \le \gamma_{\zeta}(E) + \varepsilon$.  For
each $E^{\prime} \subset W_n$, the monotonicity of the capacity shows that
$\gamma_{\zeta}(E^{\prime}) \le \gamma_{\zeta}(W_n)$.

To show that the neigbhorhood can be taken in the special form described,
fix any closed neighborhood $W$ of $E$ with 
$\gamma_{\zeta}(W) \le \gamma_{\zeta}(E) + \varepsilon$.
Let $r_{\zeta} : E \rightarrow \partial W_{\zeta}$ be the retraction map
which takes each $x$ in $E$ to the last point on the path from $x$ to $\zeta$
which belongs to $W$.

Since $E$ is contained in the interior of $W$,
none of the points $a \in r_{\zeta}(E)$ belongs to $E$,
and each lies on the interior of the path from some
point $x$ in $E$ to $\zeta$.
Such a point $a$ is necessarily of type II or III.
If $a = r_{\zeta}(x)$ and $r = \diam_{\zeta}(a)$, then
$\cB(a,r)_{\zeta} = \cB(x,r)_{\zeta}$, and $x \in \cB(x,r)_{\zeta}^-$.
Since $E$ is compact, a finite number of the discs $\cB(x,r)_{\zeta}^-$
cover $E$.  It follows that $r_{\zeta}(E)$ is a finite set of points
$\{a_1, \ldots, a_m\}$, each of which is of type II or III.  
For each $a_k$, put $r_k = \diam_{\zeta}(a_k)$, 
and put $\tilde{W} = \cup_{i=1}^m \cB(a_k,r_k)_{\zeta}$. 
Clearly $\tilde{W}$ is compact, and has $E$ in its interior.
Since $\tilde{W}$ contains points of type II or III, 
$\gamma_{\zeta}(\tilde{W}) > 0$.

We claim that 
$\gamma_{\zeta}(\tilde{W}) \le \gamma_{\zeta}(E) + \varepsilon$.  
To see this, let 
$\tilde{\mu}_)$ be any equilibrium distribution for $\tilde{W}$, and put 
$\tilde{\mu} = r_{\zeta}(\tilde{\mu}_0)$.  By Proposition \ref{PropD5},     
$\tilde{\mu}_0$ is another equilibrium distribution for $\tilde{W}$. 
It is supported on $\{a_1, \ldots, a_m\} \subset W$, 
so $V_{\zeta}(\tilde{W}) = I_{\zeta}(\tilde{\mu}) \ge V_{\zeta}(W)$.  
This is equivalent to  $\gamma_{\zeta}(\tilde{W}) \le \gamma_{\zeta}(W)$.
Since $\gamma_{\zeta}(W) \le \gamma_{\zeta}(E) + \varepsilon$, our claim
follows.   
\end{proof}

\vskip .1 in
In what follows we will speak of `the' equilibrium measure $\mu_{\zeta}$
of a compact set $E$ of positive capacity, anticipating the uniqueness
of the equilibrium measure.  However, the arguments below would apply to
any equilibrium measure.

\subsection{Potential functions.} 

For each probability measure $\nu$ supported on
$\PP^1_{\Berk} \backslash \{\zeta\}$, define the potential function
\begin{equation*}
u_{\nu}(z,\zeta) \ = \ \int -\log_v(\delta(z,w)_{\zeta}) \, d\nu(w) \ .
\end{equation*}
Recall that a real-valued function $f(z)$ is lower semi-continuous if
\begin{equation*}
\liminf_{x \rightarrow z} f(x) \ \ge \ f(z) 
\end{equation*}
for each $z$.  This is equivalent to requiring that $f^{-1}((b,\infty))$ 
be open, for each $b \in \RR$.  We will say that $f(z)$ is strongly
lower semi-continuous if for each $z$
\begin{equation*}
\liminf_{x \rightarrow z} f(x) \ = \ f(z) \ .
\end{equation*}

\begin{proposition} \label{PropD7}
Let $\nu$ be a probability measure on $\PP^1_{\Berk}$, and take
$\zeta \notin \supp(\nu)$. 
Then $u_{\nu}(z,\zeta)$ is strongly lower semi-continuous,
and is continuous at each $z \notin \supp(\nu)$ $($including $\zeta$,
if continuity is understood relative 
the extended reals, $\RR \cup \{- \infty\}$$)$.  
Moreover, for each $p \in \PP^1(\CC_v)$, if $x$ approaches $p$ along a 
path $[y,p]$, then 
\begin{equation*}
\lim \begin{Sb} x \rightarrow p \\ x \in [y,p) \end{Sb} u_{\nu}(x,\zeta) 
\ = \ u_{\nu}(p,\zeta) \ .
\end{equation*} 
\end{proposition}

\begin{proof}
We will first show $u_{\nu}(z,\zeta)$ is lower semi-continuous.  
Lower semi-continuity is automatic at $z = \zeta$, since for all $z, w$
\begin{equation*}
-\log_v(\delta(z,w)_{\zeta}) \ \ge \ -\log_v(\delta(\zeta,w)_{\zeta}) \ .
\end{equation*}
Let $K$ be any closed neighborhood of $\supp(\nu)$ which does not contain
$\zeta$.  By Lemma \ref{LemD2}, there is a sequence of continuous functions
$\{f_{k}(z,w)\}$ which increase monotonically to
$-\log_v(\delta(z,w)_{\zeta})$ on $K \times K$, so
\begin{equation*}
u_{\nu}(z,\zeta) \ = \
\lim_{k \rightarrow \infty} \int f_k(z,w) \, d\nu(w) \ .
\end{equation*}
on $K$.  Since $K$ is compact,
each $u_k(z) = \int f_k(z,w) \, d\nu(w)$ is continuous on $K$.  Thus,
$u_{\nu}(z,\zeta)$ is an increasing limit of continuous functions,
hence is lower semi-continuous on $K$.  For any $z \ne \zeta$, 
we can choose $K$ so as to contain $z$ in its interior.  
Thus $u_{\nu}(z,\zeta)$ is lower semi-continuous everywhere.

Put $E = \supp(\nu)$, and consider the decomposition
\begin{equation*}
\delta(z,w)_{\zeta} \ = \  \frac{\|z,w\|_v}{\|z,\zeta\|_v \, \|w,\zeta\|_v} \ .
\end{equation*}
Inserting this in the definition of $u_{\nu}(z,\zeta)$ we see that
\begin{eqnarray*}
u_{\nu}(z,\zeta) & = &
    \int_E -\log_v(\|z,w\|_v) \, d\nu(w)
        + \int_E \log_v(\|z,\zeta\|_v) \, d\nu(w) \\
& & \qquad + \int_{E} \log_v(\|w,\zeta\|_v) \, d\nu(w) \ .
\end{eqnarray*}
The second integral is $\log_v(\|z,\zeta\|_v)$ since the integrand
does not involve $w$.  The third integral is a finite constant $C_{\zeta}$,
because $\zeta \notin E$, so $\log_v(\|w,\zeta\|_v)$
is a bounded, continuous function of $w \in E$.  Thus
\begin{equation*}
u_{\nu}(z,\zeta) \ = \
\int_E -\log_v(\|z,w\|_v) \, d\nu(w)
        +  \log_v(\|z,\zeta\|_v) + C_{\zeta} \ .
\end{equation*}

For each $z \notin E$, the first integral is continuous at $z$
since $\|x,y\|_v$ is continuous off the diagonal and $z$ has
a closed neighborhood disjoint from $E$.
The function $\log_v(\|z,\zeta\|_v)$
is continous as a function to the extended reals,
since $\|x,y\|_v$ is continuous as a function of each variable separately.
Hence $u_{\nu}(z,\zeta)$ is continuous (and in particular,
strongly lower semi-continuous) off $\supp(\nu)$.

Next suppose $z \in E$.
Recalling that $-\log_v(\|z,w\|_v) = j_{\zeta_0}(z,w)$,
where $\zeta_0 \in \PP^1_{\Berk}$ is the Gauss point,
(the point corresponding to $B(0,1)$),
let $t$ vary along the path from $\zeta_0$ to $z$.  For each $w \in E$,
$j_{\zeta_0}(t,w)$ increases monotonically to $j_{\zeta_0}(z,w)$.
By the monotone convergence theorem, $\int_E -\log_v(\|t,w\|_v) \, d\nu(w)$
increases monotonically to $\int_E -\log_v(\|z,w\|_v) \, d\nu(w)$.
Hence as $t$ approaches $z$ along this path,
\begin{equation*}
\lim_{t \rightarrow z} u_{\nu}(t,\zeta) \ = \ u_{\nu}(z,\zeta) \ .
\end{equation*}
Thus
$\liminf_{x \rightarrow z} u_{\nu}(x,\zeta) \le u_{\nu}(z,\zeta)$.
Combined with the opposite inequality coming from lower semi-continuity, 
this gives strong lower semi-continuity.

Finally, suppose $p \in \PP^1(\CC_v)$.  Let $y \in \PP^1_{\Berk}$ be
arbitrary, and let $x$ approach $p$ along the path $[y,p]$.  
(If another point $y^{\prime}$ had been chosen,  $[y,p]$ and $[y^{\prime},p]$ 
have a terminal segment $[y^{\prime \prime},p]$ in common).  The same 
reasoning as above shows that 
\begin{equation*}
\lim \begin{Sb} x \rightarrow p \\ x \in [y,p) \end{Sb} u_{\nu}(x,\zeta) 
\ = \ u_{\nu}(p,\zeta) \ .
\end{equation*} 

\end{proof}

\vskip .1 in
In the classical theory, the two main facts about potential functions are
Maria's theorem and Frostman's theorem.  We will now establish their
analogues for the Berkovich line.

\begin{theorem} \label{ThmD8} \text{\rm (Maria)}
Let $\nu$ be a probability measure supported on 
$\PP^1_{\Berk} \backslash \{\zeta\}$.
If there is a constant $M < \infty$ such that $u_{\nu}(z,\zeta) \le M$
on $\supp(\nu)$, then $u_{\nu}(z,\zeta) \le M$ 
for all $z \in \PP^1_{\Berk} \backslash \{\zeta\}$.
\end{theorem}

\begin{proof}
Put $E = \supp(\nu)$ and
fix $z \in \PP^1_{\Berk} \backslash (E \cup \{\zeta\})$.
Since $\delta(z,w)_{\zeta}$ is continuous off the diagonal,
there is a point $\overline{z} \in E$ such that
$\delta(z,\overline{z})_{\zeta} \le \delta(z,w)_{\zeta}$ for all $w \in E$.
By the ultrametric inequality, for each $w \in E$
\begin{equation*}
\delta(\overline{z},w)_{\zeta} \ \le \
\max(\delta(z,\overline{z})_{\zeta}, \delta(z,w)_{\zeta})
\ = \ \delta(z,w)_{\zeta} \ .
\end{equation*}
Hence
\begin{eqnarray*}
u_{\nu}(z,\zeta) & = & \int -\log_v(\delta(z,w)_{\zeta}) \, d\nu(w) \\
& \le & \int -\log_v(\delta(\overline{z},w)_{\zeta}) \, d\nu(w) \\
& = & u_{\nu}(\overline{z},\zeta) \ \le \ M \ .
\end{eqnarray*}
\end{proof}

The following lemma asserts that sets of capacity $0$ are `small' in a
measure-theoretic sense.

\begin{lemma} \label{LemD9}
If $f \subset \PP^1_{\Berk} \backslash \{\zeta\}$ is a set of
capacity $0$, then $\nu(f) = 0$ for any probability measure $\nu$
supported on $\PP^1_{\Berk} \backslash \{\zeta\}$ with
$I_{\zeta}(\nu) < \infty$.
\end{lemma}

\begin{proof}
Recall that $\supp(\nu)$ is compact.
After scaling $\delta(x,y)_{\zeta}$ if necessary, we can assume
that $\delta(x,y)_{\zeta} \le 1$ on $\supp(\nu)$.
It $\nu(f) > 0$, then $\nu(e) > 0$ for some compact subset $e$ of $f$,
and $\eta := (1/\nu(e)) \cdot \nu$ is probability measure $\eta$ on $e$.
It follows that
\begin{eqnarray*}
I_{\zeta}(\eta)
& = & \iint_{e \times e} -\log_v(\delta(x,y)_{\zeta}) \, d\nu(x) d\nu(y) \\
& \le & \frac{1}{\nu(e)^2}
         \iint -\log_v(\delta(x,y)_{\zeta}) \, d\nu(x) d\nu(y)
\ < \ \infty \ .
\end{eqnarray*}
\end{proof}

\vskip .1 in
\begin{corollary} \label{CorD10}
Let $\{f_n\}_{n \ge 1}$ be a countable collection of Borel sets in
$\PP^1_{\Berk} \backslash \{\zeta\}$ such that each $f_n$ has capacity $0$.
Put $f = \cup_{n=1}^{\infty} f_n$.  Then $f$ has capacity $0$.
\end{corollary}

\begin{proof}
If  $\gamma_{\zeta}(f) > 0$, there is a probability measure $\nu$ supported
on $f$ such that $I_{\zeta}(\nu) < \infty$.  Since each $f_n$ is
Borel measureable, so is $f$, and
\begin{equation*}
\sum_{n=1}^\infty \nu(f_n) \ \ge \ \nu(f) \ = \ 1 \ .
\end{equation*}
Hence $\nu(f_n) > 0$ for some $n$, contradicting Lemma \ref{LemD9}.
\end{proof}

\vskip .1 in
Recall that an $F_{\sigma}$ set is a countable union of compact sets.

\begin{theorem} \label{ThmD11} \text{\rm (Frostman)}
Let $E \subset \PP^1_{\Berk} \backslash \{\zeta\}$
be a compact set of positive capacity.
Then the equilbrium potential $u_E(z,\zeta)$ satisfies

$A)$  $u_E(z,\zeta) \le V_{\zeta}(E)$ for all
        $z \in \PP^1_{\Berk} \backslash \{\zeta\}$.

$B)$  $u_E(z,\zeta) = V_{\zeta}(E)$ for all $z \in E$,
        except possibly on an $F_{\sigma}$ set $f \subset E$ of
        capacity $0$.
        
$C)$  $u_E(z,\zeta)$ is continuous at each point $z_0$
         where $u_E(z_0,\zeta) = V_{\zeta}(E)$.       
\end{theorem}

\begin{proof}
First, using a quadraticity argument,
we will show that $u_E(z,\zeta) \ge V_{\zeta}(E)$ on $E$ 
except on the (possibly empty) set $f$.  Then, we will show that
$u_E(z,\zeta) \le V_{\zeta}(E)$ on the support $\supp(\mu)$
of the equilibrium measure $\mu = \mu_{\zeta}$.
Since $u_E(z,\zeta) = u_{\mu}(z,\zeta)$,
it follows from Maria's theorem that $u_E(z,\zeta) \le V_{\zeta}(E)$
for all $z$.

Put
\begin{eqnarray*}
  f & = & \{ z \in E : u_E(z,\zeta) < V_{\zeta}(E) \} \ , \\
f_n & = & \{ z \in E : u_E(z,\zeta) \le V_{\zeta}(E) - 1/n \} \ , \quad
       \text{for $n = 1, 2, 3, \ldots$.}
\end{eqnarray*}
Since $u_E(z,\zeta)$ is lower semicontinuous, each $f_n$ is closed, hence
compact, so $f$ is an $F_{\sigma}$ set.  By Corollary \ref{CorD10},
$\gamma_{\zeta}(f) = 0$ if and only if $\gamma_{\zeta}(f_n) = 0$ for
each $n$.

Suppose $\gamma_{\zeta}(f_n) > 0$ for some $n$;  then there is
a probability measure $\sigma$ supported on $f_n$ such that
$I_{\zeta}(\sigma) < \infty$.  On the other hand, since
\begin{equation*}
V_{\zeta}(E)
\ = \ \iint -\log_v(\delta(z,w)_{\zeta}) \, d\mu(w) d\mu(z)
\ = \ \int u_E(z,\zeta) \, d\mu(z) \ ,
\end{equation*}
there is a point $q \in \supp(\mu)$ with $u_E(q,\zeta) \ge V_{\zeta}(E)$.
Since $u_E(z,\zeta)$ is lower semicontinuous, there is a neighborhood
$U$ of $q$ on which $u_E(z,\zeta) > V_{\zeta}(E) - 1/(2n)$.
After shrinking $U$ if necessary, we can assume that its closure
$\Ubar$ is disjoint from $f_n$. Put $e_n = E \cap U$.
By the definition  of $\supp(\mu)$,
it follows that $M := \mu(U) = \mu(e_n) > 0$.
Define a measure $\sigma_1$ of total mass $0$ on $E$ by
\begin{equation*}
\sigma_1 \ = \ \left\{ \begin{array}{ll}
                  M \cdot \sigma & \text{on $f_n$,} \\
                  -\mu & \text{on $e_n$,} \\
                  0 & \text{elsewhere.}
                       \end{array} \right.
\end{equation*}
We claim that $I_{\zeta}(\sigma_1)$ is finite.  Indeed
\begin{eqnarray*}
I_{\zeta}(\sigma_1) & = &
   M^2 \cdot \iint_{f_n \times f_n}
       -\log_v(\delta(z,w)_{\zeta}) \, d\sigma(z) d\sigma(w)  \\
   & & \qquad
   - 2M \cdot \iint_{f_n \times e_n}
       -\log_v(\delta(z,w)_{\zeta}) \, d\sigma(z) d\mu(w) \\
   & & \qquad \qquad   + \iint_{e_n \times e_n}
       -\log_v(\delta(z,w)_{\zeta}) \, d\mu(z) d\mu(w) \ .
\end{eqnarray*}
The first integral is finite by hypothesis.  The second is finite because
$e_n$ and $\Ubar$ are disjoint, so $-\log_v(\delta(z,w)_{\zeta})$ is bounded
on $e_n \times f_n$.  The third is finite because $I_{\zeta}(\mu)$ is finite.

For each $0 \le t \le 1$, the measure $\mu_t := \mu + t \sigma_1$
is a probability measure on $E$.  By an expansion like the one above,
\begin{eqnarray*}
I_{\zeta}(\mu_t) - I_{\zeta}(\mu) & = &
   2t \cdot \int_E u_E(z,\zeta) \, d\sigma_1(z)
           + t^2 \cdot I_{\zeta}(\sigma_1) \\
   & \le & 2t \cdot( (V_{\zeta}(E)-1/n) - (V_{\zeta}(E) - 1/(2n)) \cdot M
            + t^2 \cdot I_{\zeta}(\sigma_1) \\
   & = & (-M/n) \cdot t + I_{\zeta}(\sigma_1) \cdot t^2 \ .
\end{eqnarray*}
For sufficiently small $t > 0$, the right side is negative.
This contradicts the fact that $\mu$ minimizes $I_{\zeta}(\nu)$
for all probability measures $\nu$ supported on $E$.  It follows that
$\gamma_{\zeta}(f_n)$ = 0, and hence that $\gamma_{\zeta}(f) = 0$.

The second part requires showing that $u_E(z,\zeta) \le V_{\zeta}(E)$
for all $z \in \supp(\mu)$.  If $u_E(q,\zeta) > V_{\zeta}(E)$ for
some $q \in \supp(\mu)$, take $\varepsilon > 0$ such that
$u_E(q,\zeta) > V_{\zeta}(E) + \varepsilon$.
The lower semicontinuity of $u_E(z,\zeta)$ shows there is a
neighborhood $U$ of $q$ on which $u_E(z,\zeta) > V_{\zeta}(E) + \varepsilon$.
Put $e = U \cap E$.  Since $q \in \supp(\mu)$, $T := \mu(e) = \mu(U) > 0$.
By Lemma \ref{LemD9}, $\mu(f) = 0$.  Therefore, since
$u_E(z,\zeta) \ge V_{\zeta}(E)$ on $E \backslash f$,
\begin{eqnarray*}
V_{\zeta}(E)
& = & \int_{e} u_E(z,\zeta) \, d\mu(z)
             + \int_{E \backslash e} u_E(z,\zeta) \, d\mu(z) \\
& \ge & T \cdot (V_{\zeta}(E) + \varepsilon) + (1-T) \cdot V_{\zeta}(E)
\ = \ V_{\zeta}(E) + T \varepsilon
\end{eqnarray*}
which is obviously false.

Thus $u_E(z,\zeta) \le V_{\zeta}(E)$ on $\supp(\mu)$, and Maria's
theorem implies that $u_E(z,\zeta) \le V_{\zeta}(E)$ for all $z$.

The final assertion, that $u_E(z,\zeta)$ is continuous at each point $z_0$
where $u_E(z_0,\zeta) = V_{\zeta}(E)$, is now trivial.  
By lower semicontinuity,
\begin{equation*}
\liminf_{z \rightarrow z_0} u_E(z,\zeta) \ \ge \ u_E(z_0,\zeta) 
\ = \ V_{\zeta}(E) \ .
\end{equation*}
On the other hand, since $u_E(z,\zeta) \le V_{\zeta}(E)$ for all $z$,  
\begin{equation*}
\limsup_{z \rightarrow z_0} u_E(z,\zeta) \ \le \ V_{\zeta}(E) 
\ = \ u_E(z_0,\zeta) \ .
\end{equation*}
\end{proof}

\vskip .1 in
\begin{corollary} \label{CorD12}
For any compact set $E \subset \PP^1_{\Berk} \backslash \{\zeta\}$
of positive capacity, and any probability measure $\nu$ supported on $E$,
\begin{equation*}
\inf_{z \in E} u_\nu(z,\zeta) \ \le \ V_{\zeta}(E)
       \ \le \ \sup_{z \in E} u_{\nu}(z,\zeta) \ .
\end{equation*}
\end{corollary}

\begin{proof}
This follows immediately from the identity
\begin{eqnarray*}
\int_E u_{\nu}(z,\zeta) \, d\mu(z) & = &
\iint_{E \times E} -\log_v(\delta(z,w)_{\zeta}) \, d\nu(w) d\mu(z) \\
& = & \int_E u_E(w,\zeta) d\nu(w) \ = \ V_{\zeta}(E) \ .
\end{eqnarray*}
Here the second equality follows from Tonelli's theorem, and the third
from Frostman's theorem and the fact that if $f$ is the exceptional
set on which $u_E(z,\zeta) < V_{\zeta}(E)$, then $\mu(f) = 0$
by Lemma (\ref{LemD9}).
\end{proof}

\vskip .1 in
We will now show that adjoining or removing a set of capacity $0$ from
a given set $F$ does not change its capacity.  This is a consequence of
a quantitative bound which we prove first.

\begin{proposition}  \label{PropD13}
Let $\{F_m\}_{m \ge 1}$ be a countable collection of sets contained in
$\PP^1_{\Berk} \backslash \zeta$, and put $F = \bigcup_{m=1}^{\infty} F_m$.
Suppose there is an $R < \infty$ such that $\delta(x,y)_{\zeta} \le R$
for all $x, y \in F$ $($or equivalently, that there is a disc 
$\cB(a,R)_{\zeta}$ containing $F$$)$.  Then
\begin{equation} \label{FMV2}
\frac{1}{V_{\zeta}(F) + \log_v(R) + 1} \ \le \
\sum_{m=1}^{\infty} \frac{1}{V_{\zeta}(F_m) + \log_v(R) + 1} \ .
\end{equation}
\end{proposition}

\begin{proof}
First suppose  $R < 1$.  In this case we will show that
\begin{equation} \label{FMV3}
\frac{1}{V_{\zeta}(F)} \ \le \
\sum_{m=1}^{\infty} \frac{1}{V_{\zeta}(F_m)} \ .
\end{equation}
Note that our hypothesis implies that
$V_{\zeta}(F) \ge -\log(R) > 0$.

If $\gamma_{\zeta}(F) = 0$, then $\gamma_{\zeta}(F_m) = 0$ for each $m$,
so $V_{\zeta}(F) = V_{\zeta}(F_m) = \infty$ and (\ref{FMV3}) is trivial.
Hence we can assume that $\gamma_{\zeta}(F) > 0$.  Let $E \subset F$
be a compact set with $\gamma_{\zeta}(E) > 0$,
and let $\mu$ be its equilbrium distribution.
For each $m$, put $E_m = E \cap F_m$.  Then
\begin{equation}  \label{FJK0}
\sum_{m=1}^{\infty} \mu(E) \ \ge \ \mu(E) \ = \ 1 \ .
\end{equation}
We claim that for each $m$, $V_{\zeta}(E)/V_{\zeta}(E_m) \ge \mu(E_m)$.
If $\mu(E_m) = 0$ there is nothing to prove, so suppose $\mu(E_m) > 0$.
There are compact sets $e_{m,1} \subset e_{m,2} \subset \cdots \subset E_m$
such that $\lim_{i \rightarrow \infty} \mu(e_{m,i}) = \mu(E_m)$ and
$\lim_{i \rightarrow \infty} V_{\zeta}(e_{m,i}) = V_{\zeta}(E_m)$.
Without loss we can assume $\mu(e_{m,i}) > 0$ for each $i$.
Put $\nu_{m,i} = (1/\mu(e_{m,i})) \cdot \mu|_{e_{m,i}}$.  Appling Corollary
\ref{CorD12}, we see that
\begin{equation} \label{FJK1}
\sup_{z \in e_{m,i}} u_{\nu_{m,i}}(z,\zeta) \ \ge \ V_{\zeta}(e_{m,i}) \ .
\end{equation}
But also, since $-\log_v(\delta(z,w)_{\zeta}) \ge -\log_v(R) > 0$ on $E$,
for each $z \in e_{m,i}$
\begin{eqnarray}
\frac{1}{\mu(e_{m,i})} u_{E}(z,\zeta)
& = & \frac{1}{\mu(e_{m,i})}
     \int_{E} -\log_v(\delta(z,w)_{\zeta}) \, d\mu(w) \label{FJK2} \\
& \ge & \int_{e_{m,i}} -\log_v(\delta(z,w)_{\zeta}) \, d\nu_{m,i}(w)
\ = \ u_{\nu_{m,i}}(z,\zeta) \notag
\end{eqnarray}
Because $\mu(e_{m,i}) > 0$, Theorem \ref{ThmD11} shows there exist points
$z \in e_{m,i}$ where $u_E(z,\zeta) = V_{\zeta}(E)$.  Combining (\ref{FJK1})
and (\ref{FJK2}) gives 
\begin{equation*}
V_{\zeta}(E)/\mu(e_{m,i}) \ \ge \ V_{\zeta}(e_{m,i}) \ .
\end{equation*}
Transposing terms and letting $i \rightarrow \infty$ shows
$V_{\zeta}(E)/V_{\zeta}(E_m) \ge \mu(E_m)$.

By (\ref{FJK0}) 
\begin{equation*}
\frac{1}{V_{\zeta}(E)} \ \le \ \sum_{m=1}^{\infty} \frac{1}{V_{\zeta}(E_m)} \ .
\end{equation*}
But $E$ can be chosen so that $V_{\zeta}(E)$ is arbitrarily close to
$V_{\zeta}(F)$, and also such that as many of the $V_{\zeta}(E_m)$ as we
wish are arbitrarily close to $V_{\zeta}(F_m)$.  Taking a limit over
such $E$, we obtain (\ref{FMV3}).

The general case follows by scaling the Hsia kernel.
Replace $\delta(x,y)_{\zeta}$ by
$\delta^{\prime}(x,y)_{\zeta} = 1/(q_v R) \cdot \delta(x,y)_{\zeta}$,
so $\delta^{\prime}(x,y)_{\zeta} \le 1/q_v$ for all $x,y \in F$.
For each $E \subset F$, this changes $V_{\zeta}(E)$ to
$V_{\zeta}^{\prime}(E) = V_{\zeta}(E) + \log_v(R) + 1$.
Applying (\ref{FMV3}) to $V_{\zeta}^{\prime}(F)$
and the $V_{\zeta}^{\prime}(F_m)$ gives (\ref{FMV2}).
\end{proof}

\vskip .1 in
\begin{corollary} \label{CorD14}
Let $e \subset \PP^1_{\Berk} \backslash \{\zeta\}$ have capacity $0$.
Then for any $F \subset \PP^1_{\Berk} \backslash \{\zeta\}$,
\begin{equation*}
\gamma_{\zeta}(F \cup e) \ = \ \gamma_{\zeta}(F \backslash e)
                         \ = \ \gamma_{\zeta}(F) \ .
\end{equation*}
\end{corollary}

\begin{proof}
By Proposition \ref{PropD13}, if $F$ and $e$ are contained in a ball
$\cB(a,R)_{\zeta}$ with $R < \infty$, then
\begin{eqnarray*}
\frac{1}{V_{\zeta}(F \cup e) + \log(R) + 1}
& \le & \frac{1}{V_{\zeta}(F) + \log(R) + 1} + \frac{1}{V_{\zeta}(e) + R + 1} \\
& = & \frac{1}{V_{\zeta}(F) + \log(R) + 1} \ ,
\end{eqnarray*}
giving $V_{\zeta}(F) \le V_{\zeta}(F \cup e)$.
The opposite inequality is trivial, so $V_{\zeta}(F \cup e) = V_{\zeta}(F)$.
The general case follows from this, by replacing $F$ and $e$ by
$F \cap \cB(a,R)_{\zeta}$ and $e \cap \cB(a,R)_{\zeta}$, 
fixing a center $a$ and letting
$R \rightarrow \infty$.
Since any compact subset of $\PP^1_{\Berk} \backslash \{\zeta\}$ is
contained in $\cB(a,r)_{\zeta}$ for some $R$,
\begin{eqnarray*}
\gamma_{\zeta}(F \cup e)
& = & \lim_{R \rightarrow \infty}
   \gamma_{\zeta}((F \cup e) \cap \cB(a,R)_{\zeta}) \\
& = & \lim_{R \rightarrow \infty}
   \gamma_{\zeta}(F \cap \cB(a,R)_{\zeta}) \ = \ \gamma_{\zeta}(F) \ .
\end{eqnarray*}

For the equality $\gamma_{\zeta}(F \backslash e) = \gamma_{\zeta}(F)$,
apply what has been just shown to $e$ and $F \backslash e$, noting that
$(F \backslash e) \cup e = F \cup e$:
\begin{equation*}
\gamma_{\zeta}(F \backslash e) \ = \ \gamma_{\zeta} (F \cup e)
\ = \ \gamma_{\zeta} (F) \ .
\end{equation*}
\end{proof}

\vskip .1 in
As in the classical theory, there are two other important capacitary functions:  
the {\it transfinite diameter} and the {\it Chebyshev constant}.
The existence of these quantities, and the fact that for compact sets they are
equal to the logarithmic capacity, will be established below.

\subsection{The transfinite diameter $d_{\infty}(E)_{\zeta}$.}

For each $n = 2, 3, 4, \ldots$ put
\begin{equation*}
d_n(E)_{\zeta} \ = \ sup_{x_1, \ldots, x_n \in E} \ 
     (\prod_{i \ne j} \delta(x_i,x_j)_{\zeta})^{1/(n(n-1))} \ .
\end{equation*}
Note that unlike the classical case, $d_n(E)_{\zeta}$ can be nonzero
even if some of the $x_i$ coincide;  this will happen, for example,
if $E = \{a\}$ for a point of type II, III, or IV.  We claim that the
sequence $\{d_n(E)_{\zeta}\}$ is monotone decreasing.  
Fix $n$.  If $d_{n+1}(E)_{\zeta} = 0$ then certainly
$d_n(E)_{\zeta} \ge d_{n+1}(E)_{\zeta}$. Otherwise, take $\varepsilon$ 
with  $0 < \varepsilon < d_{n+1}(E)_{\zeta}$.  There are points
$z_1, \ldots, z_{n+1} \in E$ such that
\begin{equation*}
\prod^{n+1} \begin{Sb} i, j = 1 \\ i \ne j \end{Sb} \delta(z_i,z_j)_{\zeta}
\ \ge \ (d_{n+1}(E)_{\zeta} - \varepsilon)^{(n+1)n} \ .
\end{equation*}
By the definition of $d_n(E)_{\zeta}$, for each $k = 1, \ldots, n+1$
\begin{equation*}
(d_n(E)_{\zeta})^{n(n-1)} \ \ge \
  \prod \begin{Sb} i \ne j \\ i, j \ne k \end{Sb} \delta(z_i,z_j)_{\zeta} \ .
\end{equation*}
Taking the product over all $k$, we see that
\begin{equation*}
(d_n(E)_{\zeta})^{(n+1)n(n-1)} \ \ge \
  (\prod^{n+1} \begin{Sb} i, j = 1 \\ i \ne j  \end{Sb}
          \delta(z_i,z_j)_{\zeta} )^{n-1}
  \ = \ (d_{n+1}(E)_{\zeta} - \varepsilon)^{(n+1)n(n-1)} \ .
\end{equation*}
Since $\varepsilon > 0$ is arbitrary, $d_n(E)_{\zeta} \ge d_{n+1}(E)_{\zeta}$.
Put
\begin{equation*}
d_{\infty}(E)_{\zeta} \ = \ \lim_{n \rightarrow \infty} d_n(E)_{\zeta} \ .
\end{equation*}

\subsection{The Chebyshev Constant $CH(E)_{\zeta}$.}

We will define three variants of the Chebyshev constant.
For each positive integer $n$ and 
$a_1, \ldots, a_n \in \PP^1_{\Berk} \backslash \{\zeta\}$
(which need not be distinct), define the `pseudo-polynomial'
\begin{equation*}
P_n(x;a_1, \ldots, a_n) \ = \ \prod_{i=1}^n \delta(x,a_i)_{\zeta} \ .
\end{equation*}
Put $\|P_n(x;a_1, \ldots, a_n)\|_E = \sup_{x \in E} P_n(x;a_1, \ldots, a_n)$ 
and let
\begin{eqnarray*}
\CH_n^*(E)_{\zeta} & = &
\inf_{a_1, \ldots, a_n \in E} \
      (\|P_n(x;a_1, \ldots, a_n)\|_E)^{1/n} \ , \\
\CH_n^a(E)_{\zeta} & = &
\inf_{a_1, \ldots, a_n \in \PP^(\CC_v) \backslash \{\zeta\}} \
     (\|P_n(x;a_1, \ldots, a_n)\|_E)^{1/n} \ , \\
\CH_n(E)_{\zeta} & = &
\inf_{a_1, \ldots, a_n \in \PP^1_{\Berk} \backslash \{\zeta\}} \
       (\|P_n(x;a_1, \ldots, a_n)\|_E)^{1/n} \ .         
\end{eqnarray*}
We will show the three numbers
\begin{eqnarray*}
\CH^*(E)_{\zeta} & = & \lim_{n \rightarrow \infty} \CH_n^*(E)_{\zeta} \ , \\
\CH^a(E)_{\zeta} & = & \lim_{n \rightarrow \infty} \CH_n^a(E)_{\zeta} \ , \\
\CH(E)_{\zeta} & = & \lim_{n \rightarrow \infty} \CH_n(E)_{\zeta} 
\end{eqnarray*}
exist, and that $\CH^*(E)_{\zeta} \ge \CH^a(E)_{\zeta} = \CH(E)_{\zeta}$.  
They will be called the restricted Chebyshev constant, 
the algebraic Chebyshev constant, and the unrestricted Chebyshev constant
respectively.

The proofs that $\CH^*(E)_{\zeta}$, $\CH^a(E)_{\zeta}$ and $\CH(E)_{\zeta}$ 
exist are similar;  we give the argument only for $\CH(E)_{\zeta}$.
Put $\alpha = \inf_n \CH_n(E)$.  If $\alpha = \infty$,
then $\CH(E)_{\zeta} = \infty$ and there is nothing to prove.
Otherwise, fix $\varepsilon > 0$.  Then there are
an $N$ and $a_1, \ldots, a_N \in \PP^1_{\Berk} \backslash \{\zeta\}$
such that
\begin{equation*}
\sup_{x \in E} P_N(x;a_1, \ldots, a_N) \ \le\  (\alpha + \varepsilon)^N \ .
\end{equation*}
Let $\{b_j\} = \{a_1, \ldots, a_N, a_1, \ldots, a_N, \ldots \}$
be the sequence which cyclicly repeats $\{a_1, \ldots, a_N\}$.
Put $M_0 = 1$ and for each $r = 1, \ldots, N-1$,
put $M_r = \|P_r(x;a_1, \ldots, a_r)\|_E$.  Since $M_N$ is finite,
each $M_r$ is finite as well.  
Put $M = \max_{0 \le r < N} M_r/(\alpha+\varepsilon)^r$.
For each $n$ we can write $n = qN+r$
with $q, r \in \ZZ$ and $0 \le r < N$. Then
\begin{eqnarray*}
CH_n(E)_{\zeta} & \le & (\|P_n(x; b_1, \ldots, b_n)\|)^{1/N} \\
                & \le & (M_N^q \cdot M_r)^{1/n}
                \ \le \ M^{1/n} \cdot (\alpha + \varepsilon) \ .
\end{eqnarray*}
It follows that
$\limsup_{n \rightarrow \infty} CH_n(E)_{\zeta} \le \alpha + \varepsilon$.
Since $\varepsilon > 0$ is arbitrary,
\begin{equation*}
\limsup_{n \rightarrow \infty} CH_n(E)_{\zeta} \ \le \ \alpha \ = \
\liminf_{n \rightarrow \infty} CH_n(E)_{\zeta} 
\end{equation*}
and $\CH(E)_{\zeta} = \lim_{n \rightarrow \infty} CH_n(E)_{\zeta}$ exists.

\vskip .1 in
Clearly $\CH(E)^*_{\zeta} \ge \CH(E)_{\zeta}$ and 
$\CH(E)^a_{\zeta} \ge CH(E)_{\zeta}$.  We will now show that
$\CH^a(E)_{\zeta} = CH(E)_{\zeta}$.  If $CH(E)_{\zeta} = \infty$ there is
nothing to prove.  Otherwise put $\alpha = CH(E)_{\zeta}$ and 
fix $\varepsilon > 0$.  Then there are an $N$ and
$a_1, \ldots, a_N \in \PP^1_{\Berk}$ such that 
$\|P_N(x;a_1, \ldots, a_N)\|_E \le (\alpha + \varepsilon)^N$.  
Take $\eta > 0$.
We claim that for each $i$, there is an $a_i^{\prime}$ of type I such that
$\delta(x,a_i^{\prime})_{\zeta} \le (1+\eta)^2 \delta(x,a_i)_{\zeta}$ 
for all $x \in \PP^1_{\Berk} \backslash \{\zeta\}$.  
If $a_i$ of type I, put $a_i^{\prime} = a_i$.  If $a_i$ is not of type I,
let $r_i = \diam_{\zeta}(a_i)$;  then there is an $a_i^{\prime}$ of type I
in the ball $\cB(a_i,r_i(1 + \eta))_{\zeta}$.
If $x \notin \cB(a_i,r_i(1+\eta))_{\zeta}$ then
$\delta(x,a_i^{\prime})_{\zeta} = \delta(x,a_i)_{\zeta}$ by the
ultrametric inequality.  If $x \in \cB(a_i,r_i(1+\eta))_{\zeta}$ then
again by the ultrametric inequality
\begin{equation*}
\delta(x,a_i^{\prime})_{\zeta} \ \le \ (1+\eta) r_i \ \le \ 
(1+\eta) \delta(x,a_i)_{\zeta} \ ,
\end{equation*}
and in either case the claim is true.  It follows that
\begin{eqnarray*}
\|P_N(x; a_1^{\prime}, \ldots, a_N^{\prime})\|_E  & \le &
\|P_N(x; a_1, \ldots, a_N)\|_E \cdot (1 + \eta)^N \\
& \le & (\alpha+\varepsilon)^N \cdot (1 + \eta)^N \ .
\end{eqnarray*}
Since $\eta$ is arbitrary, 
$\CH_N^a(E)_{\zeta} \le \alpha + \varepsilon$,
so $\CH^a(E)_{\zeta} = \CH(E)_{\zeta}$.

\vskip .1 in
\begin{theorem} \label{ThmD15}
For any compact set $E \subset \PP^1_{\Berk} \backslash \{\zeta\}$
\begin{equation*}
\gamma_{\zeta}(E) = d_{\infty}(E)_{\zeta}
= \CH(E)_{\zeta} = \CH^*(E)_{\zeta} = \CH^a(E)_{\zeta} \ .
\end{equation*}
\end{theorem}

\begin{proof}
We have seen that
$\CH^a_{\zeta}(E) = \CH(E)_{\zeta} \le CH^*(E)_{\zeta}$.
We now show that $\gamma_{\zeta}(E) \le \CH^a(E)_{\zeta}$,
$\CH^*(E)_{\zeta} \le d_{\infty}(E)$,
and $d_{\infty}(E)_{\zeta} \le \gamma_{\zeta}(E)$.

\vskip .1in
{\noindent I. $\gamma_{\zeta}(E) \le \CH^a(E)_{\zeta}$.}

If $\gamma_{\zeta}(E) = 0$ there is nothing to prove.
Suppose $\gamma_{\zeta}(E) > 0$.
We will show that $\gamma_{\zeta}(E) \le \CH^a_n(E)_{\zeta}$ for each $n$.
Fix $n$ and take $\varepsilon > 0$;
let $a_1, \ldots, a_n \in \PP^1(\CC_v) \backslash \{\zeta\}$ be points such
that
\begin{equation*}
\sup_{z \in E} \prod_{i=1}^n \delta(z,a_i)_{\zeta} \ \le \
(\CH^a_n(E)_{\zeta} + \varepsilon)^n \ .
\end{equation*}
Any finite set of Type I points has capacity $0$.
By Corollary \ref{CorD14}, replacing $E$ by
$E \cup \{a_1, \ldots, a_n\}$ does not change its capacity.
Let $\nu$ be the probability measure on $E$ with a point mass $1/n$
at each $a_i$.  Then for each $z \in E$,
\begin{equation*}
-\log_v(\CH^a_n(E)_{\zeta} + \varepsilon)
\ \le \ (1/n) \sum_{i=1}^n -\log_v(\delta(z,a_i)_{\zeta})
\ = \ u_{\nu}(z,\zeta) \ .
\end{equation*}
By Corollary \ref{CorD12} $\inf_{z \in E} u_{\nu}(z,\zeta) \le V_{\zeta}(E)$.
Thus $-\log_v(\CH^a_n(E)_{\zeta} + \varepsilon) \le V_{\zeta}(E)$,
which gives $\gamma_{\zeta}(E) \le \CH^a_n(E)_{\zeta} + \varepsilon$.

\vskip .1 in
{\noindent II. $\CH^*(E)_{\zeta} \le d_{\infty}(E)_{\zeta}$.}

We will show that $d_{n+1}(E)_{\zeta} \ge \CH^*_n(E)_{\zeta}$ for each $n$.
Since $E$ is compact, the $\sup$ defining $d_{n+1}(E)_{\zeta}$ is achieved:
there are points $x_1, \ldots, x_{n+1} \in E$ such that
\begin{equation*}
d_{n+1}(E)^{(n+1)n} \ = \ \prod_{i \ne j} \delta(x_i,x_j)_{\zeta} \ .
\end{equation*}
If $x_1, \ldots, x_n$ are fixed, then by the definition of
$d_{n+1}(E)_{\zeta}$,
\begin{equation*}
H(z; x_1, \ldots, x_n)
:= \prod_{1 \le i < j \le n} \delta(x_i,x_j)_{\zeta}
\cdot \prod_{i=1}^n \delta(z,x_i)_{\zeta}
\end{equation*}
achieves its maximum for $z \in E$ at $z = x_{n+1}$.  But by the definition
of $\CH_n(E)_{\zeta}$,
\begin{equation*}
\max_{z \in E} \prod_{i=1}^n \delta(z,x_i)_{\zeta}
\ \ge \ (\CH^*_n(E)_{\zeta})^n \ .
\end{equation*}
Thus,
\begin{equation*}
d_{n+1}(E)_{\zeta}^{(n+1)n/2} \ \ge \
\prod_{1 \le i < j \le n} \delta(x_i,x_j)_{\zeta} \cdot
(\CH^*_n(E)_{\zeta})^n \ .
\end{equation*}
A similar inequality holds if $\{x_1, \ldots, x_n\}$ is replaced by
$\{x_1, \ldots, x_{n+1}\} \backslash \{x_k\}$ for each $k = 1, \ldots, n+1$.
Multiplying these inequalities together, we find
\begin{eqnarray*}
(d_{n+1}(E)_{\zeta}^{(n+1)n/2})^{n+1}
& \ge & (\CH^*_n(E)_{\zeta})^{(n+1)n} \cdot
\prod_{k=1}^{n+1} (\prod \begin{Sb} 1 \le i < j \le n+1 \\ i, j \ne k \end{Sb}
\delta(x_i,x_j)_{\zeta}) \\
& = & (\CH^*_n(E)_{\zeta})^{(n+1)n} \cdot
(\prod_{1 \le i < j \le n+1} \delta(x_i,x_j)_{\zeta})^{n-1} \\
& = & (\CH^*_n(E)_{\zeta})^{(n+1)n} \cdot
(d_{n+1}(E)_{\zeta}^{(n+1)n/2})^{n-1} \ .
\end{eqnarray*}
Cancelling terms and taking roots gives
$d_{n+1}(E)_{\zeta} \ge \CH^*_n(E)_{\zeta}$.

\vskip .1 in
{\noindent III. $d_{\infty}(E)_{\zeta} \le \gamma_{\zeta}(E)$.}

We will show that $d_{\infty}(E)_{\zeta} \le \gamma_{\zeta}(E) + \varepsilon$
for each $\varepsilon > 0$.

Fix $\varepsilon$.
By Corollary \ref{CorD6} there is a closed neighborhood $W$ of $E$ of the form
$W = \cup_{k=1}^m \cB(a_k,r_k)_{\zeta}$, where each $a_k$ is of type II or 
type III and $r_k = \diam_{\zeta}(a_k)$, such that
$\gamma_{\zeta}(W) \le \gamma_{\zeta}(E) + \varepsilon$.
Let $r_{\zeta} : E \rightarrow \partial W_{\zeta}$ be the retraction map
which takes each $x$ in $E$ to the last point on the path from $x$ to $\zeta$
which belongs to $W$;  thus $r_{\zeta}(E) \subset \{a_1, \ldots, a_m\}$.
Put $r = \min(r_k) > 0$.

Note that if
$x, y \in E$ and $r_{\zeta}(x) = a_k$, $r_{\zeta}(y) = a_{\ell}$,
then $\delta(x,y)_{\zeta} \le \delta(a_k,a_{\ell})_{\zeta}$.
Indeed, if $k = \ell$ then $x,y \in B(a_k,r_k)_{\zeta}$ so
$\delta(x,y)_{\zeta} \le d_k = \delta(a_k,a_{\ell})_{\zeta}$.
If $k \ne \ell$ then $\delta(x,y)_{\zeta} = \delta(a_k,a_{\ell})_{\zeta}$
since the paths from $x$ and $y$ to $\zeta$ encounter $a_k$, $a_{\ell}$
before they meet, so $\delta(x,y)_{\zeta} = \delta(a_k,a_{\ell})_{\zeta}$.

Fix $n$.
Since $E$ is compact, there are points $x_1, \ldots, x_n \in E$ such that
\begin{equation*}
(d_n(E)_{\zeta})^{n(n-1)} \ = \ \prod^n \begin{Sb} i, j = 1 \\ i \ne j \end{Sb}
                           \delta(x_i,x_j)_{\zeta}
\end{equation*}
For each $a_k$, let $m_k$ be the number of points $x_i$ for
which $r_{\zeta}(x_i) = a_k$, and let $\nu$ be the probability measure
on $W$ given by $\nu = \sum_{k=1}^m (m_k/n) \delta_{a_k}(z)$.
Then $I_{\zeta}(\nu) \ge V_{\zeta}(W)$.  It follows that
\begin{eqnarray*}
-(1-\frac{1}{n}) \log_v(d_n(E)_{\zeta})
& = & \frac{1}{n^2} \sum_{i \ne j} -\log_v(\delta(x_i,x_j)_{\zeta}) \\
& \ge & \frac{1}{n^2} \sum_{i \ne j}
        -\log_v(\delta(r_{\zeta}(x_i),r_{\zeta}(x_j))_{\zeta}) \\
& = & I_{\zeta}(\nu)
        + \frac{1}{n^2} \sum_{k=1}^m m_k \log_v(r_k))) \\
& \ge & V_{\zeta}(W) + (1/n) \log_v(r) \ .
\end{eqnarray*}
Letting $n \rightarrow \infty$,
we see that $-\log_v(d_{\infty}(E)_{\zeta}) \ge V_{\zeta}(W)$.
Equivalently, $d_{\infty}(E)_{\zeta} \le \gamma_{\zeta}(W)
                        \le \gamma_{\zeta}(E) + \varepsilon$
as claimed.
\end{proof}

\section{Harmonic functions on the Berkovich Line.}
\label{Section E}

In this section we develop a theory of harmonic functions
on subdomains of the Berkovich line.  We define a Laplacian operator
and prove analogues of the Maximum Principle, Poisson's formula, 
Harnack's theorem, the Riemann Extension theorem, 
and the theory of Green's functions.
We characterize harmonic functions
as limits of logarithms of norms of rational functions, and show their
stability under uniform limits and pullbacks by meromorphic functions.    
As a byproduct, we obtain the uniqueness of the equilibrium measure,
asserted in Section~\ref{Section D}.

\vskip .1 in
By a {\it domain} in Berkovich space we mean a nonempty, 
connected, open subset of $\PP^1_{\Berk}$. 

Our idea for constructing the Laplacian on a domain in Berkovich space is to use 
the Laplacian on metrized graphs, and to take a limit over all metrized graphs
contained in the domain.  For a function $f$ satisfying suitable hypotheses,
the Laplacian of the restriction of $f$ to each subgraph of the domain 
is a measure, so it defines a functional on continuous functions on graph.  
The Laplacians on the graphs satisfy a compatibility condition called 
{\it coherence}, which will be explained below.   Passing to the limit, 
we obtain a canonical functional on the space of continuous functions on 
the domain.  By the Riesz representation theorem, this corresponds to a 
measure;  we define it to be the Laplacian of $f$.  

\subsection{Continuous functions.}   

Given a domain $U \subset \PP^1_{\Berk}$, let $\Ubar$ be its closure.  
Our first task is to understand the space
of continuous functions $f : \Ubar \rightarrow \RR$.  Write $\cC(\Ubar)$
for this space.  

\vskip .1 in 
When $\Ubar$ is a Berkovich disc, there is a direct description of $\cC(\Ubar)$
using the definition of the Berkovich topology.  

\begin{proposition} \label{PropE1}
Let $\cB(0,R)$ be a closed Berkovich disc.
Then linear combinations of the functions $x \rightarrow [F]_x$, 
for polynomials $F(T) \in \CC_v[T]$, are dense in $\cC(\cB(0,R))$.
\end{proposition}

\begin{proof}
By the definition of the Berkovich topology, sets of the form 
$U_f(a,b) = \{ x \in \cB(0,R) : [f]_x \in (a,b) \}$ for 
$f \in \CC_v\{\{R^{-1}T\}\}$ give a basis for the open sets of $\cB(0,R)$.
Since $\cB(0,R)$ is compact, the Stone-Weierstrass theorem tells us that  
the set of polynomials in the functions $[f]_x$ is
dense in $\cC(\cB(0,R))$.  The multiplicativity of the seminorms $[ \cdot ]_x$
shows that each monomial  $[f_1]_x^{k_1} \cdots [f_n]_x^{k_n}$  
reduces to a single term $[f_1^{k_1} \cdots f_n^{k_n}]_x$.  
Finally, the Weierstrass Preparation theorem shows that
for each $f \in \CC_v\{\{R^{-1}T\}\}$ there is a polynomial $F \in \CC_v[T]$
with $[f]_x = [F]_x$ for all $x \in \cB(0,r)$.  
\end{proof}

\vskip .1 in
We now seek a description of $\cC(\Ubar)$ which works in general.

\vskip .05 in
As in Section~\ref{Section C}, 
by a {\it subgraph} $\Gamma$ of $\PP^1_{\Berk}$ we mean
a connected, finitely branched metrized subgraph, whose metric is given by 
the logarithmic path distance $\rho(x,y)$ induced from the big model,
and having finite total length.  Such a graph $\Gamma$  is necessarily a tree.   

\begin{lemma} \label{LemE2}
Let $\Gamma$ be a subgraph of $\PP^1_{\Berk}$.   Then 

$A)$ $\Gamma$ is a closed subset of $\PP^1_{\Berk}$. 

$B)$  The metric topology on $\Gamma$ coincides with the relative
topology induced from $\PP^1_{\Berk}$.
 
\end{lemma}

\begin{proof}  Part (A) follows from part (B):  since $\Gamma$ is compact
in the metric topology, it is compact in the relative topology, hence closed.    

We will now prove part (B).       
Fix a system of coordinates on $\PP^1_{\Berk}$, and regard $\Gamma$ as a 
subset of $\AA^1_{\Berk}$.  Thus, we can view $\Gamma$ as a tree
of discs, partially ordered by the function $x \rightarrow \radius(x)$:
\begin{equation*}
\Gamma \ = \ \bigcup_{i=1}^m [r_i,R]_{x_i}
\end{equation*}
where the $x_i$ are the points of minimal radius in $\Gamma$, 
and $r_i = \radius(x_i)$ for each $i$.  
  
\vskip .05 in
We will first show that each $x \in \Gamma$ has a basis of neighborhoods
in the metric topology which are open in the relative topology.
There are several cases to consider.  

(1)  $x = x_i$ is an endpoint of $\Gamma$, and is a point of minimal radius in 
its branch.  In this case, a basis of the open neighborhoods of $x$ 
given by the half-open segments $[r_i,r_i+\varepsilon)_{x_i}$ 
for sufficiently small $\varepsilon > 0$.  
Such a segment is the intersection of the disc 
$\cB(x_i,r_i+\varepsilon)^- = 
\{z \in \AA^1_{\Berk} : \delta(z,x_i)_{\infty} < r_i + \varepsilon\}$
with $\Gamma$.

(2)  $x$ is an interior point of an edge of $\Gamma$.  Suppose
$x \in [r_i,R]_{x_i}$, and put $r = \radius(x)$.  Then a basis of the open
neighborhoods of $x$ is given by the open segments 
$(r-\varepsilon,r+\varepsilon)_{x_i}$ for sufficiently small $\varepsilon > 0$.
Such a segment is the intersection of the open annulus 
$\cB(x_i,r+ \varepsilon)^-\backslash \cB(x_i,r-\varepsilon)$ 
with $\Gamma$.

(3)  $x$ is a branch point of $\Gamma$, but is not the point of maximal
radius.  In this case, after relabeling the branches if necessary, we can
assume that $[r_1,R]_{x_1}, \ldots, [r_k,R]_{x_k}$ come together at $x$.  
Put $r = \radius(x)$.  Then a basis of the open neighborhoods of $x$ is given 
by star-shaped sets of the form 
$\cup_{i=1}^k (r-\varepsilon,r+\varepsilon)_{x_i}$ for sufficiently small
$\varepsilon > 0$.  Such a set is the intersection of the punctured disc
$\cB(x_1,r+\varepsilon)^- \backslash (\cup_{i=1}^k \cB(x_i,r-\varepsilon))$
with $\Gamma$.

(4)  $x$ is the point of maximal radius $R$ in $\Gamma$.  
Some of the branches of $\Gamma$ may come together at points below $x$;  
after relabeling the branches if necessary, we can assume that 
$[r_1,R]_{x_1}, \ldots, [r_k,R]_{x_{\ell}}$ come together at $x$.
Then a basis of the open neighborhoods of $x$ is given by star-shaped 
sets of the form $\cup_{i=1}^{\ell} (R-\varepsilon,R]_{x_i}$ for 
sufficiently small $\varepsilon > 0$.  Such a set is the intersection of 
$\PP^1_{\Berk} \backslash (\cup_{i=1}^{\ell} \cB(x_i,R-\varepsilon))$
with $\Gamma$.

\vskip .05 in  
Next we will show that each subset of $\Gamma$ which is open in the
relative topology is also open in the metric topology.  It suffices to 
consider the intersection of $\Gamma$ with a basic open set of the form
$\cB(a,s)^- \backslash (\cup{j=1}^k \cB(a_j,s_j))$.   
Using the notation introduced above,  
for each $j = 1, \ldots, k$ let $T_j$ be the set of $x_i$ for
which $[r_i,R]_{x_i} \cap \cB(a_j,r_j)$ is nonempty, 
and let $T$ be the set of $x_i$, not belonging to any $T_j$, 
for which $[r_i,R]_{x_i} \cap \cB(a,r)^-$ is nonempty.  
Then the intersection
of $\cB(a,r)^- \backslash (\cup_{j=1}^k \cB(a_j,r_j))$
with $\Gamma$ is the union of the intervals of discs  
\begin{equation*}
\bigcup_{j=1}^k (\bigcup_{i \in T_j} (s_i,s)_{x_i}) 
                      \cup (\bigcup_{i \in T} [r_i,s)_{x_i})
\end{equation*}   
which is open in $\Gamma$.  
\end{proof}   

\vskip .1 in
Now let $K \subset \PP^1_{\Berk}$ be an arbitrary nonempty connected closed set,
equipped with the relative topology.  Recall that a set is connected
under the topology of $\PP^1_{\Berk}$ if and only if it is path-connected,
and there is unique path between any two points of $\PP^1_{\Berk}$.

There is a retraction map $r_K : \PP^1_{\Berk} \rightarrow K$
defined as follows.  Fix a point $k_0 \in K$.  Given $x \in \PP^1_{\Berk}$,
let $r_K(x)$ be the first point $k$ in $K$ on the path from $x$ to $k_0$.
Clearly $r_K(x) = x$ if $x \in K$.  To see that $r_K(x)$ is independent of 
the choice of $k_0$ when $x \notin K$, suppose
$k_0^{\prime} \in K$ is another point, and $k^{\prime}$ is the first point
in $K$ on the path from $x$ to $k_0^{\prime}$.  If $k^{\prime} \ne k$, then
since $K$ is connected there is a path from $k^{\prime}$ to $k$ contained
in $K$.  There is another path from $k^{\prime}$ to $k$ gotten by
concatenating the path from $k^{\prime}$ to $x$ with the path from $x$ to
$k$, and eliminating backtracking.  This second path lies outside
of $K$, apart from its endpoints.  This contradicts the fact that there
is a unique path between any two points of $\PP^1_{\Berk}$, so it must
be that $k^{\prime} = k$.

\begin{lemma} \label{LemE3}
For each nonempty closed connected subset  $K \subset \PP^1_{\Berk}$, 
the retraction map $r_K : \PP^1_{\Berk} \rightarrow K$  is continuous.
\end{lemma}

\begin{proof}  The connected open sets form a basis for the topology of
$\PP^1_{\Berk}$.  If $U \subset \PP^1_{\Berk}$ is a connected
open set, then $U \cap K$ is also connected, since if $x, y \in U \cap K$
the path $P$ from $x$ to $y$ is contained in $U$ since $U$ is connected,
and it is contained in $K$ since $K$ is connected.  Hence it is contained
in $U \cap K$.  It follows that the relative topology on $K$ has a basis
consisting of connected open sets.

Let $V \subset K$ be a connected open set, and let $U \subset \PP^1_{\Berk}$
be an open set with $V = U \cap K$.  If $U_0$ is the connected component of
$U$ containing $V$, then $U_0$ is also open;  so we can assume without loss
that $U$ is connected.  Let $\tilde{U}$ be the union of all the connected open
subsets $U \subset \PP^1_{\Berk}$ with $U \cap K = V$.  Then $\tilde{U}$
is itself connected and open, and $\tilde{U} \cap K = V$, so it is the maximal
set with these properties.

We claim that $r_K^{-1}(V) = \tilde{U}$.  First, we will show that
$r_K^{-1}(V) \subset \tilde(U)$.  Suppose $x \in r_K^{-1}(V)$ and put
$k = r_K(x) \in V$.  Consider the path $P$ from $x$ to $k$.
If $x \notin \tilde{U}$, let $\xbar$ be the first point in the closure of
$\tilde{U}$ in $P$.  Then $\xbar \notin \tilde{U}$, since $\tilde{U}$ is
open;  also, $\xbar \notin K$, since $k$ is the only point of $P$ in $K$
and $k \in \tilde{U}$. Since $K$ is compact and $\PP^1_{\Berk}$ is
Hausdorff, there is a neighborhood $W$ of $\xbar$ which is disjoint from $K$.
Without loss, we can assume $W$ is connected.  Since $\xbar$ is in the closure
of $\tilde{U}$, $W \cap \tilde{U}$ is nonempty.  It follows that
$W \cup \tilde{U}$ is connected and open, and $(W \cup \tilde{U}) \cap K = V$.
By the maximality of $\tilde{U}$ we must have $W \subset \tilde{U}$.
This contradicts the fact that $\xbar \notin \tilde{U}$.
Hence $x \in \tilde{U}$.

Next, we will show that $\tilde{U} \subset r_K^{-1}(V)$.  Fix a point
$k_0 \in V$, and let $x \in \tilde{U}$ be arbitrary.  Since $\tilde{U}$
is connected, the path from $x$ to $k_0$ must be entirely contained in
$\tilde{U}$.  Hence, the point $k = r_K(x)$, which is the first point in $K$
along that path, belongs to $\tilde{U}$.
But then $k \in \tilde{U} \cap K = V$, so $x \in r_K^{-1}(V)$.
\end{proof}

\vskip .1 in
If $K_1 \subset K_2 \subset \PP^1_{\Berk}$ are two nonempty 
connected closed subsets, the retraction map 
$r_{K_1} : \PP^1_{\Berk} \rightarrow K_1$
induces  a retraction map $r_{K_2,K_1} : K_2 \rightarrow K_1$.  
Clearly 
\begin{equation*}
r_{K_1}(x) \ = \ r_{K_2,K_1}(r_{K_2}(x))
\end{equation*}
for all $x$.  If $K_1$ and $K_2$ both have the relative topology, 
it follows from Lemma \ref{LemE3} that $r_{K_2,K_1}$ is continuous.      

\begin{proposition} \label{PropE4}
Let $U \subset \PP^1_{\Berk}$ be a domain.  

$A)$ As $\Gamma$ ranges over all subgraphs of $U$, 
and as $f$ ranges over $\cC(\Gamma)$, 
the functions of the form $f \circ r_{\Ubar,\Gamma}(x)$ 
are dense in $\cC(\Ubar)$.

$B)$ As $\Gamma$ ranges over all subgraphs of $U$, 
and as $f$ ranges over $\CPA(\Gamma)$, the functions of the form 
$f \circ r_{\Ubar,\Gamma}(x)$ are dense in $\cC(\Ubar)$.  
\end{proposition} 

\begin{proof} 
For part (A), we will apply the Stone-Weierstrass theorem. 
Functions of the form $f \circ r_{\Ubar,\Gamma}(x)$ separate
points in $\Ubar$:  given distinct points $x, y \in \Ubar$, take $\Gamma$
to be a closed segment $[p,q]$ in the path from $x$ to $y$, 
and let $f \in \cC([p,q])$ be any function with $f(p) \ne f(q)$.  Then
$f \circ r_{\Ubar,\Gamma}(x) \ne f \circ r_{\Ubar,\Gamma}(y)$.  
Likewise, each constant function on $\Ubar$ has the form 
$f \circ r_{\Ubar,\Gamma}(x)$, where $\Gamma \subset U$ is arbitrary and 
$f$ is the corresponding constant function on $\Gamma$.  

Take $F \in \cC(\Ubar)$.  By the Stone-Weierstrass theorem 
(see \cite{Kelley}, p. 244), since $\Ubar$ is compact 
the algebra of functions generated by the $f \circ r_{\Gamma}(x)$ 
is dense in $\cC(\Ubar)$.  Hence, for any $\varepsilon > 0$, there are
a finite number of subgraphs $\Gamma_{ij} \subset U$, 
$i = 1, \ldots, m$, $j = 1, \ldots, n_i$,
and functions $f_{ij} \in \cC(\Gamma_{ij})$, for which
\begin{equation*}
|\, F(x) - (\sum_{i=1}^m \prod_{j=1}^{n_i} f_{ij} \circ r_{\Ubar,\Gamma_{ij}}(x)) \, | 
\ < \ \varepsilon 
\end{equation*}
for all $x \in \Ubar$.  Let $\Gamma$ be the smallest connected 
subgraph containing all the $\Gamma_{ij}$:  then $\Gamma \subset U$.  
For each $i$, $j$, put $g_{ij} = f_{ij} \circ r_{\Gamma,\Gamma_{ij}}(x)$;  
then $g_{ij} \in \cC(\Gamma)$ and 
$g_{ij} \circ r_{\Ubar,\Gamma}(x) = f_{ij} \circ r_{\Ubar,\Gamma_{ij}}(x)$ 
for all $x \in \Ubar$.  Put 
\begin{equation*}
g  \ = \ \sum_{i=1}^m (\prod_{j=1}^{n_i} g_{ij})  \in \ \cC(\Gamma) \ .
\end{equation*}
Then $|F(x) - g \circ r_{\Ubar,\Gamma}(x)| < \varepsilon$ for all $x \in \Ubar$.  

Part (B) follows from part (A), since $\CPA(\Gamma)$ is dense in $\cC(\Gamma)$
under the $\sup$ norm, for each $\Gamma$. 
\end{proof}

\subsection{Coherent systems of measures.}

For each subgraph $\Gamma \subset U$, 
let $\mu_{\Gamma}$ be a bounded Borel measure on $\Gamma$.

\begin{definition} \label{DefE1}
A system of measures $\{\mu_{\Gamma}\}$ on the subgraphs of \,  $U$ 
is called {\rm coherent} if

$A)$  For each pair of subgraphs with $\Gamma_1 \subset \Gamma_2$, 
($r_{\Gamma_2,\Gamma_1})_*(\mu_{\Gamma_2}) = \mu_{\Gamma_1}$, 
that is, for each Borel subset \, $T \subset \Gamma_1$,
\begin{equation*}
\mu_{\Gamma_2}(r_{\Gamma_2,\Gamma_1}^{-1}(T)) \ = \ \mu_{\Gamma_1}(T) \ .
\end{equation*}

$B)$  There is a constant $B$ such that
$|\mu_{\Gamma}|(\Gamma) \le B$ for each $\Gamma$.
\end{definition}

\vskip .1 in
The collection of subgraphs $\Gamma \subset U$ forms a directed set
under containment:  for any two graphs $\Gamma_1$, $\Gamma_2$, there is
a unique minimal subgraph $\Gamma_3$ containing $\Gamma_1$ and $\Gamma_2$.

If $\mu$ is a bounded Borel measure on $\Ubar$, 
and $\mu_{\Gamma} = (r_{\Ubar,\Gamma})_*(\mu)$ for each $\Gamma \subset U$,
then $\{\mu_{\Gamma}\}$ is clearly a coherent system of measures on the
subgraphs of $U$.  The following proposition shows that every coherent 
system arises in this way:  there is a 1-1 correspondence between 
measures $\mu$ on $\Ubar$ and coherent systems of measures on subgraphs of $U$.

\begin{proposition} \label{PropE5}
If $\{\mu_{\Gamma}\}$ is a coherent system of measures in $U$, the map
\begin{equation*}
 \Lambda(F) \ = \
    \lim \begin{Sb} \longrightarrow \\ \Gamma \end{Sb}
                    \int_{\Gamma} F(x) \, d\mu_{\Gamma}(x)
\end{equation*}
defines a bounded linear functional on $\cC(\Ubar)$, so there is
a unique bounded Borel measure $\mu$ on $\Ubar$ such that 
\begin{equation*}
\Lambda(F) \ = \ \int_{\Ubar} F(x) \, d\mu(x) 
\end{equation*}
for each $F \in \cC(\Ubar)$.  
This measure is characterized by the property that 
$(r_{\Ubar,\Gamma})_*(\mu) = \mu_{\Gamma}$ for each subgraph $\Gamma \subset U$.   

\end{proposition}

\begin{proof}
Fix $\Gamma_0 \subset U$ and $f_0 \in \cC(\Gamma_0)$.
Put $F_0 = f_0 \circ r_{\Ubar,\Gamma_0}(x)$.

Since $\{\mu_{\Gamma}\}$ is a coherent system of measures,
for each subgraph $\Gamma \subset U$ containing $\Gamma_0$
\begin{equation*}
\int_{\Gamma} f_0(r_{\Gamma,\Gamma_0}(x)) \, d\mu_{\Gamma}(x)
\ = \ \int_{\Gamma_0} f_0(x) \, d\mu_{\Gamma_0}(x) \ .
\end{equation*}
Then, since $F_0|_{\Gamma} = f_{0} \circ r_{\Gamma,\Gamma_0}$ 
for each $\Gamma$ containing $\Gamma_0$,  
\begin{equation*}
\Lambda(F_0) \ = \  \lim \begin{Sb} \longrightarrow \\ \Gamma \end{Sb}
                    \int_{\Gamma} F_0(x) \, d\mu_{\Gamma}(x)
          \ = \ \int_{\Gamma_0} f_0(x) \, d\mu_{\Gamma_0}(x)
\end{equation*}
exists.  

By Proposition \ref{PropE4},
functions of the form $F = f \circ r_{\Ubar,\Gamma}$ for
$f \in \cC(\Gamma)$ and $\Gamma \subset U$ are dense in $\cC(\Ubar)$
under the $\sup$ norm, so $\Lambda$ is defined on a dense subset
of $\cC(\Ubar)$.  

On the other hand, since $|\mu_{\Gamma}|$ has total mass at most $B$ 
for each $\Gamma$, for each $G \in \cC(\Ubar)$  
\begin{equation*}
\limsup_{\Gamma} |\int_{\Gamma} G(x) \, d\mu_{\Gamma}(x)| 
\ \le \ B \cdot \|G\|_{\Ubar}
\end{equation*}
It follows that $\Lambda$ extends to a bounded linear   
functional on $\cC(\Ubar)$.  The Riesz representation theorem
(see Rudin (\cite{Rud}, Theorem 6.19, p. 139) 
tells us there is a Borel measure $\mu$ on
$\Ubar$, with $|\mu|(\Ubar) \le B$, such that
\begin{equation*}
\Lambda(F) \ = \ \int_{\Ubar} F(x) \, d\mu(x)
\end{equation*}
for each $F \in \cC(\Ubar)$.

For the final assertion, it suffices to show that 
$(r_{\Gamma})_*(\mu) = \mu_{\Gamma}$ for each $\Gamma$, since this
means that $\mu$ determines the measures $\mu_{\Gamma}$, which in turn
determine the functional $\Lambda(F)$.  

Fix $\Gamma \subset U$ and let $T$ be a connected subset of $\Gamma$.  
Let $\chi_{T}$ be the characteristic function of $T$.  
Choose a sequence of continuous, non-negative functions
$f_1, f_2, \ldots \in \cC(\Gamma)$ decreasing pointwise
to $\chi_T$.  By Lebesgue's monotone convergence theorem,
\begin{equation*}
\lim_{i \rightarrow \infty} \int_{\Gamma} f_i(x) \, d\mu_{\Gamma}(x)
\ = \ \int_{\Gamma} \chi_T(x) \, d\mu_{\Gamma}(x) \ = \ \mu_{\Gamma}(T) \ .
\end{equation*}
The functions $f_i \circ r_{\Ubar,\Gamma}$ decrease pointwise
to $\chi_{T} \circ r_{\Ubar,\Gamma}$ so by another application of
the monotone convergence theorem
\begin{equation*}
\lim_{i \rightarrow \infty}
\int_{\Ubar} f_i \circ r_{\Ubar,\Gamma}(x) \, d\mu(x)
  \ = \ \int_{\Ubar} \chi_T \circ r_{\Ubar,\Gamma}(x) \, d\mu(x)
  \ = \ \mu(r_{\Ubar,\Gamma}^{-1}(T)) \ .
\end{equation*}
Since $\int_{\Ubar} f_i \circ r_{\Ubar,\Gamma}(x) \, d\mu(x)
   = \int_{\Gamma} f_i(x) \, d\mu_{\Gamma}(x)$ for each $i$, 
$\mu(r_{\Ubar,\Gamma}^{-1}(T)) = \mu_{\Gamma}(T)$.  Because connected
sets generate the $\sigma$-algebra of Borel measurable sets on $\Gamma$,
it follows that  $\mu(r_{\Ubar,\Gamma}^{-1}(T)) = \mu_{\Gamma}(T)$ for
all Borel measurable $T \subset \Gamma$.
\end{proof}

\subsection{The Laplacian.}  

In Section~\ref{Section C}, for any metrized graph $\Gamma$ 
we have introduced a Laplacian on the space of continuous, 
piecewise affine functions $\CPA(\Gamma)$.  
Recall that if $f \in \CPA(\Gamma)$ then 
\begin{equation} \label{FFF1}
\Delta(f) \ = \ \sum_{p \in \Gamma} 
         (-\sum_{\text{$\vec{v}$ at $p$}} d_{\vec{v}}f(p)) \, \delta_p(x) \ .
\end{equation}

It is possible to define a Laplacian on larger classes of functions. 
For example, Zhang (\cite{Zh}) defined a Laplacian on the space of 
continuous, piecewise $\cC^2$ functions whose one-sided directional derivatives 
$d_{\vec{v}}f(p)$ exist for all $p \in \Gamma$.  We will write $\Zh(\Gamma)$
for this space.  Zhang put 
\begin{equation} \label{FFF2}
\Delta(f) \ = \ -f^{\prime \prime}(x) \, dx +  \sum_{p \in \Gamma} 
         (-\sum_{\text{$\vec{v}$ at $p$}} d_{\vec{v}}f(p)) \, \delta_p(x) 
\end{equation} 
where $f^{\prime \prime}(x)$ is taken relative to the arclength parametrization,  
on each segment in the complement of an appropriate vertex set for $\Gamma$. 
Again, $\Delta(f)$ is a measure on $\Gamma$.   
The condition that the directional derivatives exist for all $p$ 
is easily seen to be equivalent to requiring that 
$f^{\prime \prime} \in L^1(\Gamma)$.  For a function $f \in \CPA(\Gamma)$,
$f^{\prime \prime} = 0$ in the complement of a vertex set, 
so Zhang's Laplacian is compatible with the one in (\ref{FFF1}).  

However, one can extend the Laplacian to a still broader class of functions.
To motivate this, note that Zhang's Laplacian has the following 
property.

\begin{lemma} \label{LemE6}  {\rm (Mass Formula)} 
Let $K \subset \Gamma$ be a set which is a finite union of points and 
closed intervals.  Then for each $f \in \Zh(\Gamma)$, 
\begin{equation} \label{FFH1} 
\Delta(f)(K) \ = \ - \sum_{p \in \partial K} 
     \sum \begin{Sb} \text{$\vec{v}$ at $p$} \\ \text{outward} \end{Sb} 
                          d_{\vec{v}}f(p) 
\end{equation}
where $\partial K$ is the boundary of $K$, 
and a direction $\vec{v}$ at $p$ is called `outward' if the edge it
corresponds to leads away from $K$.  
\end{lemma}     

\begin{proof}  Choose a vertex set $S$ for $\Gamma$ which includes the
endpoints of all segments in $K$, all isolated points of $K$,
and all points where $f^{\prime \prime}(x)$ is not defined.  
Put $S_{K} = S \cap K$;  then $K \backslash S_K$ is 
a finite union of open intervals $(a_i,b_i)$.  By definition
\begin{equation} \label{FFG1}
\Delta(f)(K) \ = \ - \sum_i \int_{a_i}^{b_i} f^{\prime \prime}(x) \, dx 
    +  \sum_{p \in S_K} 
         (-\sum_{\text{$\vec{v}$ at $p$}} d_{\vec{v}}f(p)) \, \delta_p(x) \ . 
\end{equation}
On the other hand, for each interval $(a_i,b_i)$ 
\begin{equation*}
\int_{a_i}^{b_i} f^{\prime \prime}(x) dx \ = \ 
f^{\prime}(b_i) - f^{\prime}(a_i)
\end{equation*}
where $f^{\prime}(a_i)$ and $f^{\prime}(b_i)$ are one-sided derivatives
at the endpoints.  If  $a_i$, $b_i$ correspond to $p, q \in S_K$, then for 
the direction vectors $\vec{v}$, $\vec{w}$ at $p$, $q$ leading into
$(a_i,b_i)$ we have $d_{\vec{v}}(f)(p) = f^{\prime}(a_i)$, 
$d_{\vec{w}}f(q) = - f^{\prime}(b_i)$.  It follows that 
\begin{equation*}
- \int_{a_i}^{b_i} f^{\prime \prime}(x) dx \ = \ 
d_{\vec{w}}f(q) + d_{\vec{v}}f(p) \ .
\end{equation*}
Thus, in (\ref{FFG1}) the sum over the integrals cancels all the 
directional derivatives for directions leading into $K$.  
What remains is the sum over directional derivatives  
in directions outward from $K$.  
\end{proof}

\vskip .1 in
The Mass Formula says that for a closed set $K$ which is a finite union of 
segments and isolated points, $\Delta(f)(K)$ is determined simply by
knowing the directional derivatives of $f$ at the points of $\partial K$.  
In particular, taking $K = \Gamma$, we see that $\Delta(f)(\Gamma) = 0$
since $\partial \Gamma$ is empty.  Since the complement of an open interval $U$
is a closed set $K$ of the type we have been discussing, 
\begin{equation}  \label{FFH2}
\Delta(f)(U) \ = \ -\delta(f)(K) \ = \ 
 \sum_{p \in \partial U} 
     \sum \begin{Sb} \text{$\vec{v}$ at $p$} \\ \text{inward} \end{Sb} 
                          d_{\vec{v}}f(p)  \ .
\end{equation}
Taking unions of open and closed sets, it follows that 
for any $T \subset \Gamma$ which is a finite union of points and 
intervals (open, closed, or half-open) there is  
a formula for $\Delta(f)(T)$ in terms of the 
directional derivatives $d_{\vec{f}}(f)$  at points of $\partial T$.
(This formula is easy to write out, but messy.)  
The collection of such sets forms a Boolean algebra $\cA(\Gamma)$.    

Conversely, suppose $f : \Gamma \rightarrow \RR$ is any function whose
directional derivatives $d_{\vec{v}}f(p)$ exist for all $p$ and $\vec{v}$.  
Taking Lemma \ref{LemE6} as the starting point, one could hope to 
{\it define} the Laplacian $\Delta(f)$ by (\ref{FFH1}), (\ref{FFH2}).  
This does indeed give a finitely additive set function $\Delta(f)$ 
on the boolean algebra $\cA(\Gamma)$.   
However, $\Delta(f)$ does not necessarily
extend to a measure;  there are pathological examples (similar to Weierstrass's 
functions which are continuous everwhere but differentiable nowhere) 
which oscillate so much that $\Delta(f)$ is not countably additive. 
 
An additional condition is necessary:  we will say that a function 
$f : \Gamma \rightarrow \RR$, whose
directional derivatives $d_{\vec{v}}f(p)$ exist for all $p$ and $\vec{v}$,
is of {\it bounded differential variation} if there is a number $B$ 
such that for any countable collection $\{T_i\}$ of pairwise disjoint 
sets in $\cA(\Gamma)$, 
\begin{equation*}
\sum_{i=1}^{\infty} |\Delta(f)(T_i)| \ < \ B \ .
\end{equation*}
We will write $\BDV(\Gamma)$ for the space of functions of bounded
differential variation.  It is easy to see that 
$\CPA(\Gamma) \subset \Zh(\Gamma) \subset \BDV(\Gamma)$.    

In (\cite{B-R}) it is shown that if $f \in \BDV(\Gamma)$ then the 
finitely additive function $\Delta(f)$ defined on $\cA(\Gamma)$ 
extends to a bounded Borel measure on $\Gamma$, 
and conversely if $\mu$ is a bounded Borel
measure on $\Gamma$ with total mass $0$, then there is an $f \in \BDV(\Gamma)$ 
with $\Delta(f) = \mu$. (Such an $f$ is unique up to addition of a constant.)  
Thus, $\BDV(\Gamma)$ is the largest space of functions on which $\Delta(f)$
can be defined as a bounded measure.  
For future reference, we note the following facts shown in (\cite{B-R}):

\begin{lemma} \label{LemE7}  If $f, g \in \BDV(\Gamma)$, then 

$A)$  $f$ is continuous on $\Gamma$.

$B)$  $\Delta(f)(\Gamma) = 0$.

$C)$  $\Delta(f) \equiv 0$ if and only if $f = C$ is constant.  

$D)$  If $f \in \CPA(\Gamma)$ then $\Delta(f)$ is the measure defined by
$(\ref{FFF1})$;  if $f \in \Zh(\Gamma)$ then $\Delta(f)$ is the measure 
defined by $(\ref{FFF2})$.

$E)$ If $\nu$ is a probability measure on $\Gamma$,
and if $f(x) = \int_{\Gamma} j_z(x,y) \, d\nu(y)$,
then $\Delta(f) = \nu - \delta_z(x)$.

$F)$ $\int_{\Gamma} f \, \Delta(g) = \int_{\Gamma} g \, \Delta(f)$.  
\end{lemma} 

\vskip .1 in
However, we are interested in functions not just on a graph $\Gamma$, but 
on a domain $U$.  If $f|_{\Gamma} \in \BDV(\Gamma)$,
write $\Delta_{\Gamma}(f)$ for $\Delta(f|_{\Gamma})$.

\begin{definition} \label{DefE2} Let $U \subset \PP^1_{\Berk}$ be a domain.
We will say that $f : U \rightarrow \RR \cup \{ \pm \infty\}$
is of bounded differential variation on $U$, and write $f \in \BDV(U)$, if

$A)$ $f|_{\Gamma} \in \BDV(\Gamma)$ for each subgraph $\Gamma \subset U$,
and 

$B)$ there is a constant $B = B(f)$ such that for all \, $\Gamma$
\begin{equation*}
|\Delta_{\Gamma}(f)|(\Gamma) \ \le \ B \ .
\end{equation*}
\end{definition} 

\vskip .1 in
\begin{proposition} \label{PropE8}  If $f \in \BDV(U)$,
the system of measures $\{\Delta_{\Gamma}(f)\}_{\Gamma \subset U}$ is coherent.
\end{proposition}

\begin{proof}  The boundedness condition in the definition of coherence 
is built into the definition of $\BDV(U)$;  it suffices to check
the compatibility under pushforwards.

Let subgraphs $\Gamma_1 \subset \Gamma_2 \subset U$ be given.   
Since $\Gamma_2$ can be obtained by sequentially attaching a finite number of
edges to $\Gamma_1$, it suffices to consider the case where 
$\Gamma_2 = \Gamma_1 \cup T$, and $T$ is a segment attached to $\Gamma_1$ at
a point $p$.  By the definition of the Laplacian on a graph, for any Borel
set $e$ contained in $T \backslash \{p\}$ we have 
$\Delta_{\Gamma_2}(f)(e) = \Delta_T(f)(e)$.  For the point $p$,
the set of direction vectors at $p$ in $\Gamma_2$ is the union of
the corresponding sets in $\Gamma_1$ and $T$.  It follows that 
$\Delta_{\Gamma_2}(f)(\{p\}) =
\Delta_T(f)(\{p\}) + \Delta_{\Gamma_1}(f)(\{p\})$.
Since $(r_{\Gamma_2,\Gamma_1})^{-1}(\{p\}) = T$ and $\Delta_T(f)(T) = 0$,
\begin{equation*}
\Delta_{\Gamma_2}(f)(r_{\Gamma_2,\Gamma_1})^{-1}(\{p\}))
\ = \ \Delta_T(f)(T) + \Delta_{\Gamma_1}(f)(\{p\})
\ = \ \Delta_{\Gamma_1}(f)(\{p\}) \ .
\end{equation*}
Trivially  $\Delta_{\Gamma_2}(f)(e) = \Delta_{\Gamma_1}(f)(e)$
for any Borel set $e$ contained in $\Gamma_1 \backslash \{p\}$.  Hence
$(r_{\Gamma_2,\Gamma_1})_*(\Delta_{\Gamma_2}(f)) = \Delta_{\Gamma_1}(f)$.
\end{proof}

\begin{definition} \label{DefE3} 
If $f \in \BDV(U)$, the Laplacian
\begin{equation*}
\Delta(f) \ = \ \Delta_{\Ubar}(f) 
\end{equation*}
is the measure on $\Ubar$ associated to the coherent system 
$\{\Delta_{\Gamma}(f)\}_{\Gamma \subset U}$ by Proposition \ref{PropE5}, 
characterized by the property that for each subgraph $\Gamma \subset U$
\begin{equation*}
(r_{\Ubar,\Gamma})_*(\Delta_{\Ubar}(f)) \ = \ \Delta_{\Gamma}(f) \ .
\end{equation*}
\end{definition}

\begin{proposition} \label{PropE9}
For each $f \in \BDV(U)$, $\Delta_{\Ubar}(f)$ has total mass $0$.
\end{proposition}

\begin{proof} For any subgraph $\Gamma \subset U$,
\begin{equation*}  
\Delta_{\Ubar}(f)(\Ubar) 
\ = \ (r_{\Ubar,\Gamma})_*(\Delta_{\Ubar}(f))(\Gamma) 
\ = \ \Delta_{\Gamma}(f)(\Gamma) \ = \ 0 \ .
\end{equation*}
\end{proof}

\noindent{Here} are some examples of functions in $\BDV(\PP^1_{\Berk})$
and their Laplacians.

\begin{example}
\label{Example E.1}
If $f(x) = C$ is constant then
$f \in \BDV(\PP^1_{\Berk})$ and $\Delta(f) \equiv 0$.

Indeed $\Delta_{\Gamma}(f) = 0$ for each subgraph $\gamma$, so the
associated functional $\Lambda(F)$ is the $0$ functional.
\end{example}

\begin{example}
\label{Example E.2}
  If $f(x) = -\log_v(\delta(x,y)_{\zeta})$, then $f \in \BDV(\PP^1_{\Berk})$
    and 
\begin{equation*}    
    \Delta(f) \ = \ \delta_y(x) - \delta_{\zeta}(x) \ .
\end{equation*}

    To see this, let $\zeta_0 \in \PP^1_{\Berk}$ be the Gauss point,
the point corresponding to $B(0,1)$.  
Given a subgraph $\Gamma$ containing $\zeta_0$,
put $\tilde{y} = r_{\Gamma}(y)$ and $\tilde{\zeta} = r_{\Gamma}(\zeta)$.
For $x \in \Gamma$, it follows from the definition of the Hsia kernel and the
restriction formula for $j_{\zeta_0}(x,y)$ in Section~\ref{Section C} that
\begin{eqnarray}
-\log_v(\delta(x,y)_{\zeta})
& = & j_{\zeta_0}(x,y) - j_{\zeta_0}(x,\zeta) - j_{\zeta_0}(y,\zeta) \notag \\
& = & j_{\zeta_0}(x,\tilde{y}) 
         - j_{\zeta_0}(x,\tilde{\zeta}) - j_{\zeta_0}(y,\zeta)
         \label{FNX1}
\end{eqnarray}
Hence $\Delta_{\Gamma}(f) = \delta_{\tilde{y}}(x) - \delta_{\tilde{\zeta}}(x) 
= (r_{\Gamma})_*(\delta_y(x)-\delta_{\zeta}(x))$.
The only measure on $\PP^1_{\Berk}$ with this property for all $\Gamma$ 
is $\mu = \delta_y(x) - \delta_{\zeta}(x)$.
\end{example}

\begin{example}
\label{Example E.3}
Take $0 \ne g \in \CC_v(T)$;  suppose $\div(g) = \sum_{i=1}^m n_i(a_i)$.
Then $f(x) = -\log_v([g]_x) \in \BDV(\PP^1_{\Berk})$ and
\begin{equation*}
\Delta(-\log_v([g]_x)) \ = \ \sum_{i=1}^m n_i \delta_{a_i}(x) \ .
\end{equation*}

This follows from Example~\ref{Example E.2}.  Let $\zeta \in \PP^1_{\Berk}$ be arbitrary.
By the decomposition formula for the Hsia kernel, 
there is a constant $C_{\zeta}$ such that
\begin{equation*}
[g]_x \ = \ C_{\zeta} \cdot \prod_{a_i \ne \zeta} \delta(x,a_i)_{\zeta} \ .
\end{equation*}
Taking logarithms, applying Example~\ref{Example E.2}, and using $\sum_{i=1}^m n_i = 0$
gives the result.
\end{example}

\begin{example}
\label{Example E.4}
If $\nu$ is any probability measure on $\PP^1_{\Berk}$ and
$\zeta \notin \supp(\nu)$, consider the potential function
\begin{equation*}
u_{\nu}(x,\zeta) \ = \ \int -\log_v(\delta(x,y)_{\zeta}) \, d\nu(y) \ .
\end{equation*}
Then $u_{\nu}(x,\zeta) \in \BDV(\PP^1_{\Berk})$ and
\begin{equation} \label{FNX2}
\Delta(u_{\nu}(x,\zeta)) \ = \ \nu - \delta_{\zeta}(x) \ .
\end{equation}

To see this, let $\zeta_0$ be the Gauss point, 
and let $\Gamma$ be any subgraph containing $\zeta_0$.
For $x \in \Gamma$, using the definition of $\delta(x,y)_{\zeta}$,
the restriction formula $j_{\zeta_0}(x,y) = j_{\zeta_0}(x,r_{\Gamma}(y))$, 
and the fact that $\nu$ has total mass $1$, we have
\begin{eqnarray*}
u_{\nu}(x,\zeta)
& = & \int j_{\zeta_0}(x,y) - j_{\zeta_0}(x,\zeta) 
         - j_{\zeta_0}(y,\zeta) \, d\nu(y) \\
& = & \int j_{\zeta_0}(x,r_{\Gamma}(y)) \, d\nu(y)
               - \int j_{\zeta_0}(x,r_{\Gamma}(\zeta)) d\nu(y) - C_{\zeta} \\
& = & \int_{\Gamma} j_{\zeta_0}(x,t) \, d((r_{\Gamma})_*(\nu))(t)
               - j_{\zeta_0}(x,r_{\Gamma}(\zeta)) - C_{\zeta}
\end{eqnarray*}
where $C_{\zeta}$ is a finite constant.
By Lemmma \ref{LemE7}(E), $u_{\nu}(x,\zeta)|_{\Gamma} \in \BDV(\Gamma)$ and
\begin{equation*}
\Delta_{\Gamma}(u_{\nu}(z,\zeta))
\ = \ (r_{\Gamma})_*(\nu) - \delta_{r_{\Gamma}(\zeta)}(x) 
\ = \ (r_{\Gamma})_*(\nu - \delta_{\zeta}(x)) \ .
\end{equation*}
It follows that $\Delta(u_{\nu}(x,\zeta)) = \nu - \delta_{\zeta}(x)$.
\end{example}

\begin{remark}
The definition of
$f \in \BDV(U)$ only concerns the restrictions $f|_{\Gamma}$
to subgraphs $\Gamma \subset U$.  Such subgraphs can only contain
points of type II, III, or IV  (in fact, as the subgraphs vary
they exhaust $U \backslash \PP^1(\CC_v)$), so requiring that $f \in \BDV(U)$
imposes no conditions at all on its behavior on $\PP^1(\CC_v)$.
That behavior must be deduced from auxiliary hypotheses, such as continuity
or upper-semicontinuity.
\end{remark}

\vskip .1 in
The Laplacian is compatible with restriction to subdomains.

\begin{proposition} \label{PropE10}  
Suppose $U_1 \subset U_2 \subset \PP^1_{\Berk}$ are domains, and 
$f \in \BDV(U_2)$.  Then $f \in \BDV(U_1)$ and 
\begin{equation*}
\Delta_{\Ubar_1}(f) \ = \ (r_{\Ubar_2,\Ubar_1})_*(\Delta_{\Ubar_2}(f)) \ .
\end{equation*}
\end{proposition}

\begin{proof}
If $f \in \BDV(U_2)$, then trivially $f \in \BDV(U_1)$.  By the compatibility
of the retraction maps, for each $\Gamma \subset U_1$
\begin{equation*}
(r_{\Ubar_1,\Gamma})_*((r_{\Ubar_2,\Ubar_1})_*(\Delta_{\Ubar_2}(f))) 
\ = \ (r_{\Ubar_2,\Gamma})_*(\Delta_{\Ubar_2}(f)) \ = \ \Delta_{\Gamma}(f) \ .
\end{equation*}
By the characterization of $\Delta_{\Ubar_1}(f)$, it follows that
$\Delta_{\Ubar_1}(f) = (r_{\Ubar_2,\Ubar_1})_*(\Delta_{\Ubar_2}(f))$. 
\end{proof}

\subsection{Harmonic Functions.}

If $U$ is a domain and $f \in \BDV(U)$, the measure $\Delta_{\Ubar}(f)$ 
can be decomposed into two parts:  a part supported on $\partial U$, 
and a part supported on $U$. 
As will be seen, the part supported on $\partial U$ is analogous to the
classical normal derivative on the boundary, 
while the part supported on $U$ is analogous to the classical Laplacian.

\begin{definition} \label{DefE4}  If \,$U$ is a domain,
$f:U \rightarrow \RR$ is {\rm strongly harmonic} on $U$ 
if it is continuous on $U$, belongs to $\BDV(U)$, 
and satisfies $\Delta_{\Ubar}(f)|_U \equiv 0$.

If \,$V$ is an arbitrary open set, $f : V \rightarrow \RR$
is {\rm harmonic} in $V$ if for each $x \in V$ there is a 
$($connected$)$ neighborhood $U_x$ of $x$ on which $f$ is strongly harmonic.  
\end{definition}

If $f$ and $g$ are harmonic (or strongly harmonic) in $U$, 
then so are $-f$ and $a \cdot f + b \cdot g$, for any $a, b \in \RR$.  
In view of the computations of Laplacians given above, we have the following
examples of strongly harmonic functions:  

\begin{example}
\label{Example E.5}
Each constant function $f(x) = C$ on a domain $U$ is strongly harmonic.
\end{example}

\begin{example}
\label{Example E.6}
For fixed $y, \zeta$, $-\log_v(\delta(x,y)_{\zeta})$
is strongly harmonic in $\PP^1_{\Berk} \backslash \{y, \zeta\}$.
\end{example}

\begin{example}
\label{Example E.7}
If $0 \ne g \in \CC_v(T)$ has divisor
$\div(g) = \sum_{i=1}^m n_i(a_i)$, then $f(x) = -\log_v([g]_x)$
is strongly harmonic in $\PP^1_{\Berk} \backslash \{a_1, \ldots, a_m\}$.
\end{example}

\begin{example}
\label{Example E.8}
Given a probability measure $\nu$ and a point $\zeta \notin \supp(\nu)$,
the potential function $u_{\nu}(z,\zeta)$ is strongly harmonic in 
$\PP^1_{\Berk} \backslash (\supp(\nu) \cup \{\zeta\})$.
\end{example}

\begin{lemma} \label{LemE11}
\ \ 

$A)$  If $U_1 \subset U_2$ are domains, and $f$ is strongly harmonic in $U_2$,
then $f$ is strongly harmonic in $U_1$.  

$B)$  If $f$ is harmonic in an open set $V$, and $U$ is a subdomain of $V$ 
with $\Ubar \subset V$, then $f$ is strongly harmonic in $U$. 

$C)$  If $f$ is harmonic on $V$ and $K \subset V$ is compact and connected,
there is a subdomain $U \subset V$ containing $K$ such that $f$ is
strongly harmonic on $U$.  In particular, if $V$ is a domain, there is
an exhaustion of $V$ by subdomains  $U$ on which $f$ is strongly harmonic.
\end{lemma}

\begin{proof}  For (A) note that $f \in \BDV(U_1)$ and 
$(r_{\Ubar_2,\Ubar_1})_*(\Delta_{\Ubar_2}(f)) = \Delta_{\Ubar_1}(f)$
by Proposition \ref{PropE10}.  Since $f_{\Ubar_2,\Ubar_1}$ fixes $U_1$ and
maps $\Ubar_2 \backslash U_1$ to $\partial U_1$,
if $\Delta_{\Ubar_2}(f)|_{U_2} = 0$, 
then also $\Delta_{\Ubar_1}(f)|_{\Ubar_1} = 0$.

For (B), note that for each $x \in \Ubar$ there is a connected 
neighborhood $U_x$ such that $f$ is strongly harmonic in $U_x$.  
Since $\Ubar$ is compact, a finite number of the $U_x$ cover $\Ubar$,
say $U_{x_1}, \ldots, U_{x_m}$.  For each $U_{x_i}$, there is a constant $B_i$
such that for every subgraph $\Gamma_i \subset U_{x_i}$, 
\begin{equation*}
|\Delta_{\Gamma_i}(f)|(\Gamma_i) \ \le \ B_i \ .
\end{equation*}
If $\Gamma \subset U$ is a subgraph, then there are subgraphs 
$\Gamma_1, \ldots, \Gamma_m$ with $\Gamma_i \subset U_{x_i}$ such that
$\Gamma \subset \bigcup_{i=1}^m \Gamma_i$. 
(These subgraphs are not in general disjoint). 
 We have $f|_{\Gamma} \in \BDV(\Gamma)$ 
 since $f|_{\Gamma_i} \in \BDV(\Gamma_i)$ for each $i$, and 
\begin{equation*}
|\Delta_{\Gamma}(f)|(\Gamma) 
\ \le \ \sum_{i=1}^m |\Delta_{\Gamma_i}(f)|(\Gamma_i) 
\ \le \ \sum_{i=1}^m B_i \ .
\end{equation*}
Thus $f|_U \in \BDV(U)$.  

For each $x_i$, $U \cap U_{x_i}$ is a subdomain of $U_{x_i}$.  
(It is connected because there is a unique path between any two points 
of $\PP^1_{\Berk}$, and $U$, $U_{x_i}$ are both path-connected.)  
By part (A), 
\begin{equation*}
\Delta_{\Ubar \cap \Ubar_{x_i}}(f)|_{U \cap U_{x_i}} 
\ \equiv \ 0 \ . 
\end{equation*}
However, the $U \cap U_{x_i}$ cover $U$, and 
\begin{equation*}
(r_{\Ubar, \Ubar \cap \Ubar_{x_i}})_*(\Delta_{\Ubar}(f))
\ = \ \Delta_{\Ubar \cap \Ubar_{x_i}}(f) \ .
\end{equation*}
It follows that $\Delta_{\Ubar}(f)|_{U \cap U_{x_i}} \equiv 0$.  This holds 
for each $i$, so $\Delta_{\Ubar}(f)|_{U} \equiv 0$.

For (C), note that for each $x \in K$ there is a connected neighborhood
$U_x \subset V$ on which $f$ is strongly harmonic.  Since $K$ is compact,
a finite number of the $U_x$ cover $\Gamma$.  Let $U$ be the union of those
neighborhoods.  Then $U$ is connected, and by the same argument as above, 
$f$ is strongly harmonic on $U$.
\end{proof}

\vskip .1 in
Note that if $U$ is open but is not connected, and $f$ is harmonic in $U$, 
we have not defined $\Delta_{\Ubar}(f)$ even if $f$ is strongly harmonic
in each component of $U$.
Different components may share boundary points, creating unwanted
interaction between the Laplacians on the components.  For example,
let $\zeta_0 \in \PP^1_{\Berk}$ be the Gauss point,
and consider the open set $U = \PP^1_{\Berk} \backslash \{0,\zeta_0\}$.
It has one distinguished component $U_0 = \cB(0,1)^- \backslash \{0\}$
and infinitely many other components $U_i = \cB(a_i,1)_0^-$.  They all 
have $\zeta_0$ as a boundary point.  The function $f(z)$ which is
$-\log_v(\delta(x,0)_{\infty})$ in $U_0$ and is $0$ elsewhere has Laplacian
$\Delta_{\Ubar_0}(f) = \delta_0(x) - \delta_{\zeta_0}(x)$, while 
$\Delta_{\Ubar_i}(f) = 0$ for $i \ne 0$.

\vskip .1 in
A key observation is that the behavior of a harmonic function on a domain 
$U$ is controlled by its behavior on a special subset.    

\begin{definition} \label{DefE5}
If $U$ is a domain, the {\rm main dendrite} $D \subset U$ is set of all 
$x \in U$ belonging to paths between boundary points $y, z \in \partial U$.  
\end{definition}

If the main dendrite is nonempty, 
it is a tree which is finitely branched at each point.

To see this, fix $\zeta \in U$, and take $\varepsilon > 0$.  
For each $x \in \partial U$, put $r_x = \diam_{\zeta}(x)+\varepsilon$.
Since $\partial U = \Ubar \backslash U$ is compact, it can be covered by
finitely many open balls $\cB(x,r_x)_{\zeta}^-$, 
say the balls corresponding to 
$x_1, \ldots, x_m$.  Without loss we can assume these balls are pairwise
disjoint.  Each $\cB(x_i,r_{x_i})_{\zeta}^-$ has a single
boundary point $z_i$, which necessarily belongs to $U$.  Let $D_{\varepsilon}$
be the subgraph connecting the points $z_1, \ldots, z_m$;  it is a finite
tree, contained in $U$.  If $y, z \in \partial U$ belong to the same ball 
$B(x_i,r_{x_i})_{\zeta}^-$ then the path between them lies entirely in 
that ball, because a subset of $\PP^1_{\Berk}$ is connected if and only if
it is path connected, and there is a unique path between any two points.  
If $y$ and $z$  belong to distinct balls 
$B(x_i,r_{x_i})_{\zeta}^-$, $B(x_j,r_{x_j})_{\zeta}^-$, the path between 
them is contained in those balls together with $D_{\varepsilon}$.  
If we let $\varepsilon \rightarrow 0$, the balls $B(x_i,r_{x_i})_{\zeta}^-$
will eventually omit any given point $x \in U$.  
If $\varepsilon_1 > \varepsilon_2$, then $D_{\varepsilon_1}$ is 
contained in the interior of $D_{\varepsilon_2}$,
and as $\varepsilon \rightarrow 0$ 
the subgraphs $D_{\varepsilon}$ exhaust $D$. 
Thus
\begin{equation*}
D \ = \ \bigcup_{n=1}^{\infty} D_{1/n}
\end{equation*}
is a tree, which is finitely branched at every point.   

Note that the main dendrite can be empty:  
that will be the case iff $U$ has one or no boundary points, 
namely if $U$ is an open disc, if $U \cong \AA^1_{\Berk}$, 
or if $U = \PP^1_{\Berk}$.  

\begin{lemma}  \label{LemE12}
Let $f(x)$ be harmonic in a domain $U$.  
If the main dendrite is empty, then $f(x)$ is constant;  
otherwise, $f(x)$ is constant on branches off the main dendrite.  
\end{lemma}

\begin{proof}  First suppose the main dendrite is empty, 
so $\partial U$ is empty or consists of a single point.  
Then $U = \PP^1_{\Berk}$, $U \cong \AA^1_{\Berk}$, or
$U$ is an open disc.  Let $x, y \in U$ be arbitrary points not of type I, 
and let $\Gamma$ be the path connecting them.  By the description of 
$U$, there is a disc $V$ whose closure is contained in $U$, with
$\Gamma \subset V$.  

By Lemma \ref{LemE11}, $f$ is strongly harmonic in $V$.
We claim that $\Delta_{\Vbar}(f) \equiv 0$.  
By definition, $\Delta_{\Vbar}(f)|_V \equiv 0$.  
Also, $\Delta_{\Vbar}(f)(\Vbar) = 0$ by Proposition \ref{PropE9}, so
$\Delta_{\Vbar}(f)(\partial V) = 
\Delta_{\Vbar}(f)(\Vbar)-\Delta_{\Vbar}(f)(V) = 0$.  Since $\partial V$
consists of one point, our claim is established.  

It follows that 
$\Delta_{\Gamma}(f) = (r_{\Vbar,\Gamma})_*(\Delta_{\Vbar}(f)) = 0$.
By Lemma \ref{LemE7}(C), $f|_{\Gamma}$ is constant;  thus $f(x) = f(y)$.  
Since the non-type I points are dense in $U$ and $f$ is continuous, 
$f$ must be constant.
\vskip .05 in

Next suppose the main dendrite is nonempty.  Let $x \in U \backslash D$ be
arbitrary.  If $x$ is not of type I, fix a point $y_0 \in D$, 
and let $w$ be the first point in $D$ along the path from $x$ to $y_0$;  
we claim that $f(x) = f(w)$.   
Let $\Gamma$ be the path from $x$ to $w$;  
then there is a subdomain $V$ of $U$ containing $\Gamma$, 
on which $f$ is strongly harmonic.  
Since $\Delta_{\Vbar}(f)$ is supported on $\partial V$
and the retraction of $\partial V$ to $\Gamma$ is the point $w$,
it follows that 
$\Delta_{\Gamma}(f) = (r_{\Vbar,\Gamma})_*(\Delta_{\Vbar}(f))$ 
is supported on $w$.  However, $\Delta_{\Gamma}(f)(\Gamma) = 0$, 
so $\Delta_{\Gamma}(\{w\}) = 0$.
Thus $\Delta_{\Gamma}(f) \equiv 0$, and by Lemma \ref{LemE7}(C),  
$f(x) = f(w)$.  If $x$ is of type I, then it is a limit of non-type I points
along the path from $x$ to $w$, so by continuity $f(x) = f(w)$.  
\end{proof}

\vskip .1 in
Here is an example of a domain $U$ and a function $f : U \rightarrow \RR$
which is harmonic on $U$ but not strongly harmonic.  

Fix coordinates so that $\PP^1_{\Berk} = \AA^1_{\Berk} \cup \{\infty\}$,
and take $U = \PP^1_{\Berk} \backslash \ZZ_p$.  Then the main dendrite 
$D$ is a tree whose root node is the Gauss point $\zeta_0$, 
and whose other nodes are at the points corresponding to balls 
$B(a,p^{-n})$ for $a, n \in \ZZ$, with $n \ge 1$ and $0 \le a \le p^n-1$.  
It has $p$ branches extending down from each node,  
and each edge has length $1$.  

To describe $f$, it suffices to give its values on the main dendrite.  
It will have constant slope on each edge.  Let $f(\zeta_0) = 0$;    
we will recursively give its values on the other nodes.  
Suppose $z_a$ is a node on which $f(z_a)$ has already been defined.  
If $z_a$ is the root node, put $N_a = 0$;  
otherwise, let $N_a$ be the slope of $f$
on the edge entering $z_a$ from above.  Of the $p$ edges extending
down from $z_a$, choose two distinguished ones, and let $f(z)$ have slope
$N_a + 1$ on one and slope $-1$ on the other, until the next node. On the
$p-2$ other edges, let $f(z)$ have the constant value $f(z_a)$ until the
next node.  

By construction, the sum of the slopes of $f$ on the edges leading away 
from each node is $0$, so $f$ is harmonic in $U$.  
However, there are edges of $D$ on which $f$ has arbitrarily large slope;  
if $\Gamma$ is an edge of $D$  on which $f$ has slope $m_{\Gamma}$, 
then $|\Delta_{\Gamma}(f)|(\Gamma) = 2 |m_{\Gamma}|$.  
Thus, $f \notin \BDV(U)$, and $f$ is not strongly harmonic on $U$.  

\subsection{The Maximum Principle.}

We will now show that harmonic functions on Berkovich space share many  
properties of classical harmonic functions.

\begin{proposition} \label{PropE13} {\rm (Maximum Principle)} 
 
If $f$ is a nonconstant harmonic function on a domain 
$U \subset \PP^1_{\Berk}$, then $f$ does not achieve a maximum or a
minimum value on $U$.  

In particular, if $f$ is harmonic on $U$ 
and $\limsup_{x \rightarrow \partial U} f(x) \le M$,  
then $f(x) \le M$ for all $x \in U$;  
if $\liminf_{x \rightarrow \partial U} f(x) \ge m$,
then $f(x) \ge m$ for all $x \in U$.
\end{proposition} 

\begin{proof}  It suffices to deal with the case of a maximum.    

Suppose  $f$ takes on its maximum at a point $x \in U$.  Let $V \subset U$
be a subdomain containing $x$ on which $f$ is strongly subharmonic.  
We will show that $f$ is constant on $V$.  It follows that 
$\{z \in U : f(z) = f(x) \}$ is both open and closed, and hence equals $U$
since $U$ is connected, so $f$ is constant on $U$.  

By Lemma \ref{LemE12}, we can assume the main dendrite $D$ for $V$ 
is nonempty.  Let $T$ be the
branch off the main dendrite containing $x$, and let $w$ be the point
where $T$ attaches to $D$.  By Lemma \ref{LemE12}, $f(w) = f(x)$.   

Let $\Gamma \subset D$ be an arbitrary subgraph with $w$ in its interior.   
By the definition of the main dendrite, $r_{\Vbar,\Gamma}(\partial V)$ 
consists of the endpoints of $\Gamma$.   
Since $\Delta_{\Vbar}(f)$ is supported on $\partial V$, 
$\Delta_{\Gamma}(f) = (r_{\Vbar,\Gamma})_*(\Delta_{\Vbar}(f))$
is supported on the endpoints of $\Gamma$.  By Lemma \ref{LemE7}(E), 
$f|_{\Gamma} \in \CPA(\Gamma)$.  Let $\Gamma_0$ be the connected component of 
$\{ z \in \Gamma : f(z) = f(w) \}$  containing $w$.
If $\Gamma_0 \ne \Gamma$, let $p$ be a boundary point of $\Gamma_0$.  
Then $p$ cannot be an endpoint of $\Gamma$, because $\Gamma_0$ is closed. 
Since $f(p)$ is maximal, we must have $d_{\vec{v}}f(p) \le 0$ for all $\vec{v}$
at $p$.  Since $f|_{\Gamma} \in \CPA(\Gamma)$, 
and $f|_{\Gamma}$ is nonconstant near $p$,
there must be some $\vec{v}$ with $d_{\vec{v}}(f)(p) < 0$.  It follows that 
\begin{equation*}
\Delta_{\Gamma}(f)(p) 
\ = \ - \sum_{\text{$\vec{v}$ at $p$}} d_{\vec{v}}f(p) > 0 \ .
\end{equation*}
This contradicts the fact that $\Delta_{\Gamma}(f)$ is supported on the
endpoints of $\Gamma$, and we conclude that $\Gamma_0 = \Gamma$.
    
Since $\Gamma$ can be taken arbitrarily large, $f$ must be constant on the 
main dendrite.  By Lemma \ref{LemE12}, it is constant everywhere on $V$.
\end{proof}  

\vskip .1 in
There is an important strengthening of the Maximum Principle  
which allows us to ignore sets of capacity $0$ in $\partial U$.

\begin{proposition} \label{PropE14}   
{\rm (Strong Maximum Principle)} 

Let $f$ be harmonic on a domain $U \subset \PP^1_{\Berk}$.
Assume either that $f$ is nonconstant, 
or that $\partial U$ has positive capacity.
 
If $f$ is bounded above on $U$, and 
there are a number $M$ and a set \, $e \subset \partial U$ of capacity $0$ 
such that $\limsup_{x \rightarrow z} f(x) \le M$ for all 
$z \in \partial U \backslash e$, then $f(x) \le M$ for all $x \in U$.

Similarly, if $f$ is bounded below on $U$,
and there are a number $m$ and a set \, $e \subset \partial U$ of capacity $0$ 
such that $\liminf_{x \rightarrow z} f(x) \ge m$ for all 
$z \in \partial U \backslash e$, then $f(x) \ge m$ for all $x \in U$.
\end{proposition}

\vskip .1 in
For the proof, we will need  
the existence of an `Evans function'.

\begin{lemma} \label{LemE15}  Let $e \subset \PP^1_{\Berk}$ be a
compact set of capacity $0$, and let $\zeta \notin e$.  Then there is 
a probability measure $\nu$ supported on $e$ for which the potential function
$u_{\nu}(z,\zeta)$ $($which is defined on all of $\PP^1_{\Berk}$ as a function
to $\RR \cup \{\pm \infty\}$$)$ satisfies
\begin{equation*}
\lim_{z \rightarrow x} u_{\nu}(z,\zeta) \ = \ \infty
\end{equation*}
for all $x \in e$.  A function with this property is called an Evans function.
\end{lemma}

\begin{proof}  By Theorem~\ref{ThmD15},
the restricted Chebyshev 
constant  $\CH^*(E)_{\zeta}$ is $0$.  
This means that for each $\epsilon > 0$, there 
are points $a_1, \ldots, a_N \in e$ for which the pseudo-polynomial 
\begin{equation*}
P_N(z;a_1,\ldots,a_N) \ = \  \prod_{i=1}^N \delta(z,a_i)_{\zeta}
\end{equation*}
satisfies $(\|P_N(z;a_1,\ldots,a_N)\|_e)^{1/N} < \varepsilon$.  

Let $q = q_v$ be the base of the logarithm $\log_v(t)$, and 
take $\varepsilon = 1/q^{2^k}$ for $k = 1, 2, 3, \ldots$.  
For each $k$ we obtain points $a_{k,1}, \ldots, a_{k,N_k}$ 
such that the corresponding
pseudopolynomial $P_{N_k}(z)$ satisfies $(\|P_{N_k}\|_e)^{1/N_k} < 1/q^{2^k}$.   
Put
\begin{equation*}
\nu \ = \ \sum_{k=1}^{\infty} \frac{1}{2^k} 
            (\frac{1}{N_k} \sum_{i=1}^{N_k} \delta_{a_{k,i}}(x) ) \ .
\end{equation*}
Then $\nu$ is a probability measure supported on $e$, and 
\begin{equation*}
u_{\nu}(z,\zeta) \ = \ \sum_{k=1}^{\infty}  \frac{1}{2^k} 
            (- \frac{1}{N_k} \log_v(P_{N_k}(z) ))  \ .
\end{equation*}

Since $P_{N_k}(z)$ is continuous, there is a neighborhood $U_k$ where 
$P_{N_k}(z)^{1/N_k} < 1/q^{2^k}$.  Without loss we can assume 
$U_1 \supset U_2 \supset U_3 \ldots$. 
On $U_m$, for each $k = 1, \ldots, m$ we have 
\begin{equation*}
-\frac{1}{2^k} \cdot \frac{1}{N_k} \log_v(P_{N_k}(z)) 
\ \ge \ - \frac{1}{2^k} \log_v(1/q^{2^k}) \ = \ 1 \ .
\end{equation*}
Hence $u_{\nu}(z,\zeta) \ge m$ on $U_m$.  Letting $m \rightarrow \infty$,
we see that $\lim_{z \rightarrow x} u_{\nu}(z,\zeta) = \infty$ 
for each $x \in e$.
\end{proof}

\vskip .1 in                    
We can now prove Proposition \ref{PropE14}.  
\vskip .1 in

\begin{proof} 
It suffices to deal with the case of a maximum.
Suppose $f(x_0) > M$ for some $x_0 \in U$.  If $f$ is nonconstant, there
is a point $\zeta \in U$ with $f(\zeta) \ne f(x_0)$, and after 
interchanging $x_0$ and $\zeta$ if necessary we can assume $f(x_0) > f(\zeta)$.
If $\partial U$ has positive 
capacity, then for each $z \in \partial U \backslash e$ we have  
$\limsup_{x \rightarrow z} f(x) \le M$, so there is a point $\zeta \in U$
with $f(\zeta) < f(x_0)$.   

Fix $M_1 > M$ with $f(x_0) > M_1 > f(\zeta)$, 
and let $W = \{x \in U : f(x) > M_1\}$.
Let $V$ be the connected component of $W$ containing $x_0$.  
Since $f$ is continuous, $V$ is open, and is itself a 
domain in $\PP^1_{\Berk}$.  Consider its closure $\Vbar$. 
If $\Vbar \subset U$, then each $y \in \partial V$ would have
a neighborhood containing points of $V$ and points of $U \backslash W$.  
By the continuity of $f$, this implies $f(y) = M_1$.
However, $f(x_0) > M_1$.  This violates the Maximum Modulus
Principle for $V$.  

Hence $e_1 = \Vbar \cap \partial U$ is nonempty;  clearly it is compact.  
By the definition of $W$, $e_1$ is contained in $e$, the exceptional set where 
$\limsup_{x \rightarrow y} f(x) > M$.  
Since $e$ has capacity $0$, $e_1$ has capacity $0$ as well.  

Let $h(x)$ be an Evans function 
for $e_1$ with respect to $\zeta$.   Then $h(x)$ is harmonic in $V$, 
and $\lim_{x \rightarrow y} h(y) = \infty$ for each $y \in e_1$.
Since $\zeta \notin \Vbar$, $h(x)$ is bounded below on $\Vbar$.  
Let $-B$ be a lower bound for $h(x)$ on $V$.  

For each $\eta > 0$, put
\begin{equation*}
f_{\eta}(x) \ = \ f(x) - \eta \cdot h(x) 
\end{equation*}
on $V$.  Then $f_{\eta}(x)$ is harmonic in $V$.  Since $f(x)$ is bounded above, 
for each $y \in e_1$
\begin{equation*}
\limsup_{x \rightarrow y} f_{\eta}(y) \ = \ -\infty \ .
\end{equation*}
On the other hand, for each $y \in \partial V \cap U$, 
$f$ is continuous at $y$ and $f(y) = M_1$, as before.  Hence
\begin{equation*}
\limsup_{x \rightarrow y}  f_{\eta}(x) \ \le \ M_1 + \eta \cdot B \ .
\end{equation*} 
Since each $y \in \partial V$ either belongs to $e_1$ or to $\partial V \cap U$,
the Maximum Modulus principle shows that 
$f_{\eta}(x) \le M_1 + \eta \cdot B$ on $V$, hence that
\begin{equation*}
f(x) \ \le \ M_1 + \eta \cdot (B + h(x)) \ .  
\end{equation*}
Fixing $x$ and letting 
$\eta \rightarrow 0$, we see that $f(x) \le M_1$ for each $x \in V$.  
However, this contradicts $f(x_0) > M_1$.  

Thus, $f(x) \le M$ for all $x \in U$, as was to be shown.
\end{proof}  

\begin{proposition} \label{PropE16} {\rm (Riemann Extension Theorem)}
Let $U$ be a domain, and let $e \subset U$ be a compact set of 
capacity $0$.  Suppose $f(x)$ is harmonic and bounded in $U \backslash e$.
Then $f$ extends to a harmonic function on $U$.
\end{proposition}

\begin{proof}  Note that $U \backslash e$ is indeed a domain:  it is open
since $e$ is closed, and it is connected since a set of capacity $0$ is
necessarily contained in $\PP^1(\CC_v)$, and removing Type I points cannot 
disconnect any connected set.

We claim that $f(x)$ is locally constant in a neighborhood of each point 
$a \in e$.  Given $a$, take a ball $\cB(a,r)_0^-$ small enough that its
closure is contained in $U$, and consider the restriction of $f$ to the  
domain $V = \cB(a,r)_0^- \backslash e$.  The boundary $\partial V$ consists
of $\cB(a,r)_0^- \cap e$ together with a point $z \in U$ with $\diam(z) = r$.  
Since $f$ is continuous at $z$, it follows that 
\begin{equation*}
\lim \begin{Sb} x \rightarrow z \\ x \in V \end{Sb} f(x) \ = \ f(z) \ .
\end{equation*}
Applying Proposition \ref{PropE14} to $f$ on $V$, we see that for each
$x \in V$, both $f(x) \le f(z)$ and $f(x) \ge f(z)$.   Thus $f(x) = f(z)$ 
for all $z \in V$.  Extend $f$ by putting $f(x) = f(z)$
for all $x \in \cB(a,r)_0^-$.
Clearly the extended function is harmonic on $U$.
\end{proof}

\vskip .1 in
\begin{corollary} \label{CorE17}
Let $\{a_1, \ldots, a_m\}$ be a finite set of points in $\PP^1(\CC_v)$,
or in $\PP^1(\CC_v) \cap \cB(a,r)_{\zeta}^-$ for some disc.  Then the only
bounded harmonic functions in $\PP^1_{\Berk} \backslash \{a_1, \ldots, a_m\}$,
or in $\cB(a,r)_{\zeta}^- \backslash \{a_1,\ldots,a_m\}$, are the constant functions.
\end{corollary}

\begin{proof}
We have seen in Lemma \ref{LemE12} 
that if $U$ is  $\PP^1_{\Berk}$ or an open disc, 
the only harmonic functions on $U$ are the constant functions.  
Since a finite set of type I points has capacity $0$, 
Proposition \ref{PropE16} shows that a bounded harmonic function 
on $U \backslash \{a_1, \ldots, a_m\}$ extends to a function harmonic on $U$. 
\end{proof} 

\subsection{The Poisson Formula.}

In the classical theory over $\CC$, every point has a neighborhood
isomorphic to a disc.  If $f$ is harmonic on the
disc and has a continuous extension to its closure,
the Poisson Formula gives $f$ on the disc 
in terms of its values on the boundary.

\vskip .1 in
In Berkovich space, the basic open neigbhorhoods are punctured discs.
A punctured disc has only a finite number of boundary points,
and its main dendrite is the interior of a graph $\Gamma$.
We will now show that every harmonic function $f$ on a punctured
disc has a continuous extension to its closure, and 
give a formula for $f(x)$ in terms of its values on the boundary.
In one sense, the Berkovich Poisson formula is simpler than the classical one:
it is a finite sum, obtained using linear algebra.
However, in another sense, it is more complicated,
because it depends on the structure of the graph $\Gamma$.

\vskip .1 in
Consider a punctured disc
$U = \cB(a_1,r_1)_{\zeta}^- \backslash \cup_{i=2}^m \cB(a_i,r_i)_{\zeta}$.
Assume $r_1 < \diam_{\zeta}(\zeta)$, and each $r_i > 0$,
so that each disc in the representation of $U$
corresponds to a boundary point $x_i \in \partial U$,
namely the point corresponding to the classical disc $B(a_i,r_i)_{\zeta}$.

If $m = 1$, $U$ is simply the open disc $\cB(a_1,r_1)_{\zeta}^-$,
and $\partial U$ is the point $x_1$ corresponding to $B(a_1,r_1)_{\zeta}$.
By Lemma \ref{LemE12} the only harmonic functions
on $\cB(a_1,r_1)_{\zeta}^-$ are constant functions, and such functions
clearly extend continuously to $\cB(a_1,r_1)_{\zeta}$.  
Thus if $f(x_1) = A_1$, the Poisson formula is trivially 
\begin{equation} \label{FLM1}
f(x) \ = \ A_1  \qquad \text{for all $x \in \cB(a_1,r_1)_{\zeta}^-$.}
\end{equation}

Now assume $m \ge 2$.
Let $\Gamma$ be the tree spanned by $\partial U = \{x_1, \ldots, x_m\}$,   
so $\Gamma_0 := \Gamma \backslash \{x_1,  \ldots, x_m\}$  
is the main dendrite of $U$.
If $f(x)$ is harmonic in $U$, 
$\Delta_{\Gamma}(f)(p) = -\sum_{\vec{v}} d_{\vec{v}}f(p) \ = 0$
for each $p \in \Gamma_0$.
In particular, the restriction of $f(x)$ to each
edge of $\Gamma_0$ is affine.  It follows that $f(x)$ has
a continuous extension to $\Ubar$.

Fix $z \in \PP^1_{\Berk}$.  We will give a formula
for $f(x)$ on $U$ in terms of the values $A_i = f(x_i)$, 
using the Hsia kernel $\delta(x,y)_z$.
By Lemma \ref{LemE12}, $f(x)$ is determined by its restriction to $\Gamma$.
Fix $z_1 \in \Gamma$ and let
$\zbar = r_{\Gamma}(z)$ be the retraction of $z$ to $\Gamma$.
As shown in Section~\ref{Section C},
there are constants $C, C_z$ such that for all $x, y \in \Gamma$,
\begin{eqnarray}
-\log_v(\delta(x,y)_z)
& = & j_{z_1}(x,y) - j_{z_1}(x,z) - j_{z_1}(y,z) - C \notag \\
& = & (j_{z_1}(x,y) - j_{z_1}(x,\zbar) - j_{z_1}(y,\zbar) 
                     + j_{z_1}(\zbar,\zbar)) \notag \\
&  & \qquad - (C + j_{z_1}(\zbar,\zbar)) \notag \\
& = & j_{\zbar}(x,y) - C_z \label{FLM2}
\end{eqnarray}
where $C_z = C + j_{z_1}(\zbar,\zbar)$.

Consider functions on $\Gamma$ of the form
\begin{equation*}
f(x) \ = \ f_{\vec{c}}(x) \ = \ c_0 + \sum_{i=1}^m c_i j_{\zbar}(x,x_i)
\end{equation*}
where $\sum_{i=1}^m c_i = 0$.  If $f(x) \equiv 0$, then
\begin{equation*}
0 \ = \ \Delta_{\Gamma}(f) \ = \ \sum_{i=1}^m c_i \delta_{x_i}(x) \ ,
\end{equation*}
so $c_1 = \cdots = c_m = 0$, and then $0 \equiv f(x) = c_0$ so $c_0 = 0$
as well.  Let
$L = \{(c_0,c_1, \ldots, c_m) \in \RR^{m+1} : \sum_{i=1}^m c_i = 0\}$.
It follows that the map $F : L \rightarrow \RR^m$ given by
\begin{equation*}
F(\vec{c}) \ = \ (f_{\vec{c}}(x_1), \ldots, f_{\vec{c}}(x_m))
\end{equation*}
is injective.  By a dimension count, it is surjective.
Let $\tM$ be the $(m+1) \times (m+1)$ matrix
\begin{equation} \label{FLM3}
\tM \ = \ \left( \begin{array} {cccc}
0 & 1 & \cdots & 1 \\
1 & j_{\zbar}(x_1,x_1) & \cdots & j_{\zbar}(x_1,x_m) \\
\vdots & \vdots  & \ddots & \vdots \\
1 & j_{\zbar}(x_m,x_1) & \cdots & j_{\zbar}(x_m,x_m)
\end{array} \right)
\end{equation}
Then $\tM$ is nonsingular,
and for each $(A_1, \ldots, A_m) \in \RR^m$,
the numbers $c_0, c_1, \ldots, c_m$ for which $f_{\vec{c}}(x_i) = A_i$
are uniquely determined by the system of equations
\begin{equation} \label{FLM4}
\left( \begin{array}{c} 0 \\ A_1 \\ \vdots \\ A_m \end{array} \right)
\ = \ \tM \cdot
\left( \begin{array}{c} c_0 \\ c_1 \\ \vdots \\ c_m \end{array} \right) \ .
\end{equation}

In the matrix $\tM$, the $j_{\zbar}(x_i,x_j)$ can be replaced by
$-\log_v(\delta(x_i,x_j)_z)$ since
$-\log_v(\delta(x_i,x_j)_z) = j_{\zbar}(x_i,x_j) - C_z$
for all $i, j$ and the first row of $\tM$ asserts that $\sum_{i=1}^m c_i = 0$.
Thus, put
\begin{equation} \label{FLM5}
M(z) \ = \ \left( \begin{array} {ccccc}
0 & 1 & \cdots & 1 \\
1 & -\log_v(\delta(x_1,x_1)_z) & \cdots & -\log_v(\delta(x_1,x_m)_z) \\
\vdots & \vdots & \ddots & \vdots \\
1 & -\log_v(\delta(x_m,x_1)_z) &  \cdots & -\log_v(\delta(x_m,x_m)_z)
\end{array} \right) \ , 
\end{equation}
and for each $i = 0,1, \ldots m$ let $M_i(z,\vec{A})$ be the matrix
gotten by replacing the $i^{th}$ column of $M(z)$ by
$ ^t(0,A_1, \ldots, A_m)$.  The matrix $M(z)$ first appeared in
(Cantor, \cite{Cantor}); we will call it the Cantor matrix.

\begin{proposition} \label{PropE18} {\rm (Poisson Formula I)} \ 
Let $U = \cB(a_1,r_1)_{\zeta}^- \backslash \cup_{i=2}^m \cB(a_i,r_i)_{\zeta}$
be a punctured disc with boundary points $x_1, \ldots, x_m$.
Each harmonic function $f(x)$ on $U$ has continuous extension to $\Ubar$,
and there is a unique such function with prescribed boundary values
$A_1, \ldots, A_m$.  It is given as follows.  

Fix $z \in \PP^1_{\Berk}$, and put $c_i = \det(M_i(z,\vec{A}))/\det(M(z))$ 
for $i = 0, 1 \ldots, m$.  Then 
\begin{equation} \label{FLM6}
f(x) \ = \ c_0 + \sum_{i=1}^m c_i \cdot (-\log_v(\delta(x,x_i)_z)) \ .
\end{equation}
$($This should be understood as a limit, 
if $z$ is of type I and $x = z$.$)$  Moreover
\begin{equation*}
\Delta_{\Ubar}(f) \ = \ \sum_{i=1}^m c_i \delta_{x_i}(x) \ .
\end{equation*}
\end{proposition}

\begin{proof}
By the discussion above, (\ref{FLM6}) holds for all $x \in \Gamma$.
Now let $x \in U$ be arbitrary. 
If $x$ does not belong to the same branch off $\Gamma$ as $z$,
it is easy to see from (\ref{FLM2})
that $-\log_v(\delta(x,x_i)_z) = j_{\zbar}(r_{\Gamma}(x),x_i) - C_z$,
so (\ref{FLM6}) holds for such $x$.
If $x$ belongs to the same branch off $\Gamma$ as $z$,
then $j_{\zbar}(x,x_i)$ is independent of $x_i$.
It follows that $f(x) = c_0$ on that branch, except at possibly when $x = z$.
If $z$ is not of type I, then $f(z) = c_0$ by the continuity of 
$\delta(x,x_i)_z$ at $x = z$.  
If $z$ is of type I, then (\ref{FLM6}) is undefined at $x = z$, 
but $\lim_{t \rightarrow z} f(t) = c_0$.

The assertion regarding $\Delta_{\Ubar}(f)$ follows immediately
from(\ref{FLM6}).
\end{proof}

\vskip .1 in
Note that formula (\ref{FLM6}) holds even when $m = 1$,   
as is easily seen by inspection.

\vskip .1 in
There is another formula for $f$,
which has the merit of being independent of auxiliary variables. 
Since $f$ is unique, the expression (\ref{FLM6}) is valid for each $x$ and $z$.
Taking $x = z$ shows $f(z) = c_0 = \det(M_0(z,\vec{A}))/\det(M(z))$.
Note that by (\ref{FLM2}), (\ref{FLM3}) and (\ref{FLM5}),
the determinants are well-defined and finite,
and $\det(M(z)) \ne 0$, even when $z$ is of type I.

For each $i = 1, \ldots, m$, let $\hat{e}_i \in \RR^{m+1}$ be the vector
which is $1$ in the $(i+1)^{st}$ place and $0$ elsewhere, and put
\begin{equation*}
h_i(z) \ = \ \det(M_0(z,\hat{e}_i))/\det(M(z)) \ .
\end{equation*}
Then $h_i(z)$ is harmonic in $U$ and continuous on $\Ubar$. 
It takes the value $1$ at $x_i$, and is $0$ at each $x_j$ with $j \ne i$.  
It is the analogue of the classical {\it harmonic measure} 
for the boundary component $x_i$ of $U$.
By the maximum principle, $0 \le h_i(z) \le 1$ for each $z \in \Ubar$, and
\begin{equation} \
\sum_{i=1}^m h_i(z) \ \equiv \ 1 \ .
\end{equation}
In sum, we have

\begin{proposition} \label{PropE19} {\rm (Poisson Formula II) \ }
Let $U = \cB(a_1,r_1)_{\zeta}^- \backslash \cup_{i=2}^m \cB(a_i,r_i)_{\zeta}$
be a punctured disc with boundary points $x_1, \ldots, x_m$.
Let $A_1, \ldots, A_m \in \RR$ be given.
Then the unique function $f(z)$ which is harmonic on $U$,
continuous on $\Ubar$, and satisfies $f(x_i) = A_i$ for each 
$i = 1, \ldots, m$ is given by
\begin{eqnarray}
f(z) & = & \det(M_0(z,\vec{A}))/\det(M(z)) \label{FLM7} \\
& = & \sum_{i=1}^m A_i \cdot h_i(z) \ , \label{FLM8}
\end{eqnarray}
for all $z \in \Ubar$,
where $h_i(z) = \det(M_0(z,\hat{e}_i)/\det(M(z))$ is the harmonic
measure for $x_i \in \partial U$.
\end{proposition}

\subsection{Uniform Convergence.}

Poisson's formula implies that the limit of a sequence of harmonic functions 
is harmonic, under a much weaker condition than is required classically.

\begin{proposition} \label{PropE20}
Let $V$ be an open subset of $\PP^1_{\Berk}$.
Suppose $f_1, f_2, \ldots$ are harmonic in $V$ and 
converge pointwise to a function $f : V \rightarrow \RR$.
Then $f(z)$ is harmonic in $V$, 
and the $f_i(z)$ converge uniformly to $f(z)$ on compact subsets of $V$.
\end{proposition}

\begin{proof} Given $x \in V$, take a punctured disc
$U_x$ containing $x$, with closure $\Ubar_x \subset V$.
If $\partial U_x = \{x_1, \ldots, x_m\}$ then by Proposition \ref{PropE19}
\begin{equation*}
f_k(z) \ = \ \sum_{i=1}^m f_k(x_i) h_i(z)
\end{equation*}
for each $z \in \Ubar_x$.  It follows that
the $f_k(z)$ converge uniformly to $f(z)$ on $\Ubar_x$, and
$f(z) = \sum_{i=1} f(x_i) h_i(z)$
is strongly harmonic in $U_x$.  

Thus $f$ is harmonic in $V$.
Any compact $K \subset V$
is covered by finitely many sets $U_x$, so the $f_i(z)$ converge
uniformly to $f(z)$ on $K$.
\end{proof}

\vskip .1 in
Using this, we can characterize harmonic functions as local uniform limits
of logarithms of norms of rational functions.  

\begin{corollary} \label{CorE21}
If $U \subset \PP^1_{\Berk}$ is a domain and $f$ is harmonic in $U$,
there are rational functions $g_1(T), g_2(T), \ldots \in \CC_v(T)$
and numbers $R_1, R_2, \ldots \in \QQ$ such that
\begin{equation*}
f(x) \ = \ \lim_{k \rightarrow \infty} R_k \cdot \log_v([g_k]_x)
\end{equation*}
uniformly on compact subsets of $U$.
\end{corollary}

\begin{proof}
The assertion is trivial if the main denrite of $U$ is empty, since the
only harmonic functions on $U$ are constants.  

Hence we can assume the main dendrite is nonempty.
Choose an exhaustion of $U$ by punctured discs $\{U_k\}$
with closures $\Ubar_k \subset U$.
(For example, such an exhaustion can be gotten by exhausting 
the main dendrite $D$ by subgraphs $\Gamma_k$,
and putting $U_k = r_{\Gamma_k}^{-1}(\Gamma_{k,0})$.)

Fix $k$, and write $\partial U_k = \{x_{k,1}, \ldots, x_{k,m_k}\}$.
Taking $z = \infty$ in Proposition \ref{PropE18}, 
there are numbers $c_{k,i} \in \RR$, with $\sum_{i=1}^{m_k} c_{k,i} = 0$, 
such on $U_k$ 
\begin{equation*}
f(x) \ = \ c_{k,0} + \sum_{i=1}^{m_k}
         c_{k,i} \cdot (-\log_v(\delta(x,x_{k,i})_{\infty})) \ .
\end{equation*}
For each $i = 1, \ldots, m_k$ fix a type I point $a_{k,i}$
whose retraction to $\Ubar_k$ is $x_{k,i}$.
Then $\delta(x,x_{k,i})_{\infty} = \delta(x,a_{k,i})_{\infty} = [T-a_{k,i}]_x$
for each $x \in \Ubar_k$.

Choose rational numbers $d_{k,i}$, with $\sum_{i=1}^{m_k} d_{k,i} = 0$,
which are close enough to the $c_{k,i}$ that
\begin{equation*}
f_k(x) \ := \   d_{k,0} + \sum_{i=1}^{m_k}
         d_{k,i} \cdot (-\log_v(\delta(x,a_{k,i})_{\infty}))
\end{equation*}
satisfies $|f_k(x) - f(x)| < 1/k$ on $\Ubar_k$.
Let $N_k$ be a common denominator for the $d_{k,i}$,
and put $n_{k,i} = N_k \cdot d_{k,i}$.  Let $b_k \in \CC_v$
be a constant with $|b_k|_v = q_v^{-n_{k,0}}$, and put
\begin{equation} \label{FBJ1}
g_k(T) \ = \ b_k \cdot \prod_{i=1}^{m_k} (T-a_{k,i})^{n_{k,i}} \ .
\end{equation}
Then $f_k(x) = (-1/N_k) \log_v([g_k]_x)$ on $U_k$.  The result follows.
\end{proof}

\subsection{Pullbacks.} 

Harmonicity is preserved under pullbacks by meromorphic functions.  
To show this, we first need a lemma asserting that logarithms of norms of
meromorphic functions are harmonic. 

\begin{lemma} \label{LemE22}
Let $U$ be a domain, and suppose $g(x)$ is meromorphic in $U$.
Then $f(x) := \log_v([g]_x)$ is harmonic in $U \backslash \supp(\div(g))$.
\end{lemma}

\begin{proof}
Choose a covering of $U$ by punctured discs $U_k$ with $\Ubar_k \subset U$.
Let $\cA_k$ be a Tate algebra for which $\Ubar_k = \cM(\cA_k)$.
For each $k$, there are coprime polynomials $P_k(z)$ and $Q_k(z)$,
and a unit power series $u_k(z) \in \cA_k$,
such that $g(z) = (P_k(z)/Q_k(z)) \cdot u_k(z)$ on $\Ubar_k \cap \CC_v$.
A unit power series satisfies $[u_k]_x = 1$ for all $x$, so
\begin{equation*}
\log_v([g]_x) \ = \ \log_v([P_k]_x) - \log_v([Q_k]_x)
\end{equation*}
is strongly harmonic on $U_k \backslash (\div(P_k/Q_k))$.
\end{proof}

\vskip .1 in
\begin{corollary} \label{CorE23}
Let $U$ and $V$ be domains in $\PP^1_{\Berk}$.
Suppose $F(z)$ is meromorphic in $U$, with $F(U) \subset V$.
If $f$ is harmonic in $V$, then $f \circ F$ is harmonic in $U$.
\end{corollary}

\begin{proof} Let $x$ be the variable on $U$, and $y$ the variable on $V$.
By Lemma \ref{LemE22} there are rational functions $g_k$ and numbers $R_k$ 
such that
\begin{equation*}
\lim_{k \rightarrow \infty} R_k \cdot \log_v([g_k]_y)  \ = \  f(y) 
\end{equation*}
uniformly on compact subsets of $V$.  For each $y \in V$, let 
$V_y$ be a neighorhood of $y$ with $\Vbar_y \subset V$;  
since $f$ is harmonic there is a $K_y$
so that $g_k$ has no zeros or poles in $\Vbar_y$ if $k \ge K_y$.

Suppose $x \in U$ satisfies $F(x) = y$.  Choose a punctured disc
neighborhood $U_x$ of $x$ contained in $F^{-1}(V_y)$.
The functions $g_k \circ F$ are meromorphic in $U_x$, and if $k \ge K_y$
they have no zeros or poles in $U_x$.  The functions
$R_k \cdot \log_v([g_k \circ F]_z)$
converge to $f \circ F(z)$ for each $z \in U_x$.  

By Corollary \ref{PropE20},
$f \circ F$ is harmonic on $U_x$.  By Proposition \ref{PropE18}, 
$f \circ F$ is strongly harmonic on $U_x$.
Since $x \in U$ is arbitrary, $f \circ F$ is harmonic on $U$.
\end{proof}

\subsection{Harnack's Principle.}   

Harnack's principle holds as well:

\begin{proposition} \label{PropE24}
{\rm (Harnack's Principle)}
Let $U$ be a domain, and suppose $f_1, f_2, \ldots$ are harmonic in $U$
with $0 \le f_1 \le f_2 \le \ldots$.  Then either

$A)$  $\lim_{i \rightarrow \infty} f_i(z) = \infty$ for each $z \in U$, or

$B)$  $f(z) = \lim_{i \rightarrow \infty} f_i(z)$ is finite for all $z$, 
the $f_i(z)$ converge uniformly to $f(z)$ on compact subsets of $U$, 
and $f(z)$ is harmonic in $U$.
\end{proposition}

\begin{proof}  If the main dendrite $D$ for $U$ is empty,
then by Lemma \ref{LemE12} each $f_i(z)$ is constant, and our assertions
are trivial.  Hence we can assume $D$ is nonempty.

Suppose there is some $x \in U$ for which $\lim_{i \rightarrow \infty} f_i(x)$
is finite.  Since each $f_i$ is constant on branches off the main dendrite,
there is a point $x_0$ on $D$ where
$\lim_{i \rightarrow \infty} f_i(x_0)$ is finite. 

Let $\rho(x,y)$ be the logarithmic path distance on $\PP^1_{\Berk}$.
Since the main dendrite is everywhere finitely branched,
for each $p \in D$ there is an $\varepsilon > 0$ such that the closed
neighborhood of $p$ in $D$ defined by
$\Gamma(p,\varepsilon) = \{x \in D : \rho(x,p) \le \varepsilon\}$
is a star, e.g. is a union of $n$ closed segments of length $\varepsilon$
emanating from $p$, for some $n$.
Write $q_i = p + \varepsilon \vec{v}_i$, $i = 1, \ldots, n$, for
the endpoints of $\Gamma(p,\varepsilon)$.

Suppose $h(z)$ is harmonic and non-negative in $U$.
The restriction of $h(z)$ to the subgraph $\Gamma(p,\varepsilon)$
is linear on each of the segments $[p,q_i]$, and $\Delta_\Gamma(h)(p) = 0$.
Thus for each $i$, and $z = p+t \vec{v}_i \in [p,q_i]$
\begin{equation*}
0 \ = \ \Delta_{\Gamma}(h)(p)
\ = \ \frac{h(z)-h(p)}{t} + \sum_{j \ne i} \frac{h(q_j)-h(p)}{\varepsilon} \ .
\end{equation*}
Since $h(q_j) \ge 0$ for each $j$, it follows that $h(x) \le n \cdot h(p)$.
This holds for each $z \in \Gamma(p,\varepsilon)$.
Let $\Gamma \subset D$ be an arbitrary
subgraph containing $x_0$.  Proceeding stepwise from $x_0$, it follows
that there is a constant $C = C_{\Gamma}$ such that for each $x \in \Gamma$,
\begin{equation*}
0 \ \le \ h(x) \ \le C_{\Gamma} \cdot h(x_0) \ .
\end{equation*}

Applying this to the functions $f_i(x)$, since 
$f(x_0) = \lim_{i \rightarrow \infty} f_i(x_0)$ is finite, the
$f_i(x)$ are uniformly bounded on each subgraph $\Gamma \subset D$
containing $x_0$.  Hence the $f_i(x)$ converge uniformly
to $f(x)$ on $\Gamma$.  
If $K$ is any compact subset of $U$, then the image of $K$ under the 
retraction to the main dendrite is contained in some subgraph $\Gamma$.  
Hence the $f_i(x)$ converge locally uniformly to $f(x)$ on $U$, and $f(x)$
is harmonic by Proposition \ref{PropE20}.  
\end{proof}

\subsection{Green's functions.}

We will now prove the uniqueness of the equilibrium distribution,
as promised in Section~\ref{Section D}, 
and use this to build a theory of Green's functions.

\begin{proposition} \label{PropE25}
Let $E \subset \PP^1_{\Berk}$ be a compact set with positive capacity, 
and let $\zeta \in \PP^1_{\Berk} \backslash E$.  Then the equilibrium
distribution $\mu_{\zeta}$ of $E$ with respect to $\zeta$ is unique.
\end{proposition}

\begin{proof}
Suppose $\mu_1$ and $\mu_2$ are two equilibrium distributions
for $E$ with respect to $\zeta$, so that
\begin{equation*}
I_{\zeta}(\mu_1)  \ = \ I_{\zeta}(\mu_2) \ = \ V_{\zeta}(E) \ .
\end{equation*}
Let $u_1(x) = u_{\mu_1}(x,\zeta)$ and $u_2(x) = u_{\mu_2}(x,\zeta)$ be
the corresponding potential functions, and let $U$ be the connected
component of $\PP^1_{\Berk} \backslash E$ containing $\zeta$.

Both $u_1(x)$ and $u_2(x)$
satisfy Frostman's Theorem (Theorem~\ref{ThmD11}):
each is bounded above by $V_{\zeta}(E)$ for all $z$, each is equal to
$V_{\zeta}(E)$ for all $z \in E$ except possibly on an $F_{\sigma}$
set $f_i$ of capacity $0$, and each is continuous on $E$ except on $f_i$.
Each is bounded in a neighborhood of $E$.

Consider $u(z) = u_1(z) - u_2(z)$.  By Example~\ref{Example E.4}, $u(z)$ is in 
$\BDV(\PP^1_{\Berk})$, and
\begin{equation*}
\Delta_{\PP^1_{\Berk}}(u) \ = \
(\mu_{1}-\delta_{\zeta}(x))-(\mu_2 - \delta_{\zeta}(x)) \ = \ \mu_1 - \mu_2 \ .
\end{equation*}
By Proposition~\ref{PropD5}, $\mu_1$ and $\mu_2$ are both supported on
$\partial E_{\zeta} = \partial \Ubar \subset \Ubar$.
Since the retraction map $r_{\PP^1_{\Berk},\Ubar}$ fixes $\Ubar$,
it follows that
\begin{equation} \label{FGM1}
\Delta_{\Ubar}(u) \ = \ (r_{\PP^1_{\Berk},\Ubar})_*(\Delta_{\PP^1_{\Berk}}(u))
\ = \ \mu_1 - \mu_2.
\end{equation}
Thus $\Delta_{\Ubar}(u)$ is supported on $\partial \Ubar$.

Suppose $\zeta \notin \PP^1(\CC_v)$.
By Proposition~\ref{PropD5}, each $u_i(z)$ is defined and bounded on all of
$\PP^1_{\Berk}$, and is continuous on $U$.  Hence $u(z)$ has these
properties as well. By (\ref{FGM1}) $u(z)$ is strongly harmonic on $U$.
Put $f = f_1 \cup f_2$, which has capacity $0$ by Corollary~\ref{CorD14}.
By the discussion above, for each $x \in \partial U \backslash f$,
$u(z)$ is continuous at $x$ and $u(x) = 0$.  In particular
\begin{equation}  \label{FCJ1B}
\lim \begin{Sb} z \rightarrow x \\ z \in U \end{Sb} u(z) \ = \ 0 \ .
\end{equation}
By the Strong Maximum Principle (Proposition \ref{PropE14}), $u(z) \equiv 0$
on $U$.  Hence $\Delta_{\Ubar}(u) = 0$,  so $\mu_1 = \mu_2$.

If $\zeta \in \PP^1(\CC_v)$, there is a slight complication
because $u_1(\zeta) - u_2(\zeta) = \infty - \infty$ is undefined.
However, we claim that $u(z)$ is constant, and in fact is identically $0$,
in a neighborhood of $0$.  Assuming this, since
Proposition~\ref{PropD5} shows that
each $u_i(z)$ is bounded in a neighborhood of $E$ and is continuous in
$U \backslash \{\zeta\}$, it follows that $u(z)$ is harmonic and bounded
in $U \backslash \{\zeta\}$.  The boundary limits (\ref{FCJ1B})
continue to hold for all $x \in \partial U \backslash f$, and in addition
\begin{equation} \label{FCJ2}
\lim \begin{Sb} z \rightarrow \zeta \\ z \in U \end{Sb} u(z) \ = \ 0 \ .
\end{equation}
Applying the Strong Maximum Principle to $u(z)$ on the domain 
$U \backslash \{\zeta\}$, relative to the exceptional set $f \cup \{\zeta\}$
contained in its boundary $\partial U \cup \{\zeta\}$,
we conclude $u(z) \equiv 0$  in $U \backslash \{\zeta\}$.
Since the closure of $U \backslash \{\zeta\}$ is $\Ubar$,
it follows as before that $\Delta_{\Ubar}(u) = 0$,  and $\mu_1 = \mu_2$.

To establish our claim, fix $a \in \PP^1(\CC_v)$ with $a \ne \zeta$
and note that the balls $\cB(a,R)_{\zeta}^-$
exhaust $\PP^1_{\Berk} \backslash \{\zeta\}$.
Since $E$ is compact there is an $R$ with $E \subset \cB(a,R)_{\zeta}$.
For each $z \in \PP^1_{\Berk} \backslash \cB(a,R)_{\zeta}$ and each
$w \in \cB(a,R)_{\zeta}$, the ultrametric inequality for $\delta(x,y)_{\zeta}$
shows that $\delta(z,w)_{\zeta} = \delta(z,a)_{\zeta}$.  Since $\mu_1$
and $\mu_2$ are both supported on $E \subset \cB(a,R)_{\zeta}$,
for $z \notin \cB(a,R)_{\zeta}$
\begin{equation*}
u_i(z) \ = \ \int -\log_v(\delta(z,w)_{\zeta}) \, d\mu_i(z)
\ = \ -\log_v(\delta(z,a)_{\zeta}) \ .
\end{equation*}
Hence $u(z) \equiv 0$ in
$\PP^1_{\Berk} \backslash (\cB(a,R)_{\zeta} \cup \{\zeta\})$.
\end{proof}

\vskip .1 in
Since the equilibrium distribution $\mu_{\zeta}$ is unique, the potential 
function $u_E(z,\zeta) = u_{\mu_{\zeta}}(z,\zeta)$ is well-defined.  

\begin{definition} \label{DefE6} If \,$E \subset \PP^1_{\Berk}$ is a compact
set of positive capacity, then for each $\zeta \notin E$, the  
{\rm Green's function of $E$} is  
\begin{equation*}
G(z,\zeta;E) \ = \
     V_{\zeta}(E) - u_E(z,\zeta)  \qquad \text{for all $z \in \PP^1_{\Berk}$.}    
\end{equation*}
If \,$U \subset \PP^1_{\Berk}$ is a domain for which $\partial U$ has positive
capacity,  the {\rm Green's function of $U$} is 
\begin{equation*}
G_U(z,\zeta) \ = \
G(z,\zeta;\PP^1_{\Berk} \backslash U) \qquad \text{for $z, \zeta \in U$.}
\end{equation*} 
\end{definition}

\begin{proposition} \label{PropE26}  Let $E \subset \PP^1_{\Berk}$ be a compact
set of positive capacity.  Then $G(z,\zeta;E)$ has the following properties:

\vskip .05 in
$A)$  For each fixed $\zeta \notin E$, as a function of $z$,   

\quad $(1)$ $G(z,\zeta;E) \ge 0$ for all $z \in \PP^1_{\Berk}$.  

\quad $(2)$ $G(z,\zeta;E) > 0$ for $z \in U_{\zeta}$, 
where $U_{\zeta}$ is the connected component
of $\PP^1_{\Berk} \backslash E$ containing $\zeta$. 

\quad $(3)$ $G(z,\zeta;E)$ is finite and harmonic in 
$U_{\zeta} \backslash \{\zeta\}$.
  
\noindent{For} each $a \ne \zeta$, $G(z,\zeta;E) + \log_v(\delta(z,a)_{\zeta})$
               extends to a function harmonic on a neighborhood of $\zeta$.
            
\quad $(4)$ $G(z,\zeta;E) = 0$ on $\PP^1_{\Berk} \backslash U_{\zeta}$, 
except on a $($possibly empty$)$ $F_{\sigma}$ set 
$e \subset \partial E_{\zeta}$ of capacity $0$ 
which depends only on $U_{\zeta}$. 

\quad $(5)$ $G(z,\zeta;E)$ is continuous on 
                    $\PP^1_{\Berk} \backslash e$, 
                    and is strongly upper semi-continuous everywhere.  
                    
\quad $(6)$ $G(z,\zeta;E)$ is bounded if \,$\zeta \notin \PP^1(\CC_v)$. 

\noindent{It} is unbounded, with \,$G(\zeta,\zeta;E) = \infty$, 
            if $\zeta \in \PP^1(\CC_v)$. 
            
\quad $(7)$ $G(z,\zeta;E) \in \BDV(\PP^1_{\Berk})$, 
with $\Delta_{\PP^1_{\Berk}}(G(z,\zeta;E)) = \delta_{\zeta}(z) - \mu_{\zeta}$.                                 
 
\vskip .05 in                   
$B)$  As a function of two variables,   
$G(z_1,z_2;E) = G(z_2,z_1;E)$ for all $z_1, z_2 \notin E$.
\end{proposition}           
 
\begin{proof}  For (A), part (1) follows from Frostman's theorem 
(Theorem~\ref{ThmD11}).  Parts (4) and (5) (except for the fact that 
$e$ is independent of $\zeta$, which we will prove later) 
follow from Proposition~\ref{PropD7} and Frostman's theorem.  
Part (6) follows from the definition of the potential function 
$u_{\mu_{\zeta}}(z,\zeta)$ and properties of the Hsia kernel.  These 
same facts show $G(\zeta,\zeta;E) > 0$.     
Part (7) is a reformulation of Example~\ref{Example E.4}.  
For the first part of (3), we know that
$G(z,\zeta;E)$ is finite and harmonic in $U_{\zeta} \backslash \{\zeta\}$
by Example~\ref{Example E.8}.  
For the second part of (3), fix $a \ne \zeta$ and consider the function 
$f(z) = G(z,\zeta;E)+\log_v(\delta(z,a)_{\zeta})$.  
By Examples \ref{Example E.2} and \ref{Example E.4}, it belongs to $BDV(\PP^1_{\Berk})$ and satisfies 
\begin{equation*}
\Delta_{\PP^1_{\Berk}}(f) \ = \ \delta_a(z) - \mu_{\zeta} \ .
\end{equation*}
Let $V$ be a connected neighborhood of $\zeta$  
with $\Vbar \cap (E \cup \{a\}) = \phi$.
If $\zeta \notin \PP^1(\CC_v)$ then $f(z)$ is defined everywhere, and is 
strongly harmonic in $V$.  If $\zeta \in \PP^1(\CC_v)$ then $f(z)$ is 
harmonic and bounded in $V \backslash \{\zeta\}$, so it extends to a 
function strongly harmonic in $V$ by the Riemann Extension Theorem 
(Proposition \ref{PropE16}).  For (2), if $G(x,\zeta;E) = 0$
for some $x \in U_{\zeta} \backslash \{\zeta\}$, let $W$ be the component of
$U_{\zeta} \backslash \{\zeta\}$ containing $x$.  Then $G(z,\zeta;E)$
would achieve its minimum value $0$ at an interior point of $W$, 
and so would be identically $0$ on $W$ by the Maximum Principle 
(Proposition \ref{PropE13}).  Thus 
\begin{equation*}
\lim \begin{Sb} z \rightarrow \zeta \\ z \in W \end{Sb} G(z,\zeta;E)  \ = \ 0
\end{equation*}
which contradicts the fact that $G(\zeta,\zeta;E) > 0$ and 
$G(z,\zeta;E)$ is continuous at $\zeta$.  

To see that the set $e$ in part (4) is independent of $\zeta$, temporarily
write $e_{\zeta}$ for it, and take $\xi \ne \zeta$ in $U_{\zeta}$.  
Let $V$ be a neighborhood of $\zeta$ such that $\Vbar \subset U_{\zeta}$
and $\xi \notin \Vbar$.  
Since $\partial V \subset U_{\zeta} \backslash \{\zeta,\xi\}$
there is a constant $C >0$ such that 
$G(z,\xi;E) > C \cdot G(z,\zeta;E)$ on $\partial V$.  
Let $W$ be a connected component of $U_{\zeta} \backslash \Vbar$.  
Then $\partial W \subset \partial E_{\zeta} \cup \partial V$.  
Consider $f(z) = G(z,\xi;E)- C \cdot G(z,\zeta;E)$,  
which is harmonic on $W$.  For each $x \in \partial W \cap \partial V$ 
\begin{equation*}
\lim \begin{Sb} z \rightarrow x \\ z \in W \end{Sb} f(z) \ > \ 0 \ ,
\end{equation*}
while for each 
$x \in \partial W \cap (\partial E_{\zeta} \backslash e_{\zeta})$  
\begin{equation*}
\liminf \begin{Sb} z \rightarrow x \\ z \in W \end{Sb} f(z) \ \ge \ 0 \ .
\end{equation*} 
By the Strong Maximum Principle (Proposition \ref{PropE14}), $f(z) \ge 0$ for
all $z \in W$.  Hence, if $x \in \partial W \cap e_{\zeta}$, then
\begin{equation*}
\liminf \begin{Sb} z \rightarrow x \\ z \in W \end{Sb} G(z,\xi;E) 
\ \ge \  C \cdot \liminf \begin{Sb} z \rightarrow x \\ z \in W \end{Sb} 
                G(z,\zeta;E) \ > \ 0 \ .
\end{equation*} 
Since each $x \in e_{\zeta}$ belongs to $\partial W$ for some $W$,
it follows that $e_{\zeta} \subset e_{\xi}$.  
By symmetry, $e_{\zeta} = e_{\xi}$. 

\vskip .1 in

For part (B), we can assume without loss that $z_1, z_2$ belong to the
same connected component $U$ of $\PP^1_{\Berk}$, otherwise
$G(z_1,z_2;E) = G(z_2,z_1;E) = 0$.

First suppose neither $z_1$ nor $z_2$ is of type I.
Choose an exhaustion of $U$ by subdomains $W_n$ such that
$W_1 \subset W_2 \subset \cdots$, with $\Wbar_n \subset U$ for each $n$.
(Such an exhaustion can be gotten by fixing $\zeta \in U$ and covering
$E$ by balls $\cB(x,r_x+1/n)_{\zeta}^-$ where $r_x = \diam_{\zeta}(x)$.
Since $E$ is compact, finitely many balls cover $E$, for each $n$.
Take  $W_n = \PP^1_{\Berk} \backslash
\cup_{i=1}^{N_i} \cB(x_i,r_{x_i}+1/n)_{\zeta}$.)

Put $g_1(x) = G(x,z_1;E)$, $g_2(x) = G(x,z_2;E)$,
and write $\mu_1 = \mu_{z_1}$, $\mu_2 = \mu_{z_2}$.
Then $g_1, g_2 \in \BDV(\PP^1_{\Berk})$,
and $g_1,g_2 \in \cC(\Wbar_n)$ for each $n$
since $z_1, z_2 \notin \PP^1(\CC_v)$.

Fix $n$, and let $\Gamma \subset W_n$ be a subgraph.
Then $g_1|_{\Gamma}, g_2|_{\Gamma} \in \CPA(\Gamma)$
since each $g_i(x)$ is harmonic in $U \backslash \{z_i\}$.
By Proposition~\ref{PropC3}(D),
\begin{equation} \label{FMM1}
\int_{\Gamma} g_1 \, \Delta_{\Gamma}(g_2) \ = \
\int_{\Gamma} g_2 \, \Delta_{\Gamma}(g_1) \ .
\end{equation}
Taking the limit over subgraphs $\Gamma$, we obtain
\begin{equation} \label{FMM2}
\int_{\Wbar_n} g_1 \, \Delta_{\Wbar_n}(g_2) \ = \
\int_{\Wbar_n} g_2 \, \Delta_{\Wbar_n}(g_1) \ .
\end{equation}
Here $\Delta_{\Wbar_n}(g_i) = (r_{\Ubar,\Wbar_n})_*(\Delta_{\Ubar}(g_i))$
for each $i$.
Assume $n$ is large enough that $z_1, z_2 \in W_n$.
Then $(r_{\Ubar,\Wbar_n})_*(\Delta_{\Ubar}(g_i)) =
         \delta_{z_i}(x) - (r_{\Ubar,\Wbar_n})_*(\mu_i)$,
and $\mu_i$ is supported on $\partial U$.  Thus 
\begin{equation} \label{FMM3}
\int_{\Wbar_n} g_i(x) \, d(r_{\Ubar,\Wbar_n})_*(\mu_j)(x)
\ = \ \int_{\partial U} g_i(r_{\Ubar,\Wbar_n}(x)) \, d\mu_j(x) \ ,
\end{equation}
Writing $g_{i,n}(x) = g_i(r_{\Ubar,\Wbar_n}(x))$ and using (\ref{FMM2}),
(\ref{FMM3}) gives
\begin{equation}  \label{FMM4}
g_1(z_2) - \int_{\partial U} g_{1,n}(x) \, d\mu_2(x) \ = \
g_2(z_1) - \int_{\partial U} g_{2,n}(x) \, d\mu_1(x) \ .
\end{equation}
Here $\lim_{n \rightarrow \infty} r_{\Ubar,\Wbar_n}(x) = x$ for each
$x \in \partial U$.  Let $e$ be the exceptional set of capacity $0$
in $\partial U$ given in part (A.4).  By Lemma~\ref{LemD9}, $\mu_i(e) = 0$.
If $x \in \partial U \backslash e$, then
\begin{equation*}
\lim_{n \rightarrow \infty} g_{i,n}(x) \ = \ 
\lim_{n \rightarrow \infty} g_i(r_{\Ubar,\Wbar_n}(x)) \ = \ g_i(x) \ = \ 0 \ .
\end{equation*}

For each $n$, $g_{i,n}(x) = g_i(r_{\Ubar,\Wbar_n}(x))$ is continuous
on $\partial U$ since $r_{\Ubar,\Wbar_n}(x)$ is continuous and $g_i(x)$
is continuous on $U$.  Moreover, the $g_{i,n}(x)$ are uniformly bounded,
since $g_i(x) = G(x,z_i;E) = V_{\zeta}(E) - u_{\mu_i}(x,z_i)$ is bounded
in a neighborhood of $E$.  By Lebesgue's dominated convergence theorem,
\begin{equation*}
\lim_{n \rightarrow \infty} \int_{\partial U} g_{i,n}(x) \, d\mu_j(x) 
\ = \ 0 \ .
\end{equation*}
Thus (\ref{FMM4}) gives $g_1(z_2) = g_2(z_1)$, that is,
$G(z_1,z_2;E) = G(z_2,z_1;E)$, as claimed.

\vskip .05 in
Now suppose $z_1 \in U$ is of type I, but $z_2 \in U$ is not of type I.
Let $t \rightarrow z_1$ via points not of type I.  Since $G(x,z_2;E)$ is
continuous for $x \in U$,
\begin{equation*}
G(z_1,z_2;E) \ = \ \lim_{t \rightarrow z_1} G(t,z_2;E)
\ = \ \lim_{t \rightarrow z_1} G(z_2,t;E) \ .
\end{equation*}
We claim that $\lim_{t \rightarrow z_1} G(z_2,t;E) = G(z_2,z_1;E)$.
To see this, fix a point $z$ in the main dendrite $D$ for $U$.
There is a constant $C_t$ such that
\begin{equation*}
-\log_v(\delta(x,y)_{t}) \ = \ j_z(x,y) - j_z(x,t) - j_z(y,t) + C_t
\end{equation*}
for all $x, y \in \PP^1_{\Berk}$.
Let $w$ be the point where the path from $z_1$ to $z$ meets $D$.
If $t$\, lies on the path from $z_1$ to $w$,
then for all $x, y \in E$, $j_z(x,t) = j_z(x,w)$ and $j_z(y,t) = j_z(y,w)$.
In other words, for $x, y \in E$, $\log_v(\delta(x,y)_t)$ depends on $t$
only through the constant $C_t$.
Hence the same probability measure $\mu$ minimizes the energy integral
\begin{equation*}
V_t(E) \ = \ \iint_{E \times E} -\log_v(\delta(x,y)_t) \, d\mu(x) d\mu(y)
\end{equation*}
for all $t$ on the path from $w$ to $z_1$, that is, $\mu_t = \mu_{z_1} = \mu$.

Similarly,  $-\log_v(\delta(z_2,y)_t)$ depends on $t$ only through $C_t$,
if $y \in E$ and $t$ is sufficiently near $z_1$. Hence
\begin{equation*}
G(z_2,t;E) \ = \ V_t(E) - \int_E -\log_v(\delta(z_2,w)_t) \, d\mu(w)
\end{equation*}
is independent of $t$, for $t$ sufficiently near $z_1$.  Thus
$\lim_{t \rightarrow z_1} G(z_2,t,E) = G(z_2,z_1;E)$.

\vskip .05 in
Finally, let $z_1, z_2 \in U$ be arbitrary.  Let $t$ approach $z_1$ through 
points of type I.  Using continuity and an argument like the one above we find
\begin{eqnarray*}
G(z_1,z_2;E) & = & \lim_{t \rightarrow z_1} G(t,z_2;E) \\
& = & \lim_{t \rightarrow z_1} G(z_2,t;E) \ = \ G(z_2,z_1;E) \ .
\end{eqnarray*}
\end{proof}

\vskip .1 in
We will now show that the equilibrium measure $\mu_{\zeta}$ 
has an interpretation as the reproducing kernel for harmonic functions:

\begin{proposition} \label{PropE27}
Let $U \subset \PP^1_{\Berk}$ be a domain such that $\partial U$ has positive
capacity.  Suppose $f$ is harmonic in $U$ and extends to 
a continuous function on $\Ubar$.  Then for each $\zeta \in U$,
\begin{equation*}
f(\zeta) \ = \ \int_{\partial U} f(z) \, d\mu_{\zeta}(z) \ .
\end{equation*} 
\end{proposition}

\begin{proof}  
First suppose $\zeta \in U$ is not of type I.  As in the previous proof,
choose an exhaustion on $U$ by subdomains $W_n$
with $W_1 \subset W_2 \subset \cdots $ and $\Wbar_n \subset U$ for each $n$.
Without loss, we can assume $\zeta \in W_1$.
Put $g_n(x) = G(z,\zeta;W_n)$.  Since $\partial W_n$ consists of a finite
number of points, none of which is of type I, $g_n(x)$ is continuous on
$\Wbar_n$ and $g_n(x) = 0$ for each $x \in \partial W_n$.
Note that $g_n(x) \in \BDV(W_n)$, 
$f, g_n \in \cC(\Wbar_n)$ for each $n$,
and $f$ is strongly harmonic in $W_n$ by Lemma \ref{LemE11};  in particular
$f \in \BDV(W_n)$.  

Fix $n$, and let $\Gamma \subset W_n$ be a subgraph.
Then $f|_{\Gamma}, g|_{\Gamma} \in \CPA(\Gamma)$.
By Proposition~\ref{PropC3}(D),
\begin{equation} \label{FNN1}
\int_{\Gamma} f \, \Delta_{\Gamma}(g_n) \ = \
\int_{\Gamma} g_n \, \Delta_{\Gamma}(f) \ .
\end{equation}
Taking a limit over subgraphs $\Gamma$, we find that
\begin{equation} \label{FNN2}
\int_{\Wbar_n} f  \, \Delta_{\Wbar_n}(g_n)
\ = \ \int_{\Wbar_n} g_n \, \Delta_{\Wbar_n}(f) \ .
\end{equation}
Here $\Delta_{\Wbar_n}(g_n) = \delta_{\zeta}(x) - \mu_{\zeta,n}$
where $\mu_{\zeta,n}$ is the equilibrium measure of
$E_n = \PP^1_{\Berk} \backslash W_n$ with respect to $\zeta$.
It is supported on $\partial W_n$.
Since $f$ is strongly harmonic on $W_n$, $\Delta_{\Wbar_n}(f)$ is supported
on $\partial W_n$, where $g_n(x) = 0$.  Hence the integral on the right
side of (\ref{FNN2}) is $0$.  It follows that
\begin{equation} \label{FNN3}
f(\zeta) \ = \ \int f(x) \, d\mu_{\zeta,n}(x) \ .
\end{equation}

As $n \rightarrow \infty$ the measures $\mu_{\zeta,n}$ converge weakly
to $\mu_{\zeta}$, the equilibrium measure of $E = \PP^1_{\Berk} \backslash U$.
(This follows by the same argument as in the proof of Corollary~\ref{CorD6}.)
Since $f(x)$ is continuous on $\Ubar$,
\begin{equation*}
f(\zeta) \ = \ \lim_{n \rightarrow \infty} \int f(x) \, d\mu_{\zeta,n}(x)
\ = \ \int f(x) \, d\mu_{\zeta}(x) \ ,
\end{equation*}
yielding the result in this case.

If $\zeta$ is of type I, let $t$ approach $\zeta$ through points not of
type I.  As in the proof of Proposition \ref{PropE26}(B),
$\mu_{\zeta} = \mu_t$ for $t$ is sufficiently near $\zeta$.
Since $f$ is continuous, it follows from the previous case
that $f(\zeta) = \int_{\partial U} f(x) \, d\mu_{\zeta}(x)$.
\end{proof}

\vskip .1 in
In the classical theory over $\CC$, the reproducing kernel is 
the inward normal derivative of $G_U(z,\zeta;E)$ on $\partial E_{\zeta}$,  
\begin{equation*}
\frac{1}{2 \pi} \frac{\partial}{\partial n} G_U(z,\zeta) \ .
\end{equation*} 
On Berkovich space, taking $E = \PP^1_{\Berk} \backslash U$ and considering 
Proposition \ref{PropE26}(A.7) suggests that
$\mu_{\zeta} = \Delta(G(z,\zeta;E))|_{\partial U}$
should identified with the normal derivative of $G(z,\zeta;E)$ on 
$\partial E_{\zeta} = \partial U$.

\section{Subharmonic functions.}
\label{Section F}

\vskip .1 in
In this section we develop a theory of subharmonic functions
on the Berkovich line.

\subsection{Subharmonic and strongly subharmonic functions.}

\vskip .1 in
We will call a subset $V \subset \PP^1_{\Berk}$ a {\it simple domain}
if $V$ is an open disc or punctured disc, and is not $\PP^1_{\Berk}$ itself.  
Thus $V$ is a simple domain if and only if $V$ is connected 
and $\partial V$ is a nonempty finite set  $\{x_1, \ldots, x_m\}$,
where each $x_i$ is of type II, III, or IV.

\begin{definition} \label{DefF1}
If $U \subset \PP^1_{\Berk}$ is open, then $V \subset U$ is
a {\rm simple subdomain} of $U$ if $V$ is a simple domain
and $\Vbar \subset U$.
\end{definition}

By Proposition~\ref{PropE18}, if $V$ is a simple domain,
each harmonic function $h(x)$ on $V$ extends to a continuous
function on $\Vbar$, and the Poisson formula expresses $h(x)$ on $V$
in terms of its values on $\partial V$.  We will often use this implicitly,
and if $h(x)$ is harmonic on a simple domain, we will speak of its values
on $\partial V$.

Recall that a function $f : U \rightarrow \RR \cup \{\pm \infty\}$ is
{\it upper semicontinuous} if $\limsup_{z \rightarrow x} f(z) \le f(x)$
for all $x \in U$, and is {\it strongly upper semicontinuous} if
$\limsup_{z \rightarrow x} f(z) = f(x)$ for all $x \in U$.

\begin{definition} \label{DefF2}  Let $U \subset \PP^1_{\Berk}$ be open. 

\noindent{A} function  $f : U \rightarrow [-\infty,\infty)$ 
is {\rm strongly subharmonic} on $U$ if

$A)$  $f$ is strongly upper semicontinuous on $U$, and
        for each $x \in U \cap \PP^1(\CC_v)$,
\begin{equation*}
 f(x) \ = \ \limsup \begin{Sb} z \rightarrow x \\
                        z \in U \backslash \PP^1(\CC_v) \end{Sb} f(z) \ ;
\end{equation*}

$B)$  for each component $W$ of $U$, 
$f \in \BDV(W)$ and $\Delta_{\Wbar}(f)|_{W} \le 0$.
 
\vskip .05 in
\noindent{A} function $f : U \rightarrow [-\infty,\infty)$
is {\rm subharmonic} on $U$ if

$C)$  $f(x)$ is upper semicontinuous,

$D)$  $f(x) \not\equiv -\infty$ on any connected component of $U$, and

$E)$  for each simple subdomain $V \subset U$ and each function $h(x)$ 
harmonic on $V$, 
if $h(x) \ge f(x)$ on $\partial V$, then $h(x) \ge f(x)$ on $V$.
\end{definition}

\begin{remark}
The condition in (A) controlling $f$ 
on $U \cap \PP^1(\CC_v)$ is necessary, as shown by the following example:  
take $U = \cB(0,1)^-$, and put $f(z) = 0$ on $U \backslash \ZZ_p$, 
$f(z) = 1$ on $\ZZ_p$.  Then $f(z)$ is strongly upper semicontinuous on $U$,
since for each $x \in U \backslash \ZZ_p$ there is a neigbhorhood $V$ of $x$
with $V \cap \ZZ_p = \phi$.  It also belongs to $\BDV(U)$, 
with $\Delta_{\Ubar}(f) \equiv 0$.  However, it is not subharmonic on $U$,
in contrast to Proposition \ref{PropF1} below.   
\end{remark}

\begin{proposition} \label{PropF1}  Let $U \subset \PP^1_{\Berk}$ be open.

If $f$ is strongly subharmonic on $U$, then it is subharmonic on $U$.

If $f$ is subharmonic on $U$,
then $f|_V$ is strongly subharmonic on each subdomain $V \subset U$
with $\Vbar \subset U$.
If the measures $|\Delta_{\Vbar}(f)|(\Vbar)$ are uniformly bounded
for all simple subdomains $V \subset U$, 
then $f$ is strongly subharmonic on $U$.

\end{proposition}

Before giving the proof, we will need several lemmas.

\vskip .05 in
If $V$ is a simple domain, then either $V$ is a disc
with a single boundary point $x$,
or a punctured disc
$\cB(a_1,r_1)_{\zeta}^- \backslash \cup_{i=2}^m \cB(a_i,r_i)_{\zeta}$
with boundary points $x_1, \ldots, x_m$.  In that case the main dendrite of
$V$ is the graph $\Gamma$ whose endpoints are $x_1, \ldots, x_m$.  (By an
endpoint of $\Gamma$, we mean a point with a single edge emanating from it.)  
We will call
$\partial \Gamma = \{x_1, \ldots, x_m\}$ the {\it boundary} of $\Gamma$,
and $\Gamma_0 = \Gamma \backslash \partial \Gamma$
the {\it interior} of $\Gamma$.  
Note that $V$ can be recovered from $\Gamma$, 
indeed if $r_{\Gamma} : \PP^1_{\Berk} \rightarrow \Gamma$
is the retraction map, then 
\begin{equation*}
V \ = \ r_{\Gamma}^{-1}(\Gamma_0) \ .
\end{equation*}
Conversely, if $\Gamma \subset \PP^1_{\Berk}$ is a subgraph, we will call
$V = r_{\Gamma}^{-1}(\Gamma_0)$ the domain associated to $\Gamma$.
Clearly $\partial V = \partial \Gamma$, and $\Gamma_0$ is the main dendrite
of $V$.

\begin{lemma} \label{LemF2}
Let $U \subset \PP^1_{\Berk}$ be a domain.

If \, $\Gamma \subset U$ is a subgraph, then the associated
domain $V = r_{\Gamma}^{-1}(\Gamma_0)$ is a simple subdomain of $U$
if and only if $r_{\Gamma}(\partial U) \subset \partial \Gamma$.

A disc $V$ is a simple subdomain of $U$
and only if its boundary point $x$ belongs to $U$ and
$\partial U \cap V = \phi$.
\end{lemma}

\begin{proof} If $r_{\Gamma}(\partial U) \subset \partial \Gamma$, then
$V = r_{\Gamma}^{-1}(\Gamma_0)$ is a connected open set with
$\partial U \cap V = \phi$.  Since $U$ is connected and $\Gamma_0 \subset U$, 
necessarily $V \subset U$.  Since $\Gamma \subset U$ and
$\partial V = \partial \Gamma$, it follows that $\Vbar \subset U$.
Conversely, if $V$ is a simple subdomain of $U$ 
and $\Gamma$ is the main dendrite of $V$, 
then $r_{\Gamma}^{-1}(\Gamma_0) \subset U$.  Since $\partial U \cap U = \phi$, 
it follows that $\partial U \subset r_{\Gamma}^{-1}(\partial \Gamma)$.
Thus $r_{\Gamma}(\partial U) \subset \partial \Gamma$.

If a disc $V = \cB(a,r)_{\zeta}^-$ is a simple subdomain of $V$ then
by definition $\Vbar \subset U$ so $\{x\} = \partial V \subset U$,
and $\partial U \cap \Vbar = \phi$ so $\partial U \cap V = \phi$.
Conversely, if $V$ is a disc with $\{x\} = \partial V \subset U$,
then $V \cap U$ is nonempty. Since $V$ is one of the connected components of
$r_{\Gamma}^{-1}(x) \backslash \{x\}$, if $\partial U \cap V = \phi$, 
necessarily $V \subset U$.  It follows that $\Vbar \subset U$.
\end{proof}

\vskip .1 in
\begin{lemma} \label{LemF3}
Let $U \subset \PP^1_{\Berk}$ be open.

If $V$ is a simple subdomain of $U$, then there is a simple subdomain
$W$ of $U$ with $\Vbar \subset W \subset U$.  Moreover, $W$ can
be chosen so that for each $x_i \in \partial V$ which is not of type IV,
there is a point $w_j \in \partial W$ with $r_{\Vbar}(w_j) = x_i$.
\end{lemma}

\begin{proof}
Write $\partial V = \{x_1, \ldots, x_m\}$.  For each $x_i$, choose a simple
subdomain $V_i$ of $U$ which contains $x_i$.  Such a $V_i$ exists,
since $\partial U$ is closed and disjoint from $x_i$,
and punctured discs are cofinal in the open neighborhoods of $x_i$.
Put $W = V \cup_{i=1}^m$.  Then $W$ is connected and open,
and $\partial W \subset \cup_{i=1}^m \partial V_i \subset U$ is finite,
so $W$ is a simple subdomain of $U$.  Clearly $\Vbar \subset W$.

Write $\partial W = \{w_1, \ldots, w_M\}$.
If there is some $x_i$ which is not the retraction to $\Vbar$ of any $w_j$,
and $x_i$ is not of type IV,
then the component of $\PP^1_{\Berk} \backslash V$ containing $x_i$ is a
nonempty closed disc contained in $W$.
By removing a proper closed subdisc of that
disc from $W$, we would obtain a simple subdomain $\tW$ of $U$ with 
$\Vbar \subset \tW$, and having a
boundary point whose retraction was $x_i$.
Doing this sequentially for each $x_i$ gives the result.
\end{proof}

\vskip .1 in
\begin{lemma} \label{LemF4}
Let $g: (a,b) \rightarrow \RR$ be a function such that 
for each $t \in (a,b)$, both one-sided derivatives
$g_{-}^{\prime}(t) = \lim_{h \rightarrow 0^-} (g(t+h)-g(t))/h$
and $g_{+}^{\prime}(t) = \lim_{h \rightarrow 0^+} (g(t+h)-g(t))/h$ exist.
Suppose that

$A)$ $g_{-}^{\prime}(x) \le g_{+}^{\prime}(x)$ for each $x \in (a,b)$,  and

$B)$ $g_{+}^{\prime}(x) \le g_{-}^{\prime}(y)$ for each $x, y \in (a,b)$ 
with $x < y$.

\noindent{Then} $g(t)$ is convex up on $(a,b)$.
\end{lemma}

\begin{proof}
Note that $g(t)$ is continuous,
since the left and right derivatives exist at each point.
Fix $x, y \in (a,b)$ with $x < y$, and let 
\begin{equation*}
L_{x,y}(t) \ = \ \frac{g(y)-g(x)}{y-x} \cdot (t-x) + g(x)
\end{equation*} 
be the line through $(x,g(x))$ and $(y,g(y))$.
We claim that $g(t) \le L_{x,y}(t)$ for all $x < t < y$.
Suppose to the contrary that there is a point 
$z \in (x,y)$ where $g(z) > L_{x,y}(z)$.
Put $\alpha = (g(z)-g(x))/(z-x)$, $\beta = (g(y)-g(z))/(y-z)$.
Then $\alpha > \beta$.

We claim there are a point $u \in (x,z)$ where $g_{+}^{\prime}(u) \ge \alpha$,
and a point $v \in (z,y)$ where $g_{-}^{\prime}(v) \le \beta$.
This contradicts $g_{+}^{\prime}(u) \le g_{-}^{\prime}(v)$.

First consider the interval $(x,z)$.
Let $L_{x,z}(t) = \alpha \cdot (t-x) + g(x)$ be the line
through $(x,g(x))$ and $(z,g(z))$, and put $f(t) = g(t)-L_{x,z}(t)$.
Then $f(x) = f(z) = 0$.
We will now apply the argument in Rolle's theorem.
Since $f(t)$ is continuous, there is a point $u \in (x,y)$
where $f$ achieves a maximum or a minimum.  If $u$ is a minimum,
then $f_{-}^{\prime}(u) \le 0$ and $f_{+}^{\prime}(u) \ge 0$,
so $g_{-}^{\prime}(u) \le \alpha$ and $g_{+}^{\prime}(u) \ge \alpha$.
If $u$ is a maximum,
then $f_{-}^{\prime}(u) \ge 0$ and $f_{+}^{\prime}(u) \le 0$,
so $g_{-}^{\prime}(u) \ge \alpha$ and $g_{+}^{\prime}(u) \le \alpha$.
Since $g_{+}^{\prime}(z) \ge g_{-}^{\prime}(u)$ it must be that
$g_{-}^{\prime}(u) = g_{+}^{\prime}(u) = \alpha$.  In either case
we have a point $u$ where $g_{-}^{\prime}(u) \le \alpha$ and
$g_{+}^{\prime}(u) \ge \alpha$.  This second inequality is the one we want.

Next consider the interval $(z,y)$.  By the same reasoning as before, 
there is a point $v \in (z,y)$ where $g_{-}^{\prime}(v) \le \beta$
and $g_{+}^{\prime}(v) \ge \beta$.  This time it is the first
inequality that we want.
\end{proof}

\vskip .1 in
\begin{lemma} \label{LemF5}
Let $\Gamma$ be a metrized graph, and take $f \in \BDV(\Gamma)$.
Let $\Delta_{\Gamma}(f)^{+}$ and $\Delta_{\Gamma}(f)^{-}$
be the positive and negative measures in the Jordan decomposition of 
$\Delta_{\Gamma}(f)$.  Then either $f$ is constant, or it achieves its maximum
at a point of $\supp(\Delta_{\Gamma}(f)^{+})$ and its minimum at a point
of $\supp(\Delta_{\Gamma}(f)^{-})$. 
\end{lemma}

\begin{proof}  Recall that each $f \in \BDV(\Gamma)$ is
continuous (Lemma~\ref{LemE7}(A)).  Let $x_1 \in \Gamma$ be a point where
$f(x)$ achieves its maximum, and let $\Gamma_1$ be the connected component
of $\{x \in \Gamma : f(x) = f(x_1) \}$ containing $x_1$.
If $\Gamma_1 = \Gamma$, then $f(x)$ is constant.
If $\Gamma_1 \ne \Gamma$, let $x_0$ be a boundary point of $\Gamma_1$.

We claim that $x_0 \in \supp(\Delta_{\Gamma}(f)^{+})$.  Suppose not.
Let $\rho(x,y)$ be the path distance metric on $\Gamma$.
Then there is a neighborhood
$\Gamma_{x_0}(\varepsilon) = \{ x \in \Gamma : \rho(x,x_0) < \varepsilon \}$
with $\Delta_{\Gamma}(f)|_{\Gamma_{x_0}(\varepsilon)} \le 0$.  After shrinking
$\varepsilon$ if necessary, we can assume that $\Gamma_{x_0}(\varepsilon)$
is a star, a union of half-open segments $[x_0, x_0 + \varepsilon \vec{v}_i)$
where the $\vec{v_i}$ are the direction vectors at $x_0$.
Since $f(x_0)$ is the maximum value of $f(x)$ on $\Gamma$,
necessarily $d_{\vec{v}_i}f(x_0) \le 0$ for each $i$.
If $d_{\vec{v}_i}f(x_0) < 0$ for some $i$, then
\begin{equation*}
\Delta_{\Gamma}(f)(x_0) \ = \ - \sum_i d_{\vec{v}_i}f(x_0) \ > \ 0
\end{equation*}
contradicting  $x_0 \notin \supp(\Delta_{\Gamma}(f)^{+})$.
Hence $d_{\vec{v}_i}f(x_0) = 0$ for each $i$.

Since $x_0$ is a boundary point of $\Gamma_1$,
there is a point $y_0 \in \Gamma_{x_0}(\varepsilon)$ where $f(y_0) < f(x_0)$.
Let $i$ be such that $y_0 \in e_i = (x_0,x_0+\varepsilon \vec{v}_i)$.
As $\Delta_{\Gamma}(f)|_{\Gamma_{x_0}(\varepsilon)} \le 0$,
it follows that $\Delta_{\Gamma}(f)|_{e_i} \le 0$.
For each $y = x_0 + t \vec{v_i} \in e_i$,
there are two direction vectors at $y$:  
write $\vec{v}_{+}$ for the one that leads away from $x_0$,
and $\vec{v}_{-}$ for the one that leads towards $x_0$.

By the definition of $\Delta_{\Gamma}(f)$, for each $y \in e_i$
\begin{equation*}
0 \ \ge \ \Delta_{\Gamma}(f)(y)
  \ = \ -d_{\vec{v}_{-}}f(y) - d_{\vec{v}_{+}}f(y)
\end{equation*}
so $d_{\vec{v}_{+}}f(y) \ge -d_{\vec{v}_{-}}f(y)$.
Similarly, for each open subsegment  $(y_1,y_2) \subset e_i$,
\begin{equation*}
0 \ \ge \ \Delta_{\Gamma}(f)((y_1,y_2))
  \ = \ d_{\vec{v}_{+}}f(y_1) + d_{\vec{v}_{-}}f(y_2)
\end{equation*}
so $-d_{\vec{v}_{-}}f(y_2) \ge d_{\vec{v}_{+}}f(y_1)$.

For the function $g(t) = f(x_0 + t \vec{v_i})$, we have
$g_{-}^{\prime}(t) = -d_{\vec{v}_{-}}f(y)$ and 
$g_{+}^{\prime}(t) = d_{\vec{v}_{+}}f(y)$.  
By Lemma \ref{LemF4} $f$ is convex up on $e_i$.
Since $d_{\vec{v}_i}(f)(x_0) = 0$,
$f(x_0 + t \vec{v}_i)$ is nondecreasing in $t$.
Hence $f(y_0) \ge f(x_0)$, a contradiction.
Thus $x_0 \in \supp(\Delta_{\Gamma}(f)^{+})$.

The case of a minimum is similar.  
\end{proof}
   
\begin{lemma} \label{LemF6}
Let $f$ be subharmonic on a domain $U \subset \PP^1_{\Berk}$.
Then $f(x) \ne -\infty$ on $U \backslash \PP^1(\CC_v)$.
\end{lemma}

\begin{proof}
Suppose $f(x_1) = -\infty$ for some $x_1 \in U \backslash \PP^1(\CC_v)$.
We will show that $f(x) \equiv -\infty$ on $U$,
which contradicts the definition of subharmonicity.

Fix $x \in U$ with $x \ne x_1$,
and take a simple subdomain $V_1$ of $U$
which contains $x$ and $x_1$.
Let $V$ be the connected component of
$V_1 \backslash \{x_1\}$ which contains $x$.  Then $V$
is a simple subdomain of $U$ with $x_1$ as a boundary point, which
contains $x$ in its interior.
(Note that this uses $x_1 \notin \PP^1(\CC_v)$).

If $x_1$ is the only boundary point of $V$, then $V$ is a disc.
In that case, by Lemma~\ref{LemE12}, each harmonic function on $V$
is constant and is determined by its value on $\partial V = \{x_1\}$.
Since $f(z)$ is subharmonic and $f(x_1) = -\infty$, we have $f(z) \le C$
on $V$ for each $C \in \RR$.  Thus $f(x) = -\infty$.

If $x_1$ is not the only boundary point of $V$, let $\Gamma$ be
the main dendrite of $V$, and let $\partial \Gamma = \{x_1, \ldots, x_m\}$.
Fix numbers $A_2, \ldots, A_m$ with $A_i \ge f(x_i)$ for each $i$.
Given $A_1 \in \RR$, Poisson's formula (Proposition~\ref{PropE19})
constructs a harmonic function
\begin{equation*}
h_{A_1}(z) \ = \ \sum_{i=1}^m A_i h_i(z)
\end{equation*}
on $V$, where $h_i(z)$ is the harmonic measure with $h_i(x_i) = 1$,
 $h_i(x_j) = 0$ for each $j \ne i$.  
 Here $0 < h_i(z) < 1$ on $V$, otherwise $h_i(z)$ would
achieve its maximum or minimum value on the interior of $V$,
contradicting the Maximum Principle (Proposition~\ref{PropE13}).
Fixing $z$ and letting $A_1 \rightarrow -\infty$,
we see that $h_{A_1}(z) \rightarrow -\infty$.
Since $f$ is subharmonic on $V$, again  $f(x) = -\infty$.

Since $x \in U$ is arbitrary, 
we have shown that $f(x) \equiv -\infty$ on $U$.
\end{proof}

\begin{lemma} \label{LemF7}
Let $f$ be subharmonic on a domain $U \subset \PP^1_{\Berk}$.
If $\Gamma \subset U$ is a subgraph for which 
$V = r_{\Gamma}^{-1}(\Gamma_0)$ is a simple subdomain of $U$,
then for each $p \in \Gamma$ and each direction $\vec{v}$ at $p$, 
the directional derivative $d_{\vec{v}}f(p)$ exists and is finite.
In particular, the restriction of $f$ to $\Gamma$ is continuous.
\end{lemma}

\begin{proof} By Lemma \ref{LemF2},  
$r_{\Gamma}(\partial U) \subset \partial \Gamma$.  
In the argument below we will want to consider subgraphs
$\Gamma^{\prime} \subset \Gamma$ for which the associated domain 
$V^{\prime} = r_{\Gamma^{\prime}}^{-1}(\Gamma_0^{\prime})$ 
is also a simple subdomain of $U$.  

There are two types of $\Gamma^{\prime}$ for which this is assured.
One is if $\Gamma^{\prime} = [a,b]$ is a segment
contained in an edge of $\Gamma$ (that is, no branch point of
$\Gamma$ is contained in $\Gamma^{\prime}_0$).
The other is if $\Gamma^{\prime} = \overline{\Gamma_{x_0}(\varepsilon)}
= \{ x \in \Gamma : \rho(x,x_0) \le \varepsilon \}$ is a closed neighborhood
of a point $x_0 \in \Gamma$, for some $\varepsilon > 0$.  If $\varepsilon$
is small enough, then $\Gamma^{\prime}$ is a star, a union of closed 
segments $[x_0,x_i]$ for $i = 1, \ldots m$.  In either case,
$r_{\Gamma,\Gamma^{\prime}}(\partial \Gamma) \subset \partial \Gamma^{\prime}$,
so $r_{\Gamma^{\prime}}(\partial U) \subset \partial \Gamma^{\prime}$.

By Lemma \ref{LemF6}, $f(x) \in \RR$ for all $x \in \Gamma$.

By abuse of notation, write $f$ for $f|_{\Gamma}$.
We will first show that for each $p \in \Gamma$,
and each $\vec{v}$ at $p$,
the derivative $d_{\vec{v}}f(p)$ exists in $\RR \cup \{-\infty\}$.
Fix $p$ and $\vec{v}$.  For each sufficiently small $T > 0$,
the segment $\Gamma^{\prime} = [p,p+T \vec{v}]$ corresponds
to a simple subdomain of $U$.  The harmonic function on $U$ whose
values agree with those of $f(x)$ at $p$ and $p+T \vec{v}$
is linear on $\Gamma^{\prime}$ and is given by
\begin{equation} \label{FNF1}
h(p+t \vec{v}) \ = \ (1-\frac{t}{T}) \cdot f(p)
         + \frac{t}{T} \cdot f(p+ T \vec{v}) \ .
\end{equation}
for $0 \le t \le T$.
Since $f$ is subharmonic, $f(p+t\vec{v}) \le h(p+t\vec{v})$ for all $t$.
Using (\ref{FNF1}), this gives
\begin{equation} \label{FNF2}
\frac{f(p+t\vec{v}) - f(p)}{t} \ \le \ \frac{f(p+T\vec{v}) - f(p)}{T} \ .
\end{equation}
Hence 
\begin{equation} \label{FNF3}
d_{\vec{v}}f(p) \ = \ \lim_{t \rightarrow 0^+} \frac{f(p+t\vec{v}) - f(p)}{t}
\end{equation}
exists in $\RR \cup \{-\infty\}$, 
since the limit on the right side is non-increasing.

\vskip .05 in
Next, we will show that if $p \notin \partial \Gamma$,
then $d_{\vec{v}}f(p) \ne -\infty$.  Since $p \notin \partial \Gamma$,
for sufficiently small $T > 0$,
the star $\Gamma^{\prime} = \overline{\Gamma_{p}(T)}$
is a subgraph of $\Gamma$ with $p \in \Gamma_0^{\prime}$.
Let $x_1, \ldots, x_m$ be its endpoints, where $m \ge 2$,
and let $\vec{v}_1, \ldots, \vec{v}_m$ be the direction vectors at $p$,
so $\Gamma^{\prime} = \cup_{i=1}^m [p,p+T \vec{v_i}]$.
Without loss, suppose $\vec{v} = \vec{v}_1$.
We claim that for each $0 < t < T$,
\begin{equation} \label{FNF4}
\frac{f(p + t \vec{v}) - f(p)}{t} \ \ge \
     - \sum_{i=2}^m \frac{f(x_i) - f(p)}{T} \ .
\end{equation}
If this fails for some $t$, then
\begin{equation} \label{FNF5}
f(p) \ > \ \frac{T f(p + t \vec{v}) + t \cdot \sum_{i=2}^m f(x_i)}
              {T + (m-1) \cdot t} \ .
\end{equation}

Consider the subgraph $\Gamma^{\prime \prime}
= [p, p+t \vec{v}] \cup (\bigcup_{i=2}^m [p,p+T \vec{v}_i])$, 
which corresponds to a simple subdomain  $V^{\prime \prime} \subset U$
with boundary points $x_1^{\prime \prime} = p + t \vec{v}, x_2, \ldots, x_m$.
The harmonic function $h(z)$ on $V^{\prime \prime}$ with 
\begin{equation*}
h(x_1^{\prime \prime}) = f(x_1^{\prime \prime})\ ,
\qquad h(x_i) = f(x_i) \text{\quad for $i \ge 2$\ ,  }
\end{equation*}
satisfies
\begin{equation*}
0 \ = \ \Delta_{\Gamma}h(p) \ = \ \frac{f(x_1^{\prime \prime}) - h(p)}{t}
         + \sum_{i=2}^m \frac{f(x_i) - h(p)}{T} \ , 
\end{equation*}
which leads to
\begin{equation*}
h(p) \ = \ \frac{T f(p + t \vec{v}) + t \cdot \sum_{i=2}^m f(x_i)}
              {T + (m-1) \cdot t}
     \ < \ f(p) \ .
\end{equation*}
This contradicts the subharmonicity of $f(x)$.
By (\ref{FNF4}), the right side of
\begin{equation*}
d_{\vec{v}}f(p) \ = \ \lim_{t \rightarrow 0^+} \frac{f(p + t\vec{v})- f(p)}{t}
\end{equation*}
is bounded from below, so $d_{\vec{v}}f(p) > - \infty$.

\vskip .05 in
It remains to consider points $p \in \partial \Gamma$.

First suppose $p \in \partial \Gamma \backslash r_{\Gamma}(\partial U)$.
There is a single edge $e$ of $\Gamma$ emanating from $p$;  let $x$ approach
$p$ along this edge.  For each $x$, let $V_x$ be the connected
component of $\PP^1_{\Berk} \backslash \{x\}$ containing $p$.
Then $V_x$ is a disc with $\partial V_x = \{x\}$.  Since
$x \notin r_{\Gamma}(\partial U)$, it follows that $\Vbar_x \subset U$, and so  
$V_x$ is a simple subdomain of $U$.  Let $h$ be the harmonic
function on $V_x$ with $h(z) = f(x)$ for all $z \in V_x$.
Since $f$ is subharmonic, $f(z) \le f(x)$ for all $z \in V_x$.
Thus, $f(x)$ is decreasing as $x \rightarrow p$.  If $\vec{v}$
is the unique direction vector at $p$, it follows that
\begin{equation}  \label{FCV1}
d_{\vec{v}}f(p) \ = \ \lim_{t \rightarrow 0^+} \frac{f(p + t\vec{v})- f(p)}{t}
                \ \ge \ 0 \ .
\end{equation}

Finally, suppose $p \in r_{\Gamma}(\partial U) \subset \partial \Gamma$.  
Let $\vec{v}$ be the unique direction vector at $p$, 
and write  $\partial \Gamma = \{x_1, \ldots, x_m\}$.
Without loss, suppose $p = x_1$.
Let $W \subset U$ be a punctured disc containing $p$
but not $x_2, \ldots, x_m$, such that $\Wbar$ is contained in $U$.
Recalling that $V = r_{\Gamma}^{-1}(\Gamma_0)$, put $\tV = W \cup V$.
Then $\tV$ is a connected open set with a finite
number of boundary points whose closure is contained in $U$, 
so it is a simple subdomain of $U$.
Let $\tGamma$ be the subgraph spanned by $\partial \tV$.
Then $\Gamma \subset \tGamma$, since $\{x_2, \ldots, x_m\}$
remain boundary points of $\tV$, and if $\tx$ is a boundary
points of $\tV$ which is not contained in $V$, 
then the path from $\tx$ to each $x_i$, $i=2, \ldots, m$, passes through $p$.

>From this we see that $p$ is an interior point of $\tGamma$.
By what has been shown above, $d_{\vec{v}}f(p)$ is finite,
when $f$ is regarded as a function on $\tGamma$.
However, $d_{\vec{v}}f(p)$ depends only on the restriction of $f$ to $\Gamma$.
Thus it is finite.
\end{proof}

\begin{lemma} \label{LemF8}
Let $f$ be subharmonic on a domain $U \subset \PP^1_{\Berk}$.
If $\Gamma \subset U$ is a subgraph for which 
$V = r_{\Gamma}^{-1}(\Gamma_0)$ is a simple subdomain of $U$,
then $f|_{\Gamma} \in \BDV(\Gamma)$ and
$\Delta_{\Gamma}(f)^{+} \subset r_{\Gamma}(\partial U)
                        \subset \partial \Gamma$.
\end{lemma}

\begin{proof}
We have seen in Lemma \ref{LemF6} that $f(x) \in \RR$ for all $x \in \Gamma$,
and in Lemma \ref{LemF7} that $d_{\vec{v}}f(p)$ exists and is finite,
for all $p \in \Gamma$ and all directions $\vec{v}$ at $p$.  
Thus, $\Delta_{\Gamma}(f)$ exists as a finitely additive set function 
on the Boolean algebra $\cA(\Gamma)$ generated by the 
open, closed and half-open segments in $\Gamma$.

\vskip .1 in
We first claim that $\Delta_{\Gamma}(p) \le 0$ for each $p \in \Gamma_0$.
Suppose not, and fix $p \in \Gamma_0$ with $\Delta_{\Gamma}(p) > \eta > 0$.
Since $p \in \Gamma_0$, there are at least two edges emanating from $p$.
Let $\varepsilon > 0$ be small enough
that the closed neighborhood $\overline{\Gamma_p(\varepsilon)}$ is a star.
If $\vec{v}_1, \ldots, \vec{v}_m$ are the direction vectors at $p$, then
\begin{equation} \label{FMD1}
\sum_{i=1}^m d_{\vec{v}_i}f(p) \ = \ - \Delta_{\Gamma}(f)(p) \ < \ -\eta \ .
\end{equation}
Since $d_{\vec{v}_i} f(p) 
= \lim_{t \rightarrow 0^+} (f(p+t \vec{v_i})-f(p))/t \in \RR$, 
for each $i$ there is a number $0 < t_i < \varepsilon$ such that
for $q_i = p + t_i \vec{v}_i$
\begin{equation*}
| \frac{f(q_i)-f(p)}{t_i} - d_{\vec{v}_i}f(p) |
          \ < \ \frac{\eta}{2m} \ .
\end{equation*}
Hence $\sum_{i=1}^m (f(q_i)-f(p))/t_i < 0$, which gives 
\begin{equation} \label{FMD2}
f(p) \ > \ \frac{\sum_i f(q_i)/t_i}{\sum_i 1/t_i} \ .
\end{equation}
Put $\Gamma^{\prime} = \cup_{i=1}^m [p,q_i]$.
Then $V^{\prime} = r_{\Gamma^{\prime}}^{-1}(\Gamma_0^{\prime})$
is a simple subdomain of $U$.
Let $h(z)$ be the harmonic function on $V^{\prime}$ with $h(q_i) = f(q_i)$.
Since $\Delta_{\Gamma}(h)(p) = 0$, we have 
\begin{equation} \label{FMDD3}
h(p) \ = \ \frac{\sum_i f(q_i)/t_i}{\sum_i 1/t_i} \ .
\end{equation}
Combining (\ref{FMD2}) and (\ref{FMDD3}) contradicts the subharmonicity of $f$.

Thus for each $p \in \Gamma_0$,
\begin{equation} \label{FMD3}
\Delta_{\Gamma}(f)(p) \ = \ - \sum_i d_{\vec{v}_i}f(p) \ \le \ 0 \ .
\end{equation}
If $p$ is an endpoint of $\Gamma$ which does not belong
to $r_{\Gamma}(\partial U)$, and if $\vec{v}$ is the direction
vector at $p$, it has already been shown in (\ref{FCV1})
that $d_{\vec{v}}f(p) \ge 0$.  Hence for such points as well,
\begin{equation} \label{FMD5}
\Delta_{\Gamma}(f)(p) \ = \ -d_{\vec{v}}f(p) \ \le \ 0 \ .
\end{equation}

\vskip .05 in
Now consider an open segment $(x,y)$ contained in an edge of $\Gamma$.
Put $\Gamma^{\prime \prime} = [x,y]$.
The associated subgraph $V^{\prime \prime}$ is a simple subdomain of $U$
with boundary points $x$, $y$.
Let $h(z)$ be the harmonic function on $V^{\prime \prime}$
with $h(x) = f(x)$, $h(y) = f(y)$.  Let $T = \rho(x,y)$ be the length of
$\Gamma^{\prime \prime}$ and put $\alpha = (f(y)-f(x))/T$.
The restriction of $h(z)$ to
$\Gamma^{\prime \prime}$ is linear, and $f(z) \le h(z)$ for all $z \in \Gamma$.
Let $\vec{v}_+$ be the direction vector at $x$ pointing to $y$,
and let $\vec{v}_{-}$ be the direction vector at $y$ pointing to $x$.
Then 
\begin{eqnarray*}
d_{\vec{v}_+}f(x)
   & = & \lim_{t \rightarrow 0+} \frac{f(x + t \vec{v}_{+})-f(x)}{t}
   \ \le \ \lim_{t \rightarrow 0+} \frac{h(x + t \vec{v}_{+})-h(x)}{t}
   \ = \ \alpha \ , \\
d_{\vec{v}_-}f(y)
   & = & \lim_{t \rightarrow 0+} \frac{f(y + t \vec{v}_{-})-f(x)}{t}
   \ \le \ \lim_{t \rightarrow 0+} \frac{h(x + t \vec{v}_{+})-h(x)}{t}
   \ = \ -\alpha \ .
\end{eqnarray*}
so the segment $(x,y)$ has measure
\begin{equation} \label{FMD6}
\Delta_{\Gamma}(f)((x,y)) \ = \ d_{\vec{v}_+}f(x) + d_{\vec{v}_-}f(y)
                          \ \le \ 0 \ .
\end{equation}

\vskip .05 in
We can now show that $f \in \BDV(\Gamma)$, i.e. that
$\Delta_{\Gamma}(f)$ extends to a bounded Borel measure on $\Gamma$.
By the discussion in Section~\ref{Section E}, we must show
there is a number $B$ such that for any countable collection $\{T_i\}$
of pairwise disjoint sets in $\cA(\Gamma)$,
\begin{equation} \label{FMD7}
\sum_{i=1}^{\infty} |\Delta_{\Gamma}(f)(T_i)| \ \le \ B \ .
\end{equation}
Since each $T_i$ can be decomposed as a finite disjoint
union of points and open intervals, it suffices to prove (\ref{FMD7})
under the assumption that each $T_i$ is a point or an open interval.
Since $\Gamma$ has only finite many edges, endpoints, and branch points,
it also suffices to prove (\ref{FMD7}) assuming that all the $T_i$
are contained in the interior of an edge $e = [a,b]$.

In this case the fact that
$\Delta_{\Gamma}(f)$ is finitely additive, with 
$\Delta_{\Gamma}(f)(p) \le 0$ and $\Delta_{\Gamma}(f)((x,y)) \le 0$ for each 
point $p$ and open interval $(x,y)$ contained in $(a,b)$, means that
for any finite sum
\begin{equation*}
\sum_{i=1}^n |\Delta_{\Gamma}(T_i)| \ = \
- \sum_{i=1}^n \Delta_{\Gamma}(T_i) \ \le \ |\Delta_{\Gamma}((a,b))| \ .
\end{equation*}
Letting $n \rightarrow \infty$ gives (\ref{FMD7}).
The argument also shows that $\Delta_{\Gamma}(f)$
is $\le 0$ on $\Gamma \backslash r_{\Gamma}(\partial U)$.
Hence $\supp(\Delta_{\Gamma}(f)^+) \subset r_{\Gamma}(\partial U)$.
\end{proof}

\vskip .1 in
We can now prove Proposition \ref{PropF1}.

\vskip .1 in
\begin{proof} (of Proposition \ref{PropF1}):

First suppose $f$ is strongly subharmonic on $U$.  
We will show it is subharmonic.
Since $f$ is strongly upper semicontinuous on $U$, 
it is upper semicontinuous.   Let $W$ be a component of $U$.
Since $f \in \BDV(W)$, it follows that $f(x)$ is finite on each subgraph 
$\Gamma \subset W$, hence $f(x) \not\equiv -\infty$ on $W$.  
It remains to show that
if $V$ is a simple subdomain of $W$, and if $h(z)$ is harmonic on $V$
and satisfies $h(x) \ge f(x)$ on $\partial V$, then $h(z) \ge f(z)$ on $V$.

\vskip .05 in
First assume $V$ is a disc, and let $x$ be its boundary point.
For each $p \in V$, let $[x,p]$ be the path from $x$ to $p$.

First suppose $p \in V \backslash \PP^1(\CC_v)$.
Then $[x,p]$ is a subgraph $\Gamma$,
and $r_{\Gamma}(\partial W) = \{x\} \subset \partial \Gamma$.
By the definition of the measure $\Delta_{\Wbar}(f)$,
\begin{equation}
\Delta_{\Gamma}(f) \ = \ (r_{\Wbar,\Gamma})_*(\Delta_{\Wbar}(f)) \ .
\end{equation}
We are assuming $\Delta_{\Wbar}(f)|_{W} \le 0$, 
so $\Delta_{\Gamma}(f) \le 0$ on $\Gamma \backslash \{x\}$, and 
\begin{equation*} 
\supp(\Delta_{\Gamma}(f)^+) \ \subset \ \{x\} \ .
\end{equation*}
By Lemma \ref{LemF5}, either $f$ is constant on $\Gamma$, or $f(z)$
takes its maximum on $\Gamma$ at $x$.  In either case, $f(x) \ge f(p)$. 

Next suppose $p \in V \cap \PP^1(\CC_v)$.  By what has just been shown,
$f(x) \ge f(q)$ for each $q \in V \backslash \PP^1(\CC_v)$.
By condition (A) in the definition of strong subharmonicity,
\begin{equation*}
f(p) \ = \ \limsup
  \begin{Sb} q \rightarrow p \\ q \in V \backslash \PP^1(\CC_v) \end{Sb} f(q)
     \ \le \ f(x) \ .
\end{equation*}
Each harmonic function $h$ on $V$ is constant, so if $h(x) \ge f(x)$,
then $h(z) \ge f(z)$ for all $z \in V$.

\vskip .05 in
If $V$ is a simple subdomain corresponding to a subgraph $\Gamma \subset W$,
then  by Lemma \ref{LemF2}, $r_{\Gamma}(\partial W) \subset \partial \Gamma$.  
As before, $\Delta_{\Gamma}(f) = (r_{\Wbar,\Gamma})_*(\Delta_{\Wbar}(f))$
and $\Delta_{\Wbar}(f)|_W \le 0$, so
\begin{equation*}
\Delta_{\Gamma}(f)^+ \ \subset \ r_{\Gamma}(\partial W)
          \ \subset \ \partial \Gamma \ .
\end{equation*}
Let $h(z)$ be the harmonic function on $V$ for which
$h(x) = f(x)$ on $\partial V$.  Then $\Delta_{\Vbar}(h)|_V = 0$,
so by coherence
\begin{equation*}
\Delta_{\Gamma}(h) \ = \ r_{\Vbar,\Gamma}(\Delta_{\Vbar}(h))
\end{equation*}
is supported on $r_{\Vbar,\Gamma}(\partial V) = \partial \Gamma$.

Put $g(z) = f(z) - h(z)$.  Then $g(x) = 0$ on $\partial \Gamma$,
and $\Delta_{\Gamma}(g)|_{\Gamma_0} = \Delta_{\Gamma}(f)|_{\Gamma_0} \le 0$,
so $\Delta_{\Gamma}(g)^+ \subset \partial \Gamma$.  By Lemma \ref{LemF5},
$g(z) \le 0$ on $\Gamma$, that is $h(z) \ge f(z)$.

For each $z \in V \backslash \Gamma$, there is a unique $x \in \Gamma$
such that $z$ belongs to a branch off $\Gamma$ at $x$.  Let $V_x$ be the
connected component of $V \backslash \{x\}$ containing $z$.
Then $V_x$ is a disc.  By the same argument as before, 
$f(z) \le f(x)$, and since each harmonic function on a disc is constant, 
if $h$ is harmonic on $V_x$ and satisfies $h(x) \ge f(x)$, then
$h(z) \ge f(z)$.

Thus, if $f$ is strongly subharmonic, it is subharmonic.

\vskip .1 in
Now suppose $f$ is subharmonic on $U$.  We will show it is strongly
subharmonic on each subdomain $V \subset U$ with $\Vbar \subset U$.

It suffices to prove this for simple subdomains, since an arbitrary
subdomain $V$ with $\Vbar \subset U$ is contained
in a simple subdomain $V^{\prime}$.  (Cover $\Vbar$
with a finite number of punctured discs whose closures are contained
in $U$, and let $V^{\prime}$ be their union.  $V^{\prime}$ is connected
since $\Vbar$ is connected and each punctured disc is connnected;
$\Vbar^{\prime} \subset U$;  and $\partial V^{\prime}$ is finite,
since the boundary of each punctured disc is finite.)
If $f$ is strongly subharmonic on $V^{\prime}$,
then it is strongly subharmonic on $V$ since
\begin{equation*}
\Delta_{\Vbar}(f) \ =
\ (r_{\Vbar^{\prime},\Vbar})_*(\Delta_{\Vbar^{\prime}}(f)) \ .
\end{equation*}

So, let $V$ be a simple subdomain contained in a component $W$ of $U$.
We will first show that $f \in \BDV(V)$.
Let $\Gamma \subset V$ be an arbitrary subgraph.  We must show that
$f \in \BDV(\Gamma)$, and that there is a bound $B$ independent of $\Gamma$
such that $|\Delta_{\Gamma}(f)|(\Gamma) \le B$.
We will do this by enlarging $\Gamma$ to a graph $\tGamma$
which contains all the boundary points of $V$, and applying Lemma \ref{LemF8}.

\vskip .05 in
If $V$ is a disc, let $x$ be its boundary point.  Fix $q \in \Gamma$, 
and let $[x,q]$ be the path from $x$ to $q$.  It first meets $\Gamma$ at
a point $p$.  Put $\tGamma = \Gamma \cup [x,p]$, and 
let $\tV = r_{\tGamma}^{-1}(\tGamma_0)$ be the domain associated
to $\tGamma$.  Since $\Gamma_0 \subset V$, and $x$ is the only boundary
point of $V$, clearly $\tV \subset V$.  Hence the closure of $\tV$ is contained
in $\Vbar$, which in turn is contained in $W$, so $\tV$ is a simple subdomain
of $W$.  By Lemma \ref{LemF8}, $f|_{\tGamma} \in \BDV(\tGamma)$,
and $\Delta_{\tGamma}(f)^+ \subset r_{\tGamma}(\partial W) = \{x\}$.

Since $x$ is an endpoint of $\tGamma$, there is a single direction vector
$\vec{v}$ at $x$, and 
\begin{equation*}
\Delta_{\tGamma}(f)(x) \ = \ -d_{\vec{v}}f(x) \ .
\end{equation*}
Because $x$ is the only point 
where $\Delta_{\tGamma}(f)$ can have positive mass,
and since $\Delta_{\tGamma}(f)$ has total mass $0$, we see that
\begin{equation} \label{FZC1B}
|\Delta_{\tGamma}(f)|(\tGamma) \ = \ 2 \cdot |d_{\vec{v}}f(x)|
\end{equation}
Put $B = 2 \cdot |d_{\vec{v}}f(x)|$.  We claim that 
$B$ independent of the graph $\tGamma$.  Indeed, 
for any two points $p, p^{\prime} \in V$, the paths $[x,p]$
and $[x,p^{\prime}]$ must diverge at a point $y \in V$ 
and hence have an initial segment $[x,y]$ in common.  
This is because $V$ is connected:
if the paths diverged at $x$, then since $x \notin V$, $p$ and $p^{\prime}$
would lie in different components of $V$.

Taking the retraction to $\Gamma$, we see that 
\begin{equation*}
\Delta_{\Gamma}(f) \ = \ (r_{\tGamma,\Gamma})_*(\Delta_{\tGamma}(f))
\end{equation*}
has total mass at most $B$.
Since $\Gamma$ is arbitrary, $f \in \BDV(V)$.
Furthermore, $\Delta_{\Gamma}(f)^+$ is supported on
$r_{\Vbar,\Gamma}(x) = p$.  

Since $\Delta_{\Gamma}(f) = (r_{\Vbar,\Gamma})(\Delta_{\Vbar}(f))$ 
for each $\Gamma$, taking the limit over subgraphs $\Gamma$ shows 
$\Delta_{\Vbar}(f)^+$ is supported on $x$, that is,
$\Delta_{\Vbar}(f)|_V \le 0$.

\vskip .05 in
Next suppose $V$ is a simple subdomain of $W$
associated to a graph $\Gamma^{\prime}$.
Thus, $\Gamma_0^{\prime}$ is the main dendrite of $V$,
and $\partial V = \partial \Gamma^{\prime}$.

Let $\Gamma \subset V$ be an arbitrary subgraph.
If $\Gamma$ intersects $\Gamma^{\prime}$,
put $\tGamma = \Gamma \cup \Gamma^{\prime}$;
if not, take $x \in \Gamma^{\prime}_0$ and $p \in \Gamma$,
and put $\tGamma = \Gamma \cup [x,p] \cup \Gamma^{\prime}$.
Let $\tV$ be the simple domain associated to $\tGamma$.
Since all the boundary points of $V$ are contained in $\partial \tGamma$,
it follows that $\tV \subset V$.
Thus, $\tV$ is a simple subdomain of $W$.
By Lemma \ref{LemF8}, $f|_{\tGamma} \in \BDV(\tGamma)$, and
\begin{equation*}
\Delta_{\tGamma}(f)^+ \ \subset \ r_{\tGamma}(\partial W)
\ \subset \ \partial V \ = \ \partial \Gamma^{\prime}
\ \subset \ \partial \tGamma \ . 
\end{equation*}

Write $\partial V = \{x_1, \ldots, x_m\}$.
Since each $x_i$ is an endpoint of $\tGamma$,
there is a single direction vector $\vec{v}_i$ at $x_i$, and 
\begin{equation*}
\Delta_{\tGamma}(f)(x_i) \ = \ -d_{\vec{v}_i}f(x_i) \ .
\end{equation*}
The $x_i$ are the only points where $\Delta_{\tGamma}(f)$
can have positive mass.  Since $\Delta_{\tGamma}(f)$ has total mass $0$,
we see that
\begin{equation} \label{FWC1}
|\Delta_{\tGamma}(f)|(\tGamma) \ = \ 2 \cdot \sum_{i=1}^m |d_{\vec{v}_i}f(x_i)|
\end{equation}
The right side of (\ref{FWC1}) is a bound $B$ independent of the graph $\tGamma$.

Taking the retraction to $\Gamma$, it follows that
\begin{equation*}
\Delta_{\Gamma}(f) \ = \ (r_{\tGamma,\Gamma})_*(\Delta_{\tGamma}(f))
\end{equation*}
has total mass at most $B$.
Since $\Gamma$ is arbitrary, $f \in \BDV(V)$.
Furthermore $\Delta_{\Gamma}(f)^+$ is supported on
$r_{\Vbar,\Gamma}(\partial V)$.  Since
\begin{equation*}
\Delta_{\Gamma}(f) \ = \ (r_{\Vbar,\Gamma})(\Delta_{\Vbar}(f))
\end{equation*}
for each $\Gamma$, taking the limit over subgraphs $\Gamma$ shows 
$\Delta_{\Vbar}(f)^+$ is supported on $\partial V$, that is,
$\Delta_{\Vbar}(f)|_V \le 0$.

\vskip .05 in
It remains to show that $f$ is strongly upper semicontinuous, and that
\begin{equation} \label{FCV2}
f(p) \ = \ 
\limsup \begin{Sb} z \rightarrow p \\ z \in V \backslash \PP^1(\CC_v) \end{Sb}
          f(z) 
\end{equation}
for each $p \in V \cap \PP^1(\CC_v)$.
By assumption $f$ is upper semicontinuous,
so for each $p \in V$
\begin{equation} \label{FCV3}
f(p) \ \ge \ \limsup_{z \rightarrow p} f(z) \ .
\end{equation}

Fix $p \in V$;  first suppose $p \in \PP^1(\CC_v)$.  Fix $y \in V$,
and let $x$ approach $p$ along the path $[y,p]$.  If $x$ is close enough
to $p$, then the connected component $V_x$ of $\PP^1_{\Berk} \backslash \{x\}$
containing $p$ is a disc whose closure is contained in $W$.  Since $f$ is
subharmonic, $f(z) \le f(x)$ for all $z \in V_x$.  These discs $V_x$ form
a cofinal sequence of neighborhoods of $p$, so
\begin{equation*}
f(p) \ \le \ \limsup_{z \rightarrow p} f(z) \ = \
\lim \begin{Sb} x \rightarrow p \\ x \in [y,p) \end{Sb} f(x) \ .
\end{equation*}
Combined with (\ref{FCV3}), this shows $f$ is strongly upper semicontinuous 
at $p$. It also establishes (\ref{FCV2}), and indeed shows that
\begin{equation} \label{FBB1}
f(p) \ = \ \lim \begin{Sb} x \rightarrow p \\ x \in [y,p) \end{Sb} f(x) \ .
\end{equation}

Now suppose $p \in V \backslash \PP^1(\CC_v)$, 
but assume $p$ does not belong to the main dendrite of $V$.  
If $V$ is a disc, let $x$ be the unique boundary
point of $x$;  otherwise, let $x$ be the point on the main dendrite
of $V$ where the branch containing $p$ is attached.
Consider the path $\Gamma = [x,p]$.  It is a subgraph of $W$,
and $r_{\Gamma}(\partial V) = x$.  The domain associated to $\Gamma$
is a simple subdomain of $W$, so $f|_{\Gamma}$ is continuous by Lemma
\ref{LemF7}. Thus 
\begin{equation} \label{FCV4}
\lim \begin{Sb} y \rightarrow p \\ y \in \Gamma \end{Sb} f(y) \ = \ f(p) \ .
\end{equation}
For each $y \in (x,p)$, the connected component $V_y$ of
$\PP^1_{\Berk} \backslash \{y\}$ containing $p$ 
is a disc with closure $\Vbar_y \subset V$.
Since $f$ is subharmonic, $f(z) \le f(y)$ for all $z \in V_y$.  Combined with
(\ref{FCV4}), this shows $\limsup_{z \rightarrow p} f(z) \le f(p)$.
Hence $f$ is strongly upper semicontinuous at $p$.

Finally, suppose $p$ belongs to the main dendrite $\Gamma$ of $V$.
Since $V$ is a simple subdomain of $W$, Lemma \ref{LemF7} says 
 $f$ is continuous on $\Gamma$.
Given $\varepsilon > 0$, take a neighborhood
$\Gamma_p(\eta) = \{x \in \Gamma : \rho(x,p) < \eta \}$ where 
$f(x) \< f(p) + \epsilon$.  Put $\Gamma^{\prime} = \overline{\Gamma_p(\eta)}$.
Since $p \in \Gamma_0$,  we can assume without loss that
$\Gamma^{\prime} \subset \Gamma_0$.
Then $r_{\Gamma^{\prime}}(\partial V) \subset \partial \Gamma^{\prime}$, so
$V^{\prime} = r_{\Gamma^{\prime}}^{-1}(\Gamma_0^{\prime})$
is a simple subdomain of $V$.  As we have seen before,
$f(z)$ is non-increasing on branches off the main dendrite,
so $f(z) \le f(p) + \varepsilon$ for all $z \in V^{\prime}$.
Since $\varepsilon > 0$ is arbitrary, 
\begin{equation*}
\limsup_{z \rightarrow p} f(z) \ \le \ f(p) \ .
\end{equation*}
Combined with (\ref{FCV3}), this shows $f$ is strongly upper semicontinuous.

\vskip .1 in
The final assertion in Proposition \ref{PropF1} is that if the measures
$|\Delta_{\Vbar}(f)|(\Vbar)$ are uniformly bounded for all simple subdomains
$V$ of $U$, then $f$ is strongly subharmonic on $U$.  

This is trivial.  Fix a component $W$ of $U$.    
Since simple subdomains exhaust $W$, each subgraph $\Gamma \subset W$ 
is contained in some simple subdomain $V$, 
and $\Delta_{\Gamma}(f) = r_{\Vbar,\Gamma}(\Delta_{\Vbar}(f))$.
Thus the measures $|\Delta_{\Gamma}(f)|$ are uniformly bounded, 
so $f \in \BDV(W)$.  For each $V$, the retraction map $r_{\Wbar,\Vbar}$ 
takes $\partial W$ to $\partial V$ and fixes $V$.  Hence 
$\Delta_{\Wbar}(f)|_V = \Delta_{\Vbar}(f)|V$.  
Since $\Delta_{\Vbar}(f)|_v \le 0$, 
it follows that $\Delta_{\Wbar}(f)|_W \le 0$.  
Finally, the semicontinuity assertions for $f$ on $W$ follow 
from those on the subdomains $V$.  
\end{proof}

\vskip .1 in
We record the following facts shown in the proof of Proposition \ref{PropF1}:  

\begin{corollary} \label{CorF9}
Let $f$ be subharmonic on a simple domain $U$.  Then $f(x)$ is non-increasing
on paths off the main dendrite of $U$.  If $U$ is a disc, then $f(x)$ is
non-increasing on paths away from the boundary $\partial U = \{q\}$.

For each $p \in U \cap \PP^1(\CC_v)$, and for any path $[y,p] \subset U$ 
\begin{equation*}
\lim \begin{Sb} x \rightarrow p \\ x \in [y,p) \end{Sb} f(x) \ = \ f(p) \ .
\end{equation*}
\end{corollary}

A function $f : U \rightarrow \RR \cup \{\infty\}$
will be called {\it strongly superharmonic} on $U$
if $-f$ is strongly subharmonic.
It will be called {\it superharmonic} if $-f$ is subharmonic.
These concepts can be reformulated via lower semicontinuity
and submajorization by harmonic functions, as in Definition \ref{DefF2}.

\vskip .1 in
Here are some examples of subharmonic and superharmonic functions.

\begin{example}
\label{Example F.1}
A function $f : U \rightarrow \RR$
is harmonic in an open set $U \subset \PP^1_{\Berk}$
if and only if is both subharmonic and superharmonic.

The only functions which are subharmonic on all of $\PP^1_{\Berk}$ are
the constant functions.  Indeed, a function $f(x)$ subharmonic on all of 
$\PP^1_{\Berk}$ is strongly subharmonic by Proposition \ref{PropF1},
hence $\Delta_{\PP^1_{\Berk}}(f) \le 0$.  Since the total mass
of the Laplacian is $0$, this means $\Delta_{\PP^1_{\Berk}}(f) = 0$.
Thus, $f$ is harmonic on $\PP^1_{\Berk}$ and so is constant.
\end{example}

\begin{example}
\label{Example F.2}
For fixed $a, \zeta \in \PP^1_{\Berk}$ 
with $a \ne \zeta$, $f(x) = \log_v(\delta(x,a)_{\zeta})$ 
is strongly subharmonic in $\PP^1_{\Berk} \backslash \{\zeta\}$, 
and strongly superharmonic in $\PP^1_{\Berk} \backslash \{a\}$.  
Indeed, $\delta(x,a)_{\zeta}$ is continuous by Proposition~\ref{PropC10}, and
\begin{equation*}
\Delta_{\PP^1_{\Berk}}(\log_v(\delta(x,a)_{\zeta}) \ = \
\delta_{\zeta}(x) - \delta_{a}(x)
\end{equation*}
by Example~\ref{Example E.2}.

Correspondingly, $-\log_v(\delta(x,a)_{\zeta})$ is strongly superharmonic in
$\PP^1_{\Berk} \backslash \{\zeta\}$, and is strongly subharmonic in
$\PP^1_{\Berk} \backslash \{a\}$.
\end{example}

\begin{example}
\label{Example F.3}
If $f \in \CC_v(T)$ is a nonzero rational
function with divisor $\div(f) = \sum_{i=1}^m n_i \delta_{a_i}(x)$,
let $\supp^-(\div(f))$,  $\supp^+(\div(f))$ be its be its polar locus and
zero locus, respectively.  Then $\log_v([f]_x)$ is strongly subharmonic on
$\PP^1_{\Berk} \backslash \supp^-(\div(f))$ and strongly superharmonic
on $\PP^1_{\Berk} \backslash \supp^+(\div(f))$.  Likewise $-\log_v([f]_x)$
is strongly superharmonic on
$\PP^1_{\Berk} \backslash \supp^-(\div(f))$ and strongly subharmonic
on $\PP^1_{\Berk} \backslash \supp^+(\div(f))$.

These assertions follow from the continuity of $[f]_x$ and from 
Example~\ref{Example E.3}.
\end{example}

\begin{example}
\label{Example F.4}
If $\nu$ is a probability measure on
$\PP^1_{\Berk}$ and $\zeta \notin \supp(\nu)$,
then the potential function $u_{\nu}(x,\zeta)$
is strongly superharmonic in $\PP^1_{\Berk} \backslash \{\zeta\}$
and is strongly subharmonic in $\PP^1_{\Berk} \backslash \supp(\nu)$.
These assertions follow from Proposition~\ref{PropD7} and its proof,
and Example~\ref{Example E.4}.
\end{example}

\begin{example}
\label{Example F.5}
If $E \subset \PP^1_{\Berk}$ is a compact set of positive capacity
and $\zeta \notin E$, then the Green's function $G(z,\zeta;E)$
is strongly subharmonic on $\PP^1_{\Berk} \backslash \{\zeta\}$.
Indeed, if $\mu_{\zeta}$ is the equilibrium distribution of $E$
with respect to $\zeta$,
then $G(z,\zeta;E) = V_{\zeta}(E) - u_{\mu_{\zeta}}(z,\zeta)$.
\end{example}

\subsection{Locality.} \ 

Subharmonicity is a local property:

\begin{proposition} \label{PropF10}
Let $U \subset \PP^1_{\Berk}$ be open.
A function $f : U \rightarrow \RR \cup \{-\infty\}$
is subharmonic on $U$ if and only if for each $x \in U$,
there is a neighborhood $V_x$ of $x$ in $U$ such that $f|_{V_x}$
is subharmonic on $V_x$.
\end{proposition}

\begin{proof}  Only the direction ($\Longleftarrow$) requires attention.
Suppose that for each $x \in X$ there is a neighborhood $V_x$
such that $f$ is subharmonic on $V_x$.  Then $f$ is upper semicontinuous
on $U$, and $f(x) \not\equiv -\infty$ on any connected component of $U$,
since these properties hold for the $V_x$.

It remains to show that if $V$ is a simple subdomain of $U$,
and if $h(x)$ is a harmonic function on $V$ with $h(x) \ge f(x)$ on
$\partial V$, then $h(x) \ge f(x)$ on $V$.  For each $V_x$, there
is a simple subdomain $W_x$ of $V_x$ with $x \in W_x$.
By Proposition \ref{PropF1} $f$ is strongly subharmonic on $W_x$.
Since $\Vbar \subset U$ is compact,
we can cover $\Vbar$ with a finite number of
sets $W_{x_1}, \ldots, W_{x_M}$.  Let $W$ be their union.
We claim that $f$ is strongly subharmonic on $W$.  Property (A) in
the definition of strong subharmonicity is automatic, since it is
inherited from the $W_{x_i}$.  For property (B), note that
$f \in \BDV(W)$ since $f \in \BDV(W_{x_i})$ for each $i$.
Since $f$ is strongly subharmonic on $W_{x_i}$
\begin{equation*}
\Delta_{\Wbar}(f)|_{W_{x_i}}
\ = \ (r_{\Wbar,\Wbar_{x_i}})_*(\Delta_{\Wbar}(f))|_{W_{x_i}}
\ = \ \Delta_{\Wbar_{x_i}}(f)|_{W_{x_i}} \ \le \ 0
\end{equation*}
so $\Delta_{\Wbar}(f)|_W \le 0$ since the $W_{x_i}$ cover $W$.
Thus $f$ is strongly subharmonic on $W$.

By Proposition \ref{PropF1}, $f$ is subharmonic on $W$.
However, $\Vbar \subset W$, so $V$ is a simple subdomain of $W$.
By (E) in definition of subharmonicity, $f(x) \le h(x)$ on $V$.
\end{proof}

\subsection{Stability properties.}  

We will now show that subharmonic functions on $\PP^1_{\Berk}$
are stable under the same 
operations as classical subharmonic functions (see \cite{Kl}, p.49).  

Recall that if $f : U \rightarrow \RR$ is a
function on a topological space $U$ which is locally bounded from above, 
then the {\it upper semicontinuous regularization} $f^*(x)$ of $f(x)$ is
defined by 
\begin{equation*}
f^*(x) \ = \ \max(f(x),\limsup_{z \rightarrow x} f(z)) \ .
\end{equation*}
It is easy to check that $f^*(x)$ is upper semicontinuous, and that
if $g$ is upper semicontinuous and $g \ge f$, then $g \ge f^*$.  
In particular $(f^*)^* = f^* \ge f$.  

If $U \subset \PP^1_{\Berk}$ is open, we will write $\SH(U)$ for the
set of subharmonic functions on $U$. 

\begin{proposition} \label{PropF11}
Let $U \subset \PP^1_{\Berk}$ be open.

$A)$   $\SH(U)$ is a convex cone:  
      if $0 \le \alpha, \beta \in \RR$, and if $f, g \in \SH(U)$, 
      then $\alpha f + \beta g \in \SH(U)$.  
      
$B)$  If $U$ is connected, and if $\{f_j\}_{j \in \NN}$ is a decreasing 
      sequence of functions in $\SH(U)$, put 
      $f(x) = \lim_{j \rightarrow \infty} f_j(x)$. 
      Then either $f \in \SH(U)$, or $f(x) \equiv -\infty$ on $U$.
      
$C)$  If $\{f_j\}_{j \in \NN}$ is a sequence of functions in $\SH(U)$ which
      converge uniformly to a function $f : U \rightarrow \RR$ on 
      compact subsets of $U$, then $f \in \SH(U)$. 

$D)$  If $f, g \in \SH(U)$, then $\max(f,g) \in \SH(U)$.     
      
$E)$  If $\{f_{\alpha}\}_{\alpha \in A}$ is a family of functions in $\SH(U)$
      which is locally bounded from above, 
      and if $f(x) = \sup_{\alpha} f_{\alpha}(x)$, 
      then $f^*(x) \in \SH(U)$.  Furthermore $f^*(x) = f(x)$ for all
      $x \in U \backslash \PP^1(\CC_v)$.

$F)$  If $U$ is connected and $\{f_n\}_{n \ge 0}$
      is a sequence of functions in $\SH(U)$
      which is locally bounded from above,
      put $f(x) = \limsup_{n \rightarrow \infty} f_n(x)$.
      Then either $f(x) \equiv -\infty$ on $U$, or $f^*(x) \in \SH(U)$.
      Furthermore $f^*(x) = f(x)$ for all $x \in U \backslash \PP^1(\CC_v)$.
\end{proposition}                           

\begin{proof}  Except for the last assertion in (E), the proofs of these 
are the same as their classical counterparts, and rely on general properties 
of semicontinuity and domination by harmonic functions.

\vskip .05 in
(A)   If $f,g \in \SH(U)$ then $f$ and $g$ are
upper semicontinuous, and $\alpha f + \beta g$ is upper semicontinuous
since $\alpha, \beta \ge 0$.
Neither $f$ nor $g$ is $-\infty$ on $U \backslash \PP^1(\CC_v)$,
so $\alpha f + \beta g \ne -\infty$ on $U \backslash \PP^1(\CC_v)$,
and certainly $\alpha f + \beta g \not\equiv -\infty$ on any
component of $U$.  If $V$ is a simple subdomain of $U$, let
$\partial V = \{x_1, \ldots, x_m \}$.  Suppose $h$ is harmonic on $V$
with $h(x_i) \ge \alpha f(x_i) + \beta g(x_i)$ for each $i$.
Let $h_1$ be the harmonic function on $V$ with $h_1(x_i) = f(x_i)$
on $\partial V$, and let $h_2$ be the harmonic function on $V$
with $h_2(x_i) = g(x_i)$ on $\partial V$.  Then $h_1(z) \ge f(z)$ on $V$,
and $h_2(z) \ge g(z)$ on $V$.  Put $H = h - \alpha h_1 - \beta h_2$.
Then $H(x_i) \ge 0$ on $\partial V$, so by the Maximum Principle for
harmonic functions, $H(z) \ge 0$ on $V$.  Hence 
$\alpha f(z) + \beta g(z) \le \alpha h_1(z) + \beta h_2(z) \le h(z)$
on $V$.

\vskip .05 in
(B)  Suppose $U$ is connected, and let $\{f_j\}_{j \in \NN}$ be a decreasing
sequence of subharmonic functions on $U$.
Put $f(z) = \lim_{j \rightarrow \infty} f_j(z)$.  Then $f$ is upper
semicontinuous since each $f_j$ is.
By the same argument as in Lemma \ref{LemF6},
either $f(x) \ne -\infty$ on $U \backslash \PP^1(\CC_v)$,
or $f(x) \equiv -\infty$ on $U$.

Suppose $f \not \equiv -\infty$.  If $V$ is a simple subdomain of $U$
and $h$ is a harmonic function on $V$ with $f(x_i) \le h(x_i)$ on the
finite set $\partial V$, then for each $\varepsilon > 0$ there is an $N$
such that $f_N(x_i) \le h(x_i) + \varepsilon$ and all $x_i \in \partial V$.
It follows that $f(z) \le f_N(z) \le h(z) + \varepsilon$ on $V$.
Since $\varepsilon > 0$ is arbitrary,  $f(z) \le h(z)$.

\vskip .05 in
(C)  Suppose $\{f_j\}_{j \in \NN}$ is a sequence of subharmonic functions
converging uniformly to a function $f : U \rightarrow \RR$ on compact
subsets of $U$.  Since each $f_j$ is upper semicontinuous, $f$ is upper
semicontinuous.  By hypothesis $f \not \equiv -\infty$ on any component of $U$.

Let $V$ be a simple subdomain of $U$, and let $h$ be a harmonic function
on $V$ with $f(x_i) \le h(x_i)$ on $\partial V$.
Take $\varepsilon > 0$.  Since $\partial V$ is finite,
there is an $N_1$ so that $f_j(x_i) \le h(x_i) + \varepsilon$ on $\partial V$,
for all $j \ge N_1$.  Thus $f_j(z) \le h(z) + \varepsilon$ on $V$ for each
$j \ge N_1$.  Since the $f_j$ converge uniformly to $f$ on compact subsets,
there is an $N_2$ such that $|f(z) - f_j(z)| < \varepsilon$ on $\Vbar$
for all $j \ge N_2$.  Taking $j \ge \max(N_1,N_2)$,
we see that $f(z) \le h(z) + 2 \varepsilon$ on $V$.  Since $\varepsilon > 0$
is arbitrary, $f(z) \le h(z)$ on $V$.

\vskip .05 in
(D)  Suppose $f, g \in \SH(U)$, and put $F(z) = \max(f(z), g(z))$.
Since $f$ and $g$ are upper semicontinuous, so is $F$.  Since neither
$f$ nor $g$ is $\equiv -\infty$ on any component of $U$, the same is
true for $F$.  If $V$ is a simple subdomain of $U$ and $h$ is a harmonic
function on $U$ with $h(x_i) \ge F(x_i)$ on $\partial V$, then
$h(x_i) \ge f(x_i)$ and $h(x_i) \ge g(x_i)$ on $\partial V$,
so $h(z) \ge f(z)$ and $h(z) \ge g(z)$ on $V$, which means that
$h(z) \ge \max(f(z),g(z)) = F(z)$ on $V$.

\vskip .05 in
(E)  Let $\{f_{\alpha}\}_{\alpha \in A}$ be a family of subharmonic functions
on $U$, and put $f(z) = \sup_{\alpha} (f_{\alpha}(z))$ on $U$.
Then $f^*(z)$ is upper semicontinuous on $U$.  Since no $f_{\alpha}$ is
$\equiv -\infty$ on any component of $U$, the same is true for $f$,
and also for $f^* \ge f$.  If $V$ is a simple subdomain of $U$ and $h$
is harmonic on $U$ with $h(x_i) \ge f^*(x_i)$ on $\partial V$, then
$h(x_i) \ge f_{\alpha}(x_i)$ on $\partial V$ for each $\alpha \in A$.
It follows that $h(z) \ge f_{\alpha}(z)$ on $V$ for each $\alpha$,
so $h(z) \ge f(z)$ on $V$.  However, $h(z)$ is continuous on $V$, hence
certainly upper semicontinuous, so $h(z) \ge f^*(z)$ on $V$ by the
properties of the upper semicontinuous regularization.

We will now show that for each $x \in U \backslash \PP^1(\CC_v)$,
\begin{equation}  \label{FGG0}
f(x) \ \ge \ \limsup_{z \rightarrow x} f(z) \ .
\end{equation}
We first construct a simple subdomain $V$ of $U$ such that
$x$ lies on the main dendrite of $V$.  

If $x$ is on the main dendrite of $U$,
let $V$ be any simple subdomain of $U$ containing $x$.
Then $x$ is on the main dendrite of $V$, and in particular is an interior
point of the graph $\Gamma$ with $V = r_{\Gamma}^{-1}(\Gamma_0)$.
If $x$ is not on the main dendrite of $U$, 
let $V_0$ be an open subdisc of $U$ containing $x$.
Let $x_0$ be the boundary point of $V_0$,
and let $x_1 \in V_0 \backslash \PP^1(\CC_v)$ be a point such that
$x$ is in the interior of the path $\Gamma = [x_0,x_1]$.
Put $V = r_{\Gamma}^{-1}(\Gamma_0)$;
then $x$ lies on the main dendrite of $V$.

Next, we claim that $f|_{\Gamma_0}$ is continuous.  Indeed, for each
subgraph $\Gamma^{\prime} \subset \Gamma_0$ of the type considered
inthe proof of Lemma \ref{LemF7}, let $V^{\prime}$ be the 
corresponding simple subdomain of $V$.
Let $h(x)$ be the harmonic function on $V^{\prime}$ whose value
at each $p \in \partial V^{\prime} = \partial \Gamma^{\prime}$ is $f(p)$.
For each $\alpha \in A$ and each $p \in \partial \Gamma^{\prime}$,
$f_{\alpha}(p) \le f(p)$.  Hence $f_{\alpha}(z) \le h(z)$ for all 
$z \in V^{\prime}$, and in turn $f(z) = \sup_{\alpha} f_{\alpha}(z) \le h(z)$.
By the same argument as in proof of Lemma \ref{LemF7}, 
the directional derivative $d_{\vec{v}}f(p)$
exists for each $p \in \Gamma_0$ and each $\vec{v}$ at $p$. 
Thus $f$ is continuous on $\Gamma_0$.

Now fix $\varepsilon > 0$ and let
$\Gamma_p(\delta) = \{ q \in \Gamma : \rho(q,x) < \delta\}$ be a
neighborhood of $x$ in $\Gamma$ on which $f(q) < f(x) + \varepsilon$.
Let $V^{\prime \prime}$ be the simple subdomain of $V$ associated to
$\overline{\Gamma_p(\delta)}$.  For each $\alpha \in A$ and
each $q \in \Gamma_p(\delta)$, $f_{\alpha}(q) \le f(x) + \varepsilon$.
By Corollary \ref{CorF9}, $f_{\alpha}$ is non-increasing on paths off
the main dendrite $\Gamma_p(\delta)$ of $V^{\prime \prime}$,
so $f_{\alpha}(z) \le f(x) + \varepsilon$ on $V^{\prime \prime}$.
Hence $f(z) = \sup_{\alpha} f_{\alpha}(z) \le f(x) + \varepsilon$
for each $z \in V^{\prime \prime}$.
Since $\varepsilon$ is arbitrary, this gives (\ref{FGG0}),
and shows that $f^*(x) = f(x)$.

\vskip .05 in
(F)  For each $n$, put $F_n(z) = \sup_{m \ge n} f_m(x)$ 
and let $F_n^*(z)$ be the
upper semicontinuous regularization of $F_n$.  By the final assertion in 
part (E), $F_n^*(x) = F_n(x)$ for each $x \in U \backslash \PP^1(\CC_v)$.
Then, $F_1^*(z) \ge F_2^*(z) \ge \cdots$ is a decreasing sequence of
subharmonic functions, and $f(x) = \lim_{n \rightarrow \infty} F_n^*(x)$,
so by part (B) either $f(x) \equiv -\infty$, or $f(x) \in \SH(U)$.
\end{proof}

\vskip .1 in
Subharmonic functions are also stable under integration over
suitably bounded families on a parameter space.

\begin{proposition} \label{PropF12}
Let $\mu$ be a non-negative $\sigma$-finite measure on a measure space $T$,
and let $U \subset \PP^1_{\Berk}$ be a domain.
Suppose that $F: U \times T \rightarrow [-\infty,\infty)$ is a measurable
function such that

$A)$  For each $t \in T$, the function
$F_t(z) = F(z,t) : U \rightarrow [-\infty,\infty)$ is subharmonic in $U$; 

$B)$ there is a majorizing function $g : T \rightarrow (-\infty,\infty]$
which belongs to $L^1(\mu)$ such that $F_t(z) \le g(t)$ 
for all \ $t \in T$ and all $z \in X$.

Then the function  
\begin{equation*}
f(z) \ := \ \int_T F(z,t) \, d\mu(t)  
\end{equation*}
is either subharmonic in $U$, or is $\equiv -\infty$ on $U$.  
\end{proposition}

\begin{proof}  (See \cite{Kl}, Theorem 2.6.5, p.51.)
Suppose $f(z) \not \equiv -\infty$.

Fix $x \in U$, and let $z_1, z_2, \ldots$ be a sequence converging to $x$.
Fatou's Lemma, applied to the sequence of functions
$h_n(t) = F(z_n,t) - g(t)$ on $X$, implies that
\begin{equation*}
\int_T (\limsup_{n \rightarrow \infty} h_n(t)) \, d\mu(t)
\ \ge \ \limsup_{n \rightarrow \infty} (\int_T h_n(t) \, d\mu(t)) \ .
\end{equation*}
However, for each $t$
\begin{equation*}
\limsup_{n \rightarrow \infty} h_n(t)
\ = \ \limsup_{n \rightarrow \infty} (F(z_n,t)) - g(t) 
\ \le \ F(x,t) - g(t)
\end{equation*}
by the upper semicontinuity of $F_t(x)$.  Since $\int_T g(t) \, d\mu(t)$
is finite, this gives 
\begin{equation*}
f(x) \ \ge \ \limsup_{n \rightarrow z_n} f(z_n) \ ,
\end{equation*}
so $f(x)$ is upper semicontinuous.

Now let $V \subset U$  be a simple subdomain,
with boundary $\partial V = \{x_1, \ldots, x_m\}$.
Let $h$ be harmonic on $V$
with $h(x_i) \ge f(x_i)$ on $\partial V$.  Let $h_1(z), \ldots, h_m(z)$
be the harmonic measures for $V$ with $h_i(x_j) = \delta_{i,j}$.
Recall that for each $z \in V$,
$0 \le h_i(z) \le 1$ and $\sum_{i=1}^m h_i(z) = 1$.
Since each $F_t(x)$ is subharmonic, for each $z \in V$
\begin{equation*}
F_t(z) \ \le \ \sum_{i=1}^m h_i(z) F_t(x_i) \ .
\end{equation*}
Integrating over $\mu$ gives 
\begin{equation*}
f(z) \ \le \ \sum_{i=1}^m f(x_i) h_i(z)  
     \ \le \ \sum_{i=1}^m h(x_i) h_i(z) \ = \ h(z)
\end{equation*}
Thus $f$ is subharmonic.
\end{proof}

\vskip .1 in
Here are more ways of getting new subharmonic functions from old ones:

\begin{lemma} \label{LemF13}
Let $U \subset \PP^1_{\Berk}$ be open.
If $f$ is subharmonic on $U$,
and $\varphi : \RR \rightarrow \RR$ is convex up and non-decreasing,
then $\varphi \circ f$ is subharmonic on $U$.
$($Here, $\varphi(-\infty)$ is to be interpreted as
$\lim_{t \rightarrow -\infty} \varphi(t)$.$)$
\end{lemma}

\begin{proof}  (See \cite{Kl}, Theorem 2.6.6, p.51.)  Note that
$\varphi(t)$ can be written as
\begin{equation*}
\varphi(t) \ = \ \sup( \{ a \cdot t + b :  a \ge 0, b \in \RR,
     \ \text{and $a \cdot t + b \le \varphi(t)$ for all $t \in \RR$} \} ) \ .
\end{equation*}
Let $A$ be the corresponding set of pairs  $(a,b) \in \RR^2$;  then
\begin{equation*}
\varphi \circ f(z) \ = \ \sup_{(a,b) \in A} a \cdot f(z) + b \ .
\end{equation*}
For each $(a, b) \in A$, we have $a \cdot f(z) + b \in \SH(U)$.
By Proposition \ref{PropF11}(E), if $F(z) = \varphi \circ f(z)$,
then $F^*(z)$ is subharmonic.

A convex function on $\RR$ is automatically continuous, since its 
one-sided derivatives exist at each point.  We claim that 
$F(z) = \varphi(f(z))$ is upper semicontinuous.
This holds because $\varphi$ is continuous and nondecreasing,
and $f(z)$ is upper semicontinuous.  Hence  $F^*(z) = F(z)$.
\end{proof}

\begin{corollary} \label{CorF14}
Let $U \subset \PP^1_{\Berk}$ be open.

$A)$  If $f$ is subharmonic on $U$ and $q \ge 1$, then the function 
$F(z) = q^{f(z)}$ is subharmonic on $U$.

$B)$  If $f$ is subharmonic and non-negative on $U$, then for any
$\alpha \ge 1$, the function $F(z) = f(z)^{\alpha}$ is subharmonic on $U$.

\end{corollary}

\begin{proof} (See \cite{Kl}, Corollary 2.6.8, p.52).
Note that $t \rightarrow q^t$ and $t \rightarrow t^{\alpha}$
are convex up and non-decreasing, and apply Lemma \ref{LemF13}.
\end{proof}

\vskip .1 in
We can now give additional examples of subharmonic functions.

\begin{example}
\label{Example F.6}
Fix $\zeta \in \PP^1_{\Berk}$.
For each $\alpha > 0$, and each $a \ne \zeta$, the function  
$f(x) = \delta(x,a)_{\zeta}^{\alpha}$ is subharmonic in
$\PP^1_{\Berk} \backslash \{\zeta\}$.  In particular 
this applies to $\delta(x,a)_{\zeta}$.
More generally, for each $\alpha_1, \ldots, \alpha_n \ge 0$,
and each $a_1, \ldots, a_n \in \PP^1_{\Berk} \backslash \{\zeta\}$, 
the generalized pseudo-polynomial
\begin{equation*}
P(x,\vec{\alpha},\vec{a}) 
\ = \ \prod_{i=1}^m \delta(x,a_i)_{\zeta}^{\alpha_i}
\end{equation*}
is subharmonic in $\PP^1_{\Berk} \backslash \{\zeta\}$.

This follows from Corollary \ref{CorF14}(A) and Proposition \ref{PropF11}(A),
taking $q = q_v$, since $\alpha_i \cdot \log_v(\delta(x,a_i)_{\zeta})$
is subharmonic in $\PP^1_{\Berk} \backslash \{\zeta\}$.

\vskip .05 in
In particular, consider $f(x) = \delta(x,0)_{\infty}$ on $\AA^1_{\Berk}$.
It is constant on branches off the path $[0,\infty]$.  
Give $(0,\infty)$ the arclength parametrization,  
so that $x_t$ is the point corresponding to the disc $B(0,q_v^{t})$   
for $-\infty < t < \infty$.  Then $f(x_t) = q_v^t$.
For each disc $V_T := \cB(0,q_v^{T})^-$, 
the Laplacian $\Delta_{V_T}(f)$ is supported on $[0,x_T]$.
Relative to the arclength parametrization, 
\begin{eqnarray*}
\Delta_{V_T}(f) 
& = & f^{\prime}(x_T) \cdot \delta_T(t) - f^{\prime \prime}(x_t) dt \\
& = & q_v^{T} \log(q_v) \cdot \delta_T(t) - q_v^{t} (\log(q_v))^2 dt \ .
\end{eqnarray*}
The total variation of these measures grows to $\infty$ 
as $T \rightarrow \infty$.
Thus $f(x)$ is subharmonic, but not strongly subharmonic, on $\AA^1_{\Berk}$.
\end{example}

\begin{example}
\label{Example F.7}
Consider the function $\varphi(t) = \arcsin(q_v^t)$.    
It is  bounded, increasing and convex up on
$[-\infty,0)$, with a vertical tangent at $t = 0$.  

Put  $f(x) = \varphi(\log_v(\delta(x,0)_{\infty}))$ on $U := \cB(0,1)^-$.  
By Lemma \ref{LemF13}, $f(x)$ is bounded and subharmonic on $U$.
However, it is not strongly subharmonic,
and it cannot be extended to a subharmonic function on any larger domain.
\end{example}

\begin{example}
\label{Example F.8}
For a nonzero rational function $f \in \CC_v(T)$ with divisor
$\div(f)$, the function $F(x) = [f]_x$ is subharmonic on
$\PP^1_{\Berk} \backslash \supp(\div(f)^-)$.

This follows from Corollary \ref{CorF14}(A), taking $q = q_v$,
since $\log_v([f]_x)$ is subharmonic on the complement of $\supp(\div(f)^-)$.
\end{example}

\begin{example}
\label{Example F.9}
Let $E \subset \PP^1_{\Berk}$ be a compact set of positive capacity,
and take $\zeta \in \PP^1_{\Berk} \backslash E$.
Then for each $\alpha \ge 1$, the function $G(x,\zeta;E)^{\alpha}$
is subharmonic in $\PP^1_{\Berk} \backslash \{\zeta\}$.

This follows from Corollary \ref{CorF14}(B), taking $\varphi(t) = t^{\alpha}$,
since $G(x,\zeta;E)$ is subharmonic and non-negative
on $\PP^1_{\Berk} \backslash \{\zeta\}$.
\end{example}

\begin{example}
\label{Example F.10}
Let $f_1, \ldots, f_m \in \CC_v(T)$ be nonzero rational functions
with poles supported on $\{\zeta_1,\ldots,\zeta_d\}$, and let 
$N_1, \ldots, N_m$ be positive integers.  Then
\begin{equation*}
g(x) = \max(\frac{1}{N_1} \log_v([f_1]_x), \ldots,
                  \frac{1}{N_m} \log_v([f_m]_x))
\end{equation*}
is subharmonic in $\PP^1_{\Berk} \backslash \{\zeta_1,\ldots,\zeta_d\}$.

This follows from Example~\ref{Example F.8} above, and Proposition \ref{PropF11}(D).
\end{example}

\subsection{The Maximum Principle and the Comparison Theorem.}

\noindent{The} following maximum principle holds for subharmonic functions.

\begin{proposition} \label{PropF15} {\rm (Maximum Principle)}
Let $U \subset \PP^1_{\Berk}$ be open.  
Suppose $f(z)$ is subharmonic on $U$,
and $M$ is a bound such that for each $q \in \partial U$,
\begin{equation*} 
\limsup \begin{Sb} z \rightarrow q \\ z \in U \end{Sb} f(z) \ \le \ M \ .
\end{equation*}
Then $f(z) \le M$ for all $z \in U$.
\end{proposition}

\begin{proof}  Since the hypothesis holds for each component of $U$,
it suffices to prove the result when $U$ is a domain.

Fix $\varepsilon > 0$.  For each $q \in \partial U$, there is a
closed neighborhood $W_q$ of $q$ on which $f(z) < M+\varepsilon$.
As in the discussion of the main dendrite after
Definition~\ref{DefE5}, we can assume $W_q$ is a disc.
Using that $\partial U$ is compact,
take $W_{q_1}, \ldots, W_{q_m}$ which cover $\partial U$.
Put $V = U \backslash \cup_{i=1}^m W_{q_i}$.
Then $V$ is a simple subdomain of $U$,
and  $\partial V = \{q_1, \ldots, q_m\}$

If $\partial V = \{q_1\}$ is a single point, 
then $V$ is a disc and each harmonic function on $V$ is constant.  
Since $f$ is subharmonic, it follows that
$f(z) \le f(q_1) \le M + \varepsilon$ for all $z \in V$.

Otherwise, let $\Gamma$ be the subgraph of $U$ spanned by
$\{q_1, \ldots, q_m\}$;  then $V = r_{\Gamma}^{-1}(\Gamma_0)$
is the subdomain associated to $\Gamma$, and $\partial \Gamma = \partial V$.
By Lemma \ref{LemF8},
$f|_{\Gamma}$ belongs to $\BDV(\Gamma)$,
and $\supp(\Delta_{\Gamma}(f)^+) \subset r_{\Gamma}(\partial U) = \partial V$.
By Lemma \ref{LemF5}, $f|_{\Gamma}$ achieves its maximum at a point
of $\supp(\Delta_{\Gamma}(f)^+)$.  Since $f(q_i) \le M+\varepsilon$ for
each $i$, it follows that $f(x) \le M+\varepsilon$ for all $x \in \Gamma$.
However, $\Gamma$ is the main dendrite of $V$, and a subharmonic function
is non-increasing on paths off the main dendrite, by Corollary \ref{CorF9}.  
Hence $f(z) \le M+\varepsilon$ for all $z \in V$.

By construction, $f(z) \le M+\varepsilon$ on $U \cap W_{q_i}$ for each $i$.
Thus $f(z) \le M+\varepsilon$ for all $z \in U$.  Since $\varepsilon > 0$
is arbitrary, $f(z) \le M$ on $U$.
\end{proof}

\vskip .1 in
If $f(z)$ and $g(z)$ are subharmonic functions on an open set $U$,
one can also ask for conditions which assure $f(z) \le g(z)$ on $U$.
The desired result is called the Comparison Theorem.

\vskip .05 in
In order to formulate it, we need to generalize the Laplacian.
By Proposition \ref{PropF1}, if $f$ is subharmonic on $U$,
then it is strongly subharmonic on each simple subdomain $V \subset U$.
Hence $\Delta_{\Vbar}(f)$ is defined, and $\Delta_{\Vbar}(f)|_V \le 0$.
It makes sense to write $\Delta_{V}(f)$ for $\Delta_{\Vbar}(f)|_V$.  
However, we have not defined $\Delta_U(f)$ in general.  

Let $W$ be a component of $U$.  The simple subdomains of $W$ 
exhaust $W$, and form a directed set under containment.
If $V_1 \subset V_2$, then
\begin{equation*}
\Delta_{\Vbar_1}(f)|_{V_1}
\ = \ (r_{\Vbar_2,\Vbar_1})_*(\Delta_{\Vbar_2}(f))|_{V_1}
\ = \ \Delta_{\Vbar_2}(f)|_{V_1} \ .
\end{equation*}
Thus, the measures $\Delta_{\Vbar}(f)|_V$ cohere to give a well-defined Borel 
measure on $W$.  We will call this measure $\Delta_W(f)$.  By its construction,
$\Delta_W(f)$ is supported on $W$, and $\Delta_W(f) \le 0$.  
It is $\sigma$-finite, since $W$ can be exhausted by a countable sequence 
of simple domains. It may or may not have finite total mass;  
its mass is finite if and only if $f \in \BDV(W)$.  
We will write $\Delta_{U}(f)$ for the measure on $U$ whose restriction
to each component is $\Delta_{W}(f)$.   

\begin{proposition} \label{PropF16} {\rm (Comparison Theorem).}
Let $U \subset \PP^1_{\Berk}$ be an open set with nonempty boundary.
Suppose $f$ and $g$ are subharmonic on $U$, and 

\vskip .03 in
$A)$ for each $q \in \partial U$, \
$\displaystyle{\limsup \begin{Sb} z \rightarrow q \\ z \in U \end{Sb}
     f(z) - g(z) \ \le \ 0}$\ ;

\vskip .05 in
$B)$ $\Delta_U(f) \le \Delta_U(g)$ on $U$.

\vskip .05 in
\noindent{Then} $f(z) \le g(z)$ on $U$.
\end{proposition}

\begin{proof}  Consider the function $h(z) = f(z) - g(z)$ on $U$.

If we knew that $h(z)$ were subharmonic on $U$, the result would
follow from the maximum principle.
Indeed, $\limsup_{z \rightarrow q} h(z) \le 0$ for each $q \in \partial U$,
and for any simple subdomain $V$ of $U$, we have $f, g \in \BDV(V)$.
Thus $\Delta_{\Vbar}(h) = \Delta_{\Vbar}(f) - \Delta_{\Vbar}(g)$,
and $\Delta_{\Vbar}(h)|_V \le 0$.  Unfortunately, the difference of two
upper semicontinuous functions need not be upper semicontinuous,
so we must go back to first principles.

Fix a simple subdomain $V$ of $U$, and fix $p \in V$.
As noted above, $h \in \BDV(V)$ and $\Delta_{\Vbar}(h)|_V \le 0$.

First suppose the main dendrite of $V$ is nonempty,
and $p$ is on the main dendrite.
Thus, $V = r_{\Gamma}^{-1}(\Gamma_0)$
for a graph $\Gamma$, and $p \in \Gamma_0$.
Then $h_{\Gamma} \in \BDV(\Gamma)$,
and $\Delta_{\Gamma}(h) = (r_{\Vbar,\Gamma})_*(\Delta_{\Vbar}(h))$
is $\le 0$ on $\Gamma_0$.  By Lemma \ref{LemF5} $h|_{\Gamma}$ achieves
its maximum at a point of $\partial \Gamma = \partial V$, so
$h(p) \le \max_{q_i \in \partial V} h(q_i)$.

Next suppose that $p \in V \backslash \PP^1(\CC_v)$,
and that $p$ is not on the main dendrite
(whether or not the main dendrite is non-empty).
Write $\partial V = \{q_1, \ldots, q_m\}$
and let $\Gamma$ be the subgraph spanned by $\{p, q_1, \ldots, q_m\}$.
Let $V^{\prime}$ be the domain associated to $\Gamma$.
Then $V^{\prime}$ is a simple subdomain of $V$,
so $f|_{\Gamma}$, $g|_{\Gamma}$, and $h|_{\Gamma}$ belong to $\BDV(\Gamma)$,
hence $\Delta_{\Gamma}(h) = (r_{\Vbar,\Gamma})_*(\Delta_{\Vbar}(h)$
is $\le 0$ except on $r_{\Gamma}(\partial V) = \{q_1, \ldots, q_m\}$.
By Lemma \ref{LemF5} again, $h(p) \le \max_{q_i \in \partial V} h(q_i)$.

Finally suppose $p \in V \cap \PP^1(\CC_v)$.  If the main dendrite of $V$
is empty, let $y$ be the unique boundary point of $V$.
If the main dendrite of $V$ is nonempty, let $y$ be the point on
the main dendrite where the branch containing $p$ is attached.
In either case, consider the path $[y,p]$.  We have shown above
that $h(x) \le \max_{q_i \in \partial V} h(q_i)$ for each $x \in [y,p)$.
On the other hand, in (\ref{FBB1}) we saw that 
\begin{equation*}
f(p) \ = \ \lim \begin{Sb} x \rightarrow p \\ x \in [y,p) \end{Sb} f(x) \ ,
\qquad
g(p) \ = \ \lim \begin{Sb} x \rightarrow p \\ x \in [y,p) \end{Sb} g(x) \ .
\end{equation*}
Thus  $\displaystyle{
h(p) \ = \ \lim \begin{Sb} x \rightarrow p \\ x \in [y,p) \end{Sb} h(x)
\ \le \ \max_{q_i \in \partial V} h(q_i)}$.
\vskip .05 in
We have now shown that for each simple subdomain $V$ of\, $U$, $h(z)$
achieves its maximum on $\Vbar$ at a point of $\partial V$. Since
$\limsup_{z \rightarrow q} h(z) \le 0$ for each $q \in \partial U$,
it follows by same argument as in the proof of Proposition \ref{PropF15} 
that $h(z) \le 0$ on $V$.  Thus $f(z) \le g(z)$.
\end{proof}

\vskip .1 in
As a special case, the Comparison Theorem gives 

\begin{corollary} \label{CorF17}
Let $U \subset \PP^1_{\Berk}$ be an open set with nonempty boundary.
Suppose $f$ and $g$ are subharmonic functions on $U$ such that 

\vskip .03 in
$A)$ for each $q \in \partial U$, \
$\displaystyle{-\infty \ < \ 
      \limsup \begin{Sb} z \rightarrow q \\ z \in U \end{Sb} f(z)
      \ \le \ \liminf \begin{Sb} z \rightarrow q \\ z \in U \end{Sb} g(z)
      \ < \ \infty}$, \ and

$B)$ $\Delta_U(f) \le \Delta_U(g)$ on $U$.

\vskip .03 in
\noindent{Then} $f(z) \le g(z)$ on $U$.
\end{corollary}

\vskip .1 in
The proof of the Comparison Theorem yields the following useful criterion
for equality of subharmonic functions.

\begin{corollary} \label{CorF18}
Let $V \subset \PP^1_{\Berk}$ be a simple domain.
Suppose $f$ and $g$ are subharmonic on $V$,
with $\Delta_V(f) = \Delta_V(g)$.

$A)$  If $V$ is not a disc, 
assume $f(x) = g(x)$ on the main dendrite $\Gamma_0$ of $V$.

$B)$  If $V$ is a disc, let $x_1$ be its boundary point,
and let $\Gamma = [x_1,x_2] \subset \Vbar$ be a segment with $x_1$
as one of its endpoints.  Assume $f(x) = g(x)$ on $\Gamma_0$.

Then $f(z) \equiv g(z)$ on $V$.
\end{corollary}

\begin{proof}  By symmetry, it suffices to show that $f(z) \le g(z)$.

Put $h(z) = f(z)-g(z)$.
If $V$ is not a disc, the first step in the proof of Proposition
\ref{PropF16} was to show that when the main dendrite was nonempty,
$h(x) \le 0$ on $\Gamma_0$;  that is assumed here.
The rest of the argument showed that $h(z)$ was non-increasing on
branches off $\Gamma_0$;  that part goes through unchanged.

If $V$ is a disc, then by assumption $h(x) \le 0$ on $\Gamma_0$.
The proof of Proposition \ref{PropF16} shows that $h(z)$
is non-increasing on each path away from $x_1$.
Each such path shares an initial segment with $\Gamma_0$,
so $h(z) \le 0$ on $V$.
\end{proof}

\subsection{The Riesz Decomposition Theorem.}

Let $U$ be an open set in $\PP^1_{\Berk}$, and let $V \subset U$
be a simple subdomain.  Suppose $f$ is subharmonic in $U$.
By Proposition \ref{PropF1}, $f \in \BDV(V)$.
Put $\nu = -\Delta_{\Vbar}(f)|_V = -\Delta_V(f)$.
Fix $\zeta \notin \Vbar$, and consider the potential function
\begin{equation*}
u_{\nu}(z,\zeta) \ = \ \int -\log_v(\delta(x,y)_{\zeta}) \, d\nu(y) \ .
\end{equation*}

\begin{proposition} \label{PropF19} {\rm (Riesz Decomposition Theorem)}

Let $V$ be a simple subdomain of $U \subset \PP^1_{\Berk}$.
Fix $\zeta \in \PP^1_{\Berk} \backslash \Vbar$.

Suppose $f$ is subharmonic on $U$, and put $\nu = -\Delta_V(f)$.
Then there is a harmonic function $h(z)$ on $V$ such that
\begin{equation*}
f(z) \ = \ h(z) - u_{\nu}(z,\zeta)  \qquad \text{for all $z \in V$ \ .}
\end{equation*}
\end{proposition}

\begin{proof}
If $V$ is not a disc, let $\Gamma$ be the graph such that 
$V = r_{\Gamma}^{-1}(\Gamma_0)$, 
so $\Gamma_0$ is the main dendrite of $V$ and $\partial \Gamma = \partial V$.
Write $\partial \Gamma = \{x_1, \ldots, x_m\}$;  we can assume the points
are labelled in such a way that $r_{\Gamma}(\zeta) = x_1$.
If $V$ is a disc, let $x_1$ be its unique boundary point,
take $x_2 \in V \backslash \PP^1(\CC_v)$,
and let $\Gamma = [x_1,x_2]$.  Again $r_{\Gamma}(\zeta) = x_1$.
In either case, $V_1 := r_{\Gamma}^{-1}(\Gamma_0)$ 
is a simple subdomain of $U$, and $V_1 \subset V$.

Write $F(z) = -u_{\nu}(z,\zeta)$.
By Example~\ref{Example E.4}, $F(z) \in \BDV(\PP^1_{\Berk})$, and
\begin{equation*}
\Delta_{\PP^1_{\Berk}}(F) \ = \ \delta_{\zeta}(x) - \nu \ .
\end{equation*}
Put $\overline{\nu} = (r_{\Gamma})_*(\nu)$.
By the definition of the Laplacian, $F|_{\Gamma} \in \BDV(\Gamma)$
and $\Delta_{\Gamma}(F) = \delta_{x_1} - \overline{\nu}$.
Lemma \ref{LemF3} shows that there is a simple subdomain $W$ of $U$  with
$\Vbar \subset W$. By Proposition \ref{PropF1}, $f \in \BDV(W)$, and 
by Lemma \ref{LemF8}, $f|_{\Gamma} \in \BDV(\Gamma)$.
By the retraction property of the Laplacian,
$\Delta_{\Gamma}(f) = (r_{\Gamma})_*(\Delta_{\Wbar}(f))$.

If $V$ is not a disc, then $W = V$ and $\partial \Gamma = \partial V$.
It follows that
\begin{equation*}
\Delta_{\Gamma}(f)|_{\Gamma_0} \ = \ -\overline{\nu} \ ,
\end{equation*}
so $\sigma \ := \ \Delta_{\Gamma}(F) - \Delta_{\Gamma}(f)$
is a discrete measure supported on $\partial \Gamma$.

If $V$ is a disc, then $\sigma$ is supported on $x_1$
since $r_{\Gamma}(\Wbar \backslash V) = x_1$.
In this case $\sigma \equiv 0$,
since the only measure with total mass $0$ supported on a point
is the $0$ measure.

Put $h(x) = f(x) - F(x)$ on $\Gamma$, and extend $h(z)$ to $V$ by setting
$h(z) = h \circ r_{\Gamma}(z)$ for all $z \in V$.  
Then $h(z)$ is harmonic on $V$, since it is constant on branches off $\Gamma$ 
and $\Delta_{\Gamma}(h)|_{\Gamma \backslash \partial V} = 0$.
Now consider the functions $f(z)$ and $h(z) + F(z)$.
By construction, both are subharmonic on $V$ and
satisfy $\Delta_V(f) = \Delta_V(h+F) = -\nu$.
Both have the same restrictions to $\Gamma$.
By Corollary \ref{CorF18}, $f(z) \equiv h(z) + F(z)$.
\end{proof}

\vskip .1 in
For future applications, it is useful to know that if $f$
is continuous and subharmonic on $U$, then the potential function 
$u_{\nu}(z,\zeta)$ is continuous everywhere.

\begin{proposition} \label{PropF20}
Suppose $f$ is continuous and subharmonic on 
an open set $U \subset \PP^1_{\Berk}$.
Let $V$ be a simple subdomain of $U$,
and fix $\zeta \in \PP^1_{\Berk} \backslash \Vbar$.
Put $\nu = - \Delta_{\Vbar}(f)|_V = -\Delta_V(f)$.  
Then $u_{\nu}(z,\zeta)$ is continuous on all of $\PP^1_{\Berk}$.
\end{proposition}

\begin{proof}
By Corollary \ref{CorF18} there is a harmonic function $h(z)$ on $V$
such that $u_{\nu}(z,\zeta) = h(z)-f(z)$ on $V$.
Harmonic functions are continuous,
so $u_{\nu}(z,\zeta)$ is continuous on $V$.
Since $\supp(\nu) \subset \Vbar$, $u_{\nu}(z,\zeta)$
is also continuous on $\PP^1_{\Berk} \backslash \Vbar$
by Proposition~\ref{PropD7} (which includes continuity at $\zeta$).
It remains to show that $u_{\nu}(z,\zeta)$
is continuous on the finite set $\partial V$.

Since $V$ is a simple domain, each harmonic function on $V$ 
has a continuous extension to $\Vbar$.
By construction $h(z)$ is the unique harmonic function on $V$ such that 
$h(x_i) = f(x_i) + u_{\nu}(x_i,\zeta)$ for each $x_i \in \partial V$.

Fix $x_i \in \partial V$.  
Since $f$ and $h$ are continuous on $\Vbar$, 
\begin{equation*}
\lim \begin{Sb} z \rightarrow x_i \\ z \in V \end{Sb}  u_{\nu}(z,\zeta)
\ = \ u_{\nu}(x_i,\zeta) \ .
\end{equation*}
Now consider the behavior of
$u_{\nu}(z,\zeta)$ on each connected component of
$\PP^1_{\Berk} \backslash \{x_i\}$, as $z \rightarrow x_i$.  
Each such component is
a disc with $x_i$ as its boundary point.  One component contains $V$;
we have already dealt with it.  Suppose $W$ is a component which does not
contain $V$.  If $W$ contains $\zeta$, choose a point $p$ on the interior of
the path from $x_1$ to $\zeta$;
otherwise, let $p \in U \backslash \PP^1(\CC_v)$ be arbitrary.
Put $\Gamma = [x_1,p]$.
Since $u_{\nu}(z,\zeta) \in \BDV(\PP^1_{\Berk})$, its restriction
to $\Gamma$ belongs to $\BDV(\Gamma)$, hence is continuous on $\Gamma$.
Thus
\begin{equation} \label{FYM1}
\lim \begin{Sb} z \rightarrow x_i \\ z \in \Gamma \end{Sb} u_{\nu}(z,\zeta)
        \ = \ u_{\nu}(x_i,\zeta) \ .
\end{equation}
If $W$ is a component which does not contain $\zeta$,
then $u_{\mu}(z,\zeta)$ is harmonic on $W$, hence constant.
By (\ref{FYM1}), $u_{\nu}(z,\zeta) \equiv u_{\nu}(x_i,\zeta)$
on $W$ for such a component.

If $W$ is the component containing $\zeta$,
put $W_0 = r_{\Gamma}^{-1}(\Gamma_0)$;  then $W_0$ is a simple subdomain
of $\PP^1_{\Berk}$ contained in $W \backslash \{\zeta\}$, and
$u_{\nu}(z,\zeta)$ is harmonic on $W_0$.
The main dendrite of $W_0$ is $\Gamma_0$.
Since harmonic functions are constant on branches off the main dendrite, 
\begin{equation} \label{FYM2}
\lim \begin{Sb} z \rightarrow x_i \\ z \in U \end{Sb} u_{\nu}(z,\zeta)
\ = \ \lim \begin{Sb} z \rightarrow x_i \\ z \in \Gamma \end{Sb}
              u_{\nu}(z,\zeta) \ = \ u_{\nu}(x_i,\zeta) \ .
\end{equation}
Combining all cases, we see that $u_{\nu}(z,\zeta)$ is continuous at $x_i$.
\end{proof}

\vskip .1 in
The Riesz Decomposition Theorem also gives interesting information
about the structure of functions in $\BDV(U)$.  Definition~\ref{DefE2} asserts 
that $f \in \BDV(U)$ iff $f|_{\Gamma} \in \BDV(\Gamma)$ for
every subgraph $\Gamma \subset U$, and there is a uniform bound $B$
for the measures $|\Delta_{\Gamma}(f)|(\Gamma)$ for all $\Gamma \subset U$.
Nothing is said about continuity on $U$, 
or about the behavior of $f$ on $\PP^1(\CC_v)$.

\begin{proposition} \label{PropF21}
Let $U \subset \PP^1(\CC_v)$ be a domain, and let $f \in \BDV(U)$.
For any simple subdomain $V \subset U$,
there are subharmonic functions $g$, $h$ on $V$ such that
$f(z) = g(z)-h(z)$ for all $z \in V \backslash \PP^1(\CC_v)$.
\end{proposition}

\begin{proof}
Since $f \in \BDV(U)$, also $f \in BDV(V)$.  Put $\nu = \Delta_{\Vbar}(f)|_V$,
and let $\nu^{+}$, $\nu^{-}$ be the positive and negative measures
in the Jordan decomposition of $\nu$, so $\nu = \nu^{+} - \nu^{-}$.
Fix $\zeta \in \PP^1_{\Berk} \backslash \Vbar$, and put
\begin{equation*}
g_1(z) = -u_{\nu^{-}}(z,\zeta)\ , \qquad  h_1(z) = -u_{\nu^{+}}(z,\zeta) \ .
\end{equation*}
Then $g_1$ and $h_1$ are subharmonic in
$\PP^1_{\Berk} \backslash \{\zeta\}$ and belong to $\BDV(\PP^1_{\Berk})$,
with $\Delta_{\PP^1_{\Berk}}(g_1) = \delta_{\zeta}(x) - \nu^{-}$ and 
$\Delta_{\PP^1_{\Berk}}(h_1) = \delta_{\zeta}(x) - \nu^{+}$. Thus 
\begin{equation*}
\Delta_V(g_1-h_1) \ = \ \nu^{+}-\nu^{-} \ = \ \nu \ .
\end{equation*}

Put $F_1(z) = f(z) - (g_1(z)-h_1(z))$.  Then $F_1 \in \BDV(V)$ 
and $\Delta_{\Vbar}(F_1)|_V \equiv 0$.  Since $V$ is a simple domain,
its boundary $\partial V = \{p_1, \ldots, p_m\}$ is finite. 
Thus $\Delta_{\Vbar}(f)$ is a discrete measure supported on $\partial V$.

If $V$ is a disc, then $\partial V$ is a single point $p_1$.
Put $F(z) = F_1(z)$.  Since $\Delta_{\Vbar}(F)(\Vbar) = 0$, it follows that
$\Delta_{\Vbar}(F)(\{p_1\}) = 0$, and $\Delta_{\Vbar}(F) \equiv 0$.

If $V$ is not a disc, put $\nu_0 = \Delta_{\Vbar}(F_1)|_{\partial V}
= \sum_{i=1}^m c_i \delta_{p_i}(x)$.  Here $\sum c_i = 0$.
Fix $p_0 \in \PP^1_{\Berk} \backslash \PP^(\CC_v)$, and put
\begin{equation*}
H(z) \ = \ \sum_{i=1}^m c_i j_{p_0}(z,p_i) \ .
\end{equation*}
Then $H(z)$ is continuous on $\PP^1_{\Berk}$,
$H(z) \in \BDV(\PP^1_{\Berk})$,
and $\Delta_{\PP^1_{\Berk}}(H) = \sum_{i=1}^m c_i \delta_{p_i}(x)$.
In particular, $H(z)$ is harmonic on $V$.  Put $F(z) =  F_1(z) - H(z)$.
Then $F \in \BDV(V)$ and $\Delta_{\Vbar}(F) \equiv 0$.

In either case, for each subgraph $\Gamma \subset V$,
$\Delta_{\Gamma}(F) = (r_{\Vbar,\Gamma})_*(\Delta_{\Vbar}(F)) \equiv  0$.
Hence $F$ is constant on $\Gamma$.
Any two points $x_1, x_2 \in V \backslash \PP^1(\CC_v)$,
are connected by the path $[x_1,x_2]$.
Thus $F$ is constant on $V \backslash \PP^1(\CC_v)$;  
let $C$ be that constant.

We have now shown that on $V \backslash \PP^1(\CC_v)$
\begin{equation*}
f(z) \ = \ (g_1(z) + C) - (h_1(z) + H(z))
\end{equation*}
where $H(z) = 0$ if $V$ is a disc.  Thus the theorem holds, 
with $g(z) = g_1(z) + C$ and $h(z) = h_1(z) + H(z)$.
\end{proof}

\vskip .1 in
Given  $f, g \in \BDV(U)$, we will say that $f \cong g$
if $f(z) = g(z)$ on $U \backslash \PP^1(\CC_v)$.
This is an equivalence relation.

%


\subsection{Convergence of Laplacians.}

Let $U$ be an open subset of $\PP^1_{\Berk}$,
and let $f_1, f_2, \ldots$ be a sequence of subharmonic functions on $U$
which converge to a subharmonic function $F$.
One can ask under what conditions the measures  
$\Delta_{U}(f_1), \Delta_{U}(f_2), \ldots$ converge to $\Delta_U(F)$.

\begin{proposition} \label{PropF23}  
Let $U \subset \PP^1_{\Berk}$ be a domain, and let
$\{f_n(z)\}_{n \in \NN}$ be a sequence of subharmonic functions on $U$.
Suppose one of the following conditions holds:

\vskip .05 in
$A)$  The $f_n(z)$ converge uniformly to $F(z)$ on compact subsets of $U$$;$

\vskip .05 in
$B)$  The $f_n(z)$ decrease monotonically to $F(z)$, 
and $F(z) \not\equiv -\infty$ on $U$$;$

\vskip .05 in
$C)$  The $f_n(z)$ increase monotonically, 
and the family $\{f_n\}$ is locally bounded above.
Put $f(z) = \lim_{n \rightarrow \infty} f_n(z)$ and let $F(z) = f^*(z)$
be the upper semicontinuous regularization of $f(z)$.

\vskip .05 in
Then $F(z)$ is subharmonic on $U$, 
and the measures $\Delta_V(f_n)$ converge weakly to $\Delta_V(F)$
on each simple subdomain $V \subset U$.
\end{proposition}

\begin{proof}  In all three cases, Proposition \ref{PropF11} shows 
that $F(z)$ is subharmonic.  Only the convergence of the
measures needs to be established.

\vskip .05 in
\noindent{Proof of (A)}.  Let $V \subset U$ be a simple subdomain.

We first claim that the measures $|\Delta_{\Vbar}|(f_n)$ have uniformly
bounded total mass.  This depends on the fact that $V$ can be enlarged 
within $U$.  Write $\partial V = \{x_1, \ldots, x_m\}$.
For each $x_i \in \partial V$ which is not of type IV,
take a simple subdomain $W_i$ of $U$ which contains $x_i$
but does not contain any $x_j \ne x_i$.
Let $\tV$ be the union of $V$ and these $W_i$.
Write $\partial \tV = \{\tx_1, \ldots, \tx_M\}$.
   
If $V$ is not a disc, let $\Gamma$ be the subgraph spanned
by $\{x_1, \ldots, x_m\}$ and let $\tGamma$ be the subgraph spanned by
$\{\tx_1, \ldots, \tx_M\}$, so $V = r_{\Gamma}^{-1}(\Gamma_0)$
and $\tV = r_{\tGamma}^{-1}(\tGamma_0)$.  Then $\Gamma \subset \tGamma$,
and each $x_i \in \partial V$ which is not of type IV belongs to $\tGamma_0$.

If $V$ is a disc, fix a point $p \in V \backslash \PP^1(\CC_v)$,
and let $\Gamma$ be the segment $[x_1,p]$.  
Let $\tGamma$ be the subgraph spanned by $\{p,\tx_1, \ldots, \tx_M\}$.
Then $\tGamma$ corresponds to a simple subdomain of $U$,
and has $x_1$ in its interior.

Take $x_i \in \partial V$.  When $x_i$ is regarded as a point of $\Gamma$,
it has a unique direction vector $\vec{v}_i$ leading into $\Gamma$.  
As in (\ref{FWC1}) in the proof of Proposition \ref{PropF1}, 
the total mass $|\Delta_{\Gamma}(f)|(\Gamma)$ is  
\begin{equation*}
|\Delta_{\Gamma}(f)|(\Gamma)|  = 
         2 \cdot \sum_{i=1}^m |d_{\vec{v}_i}f_n(x_i)| \ .
\end{equation*}         
This is the same as the total mass $|\Delta_{\Vbar}(f_n)|(\Vbar)$ since
$\Delta_{V}(f) \le 0$.  
Thus, to bound the $|\Delta_{\Vbar}(f_n)|(\Vbar)$ it suffices to show
that the $|d_{\vec{v}_i}f_n(x_i)|$ are uniformly bounded, for all $i$ 
and $n$. 

On the other hand, when $x_i$ is regarded as a point of $\tGamma$, 
it may have several direction vectors $\vec{v}_{i,j}$, $j = 1, \ldots, n_i$.
Without loss, we can assume the $\vec{v}_{ij}$
are indexed so that $\vec{v}_i = \vec{v}_{i1}$ for each $i$. 
Fix $T > 0$ small enough that for each $i$, $j$, the point 
$q_{i,j} = x_i + t \vec{v}_{ij}$ lies on the edge emanating from
$x_i$ in the direction $\vec{v}_{ij}$.  

Consider the limit function $F$, and fix $\varepsilon > 0$.  Since the
$f_n$ converge uniformly to $F$ on $\tGamma$, there is an $N$ such that
$|f_n(q_{ij}) - F(q_{ij})| < \varepsilon$ and $|f_n(x_i) - F(x_i)|
< \varepsilon$ for all $n \ge N$ and all $i, j$.  It follows that
for such $i$, $j$ and $n$  
\begin{equation}  \label{FJJ1}
\frac{f_n(q_{ij}) - f_n(x_i)}{T}  
\ \le \ \frac{(F(q_{ij}) + \varepsilon) - (F(x_i)-\varepsilon)}{T} \ .
\end{equation}
Since each $f_n$ is convex up on edges of $\tGamma$, for each $i$, $j$
\begin{equation} \label{FJJ2}
d_{\vec{v}_{ij}}f_n(x_i) \ \le \ \frac{f_n(q_{ij}) - f_n(x_i)}{T} \ .
\end{equation}
Thus (\ref{FJJ1}) provides an upper bound $B_{ij}$  
for each $d_{\vec{v}_{ij}}f_n(x_i)$.

If $x_i \in \partial V$ is of type IV, then $x_i \in \partial \tV$
and there is only one direction vector at $x_i$ in both $\Gamma$ and
$\tGamma$.  As was shown in (\ref{FCV1}) in the proof of Lemma \ref{LemF7},
$d_{\vec{v}_1}f_n(x_i) \ge 0$, so
\begin{equation*}
0  \ \le \  d_{\vec{v}_1}f_n(x_i) \ \le \ B_{i1} \ .
\end{equation*}

If $x_\i \in \partial V$ is not of type IV, then $x_i$ is
an interior point of $\tGamma$, and $n_i \ge 2$. 
Lemma \ref{LemF8} gives $0 \ge \Delta_{\tGamma}(f_n)(x_i) 
= - \sum_{j=1}^{n_i} d_{\vec{v}_{ij}}f_n(x_i)$, so 
\begin{equation} \label{FJJ3}
d_{\vec{v}_{i1}}f_n(x_i)
\ \ge \ - \sum_{j=2}^{n_i} d_{\vec{v}_{ij}}f_n(x_i)
\ \ge \ - \sum_{j=2}^{n_i} B_{ij} \ .
\end{equation} 
Since $d_{\vec{v}_i}f_n(x_i) = d_{\vec{v}_{i1}}f_n(x_i)$,
we have shown that $|d_{\vec{v}_i}f_n(x_i)|$ is uniformly bounded  
for all $i$ and all $n \ge N$.  

\vskip .05 in
Let $B$ be the bound constructed above for
the masses $|\Delta_{\Vbar}(f_n)|(\Vbar)$ with $n \ge N$.
After increasing $B$ if necessary,
we can assume that $|\Delta_{\Vbar}(f_n)|(\Vbar) \le B$ for all $n$.

We will now show the sequence of measures $\Delta_{\Vbar}(f_n)$ 
converges weakly to a measure $\mu$ on $\Vbar$.  For this, drop
the meaning of $\Gamma$ used above, and recall that by 
Proposition~\ref{PropE4}(B), 
as $\Gamma$ ranges over all subgraphs of $V$, and $G$ ranges over
all functions in $\CPA(\Gamma)$, then the functions 
$g(x) = G \circ r_{\Gamma}(x)$ are dense in $\cC(\Vbar)$.  
For each subgraph $\Gamma^{\prime} \supset \Gamma$, 
each $g(x)$ and each $f_n$, Lemma~\ref{LemE7}(F) gives 
\begin{equation*}
\int_{\Gamma^{\prime}} g \, \Delta_{\Gamma^{\prime}}(f_n)
\ = \ \int_{\Gamma^{\prime}} f_n \, \Delta_{\Gamma^{\prime}}(g)
\ = \ \int_{\Gamma} f_n \, \Delta_{\Gamma}(G) \ .
\end{equation*}
where the last equality holds because
$\Delta_{\Gamma^{\prime}}(g)|_{\Gamma} = \Delta_{\Gamma}(G)$.
Taking a limit, we find that 
\begin{equation} \label{FMM1B} 
\int_{\Vbar} g \, \Delta_{\Vbar}(f_n) \ = \  
\int_{\Vbar} f_n \, \Delta_{\Vbar}(g) \ = \
\int_{\Gamma} f_n \, \Delta_{\Gamma}(G) .
\end{equation}
Since the $f_n$ converge uniformly to $F$ on $\Gamma$,
we can define a linear functional 
\begin{equation} \label{FMM2B}
\Lambda_F(g)
\ = \ \lim_{n \rightarrow \infty} \int_{\Vbar} f_n \, \Delta_{\Vbar}(g) 
\ = \ \int_{\Vbar} F \, \Delta_{\Vbar}(g)
\ = \ \int_{\Gamma} F \, \Delta_{\Gamma}(G)
\end{equation}
on the dense space of functions $g(x) = G \circ r_{\Gamma}(x)$
in $\cC(\Vbar)$.  On the other hand, by (\ref{FMM1B})
\begin{equation} \label{FMM3B}
\Lambda_F(g)
\ = \ \lim_{n \rightarrow \infty} \int_{\Vbar} g \, \Delta_{\Vbar}(f_n) \ .
\end{equation}
Since the measures $|\Delta_{\Vbar}(f_n)|(V_n)$ are uniformly bounded, 
$\Lambda_F$ extends to a bounded linear functional on $\cC(\Vbar)$.
By the Riesz Decomposition theorem, there is a unique measure $\mu$
such that 
\begin{equation*}
\Lambda_F(g) \ = \ \int_{\Vbar} g(x) \, d\mu(x) \ .
\end{equation*}
for all $g \in \cC(\Vbar)$.
This measure is the weak limit of the $\Delta_{\Vbar}(f_n)$.  

To complete the proof, we must show that $\mu = \Delta_{\Vbar}(F)$.
Since $\Delta_{\Vbar}(F)$ is a bounded measure on $\Vbar$, it suffices
to check that 
\begin{equation*}
\Lambda_F(g)
\ = \ \int_{\Vbar} g \, \Delta_{\Vbar}(F)(x)
\end{equation*}
on the dense space of functions considered above.  This follows immediately 
from (\ref{FMM1B}), with $f_n$ replaced by $F$.
This proves (A).

\vskip .1 in
Examining the proof, we see that the argument applies under a weaker
hypothesis, namely that for each subgraph $\Gamma \subset U$,
the $f_n$ converge uniformly to $F$ on $\Gamma$.  We will use this
observation to prove parts (B) and (C).

\vskip .05 in
For (B), fix a subgraph $\Gamma \subset U$.
As in the proof of Proposition \ref{PropF1},
after enlarging $\Gamma$ if necessary
we can assume that  $V = r_{\Gamma}^{-1}(\Gamma_0)$
is a simple subdomain of $U$.  Since $F$ is subharmonic in $U$,
Lemma \ref{LemF7} shows that $F|_{\Gamma}$ is continuous.

Similarly $f_n|_{\Gamma}$ is continuous for each $n$.
By Dini's theorem, a sequence of continuous functions on
a compact set which converge monotonically to a continuous function,
converges uniformly.  Hence the $f_n$ converge uniformly to $F$ on
each subgraph $\Gamma \subset U$, and the result follows.

\vskip .05 in
For (C), a similar argument shows that since $F(z) = f^*(z)$
is subharmonic on $U$, then for each subgraph
$\Gamma \subset U$ the restriction of $F(z)$ to $\Gamma$
is continuous.  However, by Proposition \ref{PropF11}(E), $F(z)$
coincides with $f(z) = \lim_{n \rightarrow \infty} f_n(z)$ on
$U \backslash \PP^1(\CC_v)$.  Hence $f|_{\Gamma} = F|_{\Gamma}$
is continuous, and the $f_n$ converge uniformly to $F$ on $\Gamma$.
\end{proof}

\vskip .1 in
We augment this with two simple results which deal with other
cases in Proposition \ref{PropF11}.  The first is immmediate
by the linearity of the Laplacian.

\begin{proposition}  \label{PropF24}
Let $U \subset \PP^1_{\Berk}$ be open.  Suppose $f, g \in \SH(U)$,
and $0 \le \alpha, \beta \in \RR$.  Then
$\Delta_U(\alpha f + \beta g)
 = \alpha \cdot \Delta_U(f) + \beta \cdot \Delta_U(g)$\ .
\end{proposition}

\begin{proposition} \label{PropF25}
Let $U \subset \PP^1_{\Berk}$ be open.
Suppose $f(z)$ and $g(z)$ are continuous subharmonic functions on $U$,
with $f(z) > g(z)$ outside a compact subset of $U$.
Put $h(z) = \max(f(z),g(z))$ and let $Z = \{z \in U : f(z) \le g(z)\}$.
Then

$A)$ \ \ $\Delta_U(h)|_{U \backslash Z} \equiv \Delta_U(f)|_{U \backslash Z}$,
\ and 

$B)$ \ \ $\Delta_U(h)(Z) = \Delta_U(f)(Z)$\ .
\end{proposition}

\begin{proof}  We can assume without loss that $U$ is a domain,
since the result holds if and only if it holds for each component of $U$.

\vskip .05 in
The proof of (A) is simple.  For each $x \in U \backslash \Vbar$,
there is a simple subdomain $W \subset U$ such that $x \in W$ and
$\Wbar \cap \Vbar = \varphi$.  We have $h(z) = f(z)$ on $W$,
so $\Delta_{\Wbar}(h) = \Delta_{\Wbar}(f)$, and
\begin{equation*}
\Delta_U(f)|_W \ = \ \Delta_{\Wbar}(f)|_W
  \ = \ \Delta_{\Wbar}(h)|_W \ = \ \Delta_U(h)|_W \ .
\end{equation*}
Since $U \backslash U$ can be covered by such $W$,
the result follows.

\vskip .05 in
For part (B), note that $Z$ is compact;  fix a simple subdomain
$V \subset U$ with $Z \subset V$.  By Proposition \ref{PropF1},
the restrictions of $f$, $g$ and $h$ to $V$ belong to $\BDV(V)$.
Note that $\partial V$ is compact. Fix $\zeta \in Z$, and cover
$\partial V$ with a finite number of balls $\cB(x_i,r_i)_{\zeta}^-$,
where the closure of each $\cB(x_i,r_i)_{\zeta}^-$ is disjoint from $Z$.
Without loss, we can assume the $\cB(x_i,r_i)_{\zeta}^-$ are pairwise disjoint.  
Let $p_i$ be the unique boundary point of $\cB(x_i,r_i)_{\zeta}^{-}$;
it belongs to $V \backslash Z$.  Put $V_i = U \cap \cB(x_i,r_i)$; then
$h(z) = f(z)$ on $V_i$, and $\partial V_i \subset \partial V \cup \{p_i\}$.

By the retraction formula for Laplacians,
\begin{equation*}
(r_{\Vbar,\Vbar_i})_*(\Delta_{\Vbar}(f)) \ = \ \Delta_{\Vbar_i}(f) \ = \
\Delta_{\Vbar_i}(h) \ = \ (r_{\Vbar,\Vbar_i})_*(\Delta_{\Vbar}(h)) \ .
\end{equation*}
Since $r_{\Vbar,\Vbar_i}(x)$
is the identity map on $\Vbar_i \backslash \{p_i\}$
and $r_{\Vbar,\Vbar_i}(\Vbar \backslash \Vbar) = \{p_i\}$,
\begin{equation*}
\Delta_{\Vbar}(f)|_{\Vbar_i \backslash \{p_i\}} \ = \
\Delta_{\Vbar}(h)|_{\Vbar_i \backslash \{p_i\}} \ .
\end{equation*}
Combined with the result from (A), this gives
$\Delta_{\Vbar}(f)|_{\Vbar \backslash Z} = 
\Delta_{\Vbar}(h)|_{\Vbar \backslash Z}$.
Since both $\Delta_{\Vbar}(f)$ and $\Delta_{\Vbar}(h)$ have total mass $0$,
\begin{equation*}
\Delta_{\Vbar}(f)(Z) \ = \ - \Delta_{\Vbar}(f)(\Vbar \backslash Z) 
\ = \ - \Delta_{\Vbar}(h)(\Vbar \backslash Z) \ = \ \Delta_{\Vbar}(f)(Z) \ .
\end{equation*}
\end{proof}

\subsection{Smoothing.}  

In the classical theory, each subharmonic function $f$
is locally a decreasing limit of $\cC^{\infty}$ subharmonic functions.
These are gotten by convolving $f$ with a smoothing kernel.

Here there seems to be no smoothing kernel, but nonetheless
each subharmonic function is locally a decreasing limit of continuous
subharmonic functions of a special form.  The ``nicest'' functions
we have encountered on $\PP^1_{\Berk}$ are those of the form
$f(z) = F \circ r_{\Gamma}(z)$, where $\Gamma$ is a subgraph
and $F \in \CPA(\Gamma)$.  We deem them to be the smooth functions on 
$\PP^1_{\Berk}$.

\begin{proposition}  \label{PropF26}
Let $U \subset \PP^1_{\Berk}$ be a domain,
and let $f$ be subharmonic in $U$.  Then for every simple subdomain
$V$ of $U$, there is a decreasing sequence of subharmonic functions
$f_1 \ge f_2 \ge \cdots$ on $V$ such that

\vskip .05 in
$A)$  $\lim_{n \rightarrow \infty} f_n(z) = f(z)$ pointwise,
 for each $z \in V$ $;$

\vskip .05 in
$B)$ for each $n$, there are a subgraph $\Gamma_n \subset V$ and a
function $F_n \in \CPA(\Gamma_n)$ such that
$f_n(z) = F_n \circ r_{\Gamma_n}(z)$.

\vskip .05 in
In particular, each $\Delta_V(f_n)$ 
is a discrete measure, and the $\Delta_V(f_n)$ converge
weakly to $\Delta_V(f)$ on $V$.
\end{proposition}

\begin{proof}
The construction below aims to err on the side of explicitness.
We first construct a sequence of graphs $\Gamma_n$ which exhaust
$\Vbar \backslash \PP^1(\CC_v)$, and then construct functions $F_n$
on them which approximate $f$.
If the assertion holds for a subdomain $\tV$ containing $V$, then
it holds for $V$.  Hence, after enlarging $V$ if necessary,
we can assume that each point in $\partial V$ is of type II.

\vskip .05 in
If $V$ is a disc, let $x_1$ be its unique boundary point, take 
$x_2 \in V \backslash \PP^1(\CC_v)$, and put $\tGamma = [x_1,x_2]$.
If $V$ is not a disc, let $\tGamma$ be the subgraph spanned by
$\partial V = \{x_1, \ldots, x_m\}$.  
There are countably many type II points in $V$;
list them as $\{p_n\}_{n \ge 1}$.  
Inductively define $\Gamma_1 = \tGamma \cup [x_1,p_1]$,
and $\Gamma_n = \Gamma_{n-1} \cup [x_1,p_n]$ for $n \ge 2$.
Then $\partial V \subset \Gamma_1 \subset \Gamma_2 \subset \cdots$,
and the graphs $\Gamma_n$ exhaust $\Vbar \backslash \PP^1(\CC_v)$.

The endpoints and branch points of each $\Gamma_n$ are type II points,
so the length of each edge of $\Gamma_n$ is a rational number.
Let $N_n$ be the least common denominator for the lengths of the
edges of $\Gamma_n$, and let $L_n$ be the total path length of $\Gamma_n$.
Thus $L_n = K_n/N_n$ for some integer $K_n$.
Since each edge of $\Gamma_{n-1}$ is a union of edges of $\Gamma_n$,
it follows that $N_{n-1}|N_n$ for each $n$.  It is easy to see that
$N_n \rightarrow \infty$ as $n \rightarrow \infty$.

For each $n$, choose points $p_{n,i}$, $i = 1, \ldots, M_n$
which partition $\Gamma_n$ into segments of length $1/N_n$,
i.e. which subdivide each edge into subsegments of length $1/N_n$.
The endpoints and branch points of $\Gamma_n$ are contained
in $\{p_{n,i}\}$ and because $N_{n-1}|N_n$, also
$\{p_{n-1,i}\} \subset \{p_{n,i}\}$.

Let $F_n(x) \in \CPA(\Gamma_n)$ be the function whose values at the
points $p_{n,i}$ are $F_n(p_{n,i}) = f(p_{n,i})$, and which
interpolates linearly  on the intervening segments.
By Lemmas \ref{LemF8} and \ref{LemF2} $f(x)$ is convex up on each edge
of $\Gamma_n$.  Hence $F_n(x) \ge f(x)$ for all $x \in \Gamma_n$.
For the same reason, $F_n(x) \ge F_{n+1}(x)$ for all $x \in \Gamma_{n}$.
Put $f_n(z) = F_n \circ r_{\Gamma_n}(z)$, so
$f_n(z)$ is constant on branches off $\Gamma_n$.
By Corollary \ref{CorF9}, $f(z)$ is non-increasing on branches off $\Gamma_n$,
so $f_n(z) \ge f(z)$ for all $z \in V$.  Since $F_n(x) \ge F_{n+1}(x)$
on $\Gamma_n$, $f_n(z) \ge f_{n+1}(z)$ for all $z \in V$ as well.

As $m \rightarrow \infty$ the functions $f_m(x)$ converge
uniformly to $f(x)$ on each fixed $\Gamma_n$, 
since $f|_{\Gamma_n}$ is continuous and $N_m \rightarrow \infty$.
Since each $z \in V \backslash \PP^1(\CC_v)$ belongs to some $\Gamma_n$,
we see that $f_n(z) \rightarrow f(z)$ pointwise on $V \backslash \PP^1(\CC_v)$.
Now fix $q \in V \cap \PP^1(\CC_v)$, and consider the path $[x_1,q]$.
There is a $y \in [x_1,q)$ such that $f(x)$
is non-increasing for $x$ in $[y,q)$.
Furthermore each $x \in [y,q)$ belongs to $\Gamma_n$
if $n$ is sufficiently large.  Hence for each $\varepsilon > 0$,
each $x \in [y,q)$, and each sufficiently large $n$
\begin{equation*}
f(x)+\varepsilon \ \ge \ f_n(x) \ \ge \ f_n(q) \ \ge \ f(q) \ .
\end{equation*}
By Corollary \ref{CorF9} 
\begin{equation*}
\lim \begin{Sb} x \rightarrow q \\ x \in [y,q) \end{Sb} f(x)
      \ = \ f(q) \ .
\end{equation*}
Hence $\lim_{n \rightarrow \infty} f_n(q) = f(q)$.  

\vskip .05 in
We will now show that each $f_n(z)$ is subharmonic on $V$.
Clearly
\begin{equation*}
\Delta_{\Gamma_n}(F_n) = \sum_{i=1} c_{n,i} \delta_{p_{n,i}}(x)
\end{equation*}
for some $c_{n,i} \in \RR$.  We claim that $c_{n,i} \le 0$
if $p_{n,i} \notin \partial V$.  By Lemma \ref{LemF8},
for each $p_{n,i} \notin \partial V$
\begin{equation*}
- \sum_{\text{$\vec{v}$ at $p_{n,i}$}} d_{\vec{v}}f(p_{n,i})
     \ = \ \Delta_{\Gamma_n}f(p_{n_i}) \ \le \ 0 \ .
\end{equation*}
As $d_{\vec{v}}f_n(p_{n,i}) \ge d_{\vec{v}}f(p_{n,i})$ for each $\vec{v}$,
this gives $c_{n,i}  = \Delta_{\Gamma_n}F_n(p_{n,i}) \le 0$.
Since $f_n(z)$ is constant on branches off $\Gamma_n$,
\begin{equation*}
\Delta_{\Vbar}(f_n) \ = \ \Delta_{\Gamma_n}(F_n)
\ = \ \sum_{i=1}^{M_n} c_{n,i} \delta_{p_{n,i}}(x) \ .
\end{equation*}
In particular $\Delta_V(f_n) \le 0$.  
By construction $f_n$ is continuous.
By Proposition \ref{PropF1}, $f_n$ is subharmonic on $V$.

\vskip .05 in
The final assertion, that the measures $\Delta_V(f_n)$ converge weakly to
$\Delta_V(f)$, follows from the proof of Proposition \ref{PropF23}(B).
We have shown that $f_n \rightarrow f$ monotonically
on each subgraph $\Gamma \subset V$,
which was the hypothesis needed for 
the proof of Proposition \ref{PropF23}(B).
\end{proof}

\section{Applications to Dynamics.}

\vskip .1 in

In this section we begin the study of the dynamics of a rational map 
$\varphi$ on $\PP^1_{\Berk}$.  We first develop a theory of multiplicities 
for $\varphi$ at points in $\PP^1_{\Berk}$,  
with properties analogous to the classical algebraic multiplicies.  
We define the pullback and pushforward measures under $\varphi$, 
and establish their functoriality properties. 

We then show there is a probability measure $\mu$
on $\PP^1_{\Berk}$ which satisfies $\varphi_*(\mu) = \mu$ and 
$\varphi^*(\mu) = d \cdot \mu$.  
By analogy with the classical case for $\PP^1(\CC)$ (see \cite{Ly},\cite{FLM}),
we call it the Lyubich measure.  We also show that the potential kernel
associated to the Lyubich measure satisfies a certain energy-minimization 
principle.  These results will be used in a forthcoming paper (\cite{B-R2})
concerning equidistribution of points of small dynamical height.

\vskip .1 in
We have been informed that independently Chambert-Loir/Thuillier and 
Favre/Rivera-Letelier
have constructed measures on $\PP^1_{\Berk}$ with the properties above.  
Presumably their measures coincide with ours.
One of the main insights of our construction is that the Lyubich measure is 
``the minus Laplacian of the Call-Silverman local height''.  
It will be interesting to learn what their constructions say about it.  
In this context, we note that Szpiro, Tucker and Piniero (\cite{STP})
have constructed a sequence of blowups 
of a rational map on $\PP_1/\Spec(\cO_v)$, 
leading to a sequence of discrete measures supported on the special fibres, 
which we expect should converge weakly to the Lyubich measure.  

\subsection{The action of a rational function on $\PP^1_{\Berk}$.}

Let $\varphi(T) \in \CC_v(T)$ be a rational function of degree $d \ge 1$.
The action of $\varphi(T)$ on $\PP^1_{\Berk}$ is defined by
\begin{equation} \label{FHL1}
[g]_{\varphi(x)} \ = \ [g \circ \varphi]_x
\end{equation}
for all $g \in \CC_v(T)$ for which $g \circ \varphi$ belongs to the
stalk of the Berkovich structure sheaf at $x$.  If $x$ is not of type I,
then $\CC_v(T) \subset \cO_{X,x}$, as was shown in
Section~\ref{Section B}.
If $x \in \PP^1(\CC_v)$, then (\ref{FHL1}) defines the usual action
of $\varphi$, since $y = \varphi(x)$ if and only if
$|g(y)|_v = |g(\varphi(x))|_v$ for all linear polynomials
$g = T-a \in \CC_v[T]$ and for $g = 1/T$.

Note that if $x \in \PP^1_{\Berk} \backslash \PP^1(\CC_v)$, 
then $\varphi(x) \in \PP^1_{\Berk} \backslash \PP^1(\CC_v)$, 
because the seminorm $[g]_{\varphi(x)}$ is defined for all $g \in \PP^1(\CC_v)$
and has kernel $0$, whereas for each $a \in \PP^1(\CC_v)$ the seminorm
$[g]_a = |g(a)|_v$ is only defined on a subring of $\CC_v(T)$ and 
has a nonzero kernel.  By considering the induced maps on the stalks 
and residue fields, and the description of the stalks and residue fields
given in Section~\ref{Section B}, one sees that $\varphi$ 
takes type I points to type I points, type II points to type II points, 
type III points to type III points, and type IV points to type IV points.

\subsection{Analytic Multiplicities.}

In this subsection we will use the theory of subharmonic functions
to define multiplicities for $\varphi$ on points of $\PP^1_{\Berk}$.

\vskip .1 in
For classical points $a, b \in \PP^1(\CC_v)$ with $\varphi(a) = b$, 
write $m_{\varphi,b}(a)$ for the multiplicity of $\varphi$ at $a$.  
For $a, b$ with $\varphi(a) \ne b$, put $m_{\varphi,b}(a) = 0$.
Thus for each $b \in \CC_v$ 
\begin{equation*}
\div(\varphi(T)-b) \ = \ \sum_{\varphi(a) = b} m_{\varphi,b}(a) \delta_a(x) 
         - \sum_{\varphi(a) = \infty} m_{\varphi,\infty}(a) \delta_a(x)\ ,
\end{equation*}
and for each $b \in \PP^1(\CC_v)$
\begin{equation*}
\sum_{a \in \PP^1(\CC_v)} m_{\varphi,b}(a) \ = \ d \ .
\end{equation*}
We seek multiplicities on $\PP^1_{\Berk}$ which share these properties.

\vskip .1 in
Given $b \in \CC_v$, put $h_{b,\infty}(x) = \log_v([\varphi(T)-b)]_x)$.
Our starting point is the fact that $h_{b,\infty}(x) \in \BDV(\PP^1_{\Berk})$  
and
\begin{equation} \label{FMJ1}
\Delta_{\PP^1_{\Berk}}(h_{b,\infty}) \ = \
    \sum_{\varphi(a) = \infty} m_{\varphi,\infty}(a) \delta_a(x)
        -  \sum_{\varphi(a) = b} m_{\varphi,b}(a) \delta_a(x) \ ,
\end{equation}
as shown in Example~\ref{Example E.3}.  On the other hand, 
by Corollary~\ref{CorC2}, 
for each $x \in \PP^1_{\Berk}$ and $b \in \CC_v$,
$\delta(x,b)_{\infty} = [T-b]_x$.  Hence 
\begin{equation*}
\delta(\varphi(x),b)_{\infty} \ = \ [T-b]_{\varphi(x)} = [\varphi(T)-b]_x \ ,
\end{equation*}
so $h_{b,\infty}(x) = \log_v(\delta(\varphi(x),b)_{\infty})$.

Now let $\zeta \in \PP^1(\CC_v)$ be arbitrary, and put 
$h_{b,\zeta}(x) = \log_v(\delta(\varphi(x),b)_{\zeta})$.  
There is a constant $C$ such that
\begin{equation*}
\delta(x,b)_{\zeta} \ = \ C \cdot \frac{\delta(x,b)_{\infty}}
       {\delta(x,\zeta)_{\infty} \, \delta(b,\zeta)_{\infty}} \ .
\end{equation*}
Taking logarithms and applying the Laplacian shows that 
$h_{b,\zeta}(x) \in \BDV(\PP^1_{\Berk})$ and
\begin{equation} \label{FMJ2}
\Delta_{\PP^1_{\Berk}}(h_{b,\zeta}) \ = \
    \sum_{\varphi(a) = \zeta} m_{\varphi,\infty}(a) \delta_a(x)
        -  \sum_{\varphi(a) = b} m_{\varphi,b}(a) \delta_a(x) \ .
\end{equation}
It is easy to see that this formula holds for $b = \infty$ as well.

\vskip .05 in
Now consider $h_{b,\zeta}(x) = \log_v(\delta(\varphi(x),b)_{\zeta})$ for an
arbitrary $b \in \PP^1_{\Berk} \backslash \{\zeta\}$, still assuming
$\zeta \in \PP^1(\CC_v)$.

First suppose $b$ is of type II or type III.
Put $r = \diam_{\zeta}(b)$ and
fix a point $b_0 \in \cB(b,r)_{\zeta}^{-} \cap \PP^1(\CC_v)$.
Then $\delta(x,b)_{\zeta} = \max(r,\delta(x,b_0)_{\zeta})$ so
\begin{equation*}
h_{b,\zeta}(x) \ = \ \max(\log_v(r), \log_v(\delta(\varphi(x),b_0)_{\zeta}) \ .
\end{equation*}
By Proposition~\ref{PropF11}(D), $h_{b,\zeta}(x)$ is subharmonic on 
$U := \PP^1_{\Berk} \backslash \{\varphi^{-1}(\{\zeta\})$.
The set 
\begin{equation*}
Z \  := \ \varphi^{-1}(\cB(b,r)_{\zeta}) \ = \ 
\{x \in \PP^1_{\Berk} : \delta(\varphi(x),b_0)_{\zeta} \le r \}
\ \subset \ U 
\end{equation*}
is compact and disjoint from $\varphi^{-1}(\{\zeta\})$,
so by Proposition~\ref{PropF25}(B),
\begin{equation*}
\Delta_U(h_{b,\zeta})(Z) \ = \ \Delta_U(h_{b_0,\zeta})(Z) \ = \ -d \ .
\end{equation*}
On the other hand, $h_{b,\zeta}(x) = h_{b_0}(x,\zeta)$
on $\PP^1_{\Berk} \backslash Z$.
Since $h_{b_0,\zeta}(x) \in \BDV(\PP^1_{\Berk})$,
it follows that $h_{b,\zeta} \in \BDV(\PP^1_{\Berk})$, and
\begin{equation} \label{FTR1}
\Delta_{\PP^1_{\Berk}}(h_{b,\zeta}) \ = \ 
       \Delta_U(h_{b,\zeta})
       + \sum_{\varphi(a) = \zeta} m_{\varphi,\zeta}(a) \delta_a(x) \ .
\end{equation}

\vskip .05 in
We will now consider the functions $h_{b_0,\zeta}(x)$ and $h_{b,\zeta}(x)$ 
from another point of view, which will enable us to see that 
$\Delta_U(h_{b,\zeta})$ is a discrete measure. 
Let $a_1, \ldots, a_d, \zeta_1, \ldots,\zeta_d \in \PP^1(\CC_v)$ 
be the points such that $\varphi(a_i) = b_0$ and $\varphi(\zeta_i) = \zeta$, 
listed with multiplicities.  Put
\begin{equation} \label{FTS1}
f(x) \ = \ \sum_{i=1}^d \log_v(\delta(x,a_i)_{\zeta_i}) \ .
\end{equation}
By Example~\ref{Example E.2}, $f(x) \in \BDV(\PP^1_{\Berk})$ and 
\begin{equation} \label{FTR2}
\Delta_{\PP^1_{\Berk}}(f) \ = \ 
     \sum_{i=1}^d \delta_{\zeta_i}(x) - \sum_{i=1}^d \delta_{a_i}(x)
     \ = \ \Delta_{\PP^1_{\Berk}}(h_{b_0,\zeta}) \ .
\end{equation}
Both  $h_{b_0,\zeta}(x)$ and $f(x)$ are continuous on 
$\PP^1_{\Berk} \backslash \{a_1, \ldots, a_d, \zeta_1, \ldots,\zeta_d\}$, 
and their difference is 
a function in $\BDV(\PP^1_{\Berk})$ whose Laplacian is identically $0$,  
so $h_{b_0,\zeta}(x) = f(x) + C$ for some constant $C$.  

Consider $h_{b_0,\zeta}(x)$ and $h_{b,\zeta}(x)$ on the domain 
\begin{equation*}
V \ = \ \varphi^{-1}(\PP^1_{\Berk} \backslash \{b_0,\zeta\}) 
\ = \ \PP^1_{\Berk} \backslash \{a_1, \ldots, a_d, \zeta_1, \ldots,\zeta_d\} 
\end{equation*}
whose main dendrite $D$ is the union of the paths $(a_i,\zeta_j)$ for 
$i, j = 1, \ldots d$.  The function $h_{b_0,\zeta}(x)$ is harmonic
on $V$, so it is constant on branches off $D$.  Hence, 
$h_{b,\zeta}(x) = \max(\log_v(r),h_{b_0,\zeta}(x))$ is also constant on 
branches off $D$.

The set 
\begin{equation*}
X \ = \ \{x \in D : |h_{b_0,\zeta}(x) - \log_v(r)| \le 1\}
\end{equation*}
is a compact subset of $D$, since $h_{b_0,\zeta}(x)$ is continuous
on $D$ and  
$h_{b_0,\zeta}(x) \rightarrow \infty$ as $x \rightarrow \zeta_i$, 
$h_{b_0,\zeta}(x) \rightarrow -\infty$ as $x \rightarrow a_i$ for each $i$.  
Note also that $D$ has only finitely many branch points.  Choose a subgraph
$\Gamma_1 \subset D$ which has $X$ and all the branch points of $D$
in its interior, and then choose an exhaustion of $D$ by subgraphs $\Gamma_N$
with $\Gamma_1 \subset \Gamma_2 \subset \cdots$.  
Thus, for all $N \ge M$, $\Gamma_N$ is gotten by
extending edges of $\Gamma_M$ towards 
$\{a_1, \ldots, a_d, \zeta_1, \ldots, \zeta_d\}$.  
  
Since the restriction to $\Gamma_N$ 
of each term $\log_v(\delta(x,a_i)_{\zeta})$ in (\ref{FTR2})
belongs to $\CPA(\Gamma_N)$, the restriction of 
$h_{b_0,\zeta}(x) = f(x) + C$ 
to $\Gamma_N$ belongs to $\CPA(\Gamma_N)$.  Hence the restriction of 
$h_{b,\zeta}(x)$ to $\Gamma_N$ belongs to $\CPA(\Gamma_N)$ as well.  
It follows that
\begin{equation} \label{FTR3}
\Delta_{\Gamma}(h_{b,\zeta}) 
\ = \  \sum c_i \delta_{p_i}(x) - \sum d_j \delta_{q_j}(x)
\end{equation}
for certain points $p_i, q_j \in \Gamma_N$, with $c_i, d_j > 0$.  

On the other hand, by the compatibility of Laplacians with retraction maps
and (\ref{FTR1}),  
\begin{eqnarray}
\Delta_{\Gamma_N}(h_{b,\zeta}) 
& = & (r_{\Gamma_N})_*(\Delta_{\PP^1_{\Berk}}(h_{b,\zeta})) \notag \\
& = & (r_{\Gamma_N})_*(\Delta_U(h_{b,\zeta}))
  + (r_{\Gamma_N})_*(\sum_i m_{\varphi,\zeta}(\zeta_i) \delta_{\zeta_i}(x)) \ .
          \label{FTR4}
\end{eqnarray}  
Here $(r_{\Gamma_N})_*(\Delta_{U}(h_{b,\zeta}))$ is a non-positive 
measure with mass $-d$ on $r_{\Gamma_N}(Z)$, while 
$(r_{\Gamma_N})_*(\sum_i m_{\varphi,\zeta}(\zeta_i) \delta_{\zeta_j}(x))$
is a positive measure with mass $d$ supported on the points 
$r_{\Gamma}(\zeta_i)$, which are disjoint from $r_{\Gamma_N}(Z)$.  
Since the total mass $\Delta_{\Gamma_N}(h_{b,\zeta})(\Gamma_N) = 0$,
it follows that the  $p_i$ in (\ref{FTR3}) are the points 
$r_{\Gamma}(\zeta_i)$, and the $q_j$ belong to 
$Z \cap \Gamma_N = \{x \in \Gamma_N : h_{b_0,\zeta}(x) \le \log_v(r)\}$.  
Moreover each $c_i = m_{\varphi,\zeta}(\zeta_i) \in \ZZ$.

We claim that the $q_j$ belong to 
$\{x \in \Gamma_N : h_{b_0,\zeta}(x) = \log_v(r)\}$, 
and that each $d_j \in \ZZ$.  The first assertion is trivial:
if $h_{b_0,\zeta}(x) < \log_v(r)$, then $h_{b,\zeta}(z) = \log_v(r)$ 
is constant in a neighborhood of $x$ on $\Gamma_N$,
so $\Delta_{\Gamma_N}(h_{b,\zeta})(x) = 0$.  For the second, note
that 
\begin{equation} \label{FGUM1}
d_j \ = \ \Delta_{\Gamma}(h_{b,\zeta})(q_j) 
\ = \ - \sum_{\text{$\vec{v}$ at $q_j$}} d_{\vec{v}}h_{b,\zeta}(q_j) \ .
\end{equation}
Fix $\vec{v}$,
and let $e$ be the edge emanating from $q_j$ in the direction $\vec{v}$.
As $h_{b,\zeta}(x) = \max(\log_v(h_{b_0,\zeta}(x),\log_v(r))$, 
either $h_{b,\zeta}(x)$ is constant on $e$ near $q_j$, 
in which case $d_{\vec{v}}h_{b,\zeta}(q_j) = 0$;  or $h_{b,\zeta}(x)$ 
coincides with $h_{b_0,\zeta}(x)$ for $x \in e$ near $q_j$.  
In that case, by (\ref{FTS1}) 
\begin{equation} \label{FGUM2}
d_{\vec{v}}h_{b,\zeta}(q_j) \ = \ d_{\vec{v}}h_{b_0,\zeta}(q_j)
\ = \ \sum_{i=1}^d d_{\vec{v}}(\log_v(\delta(x,a_i)_{\zeta_i})(q_j) \ .
\end{equation}
However, each function $\log_v(\delta(x,a_i)_{\zeta_i})$ has constant
slope $1$ on the path $(a_i,\zeta_i)$, and slope $0$ on each branch
off that path.  
Hence each directional derivative $d_{\vec{v}}h_{b,\zeta}(q_j)$ belongs
to $\ZZ$, and $d_j \in \ZZ$.   
  
Since $X$ is contained in the interior of $\Gamma_N$, each $q_j$
belongs to the interior of $\Gamma_N$.  By the compatibility of 
Laplacians with retraction maps, this means the $q_j$ are independent of $N$.  
Since $h_{b,\zeta}(x)$ is constant on branches off $D$, necessarily  
$\Delta_{\PP^1_{\Berk}}(h_{b,\zeta})$ is supported on $\overline{D}$.  
Hence the negative part of $\Delta_{\PP^1_{\Berk}}(h_{b,\zeta})$ is 
supported on the $q_j$.  

In summary, we have shown that for 
$h_{b,\zeta}(x) = \log_v(\delta(\varphi(x),b)_{\zeta})$,  
\begin{equation} \label{FTR5}   
\Delta_{\PP^1_{\Berk}}(h_{b,\zeta}) 
\ = \ \sum_i m_{\varphi,\zeta_i}(\zeta) \delta_{\zeta_i}(x) 
                    - \sum_j d_j \delta_{q_j}(x) \ ,
\end{equation}
where each $q_j$ belongs to $D \cap Z$, 
each $d_j > 0$ is an integer, and $\sum d_j = \deg(\varphi)$.  
Note that  
\begin{equation} \label{FTR6}
Z \ = \ Z_{b,\zeta} \ = \ 
\{x \in \PP^1_{\Berk} : \delta(\varphi(x),b)_{\zeta} \le \diam_{\zeta}(b) \} \ .
\end{equation} 

\vskip .05 in
We now will consider $h_{b,\zeta}(x)$ from a third point of view,
and show that the points $q_j$ belong to $\varphi^{-1}(\{b\})$.  
To see this, let $\xi \in \PP^1(\CC_v)$ be arbitrary, 
and note that there is a constant $C$
such that for all $y \in \PP^1_{\Berk}$. 
\begin{equation} \label{FMM0}
\delta(y,b)_{\xi} \ = \ C \cdot \frac{\delta(y,b)_{\zeta}}
                      {\delta(y,\xi)_{\zeta} \, \delta(b,\xi)_{\zeta}} \ .
\end{equation}
Here $\delta(b,\xi)_{\xi}$ is a finite nonzero constant 
since $b \notin \PP^1(\CC_v)$.  Replacing $y$ by $\varphi(x)$, 
taking logarithms and applying the Laplacian, we find that  
\begin{eqnarray}
\Delta_{\PP^1_{\Berk}}(h_{b,\xi}) 
& = & \Delta_{\PP^1_{\Berk}}(h_{b,\zeta}) 
              - \Delta_{\PP^1_{\Berk}}(h_{\xi,\zeta}) \notag \\
& = & (\sum_i m_{\varphi,\zeta_i}(\zeta) \delta_{\zeta_i}(x)
                         -\sum_j d_j \delta_{q_j}(x)) \notag \\
&   & \qquad \quad
         - (\sum_i m_{\varphi,\zeta_i}(\zeta) \delta_{\zeta_i}(x)  
                 - \sum_i m_{\varphi,\xi_i}(\xi) \delta_{\xi_i}(x)) \notag \\                          
& = &  \sum_i m_{\varphi,\xi_i}(\xi) \delta_{\xi_i}(x) 
                   - \sum_j d_j \delta_{q_j}(x)  \label{FMM1C}
\end{eqnarray}
where the $\xi_i$ are the points with $\varphi(\xi_i) = \xi$.    
Thus, the $q_j$ depend only on $b$, not $\zeta$.  Applying
(\ref{FTR6}) with $\zeta$ replaced by $\xi$, we see that each $q_j$
belongs to  
\begin{equation*}
\bigcap_{\xi \in \PP^1(\CC_v)}  
\{x \in \PP^1_{\Berk} : \delta(\varphi(x),b)_{\xi} \le \diam_{\xi}(b) \} \ .
\end{equation*}
We now apply following lemma to conclude that $\varphi(q_j) = b$:  

\begin{lemma} \label{LemG1}  For each $b \in \PP^1_{\Berk}$
\begin{equation*}
\bigcap_{\xi \in \PP^1(\CC_v)}  
\{y \in \PP^1_{\Berk} : \delta(y,b)_{\xi} \le \diam_{\xi}(b) \} \ = \ \{b\} \ .
\end{equation*}
\end{lemma}

\begin{proof}
Suppose $y \ne b$.  Consider the path $\Gamma = [y,b]$.
Let $w$ be a type II point in the interior of $\Gamma$,
and let $\xi \in \PP^1(\CC_v)$ be a point with $r_{\Gamma}(\xi) = w$.
Then $w$ lies on the interior of the path $[b,\xi]$ and is the point
where the paths $[y,\xi]$ and $[b,\xi]$ meet, so
$\delta(y,b)_{\xi} = \delta(w,b)_{\xi} > \delta(b,b)_{\xi} = \diam_{\xi}(b)$.
\end{proof}  

\vskip .1 in
For each $q_j \in \varphi^{-1}(b)$, we define the analytic multiplicity  
$m_{\varphi,q_j}(b) = d_j$.  For all other $q$, put $m_{\varphi,q}(b) = 0$.
Thus $0 \le m_{\varphi,q}(b) \in \ZZ$, 
$\sum_q m_{\varphi,q}(b) = \deg(\varphi)$, and (\ref{FTR5}) becomes 
\begin{equation} \label{FTR7}   
\Delta_{\PP^1_{\Berk}}(h_{b,\zeta}) 
\ = \ \sum_{\varphi(\zeta_i) = \zeta} 
        m_{\varphi,\zeta_i}(\zeta) \delta_{\zeta_i}(x) 
   - \sum_{\varphi(q_j) = b} m_{\varphi,q_j}(\zeta) \delta_{q_j}(x) \ .
\end{equation}

So far this only holds for points $b$ of type II or III (and I).
We will now extend it to points of type IV.
If $b$ is of type IV, take a sequence of type II points $b_1, b_2, \ldots$
approaching $b$, with $b \in \cB(b_i,\diam_{\zeta}(b_i)_{\zeta})$
for each $i$.  The functions $h_{b_i,\zeta}(x)$
decrease monotonically to $h_{b,\zeta}(x)$.
By Proposition~\ref{PropF11}(D), $h_{b,\zeta}(x)$ is subharmonic in $U$,
and by Proposition~\ref{PropF23}(B), the measures $\Delta_U(h_{b_i,\zeta})$
converge weakly to $\Delta_U(h_{b,\zeta})$ on simple subdomains of $U$.

Suppose $q \in \supp(\Delta_{U}(h_{b,\zeta}))$, and let $V$ be a neighborhood
of $V$ of $q$ contained in $U$.  For all sufficiently large $i$,
there are points of $\supp(\Delta_U(h_{b_i,\zeta}))$ in $V$.  Since 
each $\Delta_U(h_{b_i,\zeta})$ is a sum of point masses with positive 
integer weights, $\Delta_U(h_{b,\zeta})(V)$ is a positive integer.  
Choosing a sequence of neighborhoods $V_1 \supset V_2 \supset \cdots$
with $\cap_{k=1}^{\infty} V_k = \{b\}$, we see that 
\begin{equation*}
\Delta_U(h_{b,\zeta})(V_1) \ \ge \ \Delta_U(h_{b,\zeta})(V_2) \ \ge \ \cdots \ .
\end{equation*}
As each $\Delta_U(h_{b,\zeta})(V_k)$ is a positive integer, 
the $\Delta_U(h_{b,\zeta})(V_k)$ must stabilize at an integer $m > 0$,
and $\Delta_U(h_{b,\zeta})(\{q\}) = m$.  We define $m_{\varphi,q}(b) = m$.  
Since the total mass of $\Delta_U(h_{b,\zeta})$ is $\deg(\varphi)$, 
$\Delta_U(h_{b,\zeta})$ must be a discrete measure, and we put
$m_{\varphi,q}(b) = 0$ if $q \notin \supp(\Delta_U(h_{b,\zeta}))$.  
Each $q \in \supp(\Delta_U(h_{b,\zeta}))$ is a limit of points 
$q_i \in \supp(\Delta_U(h_{b_i,\zeta}))$ with $\varphi(q_i) = b_i$, 
so by continuity 
\begin{equation*}
\varphi(q) \ = \ \lim_{i \rightarrow \infty} \varphi(q_i) 
\ = \ \lim_{i \rightarrow \infty} \ = \ b \ .
\end{equation*}
Since $h_{b,\zeta}(x)$ coincides with $h_{b_i,\zeta}(x)$ in a neighborhood
of each $\zeta_i$, it follows that $h_{b,\zeta} \in \BDV(\PP^1_{\Berk})$ and
(\ref{FTR7}) holds.  By a computation similar to (\ref{FMM1C}), 
the $m{\varphi,q}(b)$ depend only on $b$, and not on $\zeta$.  

\vskip .05 in
Finally, we lift the restriction that $\zeta \in \PP^1(\CC_v)$.  
Formula  (\ref{FMM0}) holds for any $\xi \in \PP^1_{\Berk}$.  Repeating the
computation (\ref{FMM1C}), we see that for any $b, \xi \in \PP^1_{\Berk}$
with $b \ne \xi$, the function 
$h_{b,\xi}(x) = \log_v(\delta(\varphi(x),b)_{\xi})$ belongs to 
$\BDV(\PP^1_{\Berk})$ and its Laplacian satisfies (\ref{FTR7}).  
We summarize this with the following Proposition.

\begin{proposition} \label{PropG2}
Let $\varphi(T) \in \CC_v(T)$ be a nonconstant rational function.   
Then for each \ $b, \zeta \in \PP^1_{\Berk}$ with $b \ne \zeta$,  
the function $h_{b,\zeta}(x) = \log_v(\delta(\varphi(x),b)_{\zeta})$ 
belongs to $\BDV(\PP^1_{\Berk})$.

Furthermore, as $b$ ranges over $\PP^1_{\Berk}$, there
is a unique way to assign multiplicities $m_{\varphi,b}(q)$ 
to points $q \in \PP^1_{\Berk}$ in such a way that 

$A)$ \ $0 \le m_{\varphi,b}(q) \in \ZZ$,  
         with $m_{\varphi,b}(q) = 0$ unless $\varphi(q) = b$.          
         
$B)$ \ $\sum_{\varphi(q) = b} m_{\varphi,b}(q) = \deg(\varphi)$\ .            
         
$C)$  For points $q$ of type {\rm I}, the $m_{\varphi,b}(q)$ coincide with
the usual algebraic multiplicities. 

$D)$  For all $\zeta \ne b \in \PP^1_{\Berk}$, 
\begin{equation} \label{FTR8}   
\Delta_{\PP^1_{\Berk}}(h_{b,\zeta}) 
\ = \ \sum_{\varphi(\xi) = \zeta} 
        m_{\varphi,\zeta}(\xi) \delta_{\xi}(x) 
   - \sum_{\varphi(q) = b} m_{\varphi,b}(q) \delta_{q}(x) \ .
\end{equation}
\end{proposition}

Although the analytic multiplicities may seem quite mysterious, for points
$q$ of type II or III there is a simple formula for them.
Put $b = \varphi(q)$.  Then $\PP^1_{\Berk} \backslash \{b\}$ has
at least two components.  Let $b_0$ and $\zeta$ be arbitrary type I points
lying in different components of $\PP^1_{\Berk} \backslash \{b\}$,
and let $a_1, \ldots, a_d$, $\zeta_1, \ldots, \zeta_d$ be the points with
$\varphi(a_i) = b_0$, $\varphi(\zeta_i) = \zeta$, listed with their
usual algebraic multiplicities.  There are also at least two components
of $\PP^1_{\Berk} \backslash \{q\}$.  For each component $U$ of
$\PP^1_{\Berk} \backslash \{q\}$, let $N_a(U)$ be the number of points
$a_i \in U$ (counted with multiplicities), and let $N_{\zeta}(U)$
be the number of points $\zeta_i \in U$ (counted with multiplicities).
Put $N^+(U) = \max(0,N_{\zeta}(U)-N_{a}(U))$.  

\begin{corollary} \label{CorG3}
Let $\varphi(T) \in \CC_v(T)$ be a nonconstant rational function.
For each $q \in \PP^1_{\Berk}$ of type {\rm II} or type {\rm III}, 
if $\varphi(q) = b$, then  
\begin{equation} \label{FVU1}
m_{\varphi,b}(q) \ = \ 
\sum \begin{Sb} \text{components} \\ 
            \text{of $\PP^1_{\Berk} \backslash \{q\}$} \end{Sb}  N^+(U)  \ .
\end{equation}
\end{corollary} 

\begin{proof}
This is a reformulation of (\ref{FGUM1}), (\ref{FGUM2}) 
in the discussion leading to Proposition \ref{PropG2}.  
To obtain (\ref{FVU1}), recall that
\begin{eqnarray}
h_{b_0,\zeta}(x) & = & \sum_{i=1}^d \log_v(\delta(x,a_i)_{\zeta_i}) + C \ , 
                             \label{FBC1} \\
h_{b,\zeta}(x) & = & \max(h_{b_0,\zeta}(x),\log_v(r)) \ . \label{FBC2}
\end{eqnarray}
Given a direction vector $\vec{v}$ at $q$ in $\Gamma_N$, 
let $e$ be the edge of $\Gamma_N$ emanating from $q$ 
in the direction $\vec{v}$.
The directional derivative $d_{\vec{v}}h_{b_0,\zeta}(q)$ 
is the sum of the terms
$d_{\vec{v}}(\log_v(\delta(x,a_i)_{\zeta_i})$, for $i = 1, \ldots, d$.  
Each such term is $+1$
if the path from $a_i$ to $\zeta_i$ passes through edge $e$ in the
direction $\vec{v}$.  It is $-1$ if the path from $a_i$ to $\zeta_i$ passes
through $e$ in the direction opposite $\vec{v}$, and $0$ if the
path does not pass through edge $e$.
The edge $e$  belongs to a component $U$
of $\PP^1_{\Berk} \backslash \{q\}$;
it follows that $d_{\vec{v}}h_{b_0,\zeta}(q) = N_{\zeta}(U) - N_{a}(U)$.
If $d_{\vec{v}}h_{b_0,\zeta}(q) \ge 0$ then
$d_{\vec{v}}h_{b,\zeta}(q) = d_{\vec{v}}h_{b_0,\zeta}(q)$;
otherwise $d_{\vec{v}}h_{b,\zeta}(q) = 0$.  For the components $U$ of
$\PP^1_{\Berk} \backslash \{q\}$ not containing edges of $\Gamma_N$,
automatically $N_a(U) = N_{\zeta}(U) = 0$.  Hence
\begin{eqnarray*}
m_{\varphi,b}(q)
& = & \sum_{\text{$\vec{v}$ at $q$}} d_{\vec{v}}h_{b,\zeta}(q) \\
& = & \sum_U N^+(U) \ .
\end{eqnarray*}
\end{proof}

\vskip .05 in
Note that we have not yet shown that if $\varphi(q) = b$, 
then $m_{\varphi,b}(q) > 0$.  That fact follows from 

\begin{proposition} \label{PropG4}  
Let $\varphi(T) \in \CC_v(T)$ be a nonconstant rational function.
Let $b_1, b_2, \ldots  \in \PP^1_{\Berk}$ approach $b$
in the Berkovich topology.  
Then the measures
\begin{equation*}
\mu_{b_i} = \sum_{\varphi(q_{ij}) = b_i}
        m_{\varphi,b_i}(q_{ij}) \delta_{q_{ij}}(x)
\end{equation*}         
converge weakly to
$\mu_b = \sum_{\varphi(q_j) = b} m_{\varphi,b}(q_j) \delta_{q_j}(x)$ 
on each simple subdomain of $\PP^1_{\Berk}$.
\end{proposition}

\begin{proof}
Let $V$ be a simple subdomain of $\PP^1_{\Berk}$.  Without loss
we can assume $\PP^1_{\Berk} \backslash \Vbar$ is nonempty.
Take $\zeta \in \PP^1(\CC_v) \backslash \Vbar$.
Put $U = \PP^1_{\Berk} \backslash \varphi^{-1}(\{\zeta\}$;   
we must show the measures $\Delta_U(h_{b_i,\zeta})$ converge 
weakly to $\Delta_U(h_{b_i,\zeta})$ on $V$.

By  Proposition~\ref{PropF23}, it is enough to show that for each
subgraph $\Gamma \subset V$, the subharmonic functions $h_{b_i,\zeta}(x)$
converge uniformly to $h_{b,\zeta}(x) = \log_v(\delta(\varphi(x),b)_{\zeta}$
on $\Gamma$.  For this, it suffices show that the functions
$\log_v(\delta(y,b_i)_{\zeta})$
converge uniformly to
$\log_v(\delta(y,b)_{\zeta})$ on $\varphi(\Gamma)$. 
There are two cases to consider.

\vskip .05 in
First, suppose $b \notin \varphi(\Gamma)$.  Since $\Gamma$ is compact
and $\varphi$ is continuous, $\varphi(\Gamma)$ is compact.
Since $\PP^1_{\Berk}$ is a compact Hausdorff space, there are neighborhoods
$W$ of $\varphi(\Gamma)$ and $Z$ of $b$ with disjoint closures.
Note that $\zeta \notin \varphi(\Gamma)$ since $\Gamma \subset U$, 
and $b \ne \zeta$ by hypothesis.  
After deleting a suitably small neighborhood of $\zeta$ from $W$ and $Z$
if necessary,  we can assume that $\zeta \notin \Wbar \cup \Zbar$.  
By Proposition~\ref{PropC10}, 
$\delta(y,z)_{\zeta}$ is continuous and bounded below
on $\Wbar \times \Zbar$.  Since $\PP^1_{\Berk}$ is a metric space,
$\log_v(\delta(y,z)_{\zeta})$ is uniformly continuous on
$\Wbar \times \Zbar$.  Hence, as $b_i \rightarrow b$, the functions
$h_{b_i,\zeta}(x) = \log_v(\delta(\varphi(x),b_i)_{\zeta})$ converge
uniformly to $h_{b,\zeta}(x) = \log_v(\delta(\varphi(x),b)_{\zeta})$
on $\Gamma$.  

\vskip .05 in
Next, suppose $b \in \varphi(\Gamma)$.  
In particular, $b \in \PP^1_{\Berk} \backslash \PP^1(\CC_v)$, 
so there is a constant $C$ such that
\begin{equation*}
\log_v(\delta(y,z)_{\zeta}) \ = \ j_b(y,z) - j_b(y,\zeta) - j_b(z,\zeta) + C
\end{equation*}
for all $y, z \in \PP^1_{\Berk}$.  Noting that $j_b(y,b) = j_b(b,\zeta) = 0$,
we see that
\begin{eqnarray}
|\log_v(\delta(y,b_i)_{\zeta}) - \log_v(\delta(y,b)_{\zeta})|
& = & |j_b(y,b_i) - j_b(b_i,\zeta)| \notag \\
& \le & j_b(y,b_i) + j_b(b_i,\zeta) \ . \label{FMS1}
\end{eqnarray}

Fix $\varepsilon > 0$, and let $\rho(x,y)$ denote the path length metric
on $\PP^1_{\Berk} \backslash \PP^1(\CC_v)$.
Note that we do not know that $\varphi(\Gamma)$ is a subgraph:
it is compact and connected, but it could conceivably have infinitely
many branch points.  Nonetheless it is shown in Lemma \ref{LemG11} below
that the path metric topology on $\varphi(\Gamma)$ coincides
with the induced topology from $\PP^1_{\Berk}$.
Let $W_1 \subset \PP^1_{\Berk}$ be a connected open set containing $b$
such that
\begin{equation*}
\varphi(\Gamma) \cap W_1 \ \subset \ B_{\rho}(b,\varepsilon)
\ := \ \{ z \in \varphi(\Gamma) : \rho(z,b) < \varepsilon \} \ .
\end{equation*}
Let $c \in (b,\zeta)$ be a point such that $\rho(c,b) < \varepsilon$,
put $\tGamma = [b,c]$, and put $W_2 = r_{\tGamma}^{-1}([b,c))$.
Then $W_2$ is an open set containing $b$, and $r_{[b,\zeta]}(W_2) = [b,c)$.
Put $W = W_1 \cap W_2$.  Since the $b_i$ converge to $b$ in the Berkovich
topology, there is an $N$ such that $b_i \in W$ for all $i \ge N$.

Take $i \ge N$, and let $y$ range over $\varphi(\Gamma)$.
Recall that $j_b(y,b_i) = \rho(b,w)$,
where $w$ is the point where the paths $[y,b]$ and $[b_i,b]$ meet.
Since $b_i \in W$ and $y \in \varphi(\Gamma)$, both of which are connected,
it follows that
$w \in \varphi(\Gamma) \cap W \subset B_{\rho}(b,\varepsilon)$.
Hence $0 \le j_b(y,b_i) = \rho(b,w) < \varepsilon$.
Likewise, $j_b(b_i,\zeta)$ is the point $t$ where the paths
$[b_i,b]$ and $[\zeta,b]$ meet.  Since $b_i \in W_2$,
it follows that $t = r_{[b,\zeta]}(b_i) \in [b,c)$, so
$0 \le j_b(b_i,\zeta) \le \rho(b,t) < \varepsilon$.
By (\ref{FMS1}),
\begin{equation*}
|\log_v(\delta(y,b_i)_{\zeta}) - \log_v(\delta(y,b)_{\zeta})|
\ < \ 2 \varepsilon \ .
\end{equation*}
\end{proof}

\vskip .1 in
Before proving the assertion concerning the topology of $\varphi(\Gamma)$,
let us note some consequences of Proposition \ref{PropG4}.

\begin{corollary} \label{CorG5}
Let $\varphi(T) \in \CC_v(T)$ be a nonconstant rational function.
Then for each $q$ and $b$ in $\PP^1_{\Berk}$,
\begin{equation*}
\text{$m_{\varphi,b}(q) > 0$ \quad if and only if \quad $\varphi(q) = b$\ .}
\end{equation*}
\end{corollary}

\begin{proof} Only the direction ($\Longleftarrow$) requires proof.
Suppose $\varphi(q) = b$.  Take a sequence of type I points $q_1, q_2, \cdots$
converging to $q$, and put $b_i = \varphi(q_i)$ for each $i$.
For Type I points, the analytic multiplicites coincide with the usual
algebraic multiplicities, so for each $q_{ij} \in \varphi^{-1}(\{b_i\})$
we have $m_{\varphi,b_i}(q_{ij}) > 0$.  
In particular $m_{\varphi,b_i}(q_i) > 0$.
Since the measure $\mu_b$ in Proposition \ref{PropG4}
is the weak limit of the measures $\mu_{b_i}$ and all the multiplicities
are integers, it follows that $m_{\varphi,q}(b) > 0$.
\end{proof}

\begin{corollary} \label{CorG6}
Let $\varphi(T) \in \CC_v(T)$ be a nonconstant rational function.
Then  $\varphi : \PP^1_{\Berk} \rightarrow \PP^1_{\Berk}$
is open in the Berkovich topology.  It is surjective, and for each
$b \in \PP^1_{\Berk}$ there are at most $d = \deg(\varphi)$ points
$q$ with $\varphi(q) = b$.
\end{corollary}

\begin{proof}
If there were an open set $U \subset \PP^1_{\Berk}$ such that
$\varphi(U)$ was not open, there would be a point $q \in U$
such that $b = \varphi(q)$ was a limit of points
$b_1, b_2, \ldots \notin \varphi(U)$.  By Corollary \ref{CorG5},
$m_{\varphi,b}(q) > 0$, and by Proposition \ref{PropG4}
the measures $\mu_{b_i}$ converge weakly to $\mu_b$.
Hence for each sufficiently large $i$ there is a point
$q_{ij} \in \supp(\mu_{b_i}) = \varphi^{-1}(\{b_i\}$ which belongs to $U$.
It follows that $b_i = \varphi(q_{ij}) \in U$,
contradicting our assumption that $b_i \notin U$.

The fact that the measure $\mu_b$ in
Proposition \ref{PropG4} has total mass $d$ gives surjectivity.  
Corollary \ref{CorG5}, together with the fact that each 
$m_{\varphi,q}(b)$ is a non-negative integer,
shows there are at most $d$ points with $\varphi(q) = b$.
\end{proof}

\begin{corollary} \label{CorG7}
Let $\varphi(T) \in \CC_v(T)$ be a nonconstant rational function,
and take $q \in \PP^1_{\Berk}$.  For each neighborhood $U$ of $q$,
there is a neighborhood $V$ of $q$ such that the function
\begin{equation*}
M_U(x) \ = \ \sum \begin{Sb} y \in U \\ \varphi(y) = \varphi(x) \end{Sb}
                 m_{\varphi,\varphi(x)}(y)
\end{equation*}
is constant for $x \in V$.  
\end{corollary}

\begin{proof}  This follows immediately from weak convergence;  
it is equivalent to the assertion that if
$q_1, q_2, \ldots \in \PP^1_{\Berk}$
is a sequence of points converging to $q$, and if $b = \varphi(q)$
and $b_i = \varphi(q_i)$, then for the measures $\mu_b$ and $\mu_{b_i}$
in Proposition \ref{PropG4}, 
$\lim_{i \rightarrow \infty} \mu_{b_i}(U) = \mu_b(U)$.
\end{proof}

\vskip .1 in
Define the ``ramification function''
$R_{\varphi} : \PP^1_{\Berk} \rightarrow \ZZ_{\ge 0}$ by 
\begin{equation} \label{FER1}
R_{\varphi}(x) \ = \ m_{\varphi,\varphi(x)}(x) - 1 \ .
\end{equation}

\begin{corollary} \label{CorG8}
Let $\varphi(T) \in \CC_v(T)$ be a nonconstant rational function.
The ramification function $R_{\varphi}(x)$ is upper semicontinuous,
and takes integer values in the range $0 \le n \le \deg(\varphi)-1$.
Furthermore 
$R_{\varphi}(q) = 0$ if and only if $\varphi$ is locally injective at $q$.  
\end{corollary}

\begin{proof}
Proposition \ref{PropG2} and Corollary \ref{CorG5} show 
that $R_{\varphi}(x)$ takes integer values in the range $0 \le y \le d-1$.

For the upper semicontinuity, fix $q \in \PP^1_{\Berk}$
and let $U$ be a neighborhood of $q$ such that $q$ is the only point in 
$\varphi^{-1}(\{\varphi(q)\}) \cap U$.  
By Corollary \ref{CorG7}, there is a neighborhood $V$ of
$q$ such that for each $x \in V$,
\begin{equation*}
\sum \begin{Sb} y \in U \\ \varphi(y) = \varphi(x) \end{Sb}
     m_{\varphi,\varphi(x)}(y) \ = \ m_{\varphi,\varphi(q)}(q)
\end{equation*}
In particular, $m_{\varphi,\varphi(x)}(x) \le m_{\varphi,\varphi(q)}(q)$ 
for each $x \in V$.

If $R_{\varphi}(q) = 0$ then 
$m_{\varphi,\varphi(q)}(q) = 1$, so Corollary \ref{CorG7} shows 
there is a neigbhorhood $V$ of $q$ such that $m_{\varphi,\varphi(x)} = 1$
for all $x \in V$, and for each $x \in V$ there are no other points 
$x^{\prime}$ in $U$ with $\varphi(x^{\prime}) = \varphi(x)$.  
Conversely, if $\varphi$ is one-to-one on a neighborhood $V$ of $q$, 
the type I points all have multiplicity $1$, and they are dense in $V$.
By Proposition \ref{PropG4}, $m_{\varphi,\varphi(q)}(q) = 1$ for 
all $q \in V$.   
\end{proof}

\begin{example}
It is well known that there are only finitely 
many points in $\PP^1(\CC_v)$ where $\varphi(T)$ is ramified.  However, the 
situation is very different on Berkovich space.  

Consider $\varphi(T) = T^2$.  For a point $x \in \AA^1_{\Berk}$ 
corresponding to a disc $B(a,r)$ (i.e. for a point of type II or III), 
Proposition \ref{PropG4} shows that $\varphi$ 
is ramified at $x$ (that is, $R_{\varphi}(x) \ge 1$)
if and only if there are distinct points
$a_1, a_2 \in B(a,r)$ for which $\varphi(a_1) = \varphi(a_2)$, that is,
$a_1^2 = a_2^2$, so $a_1 = -a_2$.  It follows that $x$ is ramified
if and only if $r \ge |a_1-(-a_1)| = |2 a_1|_v$ for some $a_1 \in B(a,r)$.  

If $r \ge |2a_1|$ for some $a_1 \in B(a,r)$, 
then $r \ge |a-a_1|_v \ge |2|_v \cdot |a-a_1|_v = |2a-2a_1|_v$,
so  $r \ge |2a|_v$ by the ultrametric inequality.  Conversely, if
$r \ge |2a|_v$, then $-a \in B(a,r)$.  Thus $\varphi$ is ramified at 
$x$ if and only if $r \ge |2a|_v$.  

\vskip .05 in
If the residue characteristic of $\CC_v$ is not $2$, then $|2a|_v = |a|_v$.
so $\varphi$ is ramified at the point corresponding to $B(a,r)$ 
if and only $r \ge |a|_v$, or equivalently, if and only if $B(a,r) = B(0,r)$.
These are the points corresponding to the interior of the path $[0,\infty]$
in $\PP^1_{\Berk}$.  One sees easily that $\varphi$ is not ramified at any
points of type IV.  On the other hand its classical ramification points
in $\PP^1(\CC_v)$ are precisely $0$ and $\infty$.  Thus, the ramification
locus of $\varphi(T)$ is $[0,\infty]$.

\vskip .05 in
If the residue characteristic of $\CC_v$ is $2$, then $|2a|_v = |a|_v/2$
(assuming $|x|_v$ is normalized so that it extends the usual absolute value 
on $\QQ_2$).  In this case, $\varphi$ is ramified at a point $x$ 
corresponding to a disc $B(a,r)$ if and only if $r \ge |a|_v/2$. 
 A little thought shows this is 
a much larger set than $[0,\infty]$:  for each $a \in \CC_v^{\times}$
it contains the line of discs $\{B(a,t) : |a|_v/2 \le t \le |a|_v\}$
which is a path leading off $[0,\infty]$.  The union of these paths gives
an infinitely branched dendritic structure whose interior is open in the
path metric topology.  
       If one views $\PP^1_{\Berk}$ as $\AA^1_{\Berk} \cup \{\infty\}$, 
the ramification locus of $\varphi$ could be visualized as an inverted 
Christmas tree.  If one considers the small model of $\PP^1_{\Berk}$, 
it would look like a spindle around the axis $[0,\infty]$. 
{\ $\square$}
\end{example}

\vskip .1 in
We now return to the assertion concerning the path-length topology 
on $\varphi(\Gamma)$ which was used in the proof of
Proposition~\ref{PropG4}.  This will be a consequence of continuity properties
of $\varphi$ relative to the path-length topology 
which are of independent interest.  

\vskip .05 in
Put $\PP^1_{\Berk,0} = \PP^1_{\Berk} \backslash \PP^1(\CC_v)$, the 
set of type II, III, and IV points in $\PP^1_{\Berk}$.  We have 
encountered $\PP^1_{\Berk,0}$ as a set many times before, 
but now we will study it as a space in its own right, giving it
the path-length topology.  Let $\rho(x,y)$ be the path length metric;  
it is finite for all $x, y \in \PP^1_{\Berk,0}$.

\begin{lemma}  \label{LemG9} 
Let $\varphi(T) \in \CC_v(T)$ be a nonconstant rational function
of degree $d = \deg(\varphi)$.  
Then $\varphi$ acts on $\PP^1_{\Berk,0}$, 
and for all $x, y \in \PP^1_{\Berk,0}$
\begin{equation*}
\rho(\varphi(x),\varphi(y)) \ \le \ d \cdot \rho(x,y) \ .
\end{equation*}
\end{lemma}

\begin{proof}  We already know that $\varphi$ takes $\PP^1_{\Berk,0}$
to $\PP^1_{\Berk,0}$.  For the assertion about path distances,
we use the notation from the discussion preceding Proposition 
\ref{PropG2}.  Consider the subgraph $\Gamma = [\varphi(x),\varphi(y)]$ 
in $\PP^1_{\Berk}$.  Take $b_0, \zeta \in \PP^1(\CC_v)$ with
$r_{\Gamma}(b_0) = \varphi(x)$, $r_{\Gamma}(\zeta) = \varphi(y)$.  Thus 
$\varphi(x)$ and $\varphi(y)$ lie on the path $[b_0,\zeta]$.  It follows that
\begin{equation} \label{FXN1B} 
\rho(\varphi(x),\varphi(y)) \ = \ 
   | \log_v(\delta(\varphi(x),b_0)_{\zeta}) 
                 - \log_v(\delta(\varphi(y),b_0)_{\zeta}) |  \ .
\end{equation}
However, the function 
$h_{b_0,\zeta}(x) = \log_v(\delta(\varphi(x),b_0)_{\zeta})$ has the 
representation (\ref{FBC1}) 
\begin{equation*}                
h_{b_0,\zeta}(x) \ = \  \sum_{i=1}^d \log_v(\delta(x,a_i)_{\zeta_i}) + C 
\end{equation*}
where the $a_i$ and $\zeta_i$ are the pre-images of $b_0$ and $\zeta$, 
listed with multiplicities.  Inserting this in (\ref{FXN1B}) gives 
\begin{equation} \label{FNX2B} 
\rho(\varphi(x),\varphi(y)) \ \le \ 
   \sum_{i=1}^d | \log_v(\delta(x,a_i)_{\zeta_i}) 
                 - \log_v(\delta(y,a_i)_{\zeta_i}) |  \ .
\end{equation}
For each $q \in \PP^1_{\Berk,0}$  
the representation of the Hsia kernel in terms of $j_q(z,w)$ gives 
\begin{eqnarray*}
& & \log_v(\delta(x,a_i)_{\zeta_i}) 
                 - \log_v(\delta(y,a_i)_{\zeta_i}) \\                 
&  & \qquad \qquad \qquad = \ 
            (j_q(x,a_i) - j_q(x,\zeta_i) - j_q(a_i,\zeta_i) + C) \\ 
&  & \qquad \qquad \qquad \qquad \qquad 
             - (j_q(y,a_i) - j_q(y,\zeta_i) - j_q(a_i,\zeta_i) + C) \\ 
&  & \qquad \qquad \qquad = \ 
             j_q(x,a_i) - j_q(x,\zeta_i) + j_q(y,\zeta_i) - j_q(y,a_i) \ .
\end{eqnarray*}
Now take $q = x$;  then $j_x(x,a_i) = j_x(x,\zeta_i) = 0$, so 
\begin{equation*} 
|\log_v(\delta(x,a_i)_{\zeta_i}) 
                 - \log_v(\delta(y,a_i)_{\zeta_i})|
\ = \ |j_x(y,\zeta_i) - j_x(y,a_i) | \ .                 
\end{equation*}
Here $j_x(y,\zeta_i) = \rho(x,w)$ where $w \in [x,y]$ is the point
where the paths $[y,x]$ and $[\zeta_i,x]$ meet.  Likewise,
$j_x(y,a_i) = \rho(x,t)$ where $t \in [x,y]$ is the point
where the paths $[y,x]$ and $[a_i,x]$ meet.  Since $w, t \in [x,y]$
clearly $|\rho(x,w) - \rho(x,t)| < \rho(x,y)$.
Hence 
\begin{equation*}
|\log_v(\delta(x,a_i)_{\zeta_i}) 
                 - \log_v(\delta(y,a_i)_{\zeta_i})| \ \le \ \rho(x,y) \ .
\end{equation*}
Inserting this in (\ref{FNX2B}) gives the result.
\end{proof}

\vskip .1 in
The bound in the lemma is sharp, as shown by the example
$\varphi(T) = T^d$, which takes the point $x \in [0,\infty]$ corresponding
to  $B(0,r)$ to the point $\varphi(x) \in [0,\infty]$ 
corresponding to $B(0,r^d)$.  
Using the definition of the logarithmic path length, we see that 
$\rho(\varphi(x),\varphi(y)) = d \cdot \rho(x,y)$ 
for each $x,y \in (0,\infty)$. 
                  
\begin{corollary} \label{CorG10}
Let $\varphi(T) \in \CC_v(T)$ be a nonconstant rational function.  Then 

$A)$  $\varphi : \PP^1_{\Berk} \rightarrow \PP^1_{\Berk}$ is continuous
for the Berkovich topology.

$B)$ $\varphi : \PP^1_{\Berk,0} \rightarrow \PP^1_{\Berk,0}$ is continuous
for the path length topology.
\end{corollary}

\begin{proof} We already know $\varphi$ is continuous on $\PP^1_{\Berk}$
for the Berkovich topology;  this follows from the definition of the Berkovich
topology.  The assertion about its continuity on $\PP^1_{\Berk,0}$ for
the path length topology is immediate from Lemma \ref{LemG9}.
\end{proof}

\vskip .1 in
Write $\ell(\Gamma)$ for the path length of a subgraph 
$\Gamma \subset \PP^1_{\Berk}$.  For any connected subset 
$Z \subset \PP^1_{\Berk}$, define the path length
\begin{equation*}
\ell(Z) \ = \ \sup_{\text{subgraphs $\Gamma \subset Z$}} \ell(\Gamma) \ .
\end{equation*}

\begin{lemma} \label{LemG11}
Let $\varphi(T) \in \CC_v(T)$ be a nonconstant rational function
of degree $d = \deg(\varphi)$, 
and let $\Gamma \subset \PP^1_{\Berk}$ be a subgraph.  Then

$A)$ \ $\varphi(\Gamma)$ is compact and connected for the path length topology.

$B)$ \ $\ell(\varphi(\Gamma)) \le d \cdot \ell(\Gamma)$\ .

$C)$ \ The path length topology on $\varphi(\Gamma)$ coincides
with the induced topology from the Berkovich topology on $\PP^1_{\Berk}$. 
\end{lemma}  

\begin{proof}
Assertion (A) is trivial, since $\varphi$ is continuous for the
path length topology, and $\Gamma$ is compact and connected in the 
path length topology.

\vskip .05 in
For (B), let $\tGamma$ be an arbitrary subgraph of $\varphi(\Gamma)$.   
Choose a set of partition points $X = \{x_1, \ldots, x_m\} \subset \Gamma$ 
such that all the endpoints and branch points of $\Gamma$ are among the $x_i$,
and such that for each endpoint and branch point $p$ of $\tGamma$,
there is an $x_i \in X$ with $\varphi(x_i) = p$.  
Let $A$ be the set of pairs $(i,j)$ such that $x_i$ and $x_j$ are
adjacent in $\Gamma$, and put  
\begin{equation*}
\Gamma^{\prime} \ = \ \bigcup_{(i,j) \in A} [\varphi(x_i),\varphi(x_j)] \ .
\end{equation*}
(Note that we are taking the segments $[\varphi(x_i),\varphi(x_j)]$,
not the path images $\varphi([x_i,x_j])$).  
Then $\Gamma^{\prime}$ is a subgraph of $\varphi(\Gamma)$ which 
contains all the endpoints and branch points of $\tGamma$,
and hence contains $\tGamma$.  
For each $(i,j) \in A$, Lemma \ref{LemG9} gives 
$\rho(\varphi(x_i),\varphi(x_j)) \le d \cdot \rho(x_i,x_j)$.  Hence
\begin{equation*}
\ell(\tGamma) \ \le \ \ell(\Gamma^{\prime}) \ \le \ 
\sum_{(i,j) \in A} \rho(\varphi(x_i),\varphi(x_j)) 
\ = \ d \cdot \ell(\Gamma) \ .
\end{equation*}
Taking the $\sup$ over all subgraphs 
shows $\ell(\varphi(\Gamma)) \le d \cdot \ell(\Gamma)$.   

\vskip .05 in
For (C), we introduce the following notation.  Given a connected subset
$Z \subset \PP^1_{\Berk,0}$ and a point $b \in Z$, put
\begin{equation*}
\radius(b,Z) \ = \ \sup_{x \in Z} \ \rho(b,x) \ .
\end{equation*} 
We will show that for each $b \in \varphi(\Gamma)$, 
and for each $\varepsilon > 0$, there are a finite number of points 
$p_1, \ldots, p_m \in \varphi(\Gamma)$ distinct from $b$, such that 
if $V$ is the connected component of 
$\varphi(\Gamma) \backslash \{p_1,\ldots,p_m\}$ containing $b$, 
then $\radius(b,V) < \varepsilon$.  
To construct the points $p_i$, put $Z_0 = \varphi(\Gamma)$.    
If $\radius(b,Z_0) < \varepsilon$ there is nothing to show.  Otherwise,
let $q_1 \in Z_0$ be a point with $\rho(b,q_1) \ge \varepsilon$.  
Let $p_1$ be the midpoint of $[b,q_1]$ and let $Z_1$ be 
the connected component of $Z_0$ containing $b$.  Then 
$\ell(Z_1) \le \ell(\varphi(\Gamma)) - \varepsilon/2$.  
Inductively suppose $Z_{n}$ has been constructed 
with $\ell(Z_n) \le \ell(\varphi(\Gamma)) - n \cdot \varepsilon/2$.
If $\radius(b,Z_n) < \varepsilon$, put $V = Z_n$.
Otherwise, let $q_n \in Z_n$ be a point with $\rho(b,q_n) \ge \varepsilon$,
and let $p_n$ be the midpoint of $[b,q_n]$.  Let $Z_{n+1}$ be the
connected component of $Z_n \backslash \{p_n\}$ containing $b$. 
Then $\ell(Z_{n+1}) \le \ell(\varphi(\Gamma)) - (n+1) \cdot \varepsilon/2$. 
This process must terminate, 
since $\ell(\varphi(\Gamma)) \le d \cdot \ell(\Gamma)$ is finite by part (B).
Let $U$ be the connected component
of $\PP^1_{\Berk} \backslash \{p_1, \ldots, p_m\}$ containing $b$.
Then $U$ is open in $\PP^1_{\Berk}$, and 
\begin{equation*}
\varphi(\Gamma) \cap U \ = \ V \ \subset \ B_{\rho}(b,\varepsilon) 
               \ = \  \{x \in \varphi(\Gamma) : \rho(b,x) < \varepsilon \} \ .
\end{equation*}

From this we see that each set $V \subset \varphi(\Gamma)$ which is open for
the path-length topology, is open for the relative Berkovich topology.
Conversely, if $V \subset \varphi(\Gamma)$ is open for the relative Berkovich 
topology and $b \in V$, let $W$ be a basic open 
neighborhood of $b$ in $\PP^1_{\Berk}$ for which 
$W \cap \varphi(\Gamma) \subset V$.  Then $W$ is a simple domain;
write $\partial W = \{x_1, \ldots, x_m\}$.  
Put $\varepsilon = \min_j (\rho(b,x_j))$.  Then $B_{\rho}(b,\varepsilon)$
is an open neighbhorhood of $b$ in the path length topology 
on $\varphi(\Gamma)$ which is contained in $V$.  
\end{proof}            
   
\subsection{The Pushforward and Pullback measures.}

Let $\varphi(T) \in \CC_v(T)$ be a nonconstant rational function.

\vskip .05 in
If $\nu$ is a bounded Borel measure on $\PP^1_{\Berk}$, 
the pushforward measure $\varphi_*(\nu)$ is the Borel measure defined by
\begin{equation*}
\varphi_*(\nu)(U) \ = \ \nu(\varphi^{-1}(U)) 
\end{equation*}
for all Borel subsets $U \subset \PP^1_{\Berk}$.  Here $\varphi^{-1}(U)$
is a Borel set because $\varphi$ is continuous for the Berkovich topology.

\vskip .05 in
The pullback measure $\varphi^*(\nu)$ is more complicated to define.  
If $g(x)$ is a continuous function on $\PP^1_{\Berk}$, then $\varphi_*(g)$
is the function given by 
\begin{equation*}
\varphi_*(g)(y) \ = \ \sum_{\varphi(x) = y} m_{\varphi,x}(y) \cdot g(x) \ .
\end{equation*}
We claim that $\varphi_*(g)$ is continuous.  
To see this, fix $p \in \PP^1_{\Berk}$, and take $\varepsilon > 0$.
Let $q_1, \ldots, q_r \in \PP^1_{\Berk}$
be the points with $\varphi(q_i) = p$.  
Choose disjoint neighborhoods $U_1, \ldots, U_r$
with $q_i \in U_i$.  Since $g$ is continuous,
after shrinking $U_i$, we can assume that
\begin{equation*}
|g(x)-g(q_i)| \ < \ \varepsilon 
\end{equation*}
for all $x \in U_i$.  For each $U_i$, let $V_i$ be the neighborhood of $q_i$
given by Corollary \ref{CorG7} such that 
\begin{equation*}
M_{U_i}(x) \ = \ \sum \begin{Sb} t \in U_i \\ \varphi(t) = \varphi(x) \end{Sb}
                 m_{\varphi,\varphi(x)}(t)
\end{equation*}
is constant for $x \in V_i$.  Without loss we can assume $V_i \subset U_i$.  

Put $W = \cap_{i=1}^r \varphi(V_i)$.
By Corollary \ref{CorG6} each $\varphi(V_i)$ is open and contains $p$, 
so $W$ is a neighbhorhood of $V$.  
Fix $y \in W$, and let $y_1, \ldots, y_s$ be the
preimages of $y$ under $\varphi$.  Each $y_j$ belongs to $V_i$ 
for some $i$.  On the other hand, for each $U_i$,
\begin{equation*}
\sum_{y_j \in U_i} m_{\varphi,y}(y_j) \ = \ m_{\varphi,p}(q_i) 
\end{equation*}
and $\sum_j m_{\varphi,y}(y_j) = \sum_i m_{\varphi,p}(q_i) = \deg(\varphi)$.
Therefore 
\begin{eqnarray*}
|\varphi_*(g)(y) - \varphi_*(p)| 
& \le & \sum_{i=1}^r \sum_{y_j \in U_i} 
       m_{\varphi,y}(y_j) \cdot |g(y_j) - g(q_i)| \\
& \le &  (\sum_j m_{\varphi,y}(y_j)) \cdot \varepsilon 
\ = \ d \cdot \varepsilon \ .
\end{eqnarray*}
This shows $\varphi_*(g)$ is continuous.  
The construction also shows $\|\varphi_*(g)\| \le d \cdot \|g\|$, 
where $\|g\| = \sup_{x \in \PP^1_{\Berk}}|g(x)|$.

The pullback measure $\varphi^*(\nu)$ is the measure representing the 
bounded linear functional $\Lambda : \cC(\PP^1_{\Berk}) \rightarrow \RR$
defined by $\Lambda(g) = \int \varphi_*(g)(y) \, d\nu(y)$.  Thus
\begin{equation}
\int g(x) \, d\nu^*(x) \ = \ \int \varphi_*(g)(y) \, d\nu(y) \ . \label{FDR1}
\end{equation}

\vskip .05 in 
These measures satisfy the usual formal functorial properties.

\begin{proposition}  \label{PropG12}
Let $\varphi(T) \in \CC_v(T)$ be a nonconstant rational function
of degree $d = \deg(\varphi)$.  Suppose $\nu$ is a bounded Borel measure
on an open set $U \subset \PP^1_{\Berk}$.

$A)$ \ Let $\nu_1, \nu_2, \cdots$ be a sequence of signed Borel measures
on $U$ for which the masses of $|\nu_i|$ are uniformly bounded, 
and which converge weakly to $\nu$ on compact subsets of $U$.  
Then  $\varphi^*(\nu_1), \varphi^*(\nu_2), \cdots$
converge weakly to $\varphi^*(\nu)$ on compact subsets of $\varphi^{-1}(U)$. 

If $U = \varphi^{-1}(\varphi(U))$, then 
$\varphi_*(\nu_1), \varphi_*(\nu_2), \cdots$ converge weakly to 
$\varphi_*(\nu)$ on compact subsets of $\varphi(U)$.

$B)$ \ $\varphi_*(\varphi^*(\nu)) = d \cdot \nu$ \ .
\end{proposition}

\begin{proof}  The first assertion in (A) holds 
because for each $g \in \cC(\varphi^{-1}(\Ubar))$
\begin{eqnarray*}
\int g(x) \, d\varphi^*(\nu) & = & \int \varphi_*(g)(y) \, d\nu(y) \\
& = & \lim_{n \rightarrow \infty} \int \varphi_*(g)(y) \, d\nu_n(y) 
\ = \ \lim_{n \rightarrow \infty} \int g(x) \, d\varphi^*(\nu_n)  \ .
\end{eqnarray*}
The second holds because for each open subset 
$V \subset \varphi(U)$,
\begin{eqnarray*}
\nu_*(V) & = & \nu(\varphi^{-1}(V)) \\
& = & \lim_{n \rightarrow \infty} \nu_n(\varphi^{-1}(V)) 
\ = \ \lim_{n \rightarrow \infty} \varphi_*(\nu_n)(V) \ .
\end{eqnarray*}
and equality of measures on open sets implies their
equality for all Borel sets.

Part (B) also follows from a simple computation.  Formula (\ref{FDR1})
extends to characteristic functions of Borel sets.  
For each Borel set $E \subset \PP^1_{\Berk}$, the identity 
$\sum_{\varphi(x) = y} m_{\varphi,y}(x) = d$ 
means that $\varphi_*(\chi_{\varphi^{-1}(E)}) = d \cdot \chi_E$.  Hence
\begin{eqnarray*}
\varphi_*(\varphi^*(\nu))(E) & = & \varphi^*(\nu)(\varphi^{-1}(E)) \\
& = & \int \chi_{\varphi^{-1}(E)}(x) \, d\varphi^*(\nu)(x) 
\ = \ \int \varphi_*(\chi_{\varphi^{-1}(E)})(y) \, d\nu(y) \\
& = & \int d \cdot \chi_E(y) \, d\nu(y) \ = \ d \cdot \nu(E) \ .
\end{eqnarray*}
\end{proof}

\subsection{The Pullback Formula for subharmonic functions.}

\begin{proposition} \label{PropG13}
Let $\varphi(T) \in \CC_v(T)$ be a nonconstant rational map of degree $d$.  
Let $U \subset \PP^1_{\Berk}$ be a domain, and let $f$ be subharmonic in $U$.
Then

$A)$ The function $\varphi^*(f) = f \circ \varphi$ 
is subharmonic in $\varphi^{-1}(U)$.

$B)$ $\Delta_{\varphi^{-1}(U)}(\varphi^*(f)) = \varphi^*(\Delta_U(f))$ \ .
\end{proposition} 

\begin{proof}  

(A)  To see that $f \circ \varphi$ is subharmonic, 
we use the Riesz Decomposition Theorem (Proposition~\ref{PropF19}).  
Let $V$ be a simple subdomain of $U$, and put $\nu = -\Delta_V(f)$.  
Let $\zeta$ be an arbitrary point in $\PP^1_{\Berk} \backslash \Vbar$.  
Proposition~\ref{PropF19} shows there is a harmonic function $h$ on $V$
such that for all $x \in V$
\begin{equation*}
f(y) \ = \ h(y) + \int \log_v(\delta(y,z)_{\zeta}) \, d\nu(z) \ .
\end{equation*}
Composing with $\varphi(x)$ gives
\begin{equation*}
f(\varphi(x)) \ = \ h(\varphi(x)) 
         + \int_V \log_v(\delta(\varphi(x),z)_{\zeta}) \, d\nu(z) \ .
\end{equation*}
By Corollary~\ref{CorE23}, $h(\varphi(x))$ is harmonic in $V$.  By Proposition 
\ref{PropG2}, for each $z \in V$ the function 
$\log_v(\delta(\varphi(x),z)_{\zeta})$ is subharmonic
in $\PP^1_{\Berk} \backslash \varphi^{-1}(\{\zeta\})$, hence in particular 
in $\varphi^{-1}(V)$.  Finally, by Proposition~\ref{PropF12}, 
\begin{equation*}
F(x) \ = \ \int_V \log_v(\delta(\varphi(x),z)_{\zeta}) \, d\nu(z)
\end{equation*}
is subharmonic on $\varphi^{-1}(V)$.  (The majorizing function required
in Proposition~\ref{PropF12} 
can be taken to a constant, since $\delta(y,z)_{\zeta}$
is bounded from above on $\Vbar \times \Vbar$, and $\nu$ has finite mass.
The function $F(x)$ is not $\equiv -\infty$ on any component of 
$\varphi^{-1}(V)$, since $F(x) = f(\varphi(x))-h(\varphi(x))$.)  
It follows that $f(\varphi(x)) = h(\varphi(x)) + F(x)$ is subharmonic on 
$\varphi^{-1}(V)$.  Exhausting $U$ by simple subdomains $V$, we see that
$f(\varphi(x))$ is subharmonic on $\varphi^{-1}(U)$. 

\vskip .05 in
(B)  Let $V$ be a simple subdomain of $U$.  By the smoothing theorem 
(Proposition~\ref{PropF26}), 
there is a decreasing sequence of subharmonic functions
$f_1 \ge f_2 \ge \cdots$ on $V$ such that 
$\lim_{n \rightarrow \infty} f_n(z) = f(z)$ pointwise for all $z \in V$,
and such that each $\Delta_V(f_n)$ is a discrete measure supported on a finite 
number of points.  We will first prove the pullback formula for the $f_n$.  

Fix $\zeta \in \PP^1_{\Berk} \backslash \Vbar$, and write 
\begin{equation*}
\nu_n = -\Delta_V(f_n) \ = \ \sum_{i=1}^{M_n} c_{ni} \delta_{p_{ni}}(x)  
\end{equation*}
As in the proof of part (A), there is a harmonic function $h_n$ on $V$
such that for all $u \in V$ 
\begin{eqnarray*} 
f_n(y) & = & h_n(y) + \int \log_v(\delta(y,z)_{\zeta}) \, d\nu_n(z) \\
    & = & h_n(y) + \sum_{i=1}^{M_n} c_{ni} \log_v(\delta(y,p_{ni})_{\zeta}) \ .
\end{eqnarray*}
Composing with $\varphi(x)$ gives
\begin{equation*}
f_n(\varphi(x)) \ = \  h_n(\varphi(x)) 
    + \sum_{i=1}^{M_n} c_{ni} \log_v(\delta(\varphi(x),p_{ni})_{\zeta}) \ .
\end{equation*} 
As before, $h_n(\varphi(x))$ is harmonic in $\varphi^{-1}(V)$.  
By Proposition \ref{PropG2} each function 
$h_{ni}(x) = \log_v(\delta(\varphi(x),p_{ni})_{\zeta})$ 
belongs to $\BDV(\PP^1_{\Berk})$ and 
\begin{equation*}
\Delta_{\PP^1_{\Berk}}(h_{ni}) \ = \ 
\sum_{\varphi(\xi) = \zeta} m_{\varphi,\zeta}(\xi) \cdot \delta_{\xi}(x) 
- \sum_{\varphi(a) = p_{ni}} m_{\varphi,p_{ni}}(a) \cdot \delta_{a}(x)
\end{equation*}
In particular, $h_{ni}(x)$ is subharmonic in $\varphi^{-1}(V)$.  
We also note that for any $p$, 
the measure $\sum_{\varphi(a) = p} m_{\varphi,p}(a) \delta_a(x)$ is precisely
the pullback measure $\varphi^*(\delta_{p}(x))$.  Hence 
\begin{equation*}
\Delta_{\varphi^{-1}(V)}(h_{ni}) \ = \ -\varphi^*(\delta_{p_{ni}}(x)) \ .
\end{equation*}
Summing over all $i$, we find that
\begin{equation*}
\Delta_{\varphi^{-1}(V)}(f_n \circ \varphi)   
  \ = \ - \sum_{i = 1}^{M_n}c_{ni} \varphi^*(\delta_{p_{ni}}(x))  
  \ = \ \varphi^*(\Delta_V(f_n)) \ .
\end{equation*}

The functions $f_n \circ \varphi(x)$ decrease monontonically
to $f \circ \varphi(x)$.  By the proof of 
Proposition~\ref{PropF23}(B), the measures 
$\Delta_{\varphi^{-1}(V)}(f_n \circ \varphi)$ converge weakly to 
$\Delta_{\varphi^{-1}(V)}(f \circ \varphi)$.  Proposition~\ref{PropF23}(B)
also shows the measures $\Delta_V(f_n)$ converge weakly to $\Delta_V(f)$.  
Hence, by Proposition \ref{PropG12} the measures $\varphi^*(\Delta_V(f_n))$
converge weakly to $\varphi^*(\Delta_V(f))$.  Thus 
\begin{equation*}
\Delta_{\varphi^{-1}(V)}(f \circ \varphi) \ = \ \varphi^*(\Delta_V(f)) \ .
\end{equation*}
Exhausting $U$ by a sequence of simple subdomains $V_n$, we obtain the
result.
\end{proof} 

\subsection{Construction of the Lyubich measure.}

Let $\varphi(T) \in \CC_v(T)$ be a rational function of degree $d \ge 2$.
Thus, there are coprime polynomials 
$f_1(T), f_2(T) \in \CC_v[T]$ with $\max(\deg(f_1),\deg(f_2)) = d$
such that $\varphi(T) = f_2(T)/f_1(T)$.  

We are interested in the dynamics of $\varphi$ on $\PP^1_{\Berk}$.
Put $\varphi^{(1)}(T) = \varphi(T)$, and inductively define
\begin{equation*}
\varphi^{(n)}(T) \ = \ \varphi(\varphi^{(n-1)}(T)) 
\ = \ \varphi^{(n-1)}(\varphi(T)) \ .
\end{equation*}
Homogenizing $f_1(T)$ and $f_2(T)$, we obtain coprime homogeneous polynomials
$F_1(X,Y), F_2(X,Y) \in \CC_v[X,Y]$ of degree $d$ 
such that $f_1(T) = F_1(1,T)$, $f_2(T) = F_2(1,T)$.  Write
\begin{equation*}
F_1(X,Y) = \sum_{i+j = d} c_{1,ij} X^i Y^j\ , \quad
F_2(X,Y) = \sum_{i+j = d} c_{2,ij} X^i Y^j \ .
\end{equation*}

For the moment, regard $F_1(X,Y)$ and $F_2(X,Y)$ as functions on $\CC_v^2$
and write $F(x,y) = (F_1(x,y),F_2(x,y))$.  Let 
\begin{equation*} 
\|(x,y)\|_v \ = \ \max(|x|_v,|y|_v)
\end{equation*}
be the usual $\sup$ norm on $\CC_v^2$.
We claim there are numbers $0 < B_1 < B_2 < \infty$ such that for all
$(x,y) \in \CC_v^2$,
\begin{equation} \label{FZR1}
B_1 \cdot \|(x,y)\|_v^d \ \le \ \|F(x,y)\|_v
                       \ \le \ B_2 \cdot \|(x,y)\|_v^d \ .
\end{equation}
The upper bound is trivial, since 
if $B_2 = \max(|c_{1,ij}|_v,|c_{2,ij}|_v)$ then by the ultrametric
inequality 
\begin{equation*}
|F_1(x,y)|_v, \ |F_2(x,y)|_v \ \le \ \max_{i,j}(|c_{1,ij}|_v |x|_v^i |y|_v^j)
\ \le \ B_2 \|(x,y)\|_v^d \ .
\end{equation*}

For the lower bound,
note that since $F_1(1,T)$ and $F_2(1,T)$ are coprime,
their resultant $b_1 = \Res(F_1(1,T),F_2(1,T)) \in \CC_v$ is nonzero.
By properties of the resultant,
there are polynomials $g_1(T), g_2(T) \in \CC_v[T]$ of degree at most $d-1$
such that
\begin{equation*}
g_1(T) F_1(1,T) + g_2(T) F_2(1,T) \ = \ b_1 \ .
\end{equation*}
Homogenizing, we obtain homogenous polynomials $G_1(X,Y), G_2(X,Y)$ of
degree $d-1$ for which
\begin{equation} \label{FLE1}
G_1(X,Y) F_1(X,Y) + G_2(X,Y) F_2(X,Y) \ = \ b_1 X^{2d-1}
\end{equation}
Similarly, since $F_1(U,1)$ and $F_2(U,1)$ are coprime, their
resultant $b_2 = \Res(F_1(U,1),F_2(U,1))$ is nonzero, and there are
homogeneous polynomials $H_1(X,Y), H_2(X,Y)$ of degree $d-1$ such that
\begin{equation} \label{FLE2}
H_1(X,Y) F_1(X,Y) + H_2(X,Y) F_2(X,Y) \ = \ b_2 Y^{2d-1} \ .
\end{equation}
By the upper bound argument applied to $G = (G_1,G_2)$
and $H = (H_1,H_2)$, there is an $A_2 > 0$ 
such that $\|G(x,y)\|_v, \|H(x,y)|_v \le A_2 \|(x,y)\|_v^{d-1}$
for all $(x,y) \in \CC_v^2$.  By (\ref{FLE1}), (\ref{FLE2})
and the ultrametric inequality, 
\begin{eqnarray*}
|b_1|_v|x|_v^{2d-1} & \le & A_2 \|(x,y)\|_v^{d-1} \cdot \|F(x,y)\|_v \ , \\
|b_2|_v|y|_v^{2d-1} & \le & A_2 \|(x,y)\|_v^{d-1} \cdot \|F(x,y)\|_v \ .
\end{eqnarray*}
Writing $A_1 = \min(|b_1|_v,|b_2|_v)$, it follows that  
\begin{equation*}
A_1 \max(\|x,y\|)^{2d-1} \ \le \ A_2 \|(x,y)\|_v^{d-1} \cdot \|F(x,y)\|_v \ .
\end{equation*}  
Thus if $B_1 = A_1/A_2$, then $\|F(x,y)\|_v \ge B_1 \|(x,y)\|_v^d$.

Taking logarithms, and putting $C_1 = \max(\log_v(B_1), \log_v(B_2))$, we have 
\begin{equation} \label{FLE3}
|\frac{1}{d} \log_v(\|F(x,y)\|_v) - \log_v(\|(x,y)\|_v)|
 \ \le \ \frac{C_1}{d} \ .
\end{equation}
for all $\vec{0} \ne (x,y) \in \CC_v^2$.
Now iterate $F(X,Y)$:  put $F^{(1)}(X,Y) = F(X,Y)$, and inductively define 
\begin{equation*}
F^{(n)}(X,Y) \ = \ F(F^{(n-1)}(X,Y)) \ = \ F^{(n-1)}(F(X,Y)) \ .
\end{equation*}
Then $F^{(n)}(X,Y) = (F_1^{(n)}(X,Y), F_2^{(n)}(X,Y))$
where $F_1^{(n)}(X,Y), F_2^{(n)}(X,Y) \in \CC_v[X,Y]$ are homogeneous
polynomials of degree $d^n$.  Inserting $F(x,y)$ for $(x,y)$ in (\ref{FLE3})
and iterating, we find that for each $n$
\begin{equation} \label{FLE4}
|\frac{1}{d^n} \log_v(\|F^{(n)}(x,y)\|_v)
 - \frac{1}{d^{n-1}} \log_v(\|F^{(n-1)}(x,y)\|_v)|
       \ \le \ \frac{C_1}{d^n} \ .
\end{equation}
Put
\begin{equation*}
h^{(n)}_{\varphi,v}(x) \ = \
   \frac{1}{d^n} \log_v(\max(|F_1^{(n)}(1,x)|_v,|F_2^{(n)}(1,x)|_v)) \ .
\end{equation*}
The Call-Silverman local height for $\varphi$ on $\PP^1(\CC_v)$  
(relative to the point $\infty$ and the dehomogenization 
$F_1(1,T), F_2(1,T)$), is defined for $x \in \AA^1(\CC_v)$ by
\begin{equation} \label{FLE5}
h_{\varphi,v}(x) \ = \ \lim_{n \rightarrow \infty} h^{(n)}_{\varphi,v}(x) \ .
\end{equation}
The fact that the limit exists follows from (\ref{FLE4}), 
as does the bound
\begin{equation} \label{FLE6}
|h_{\varphi,v}(x) - \log_v(\max(1,|x|_v))|
\ \le \ \sum_{n=1}^{\infty} \frac{C_1}{d^n} 
            \ = \ C \ .
\end{equation}
Taking logarithms in the identity
\begin{equation*} 
F^{(n-1)}(1,\varphi(x)) \ = \ F^{(n-1)}(1,\frac{F_2(1,x)}{F_1(1,x)})
\ = \ \frac{F^{(n)}(1,x)}{F_1(1,x)^{d^{n-1}}} 
\end{equation*}
and letting $n \rightarrow \infty$ gives the functional equation
\begin{equation} \label{FLE9}
h_{\varphi,v}(\varphi(x))
  \ = \ d \cdot h_{\varphi,v}(x) - \log_v(|F_1(1,x)|_v) \ , 
\end{equation}
valid on $\PP^1(\CC_v) \backslash (\{\infty\} \cup \varphi^{-1}(\{\infty\}))$.  
The two properties (\ref{FLE6}), (\ref{FLE9}) characterize the Call-Silverman
local height.

\vskip .1 in
We will now ``Berkovichize'' the local height.  For each $n$, put
\begin{equation} \label{FLE10}
\hat{h}^{(n)}_{\varphi,v}(x) \ = \
 \frac{1}{d^n} \max(\log_v([F_1^{(n)}(1,T)]_x),
                    \log_v([F_2^{(n)}(1,T)]_x)))
\end{equation}
for $x \in \PP^1_{\Berk}$ 
(we regard it as having the value $\infty$ at $x = \infty$).  
Then $\hat{h}^{(n)}_{\varphi,v}(x)$
coincides with $h^{(n)}_{\varphi,v}(x)$ on $\AA^1(\CC_v)$ and is
continuous and strongly subharmonic in $\AA^1_{\Berk}$.
We claim that for all $x \in \AA^1_{\Berk}$
\begin{equation} \label{FLE11}
|\hat{h}^{(n)}_{\varphi,v}(x) - \hat{h}^{(n-1)}_{\varphi,v}(x)|
       \ \le \ \frac{C_1}{d^n} \ .
\end{equation}
Indeed, this holds for all type I points in $\AA^1(\CC_v)$;  such points
are dense in $\AA^1_{\Berk}$ and the functions involved are continuous,
so it holds for all $x \in \AA^1_{\Berk}$.  

It follows that the functions $\hat{h}^{(n)}_{\varphi,v}(x)$ 
converge uniformly to a continuous subharmonic function
$\hat{h}_{\varphi,v}(x)$ on $\AA^1_{\Berk}$ 
which extends the Call-Silverman local height $h_{\varphi,v}(x)$.
By the same arguments as before, for all $x \in \AA^1_{\Berk}$ 
\begin{equation} 
|\hat{h}_{\varphi,v}(x) - \log_v(\max(1,[T]_x))| \ \le \ C \ , 
     \label{FLE13} 
\end{equation}
and for all 
$x \in \PP^1_{\Berk} \backslash (\infty \cup \varphi^{-1}(\{\infty\}))$,
\begin{equation}     
\hat{h}_{\varphi,v}(\varphi(x))
 \ = \ d \cdot \hat{h}_{\varphi,v}(x) - \log_v([F_1(1,T)]_x) \ . 
     \label{FLE14}
\end{equation}
Actually, this can be viewed as an identity for all 
$x \in \PP^1_{\Berk}$,
if one views the right side as given by its limit as $x \rightarrow \infty$,
for $x = \infty$.  

Let $V_1 = \AA^1_{\Berk} = \PP^1_{\Berk} \backslash \{\infty\}$.
Since $\hat{h}_{\varphi,v}$ is subharmonic on $V_1$,
there is a non-negative measure $\mu_1$ on $V_1$ such that
$\Delta_{V_1}(\hat{h}_{\varphi,v}) = -\mu_1$.
On the other hand, each $\hat{h}^{(n)}_{\varphi,v}(x)$ 
belongs to $\BDV(\PP^1_{\Berk})$ and satisfies
\begin{equation*}
\Delta_{\PP^1_{\Berk}}(\hat{h}^{(n)}_{\varphi,v}) \ = \
\delta_{\infty} - \mu_1^{(n)},
\end{equation*}
where $\mu_1^{(n)} \ge 0$ has total mass $1$.  
Since the $\hat{h}^{(n)}_{\varphi,v}$ converge uniformly
to $\hat{h}_{\varphi,v}$ on $V_1$, the $\mu_1^{(n)}$ converge weakly
to $\mu_1$ on simple subdomains of $V_1$, and so $\mu_1(V_1) \le 1$.
It follows that $\hat{h}_{\varphi,v}$ belongs to $\BDV(\PP^1_{\Berk})$, 
and there is a non-negative measure $\mu$ on $\PP^1_{\Berk}$ 
of total mass $1$ such that 
\begin{equation} \label{FGM1B}
\Delta_{\PP^1_{\Berk}}(\hat{h}_{\varphi,v}) \ = \ \delta_{\infty}(x) - \mu \ .
\end{equation}
(Note that functions in $\BDV(\PP^1_{\Berk})$ are allowed take values
$\pm \infty$ on points of $\PP^1(\CC_v)$;  the definition of the Laplacian
only involves their restriction to subgraphs $\Gamma$, which are necessarily 
contained in $\PP^1_{\Berk} \backslash \PP^1(\CC_v)$.)   

In the affine patch $V_2 := \PP^1_{\Berk} \backslash \{0\}$,
relative to the coordinate function $U = 1/T$, 
the map $\varphi$ is given by $F_1(U,1)/F_2(U,1)$.
By a construction similar to the one above,
using the functions $F_1(U,1)$ and $F_2(U,1)$,
we obtain a function $\hat{g}_{\varphi,v}(x) \in \BDV(\PP^1_{\Berk})$
which is continuous and subharmonic on $V_2$  
and extends the Call-Silverman local height relative to the point $0$.
For all $x \in V_2$, it satisfies
\begin{eqnarray}  
|\hat{g}_{\varphi,v}(x) - \log_v(\max(1,[1/T]_x))| \ \le \ C \ , \qquad
\label{FLE15} \\
\hat{g}_{\varphi,v}(\varphi(x))
\ = \ d \cdot \hat{g}_{\varphi,v}(x) - \log_v([F_2(1/T,1)]_x) \ ,
             \label{FLE16}
\end{eqnarray}
where (\ref{FLE16}) holds in the same sense as (\ref{FLE14})
Using the identity $F^{(n)}(U,1) = F^{(n)}(1,T)/T^{d^n}$,
taking logarithms, and letting $n \rightarrow \infty$ gives 
\begin{equation} \label{FLE17}
\hat{g}_{\varphi,v}(x) \ = \ \hat{h}_{\varphi,v}(x) - \log_v([T]_x) \ .
\end{equation}
Applying the Laplacian and using (\ref{FGM1B}) shows that  
\begin{equation} \label{FGM2}
\Delta_{\PP^1_{\Berk}}(\hat{g}_{\varphi,v}) \ = \ \delta_{0}(x) - \mu \ .
\end{equation}

\vskip .1 in
The construction of the measure $\mu$ has been the main goal of this
section.  It will be called the Lyubich measure.

\begin{theorem} \label{ThmG14}
Let $\varphi(T) \in \CC_v(T)$ have degree $d \ge 2$.  
The Lyubich measure $\mu = \mu_{\varphi}$ on $\PP^1_{\Berk}$ is non-negative 
and has total mass $1$. It satisfies the functional equations  
$\varphi^*(\mu) = d \cdot \mu$ and $\varphi_*(\mu) = \mu$.
\end{theorem}

\begin{proof}
To show $\varphi^*(\mu) = d \cdot \mu$, 
we use (\ref{FLE14}) and (\ref{FLE16}).
Let $x_1, \ldots, x_d$ be the zeros (not necessarily distinct)
of $F_1(1,T)$, and put 
$U_1 = \varphi^{-1}(V_1) = \PP^1_{\Berk} \backslash \{x_1, \ldots, x_d\}$. 
By the pullback formula for subharmonic functions (Proposition \ref{PropG13}), 
$H(x) := \hat{h}_{\varphi,v}(\varphi(x))$ is subharmonic on $U_1$, and 
\begin{equation*}
\Delta_{U_1}(H) \ = \ \varphi^*(\Delta_{V_1}(\hat{h}_{\varphi,v})) \ .
\end{equation*}
Taking Laplacians in (\ref{FLE14}) gives 
$\varphi^*(\mu|_{V_1}) \ = \ d \cdot \mu|_{U_1}$.
Likewise, put
$U_2 = \varphi^{-1}(V_2) = \PP^1_{\Berk} \backslash \{y_1, \ldots, y_d\}$,
where $y_1, \ldots, y_d$ are the zeros of $F_2(1/T,1) = F_2(1,T)/T^d$.
Since $F_1(1,T)$ and $F_2(1,T)$ are coprime, the sets
$\{x_1, \dots, x_d\}$ and $\{y_1, \ldots, y_d\}$ are disjoint.
Then $G(x) := \hat{g}_{\varphi,v}(\varphi(x))$ is subharmonic on $U_2$,
and satisfies
\begin{equation*}
\Delta_{U_2}(G) \ = \ \varphi^*(\Delta_{V_2}(\hat{g}_{\varphi,v})) \ .
\end{equation*}
Taking Laplacians in (\ref{FLE16}) gives 
$\varphi^*(\mu|_{V_2}) = d \cdot \mu|_{U_2}$.
Since $V_1 \cup V_2 = U_1 \cup U_2 = \PP^1_{\Berk}$, it follows that
$\varphi^*(\mu) = d \cdot \mu$.

\vskip .05 in
The identity $\varphi_*(\mu) = \mu$ follows formally 
from $\varphi^*(\mu) = d \cdot \mu$.  By Proposition \ref{PropG12}(B), 
$\varphi_*(\varphi^*(\mu)) = d \cdot \mu$.
Since $\varphi^*(\mu) = d \cdot \mu$, this gives $\varphi_*(\mu) = \mu$.
\end{proof}

\vskip .1 in
The Lyubich measure is better behaved than arbitrary measures;  
in particular, it is ``log-continuous'' in the following sense:

\begin{definition} \label{DefG1}
\noindent{A} positive measure $\omega$ on $\PP^1_{\Berk}$
is {\rm log-continuous} if
\begin{equation*}
u_{\omega}(x)
\ = \ \ \int_{\PP^1_{\Berk}} -\log_v(\|x,y\|_v) \, d\omega(y)
\end{equation*} 
is continuous $($hence uniformly bounded$)$ on all of $\PP^1_{\Berk}$.
\end{definition}

In the definition,
$-\log_v(\|x,y\|_v) = -\log_v(\delta(x,y)_{\zeta_0}) = j_{\zeta_0}(x,y)$
where $\zeta_0 \in \PP^1_{\Berk}$ is the Gauss point,
the point corresponding to $B(0,1)$.
Thus $u_{\omega}(x)$ is the potential function $u_{\omega}(x,\zeta_0)$.
In the definition, $\|x,y\|_v$ could be replaced by any other kernel
which differs from it by a bounded continuous function;   
in particular, it could be replaced by $\delta(x,y)_{\zeta}$
for any $\zeta \in \PP^1_{\Berk} \backslash \PP^1(\CC_v)$.

The log-continuity of the Lyubich measure follows from

\begin{proposition} \label{PropG15}
Let $\omega$ be a positive measure on $\PP^1_{Berk}$
for which $-\omega$
is locally the Laplacian of a continous subharmonic function.
Then $\omega$ is log-continuous.
\end{proposition}

\begin{proof}   
Fix $x \in \PP^1_{\Berk}$, and let $V \subset \PP^1_{\Berk}$ be a 
neighborhood on which $-\omega|_V = \Delta_{V}(f)$ for some continuous 
subharmonic function $f$.  
After shrinking $V$ if necessary, we can that assume $V$
is a simple domain, and that $f$ is bounded and strongly subharmonic on $V$.
Put $\omega_V = \omega|_V$, and consider the potential function  
\begin{equation*}
u_{\omega_V}(z,\zeta)
\ = \ \int -\log_v(\delta(x,y)_{\zeta}) \, d\omega_V(y) \ .
\end{equation*}
By Proposition~\ref{PropF20}, 
$u_{\omega_V}(z,\zeta_0)$ is continuous on all of 
$\PP^1_{\Berk}$.  It is bounded since $\zeta_0 \notin \PP^1(\CC_v)$.

Since $\PP^1_{\Berk}$ is compact, a finite number of such simple
domains $V_i$ cover $\PP^1_{\Berk}$.  It is easy to see that
the assertions made for the restriction of $\omega$ to the $V_i$
hold also for its restriction to their intersections
$V_{i_1} \cap \cdots \cap V_{i_r}$.  By the inclusion-exclusion formula,
$u_{\omega}(z,\zeta_0)$ is continuous on $\PP^1_{\Berk}$.
\end{proof}

\begin{example}
A rational function $\varphi(T) \in \CC_v(T)$
has {\it good reduction} if it can be written as $\varphi(T) = f_2(T)/f_1(T)$
where $f_1(T), f_2(T) \in \CC_v[T]$ are coprime polynomials whose
coefficients belong to the ring of integers $\hat{\cO}_v$ of $\CC_v$,
and where the resultant $\Res(f_1,f_2)$ is a unit in $\hat{\cO}_v$.
(If $\varphi(T)$ is defined over a global field, then it has good
reduction at all but finitely many $v$.)

Suppose $\varphi(T)$ has good reduction, and has degree $d \ge 2$.
We claim that the Lyubich measure $\mu_{\varphi}$ coincides with the
Dirac measure $\delta_{\zeta_0}(x)$, where $\zeta_0$ is the point
of $\PP^1_{\Berk}$ corresponding to $B(0,1)$.

To see this, recall formula (\ref{FZR1})
and note that the proof of that formula shows that under our
hypotheses, $B_1 = B_2 = 1$.  Thus for all $(x,y) \in \CC_v^2$,
\begin{equation*}
\|F(x,y)\|_v \ = \ \|(x,y)\|_v^d \ .
\end{equation*}
By iteration, for each $n$, $\|F^{(n)}(x,y)\| = \|(x,y)\|_v^{d^n}$.
Examining the construction of $\hat{h}_{\varphi,v}$, one finds that
\begin{equation*}
\hat{h}_{\varphi,v}(x) \ = \ 
\log_v(\max(1,[T]_x)) \ = \ \log_v(\delta(x,\zeta_0)_{\infty}) \ , 
\end{equation*}
so $\Delta(\hat{h}_{\varphi,v}) = \delta_{\infty}(x) - \delta_{\zeta_0}(x)$.
Hence $\mu_{\varphi} = \delta_{\zeta_0}(x)$.
\end{example}

\vskip .1 in


\subsection{Fatou and Julia sets.}

In this subsection, we define the Fatou and Julia sets of a rational map $\varphi$ defined over
$\CC_v$; these are open and closed subsets, respectively, of $\PP^1_{\Berk}$.
Assuming that the degree of $\varphi$ is at least 2, we will show that the support of
the Lyubich measure is contained in the Julia set of $\varphi$.  In
particular, it will follow that the Julia
set in $\PP^1_{\Berk}$ of a rational map $\varphi$ of degree at
least two is always non-empty.

In order to define the Fatou and Julia sets, we need to discuss the
general notion of topological equicontinuity.
The following definition is taken from \cite[Section 14.2]{Royden}.

\begin{definition}
Let $X$ and $Y$ be topological spaces, 
let $\cF$ be a family of continuous maps from $X$ to $Y$,
and let $x \in X$ and $y \in Y$ be arbitrary points.
\begin{itemize}
\item[(a)] 
$\cF$ is {\it topologically equicontinuous at $(x,y)$} if, given any open neighborhood $O$ of
$y$, there exist open neighborhoods $U$ of $x$ and $V$ of $y$ such
that for every $f \in \cF$, we have $f(U) \subseteq O$ whenever $f(U)
\cap V$ is nonempty.
\item[(b)] $\cF$ is {\it topologically
equicontinuous at $x$} if it is 
topologically equicontinuous at $(x,y)$ for all $y \in Y$.
\item[(c)] $\cF$ is {\it topologically
equicontinuous} if it is topologically equicontinuous at all $x \in
X$.
\end{itemize}
\end{definition}

We will need the following lemma from \cite[Proposition~14.2.1]{Royden}:

\begin{lemma}
\label{LemRoydenLemma}
Let $\cF$ be a family of continuous maps from $X$ to $Y$ which is
topologically equicontinuous at $x \in X$.  Then given $K \subseteq O
\subseteq Y$ with $K$ compact and $O$ open, there is a neighborhood
$U$ of $x$ such that $f(U) \subseteq O$ whenever $f \in \cF$ and $f(U)
\cap K$ is nonempty.
\end{lemma}


\begin{definition}
Suppose that $X$ is a topological space and that $Y$ is a metric
space.
Let $\cF$ be a family of continuous maps from $X$ to $Y$, and let $x$
be a point of $X$.
We say that $\cF$ is {\it equicontinuous} at $x$ if for every
$\varepsilon > 0$, there exists an open neighborhood $U$ of $X$ such
that $d(f(x),f(x')) < \varepsilon$ for all $x' \in U$ and all $f \in
\cF$.
We say that $\cF$ is {\it equicontinuous} if it is equicontinuous at
every $x \in X$.
\end{definition}

\vskip .1 in

Topological equicontinuity is then related to equicontinuity as
follows (see \cite[Exercise 14.2.3]{Royden}:

\begin{proposition}
Suppose that $X$ is a topological space and that $Y$ is a metric
space.
Let $\cF$ be a family of continuous maps from $X$ to $Y$, and let $x$
be a point of $X$.  Then:

A) If $\cF$ is equicontinuous at $x$, then $\cF$ is
  topologically equicontinuous at $x$. 

B) If $Y$ is compact and $\cF$ is topologically equicontinuous at $x$, then $\cF$ is
  equicontinuous at $x$. 
\end{proposition}

\begin{proof}
For (A), suppose that $\cF$ is equicontinuous at $x$,
let $y \in Y$, and let $O$ be an open neighborhood of $y$ in $Y$.  Choose
$\varepsilon > 0$ such that $V := D(y,\varepsilon) \subseteq O$.  
By definition of equicontinuity, there exists an open neighborhood $U$ of $x$ such
that $d(f(x),f(u)) < \varepsilon/3$ for all $u \in U$ and all $f \in
\cF$.  Let $f \in \cF$, and suppose that $f(U) \cap V$ is nonempty, so
that there exists some element $x' \in U$ such that $f(x') \in V$.  
Then for all $u \in U$, we have 
\begin{equation*}
\begin{array}{lll}
d(f(u),y) & \leq & d(f(u),f(x)) + d(f(x),f(x')) + d(f(x'),y) \\
          & < & \varepsilon/3 + \varepsilon/3 + \varepsilon/3 \\
          & = & \varepsilon,
\end{array}
\end{equation*}
so that $f(u) \in O$.

For (B), suppose $Y$ is compact and $\cF$ is topologically
equicontinuous at $x$.  
Let $\varepsilon > 0$ be arbitrary. 
Then for each $y \in Y$,
there exist open neighborhoods $U_y$ and
$V_y$ of $x$ and $y$, respectively, such that $f(U_y) \subseteq D(y,\varepsilon/2)$
whenever $f(U_y) \cap V_y \neq \emptyset$.  

Since $Y$ is compact, there is a finite subcover 
$V_1,\ldots,V_n := V_{y_1},\ldots, V_{y_n}$ which covers $Y$.  Let $U_1\ldots,U_n$ be the
corresponding sets in $X$, and let $U = \cap U_j$, so that $U$ is an
open neighborhood of $x$.  
 
We claim that $d(f(x),f(x')) < \varepsilon$ whenever $x' \in U$ and $f
\in \cF$.  To see this, note that since $x' \in U$, we have $x'
\in U_j$ for all $j$.  Choosing an index $k$ such that 
$f(x') \in V_k$, we have $f(U_k) \cap V_k \neq \emptyset$, and therefore
$f(x') \in D(y_k,\varepsilon/2)$.  It follows that
\begin{equation*}
\begin{aligned}
d(f(x),f(x')) &\leq d(f(x),f(y_k)) + d(f(y_k),f(x')) \\
          &<  \varepsilon/2 + \varepsilon/2 \\
          &= \varepsilon \\
\end{aligned}
\end{equation*}
as desired.
\end{proof}

\vskip .1 in

In particular, if $X$ is a topological space and $Y$ is a
compact metric space, then the notions of equicontinuity and
topological equicontinuity coincide.  It follows that the equicontinuity of
a family $\cF$ of continuous maps from $X$ to $Y$ depends only on the
underlying topology on $Y$, not on the choice of a particular metric
on $Y$.

\begin{remark}
It is a consequence of the Arzela-Ascoli theorem that
if $X$ and $Y$ are metric spaces, with $Y$ compact,
then the notions of equicontinuity and normality coincide.  (For the
definition of normality, see \cite{Milnor}). 
\end{remark}

We now specialize to the case where $X=Y=\PP^1_{\Berk}$.
Let $\cF : U \to \PP^1_\Berk$ be a family of continuous maps, where $U \subseteq
\PP^1_\Berk$ is an open set.  We say that $\cF$ is equicontinuous on $U$ iff it
is topologically equicontinuous in the above sense.

\begin{remark}
\label{TwoTopologiesRemark}
Using the fact that the subspace topology on $\PP^1(\CC_v)$ coincides with the spherical topology (i.e.,
the topology defined by the spherical metric), it follows that
equicontinuity in the Berkovich topology at a point $x \in \PP^1(\CC_v)$ is equivalent to equicontinuity at $x$
with respect to the spherical metric on $\PP^1(\CC_v)$.  
However, equicontinuity at a point of $\PP^1_\Berk \backslash
\PP^1(\CC_v)$ does not seem easy to describe in terms of the spherical
metric on $\PP^1(\CC_v)$.
\end{remark}

\vskip .1 in

We now define the Fatou and Julia sets of a rational function
$\varphi(T) \in \CC_v(T)$.

\begin{definition}
The {\em Fatou set} of $\varphi$ is the open set $F_\varphi$ consisting of all $x
\in \PP^1_{\Berk}$ such that the family $\cF_\varphi = \{ \varphi^{(n)} \}_{n\geq 1}$ of
iterates of $\varphi$ is equicontinuous at $x$.  
The {\em Julia set} $J_\varphi$ of $\varphi$ is the complement of the
Fatou set, i.e., $J_\varphi := \PP^1_\Berk \backslash F_\varphi$.
\end{definition}

By definition, the Fatou set of $\varphi$ is open, and the Julia set
of $\varphi$ is compact.

\vskip .1 in

The sets $J_\varphi(\CC_v) = J_\varphi \cap \PP^1(\CC_v)$ and
$F_\varphi(\CC_v) = F_\varphi \cap \PP^1(\CC_v)$ coincide, respectively, with the Julia sets and
Fatou sets of $\varphi$ in $\CC_v$ as defined in \cite{Benedetto},
\cite{Hsia2}, and \cite{MS}.  This follows from 
Remark~\ref{TwoTopologiesRemark}.


\vskip .1 in

We now show that if $\deg(\varphi) \geq 2$, then $J_\varphi$ is always non-empty.  This is in contrast
to $J_\varphi(\CC_v)$, which is empty whenever the map $\varphi$ has good
reduction (see \cite{MS}, \cite{Benedetto}).  
We will deduce the fact that $J_\varphi$ is non-empty from the 
stronger fact that the Lyubich measure $\mu_\varphi$ is supported on $J_\varphi$.
It seems likely that the support of $\mu_\varphi$ is in fact equal to
$J_\varphi$.

\begin{theorem}
\label{ThmLyubichSupport}
The support of $\mu_\varphi$ is contained in the Julia set $J_\varphi$.
\end{theorem}

\begin{proof}
Let $x_0$ be a point of the Fatou set.  We will show that if $U$ is a
sufficiently small open neighborhood of $x_0$, then 
$\mu_\varphi$ restricted to $U$ is the zero measure.

First assume that $x_0 \in \AA^1_\Berk$.
Then it suffices to show that 
$\hat{h}_{\varphi,v}(x)$ is harmonic on some neighborhood $U \subset
\AA^1_{\Berk}$ of $x_0$.
To do this, we follow the general outline of an argument of Fornaess and Sibony \cite{FS}.
Define 
\begin{equation*}
\begin{aligned}
O = \{ x \in \AA^1_\Berk \; : \; [T]_x < 2\} \\
K = \{ x \in \AA^1_\Berk \; : \; [T]_x \leq \frac{1}{2} \}. \\
\end{aligned}
\end{equation*}

By Lemma~\ref{LemRoydenLemma}, there exists a neighborhood $U$ of $x_0$
such that $\varphi^{(n)}(U) \subseteq O$ whenever $\varphi^{(n)}(U)
\cap K \neq\emptyset$.  Therefore for all $n$ we have either
$\varphi^{(n)}(U) \subseteq O$ or 
$\varphi^{(n)}(U) \subseteq \PP^1_\Berk \backslash K.$
Passing to a subsequence $\varphi^{(n_k)}$ and replacing $\varphi$ by
$1/\varphi$ if necessary, we may therefore assume without loss of
generality that 
$\varphi^{(n_k)}(U) \subseteq O$ for all $k$, or equivalently that
$[\varphi^{(n_k)}(T)]_x < 2$ for all $z \in U$ and all $k$.  
In particular, $F_1^{(n_k)} \neq 0$ on $U$ for all $k$.

We now find that, for $x \in U$, we have
\begin{equation} 
\label{eq:decomp}
\begin{aligned}
\hat{h}^{(n_k)}_{\varphi,v}(x)  &= 
 \frac{1}{d^{n_k}} \max(\log_v([F_1^{(n_k)}(1,T)]_x),
                    \log_v([F_2^{(n_k)}(1,T)]_x)) \\
&= \frac{1}{d^{n_k}} \log_v([F_1^{(n_k)}(1,T)]_x) + 
 \frac{1}{d^{n_k}} \max(1, \log_v([\varphi^{(n_k)}(T)]_x)). \\
\end{aligned}
\end{equation}

The last term in (\ref{eq:decomp}) converges uniformly to zero, since
the quantity $[\varphi^{(n_k)}(T)]_x$ is uniformly bounded as $k$ varies.
Moreover, the term $\frac{1}{d^{n_k}} \log_v([F_1^{(n_k)}(1,T)]_x)$ is
harmonic on $U$ for all $k$.
Since $\hat{h}^{(n_k)}_{\varphi,v}$ converges uniformly to
$\hat{h}_{\varphi,v}$, it follows that 
$\frac{1}{d^{n_k}} \log_v([F_1^{(n_k)}(1,T)]_x)$ converges uniformly
to $\hat{h}_{\varphi,v}$ on $U$.  Therefore
$\hat{h}_{\varphi,v}$ is harmonic on $U$ as desired.

Finally, for the case $x_0 = \infty$, apply a similar argument 
to the function $\hat{g}_{\varphi,v}(x)$, using the coordinate function
$U = 1/T$ instead of $T$.  
\end{proof}


\vskip .1 in
\subsection{The Energy Minimization Principle.}

\vskip .05 in
Recall that for each $\zeta \in \PP^1_{\Berk}$, the Hsia kernel is
\begin{equation} \label{FKK1}
\delta(x,y)_{\zeta} \ = \ \frac{\|x,y\|_v}{\|x,\zeta\|_v \, \|y,\zeta\|_v}
\ = \ \frac{\delta(x,y)_{\zeta_0}}
        {\delta(x,\zeta)_{\zeta_0} \, \delta(y,\zeta)_{\zeta_0}}\ .
\end{equation}
Given a log-continuous probability measure $\mu$, define the potential kernel
\begin{equation} \label{FNC1}
g_{\mu}(x,y) \ = 
\ \int_{\PP^1_{\Berk}} -\log_v(\delta(x,y)_{\zeta}) \, d\mu(\zeta)  + C
\end{equation}
where the normalizing constant $C$ is chosen so that 
\begin{equation} \label{FNC2}
\iint g_{\mu}(x,y) \, d\mu(x) d\mu(y) \ = \ 0 \ .
\end{equation}
Using (\ref{FKK1}) one sees that
\begin{equation*}
g_{\mu}(x,y) \ = \ j_{\zeta_0}(x,y) - u_{\mu}(x,\zeta_0) 
                      - u_{\mu}(y,\zeta_0) + C 
\end{equation*}
and that $C = \iint j_{\zeta_0}(x,y) \, d\mu(x) d\mu(y)$.
Since $u_{\mu}(z,\zeta_0)$ is continuous, 
$g_{\mu}(x,y)$ inherits the following properties from $j_{\zeta_0}(x,y)$.  

\begin{proposition} \label{PropG16}
Let $\mu$ be a log-continuous probability measure on $\PP^1_{\Berk}$.
Then $g_{\mu}(x,y)$ is symmetric, lower semicontinuous everywhere,
continuous off the diagonal, and bounded from below.
For each fixed $y$, the function $G_y(x) = g_{\mu}(x,y)$
is continuous and belongs to $\BDV(\PP^1_{\Berk})$;
it is subharmonic in $\PP^1_{\Berk} \backslash \{y\}$
and satisfies $\Delta_{\PP^1_{\Berk}}(G_y(x)) = \delta_y(x) - \mu$.  
\end{proposition}

Given a probability measure $\nu$ on $\PP^1_{\Berk}$,
define the generalized potential function
\begin{equation} \label{FLG1}
u_{\nu}(x,\mu) \ = \ \int g_{\mu}(x,y) \, d\nu(y) \ ,
\end{equation}
and the $\mu$-energy integral 
\begin{equation} \label{FLG2}
I_{\nu}(\mu) \ = \ \iint g_{\mu}(x,y) \, d\nu(x) d\nu(y) 
\ = \ \int u_{\nu}(x,\mu) \, d\nu(x) \ .
\end{equation}

Our goal is to prove the following
energy minimization principle: 

\begin{theorem} \label{ThmG147}
Let $\mu$ be a log-continuous probability measure on $\PP^1_{\Berk}$.  Then 

$A)$  $I_{\mu}(\nu) \ge 0$ for each probability measure 
          $\nu$ on $\PP^1_{\Berk}$, and 
          
$B)$  $I_{\mu}(\nu) = 0$ if and only if $\nu = \mu$.  
\end{theorem}

The proof rests on establishing analogues of Maria's theorem and 
Frostman's theorem for the functions $u_{\nu}(x,\mu)$.
We begin by showing that generalized 
potential functions have the same properties as usual ones.  

Recall that a real-valued function $f(z)$ is lower semi-continuous if
\begin{equation*}
\liminf_{z \rightarrow z} f(x) \ \ge \ f(x)
\end{equation*}
for each $x$.  This is equivalent to requiring that $f^{-1}((b,\infty))$
be open, for each $b \in \RR$.  Recall also that $f(z)$ is strongly
lower semi-continuous if for each $x$
\begin{equation*}
\liminf_{z \rightarrow x} f(z) \ = \ f(x) \ .
\end{equation*}

\begin{proposition} \label{PropG18}
Let $\mu$ be a log-continuous probability measure on $\PP^1_{\Berk}$, 
and let $\nu$ be an arbitrary probability measure on $\PP^1_{\Berk}$.  
Then $u_{\nu}(x,\mu)$ is strongly lower semi-continuous everywhere,
and is continuous at each $z \notin \supp(\nu)$.
For each $p \in \PP^1(\CC_v)$, as $x$ approaches $p$ along any path $[y,p]$,
\begin{equation*}
\lim \begin{Sb} x \rightarrow p \\ x \in [y,p) \end{Sb} f(x) \ = \ f(p) \ .
\end{equation*}
Moreover, $u_{\nu}(z,\mu)$ belongs to $\BDV(\PP^1_{\Berk})$, and
$\Delta_{\PP^1_{\Berk}}(u_{\nu}(x,\mu)) = \nu - \mu$.
\end{proposition}

\begin{proof}  
Note that 
\begin{eqnarray*}
u_{\nu}(x,\mu) & = &
\iint j_{\zeta_0}(x,y) - j_{\zeta_0}(x,\zeta) - j_{\zeta_0}(y,\zeta)
                          \, d\mu(\zeta) d\nu(y) \\
& = & u_{\nu}(x,\zeta_0) - u_{\mu}(x,\zeta_0) - D \ .
\end{eqnarray*}
Here $D = \int u_{\mu}(y,\zeta_0) \, d\nu(y)$ is a finite constant.  
Since $u_{\mu}(x,\zeta_0)$ is bounded and continuous everywhere,  
the assertions about continuity, strong lower-semicontinuity,
and path limits follow from the corresponding facts for $u_{\nu}(x,\zeta_0)$
proved in Proposition~\ref{PropD7}.

We have also seen in Example~\ref{Example E.4} that $u_{\nu}(x,\zeta_0)$
and $u_{\mu}(x,\zeta_0)$ belong to $\BDV(\PP^1_{\Berk})$.  
By the computations of the Laplacians there, $\Delta(u_{\nu}(x,\mu)) 
= (\nu - \delta_{\zeta_0}(x)) - (\mu - \delta_{\zeta_0}(x)) = \nu - \mu$.
\end{proof}    

\vskip .1 in
The following lemma plays a key role in the proof of 
Maria's theorem for the generalized potential function $u_{\nu}(z,\mu)$. 
We have used its analogue for standard potential functions many times.  

\begin{lemma} \label{LemG19} 
Let $\Gamma$ be a subgraph of $\PP^1_{\Berk}$,
and let $\nu$ be an arbitrary probability measure on $\PP^1_{\Berk}$.  
Then for each $\zeta \in \PP^1_{\Berk} \backslash \PP^1(\CC_v)$, 
the restriction of $u_{\nu}(x,\zeta)$ to $\Gamma$ is finite and continuous.
\end{lemma} 

\begin{proof}  First note that for any  
$\xi \in \PP^1_{\Berk} \backslash \PP^1(\CC_v)$, 
there is a constant $C$ such that 
\begin{equation} \label{FAA1}
u_{\nu}(x,\zeta) 
\ = \ u_{\nu}(x,\xi) + j_{\zeta_0}(x,\xi) - j_{\zeta_0}(x,\zeta) + C \ .
\end{equation}
To see this, compute  
\begin{eqnarray*}
u_{\nu}(x,\zeta) - u_{\nu}(x,\xi) & = &         
\int j_{\zeta_0}(x,y) - j_{\zeta_0}(x,\zeta) 
                      - j_{\zeta_0}(y,\zeta) \, d\nu(y) \\
&  &         
\quad -  \int j_{\zeta_0}(x,y) - j_{\zeta_0}(x,\xi) 
                      - j_{\zeta_0}(y,\xi) \, d\nu(y) \\   
& = & j_{\zeta_0}(x,\xi) - j_{\zeta_0}(x,\zeta) \\    
&  &         
\quad + \int j_{\zeta_0}(y,\xi) \, d\nu(y) 
      - \int j_{\zeta_0}(y,\zeta) \, d\nu(y) \ .
\end{eqnarray*} 
Since $\zeta_0$, $\xi$ and $\zeta$ all belong to 
$\PP^1_{\Berk} \backslash \PP^1(\CC_v)$, the functions $j_{\zeta_0}(y,\xi)$
and $j_{\zeta_0}(y,\zeta)$ are bounded and continuous on $\PP^1_{\Berk}$,
and hence the integrals in the last line are finite.  

Now take $\xi = r_{\Gamma}(\zeta)$.  Note that there is a constant $C_1$
such that
\begin{equation*} 
j_{\xi}(x,y) + C_1 \ = \ 
j_{\zeta_0}(x,y) - j_{\zeta_0}(x,\xi) - j_{\zeta_0}(y,\xi)  \ .
\end{equation*} 
By (\ref{FAA1}), to show that $u_{\nu}(x,\zeta)$ is continuous on $\Gamma$ 
it suffices to show that $u_{\nu}(x,\xi)$ is continuous on $\Gamma$,
since $j_{\zeta_0}(x,\xi)$ and $j_{\zeta_0}(x,\zeta)$ are continuous on 
all of $\PP^1_{\Berk}$.  However, for  $x \in \Gamma$, 
\begin{eqnarray*}
u_{\nu}(x,\xi) & = & \int j_{\xi}(x,y) \, d\nu(y) + C_1  \\                              
& = & \int_{\Gamma} j_{\xi}(x,r_{\Gamma}(y)) \, d\nu(y) + C_1  
\ = \ \int_{\Gamma} j_{\xi}(x,w) \, d\overline{\nu}(w) + C_1 
\end{eqnarray*}
where $\overline{\nu} = (r_{\Gamma})_*(\nu)$.  Since $j_{\xi}(x,w)$ is
bounded and continuous for $(x,w)$ $\Gamma \times \Gamma$, it follows that
that $u_{\nu}(x,\xi)$ is finite and continuous on $\Gamma$.  
\end{proof}                 
                                 
\begin{proposition} \label{PropG150} \text{\rm (Maria)}
Let $\mu$ be a log-continuous probability measure on $\PP^1_{\Berk}$,
and let $\nu$ be a probability measure on $\PP^1_{\Berk}$.
If there is a constant $M < \infty$ such that $u_{\nu}(z,\mu) \le M$
on $\supp(\nu)$, then $u_{\nu}(z,\mu) \le M$
for all $z \in \PP^1_{\Berk}$.
\end{proposition}

\begin{proof}
Put $E = \supp(\nu)$, and let $U$ be a component of
$\PP^1_{\Berk} \backslash E$.  By Proposition \ref{PropG18}
$u_{\nu}(z,\mu)$ is continuous on $U$ and $\Delta_U(u_{\nu}(z,\mu)) \le 0$.
Hence $u_{\nu}(z,\mu)$ is  subharmonic on $U$.
We will use the maximum principle for subharmonic functions 
and the representation 
$u_{\nu}(z,\mu) = u_{\nu}(z,\zeta_0) - u_{\mu}(z,\zeta_0) - D$ 
to show that $u_{\nu}(z,\mu) \le M$ on $U$.

\vskip .05 in
First we will show that for each type I point $x \in \partial U$, 
\begin{equation} \label{FAB1} 
\limsup \begin{Sb} p \rightarrow x \\ p \in U \end{Sb} \ u_{\nu}(p,\mu) 
        \ \le \ u_{\nu}(x,\mu)  \ .
\end{equation}        
Fix $x \in \partial U \cap \PP^1(\CC_v)$, and let $p_1, p_2, \ldots$ 
be a sequence of points in $U$ which approach $x$.  Fix $\varepsilon > 0$.
Since $u_{\mu}(z,\zeta_0)$ is continuous, there is a neighborhood $V$ of $x$
on which $|u_{\mu}(z,\zeta_0) - u_{\mu}(x,\zeta_0)| < \varepsilon$.

Since $\|z,w\|_v$ is continuous off the diagonal,
for each $p_n$ there is a point $\overline{p}_n \in E$ such that
$\|p_n,\overline{p}_n\|_v \le \|p_n,w\|_v$ for all $w \in E$.
By the ultrametric inequality, for each $w \in E$
\begin{equation} \label{FAC1}
\|\overline{p}_n,w\|_v \ \le \ \max(\|p_n,\overline{p}_n\|_v,\|p_n,w\|_v)
\ = \ \|p_n,w\|_v \ .
\end{equation}
Since  $p_n \rightarrow x$ and $x$ is of type I,  
$\lim_{n \rightarrow \infty} \|p_n,x\|_v \rightarrow 0$.  
By (\ref{FAC1}), $\|\overline{p}_n,x\|_v \le \|p_n,x\|_v$.
It follows that $\overline{p}_n \rightarrow x$ as well.   

Let $n$ be large enough that $p_n, \overline{p}_n \in V$.  
Then $|u_{\mu}(p_n,\zeta_0) - u_{\mu}(x,\zeta_0)| \le \varepsilon$ and 
$|u_{\mu}(\overline{p}_n,\zeta_0) - u_{\mu}(x,\zeta_0)| \le \varepsilon$, so 
$|u_{\mu}(p_n,\zeta_0) - u_{\mu}(\overline{p}_n,\zeta_0)| \le 2 \varepsilon$.
Using this, (\ref{FAC1}), and the definition of a potential function, 
\begin{eqnarray*}
u_{\nu}(p_n,\zeta_0) 
& = & u_{\nu}(p_n,\zeta_0) - u_{\mu}(p_n,\zeta_0) - D \\
& = & \int -\log_v(\|p_n,w\|_v) \, d\nu(w) - u_{\mu}(p_n,\zeta_0) - D \\
& \le & \int -\log_v(\|\overline{p}_p,w\|_v) \, d\nu(w)
     - u_{\mu}(\overline{p}_n,\zeta_0) - D + 2 \varepsilon \\
& = & u_{\nu}(\overline{p}_n,\mu) + 2 \varepsilon
 \ \le \ M + 2 \varepsilon \ .
\end{eqnarray*}
Since $\varepsilon > 0$ is arbitrary, we obtain (\ref{FAB1}).

\vskip .05 in
Suppose there were a point $p_0 \in U$ where $u_{\nu}(p_0,\mu) > M$. 
Fix $\sigma > 0$ with $u_{\nu}(p_0,\mu) \ge M + \sigma$.
Write $f(z) = u_{\nu}(z,\mu)$.  We will
construct a sequence of points $p_1, p_2, \ldots$ approaching $\partial U$
and lying on a path, with $f(p_n) \ge M+\sigma$ for each $n$.

If the main dendrite of $U$ is empty, then $U$ is a disc
(it cannot be $\PP^1_{\Berk}$, since $\supp(\nu)$ is nonempty).
Let $x$ be its unique boundary point, and put $D = [x,p_0]$.  
Take a sequence of points $p_1, p_2, \ldots \in D$ with $p_n \rightarrow x$.
By Corollary~\ref{CorF9}, $f(z)$ is non-increasing on paths away
from the boundary point $x$, so $M+\sigma \le f(p_1) \le f(p_2) \le \cdots$.

If the main dendrite $D$ of $U$ is nonempty, then since $f(z)$
is non-increasing on branches off the main dendrite,
we can assume that $p_0$ lies on the main dendrite.
Take a sequence of subgraphs $\Gamma^{(n)} \subset D$ which exhaust $D$.
Without loss we can assume that $p_0 \in \Gamma^{(1)}_0$ 
and that $\Gamma^{(n)} \subset \ \Gamma^{(n+1)}_0$ for each $n$.
After enlarging $\Gamma^{(n)}$ slightly if necessary, we can also assume that
$r_{\Gamma^{(n)}}(\partial U) = \partial \Gamma^{(n)}$ for each $n$.

Inductively construct $p_1, p_2, \ldots$ as follows.
Since $u_{\nu}(z,\mu)$ is subharmonic on $U$, Lemmas 
\ref{LemF4}, \ref{LemF6}, and \ref{LemF7} show
that its restriction to each $\Gamma^{(n)}$ belongs to $\BDV(\Gamma^{(n)})$
and achieves its maximum on $\Gamma^{(n)}$ at a point of
$r_{\Gamma^{(n)}}(\partial U)$.  Let $p_1 \in \partial \Gamma^{(1)}$
be a point where $f(p_1) \ge M+\sigma$.
Suppose we have found $p_n \in \partial \Gamma^{(n)}$
with $f(p_n) \ge M+\sigma$,
which is a local maximum for $f(z)$ on $\Gamma^{(n)}$,
and such that $p_0, p_1, \ldots, p_n$ lie on a path.

Viewing $p_n$ as a point of $\Gamma^{(n+1)}$,
let the direction vectors at $p_n$ be $\vec{v}_1, \ldots, \vec{v}_m$.
Since $p_n$ belongs to $\Gamma^{(n+1)}_0$, necessarily $m \ge 2$.
As $p_n \in \partial \Gamma^{(n)}$, one of them, say $\vec{v}_1$,
leads into $\Gamma^{(n)}$.  The others lead into edges
$e_2, \ldots, e_m$ of $\Gamma^{(n+1)}$.
For each $k \ge 2$, the edge $e_k$ corresponds to a component $C_k$ of
$\Gamma^{(n+1)} \backslash \Gamma^{(n)}$;  let $T^{(k)} = C_k \cup \{p_n\}$
be the closure of $C_k$.  Then $\vec{v}_k$ is the unique direction vector
at $p_n$ leading into $T^{(k)}$.  Note that each $T^{(k)}$
is a subtree of $D$ with $r_{T^{(k)}}(\partial U) = \partial T^{(k)}$.

Since $p_n \in \Gamma^{(n+1)}_0$, Lemma~\ref{LemF8} shows that
\begin{equation*}
0 \ \ge \ \Delta_{\Gamma^{(n+1)}}(f)(p_n)
  \ = \ - \sum_{k=1}^m d_{\vec{v}_k}f(p_n) \ .
\end{equation*}
Here $d_{\vec{v}_1}f(p_n) \le 0$,
since $p_n$ is a local maximum for $f(z)$ on $\Gamma^{(n)}$.
Hence there must be some $k \ge 2$ with $d_{\vec{v}_k}f(p_n) \ge 0$.
Fix such a $k$ and consider $f(z)$ on $T^{(k)}$.
By Lemma~\ref{LemF8}, $\Delta_{T^{(k)}}f(x) \le 0$ on $T^{(k)}$.
By our choice of $k$, $\Delta_{T^{(k)}}f(p_n) \le 0$.  Hence
$\Delta_{T^{(k)}}(f)^+ \subset \partial T^{(k)} \backslash \{p_n\}$.

Lemmas \ref{LemF5} and \ref{LemF8} 
show that $f(z)$ is either constant on $T^{(k)}$
or achieves its maximum on $T_k$ at a point of $\Delta_{T^{(k)}}(f)^+$.
In either case there is a point $p \in \partial T^{(k)} \backslash \{p_n\}$
where $f(z)$ achieves its maximum.  Put $p_{n+1} = p$.  Then
$p_{n+1} \in \partial \Gamma^{(n+1)}$, $f(p_{n+1}) \ge M+\sigma$,
and $p_{n+1}$ is a local maximum for $f(z)$ on $\Gamma^{(n+1)}$.
By construction the points $p_0, p_1, \ldots, p_{n+1}$ lie on a path.

Since the graphs $\Gamma^{(n)}$ exhaust $D$, the points $p_n$
approach a boundary point $x \in \partial U \subset E$.
By what we have shown above, $x$ cannot be of type I,
since $\liminf_{n \rightarrow \infty} f(p_n) \ge M+\sigma$.

Let $\Gamma$ be the path $[p_0,x]$.
It has finite length, since neither $p_0$ nor $x$ is of type I,
so it is a subgraph of $\PP^1_{\Berk}$.
By Lemma \ref{LemG19}, the restriction of 
$u_{\nu}(z,\mu) = u_{\nu}(z,\zeta_0) - u_{\mu}(z,\zeta_0) - D$
to $\Gamma$ is continuous.  It follows that 
\begin{equation*}
u_{\nu}(x,\mu) \ = \ \lim_{n \rightarrow \infty} u_{\nu}(p_n,\mu)
               \ \ge \ M + \sigma \ .
\end{equation*}
This is a contradiction, since $x \in \partial U \subset E = \supp(\nu)$,
so $u_{\nu}(x,\mu) \le M$ by hypothesis.  
Hence $u_{\nu}(z,\mu) \le M$ on $U$.  
\end{proof}                 
 
\vskip .1 in
Define the `$\mu$-Robbin's constant' to be
\begin{equation*}
V(\mu) \ = \ \inf \begin{Sb} \nu \\ \text{prob meas} \end{Sb} I_\mu(\nu)
       \ = \ \inf \begin{Sb} \nu \\ \text{prob meas} \end{Sb} 
                      \iint g_{\mu}(x,y) \, d\nu(x) d\nu(y), 
\end{equation*}
where $\nu$ runs over all probability measures supported on $\PP^1_{\Berk}$.
Trivially $V(\mu) > -\infty$, since $g_{\mu}(x,y)$ is bounded below,
and $V(\mu) \le 0$, since $I_{\mu}(\mu) = 0$ 
by the normalization of $g_{\mu}(x,y)$.  
Taking the weak limit of a sequence of measures $\nu_1, \nu_2, \ldots$
for which $I_{\mu}(\nu_n) \rightarrow V(\mu)$ and applying the same argument
used to prove the existence of an equilibrium measure,
one obtains a probability measure $\omega$ for which
$I_{\mu}(\omega) = V(\mu)$.

\begin{proposition} \label{PropG151} \text{\rm (Frostman)}
Let $\mu$ be a log-continuous probability measure on $\PP^1_{\Berk}$,
and let $\omega$ be a probability measure for which $I_{\mu}(\omega) = V(\mu)$.
Then the generalized potential function $u_{\omega}(z,\mu)$ satisfies
\begin{equation*}
u_{\omega}(z,\mu) \ \equiv \ V(\mu) \qquad \text{on \ \ $\PP^1_{\Berk}$.}
\end{equation*}
\end{proposition}

\begin{proof}
First, using a quadraticity argument,
we will show that $u_{\omega}(z,\mu) \ge V(\mu)$ for all $z \in \PP^1_{\Berk}$
except possibly on a set $f$ of capacity $0$.  Then, we will show that
$u_{\omega}(z,\mu) \le V(\mu)$ on $\supp(\omega)$.  By Maria's theorem,
$u_{\omega}(z,\mu) \le V(\mu)$ for all $z \in \PP^1_{\Berk}$.  Since
a set of capacity $0$ is necessarily contained in $\PP^1(\CC_v)$,
it follows that $u_{\omega}(z,\mu) = 0$ on
$\PP^1_{\Berk} \backslash \PP^1(\CC_v)$.
Finally, if $p \in \PP^1(\CC_v)$, let $z \rightarrow p$ along
a path $[y,p]$. Proposition \ref{PropG18} gives 
$u_{\omega}(p,\mu) = V(\mu)$.

Put
\begin{eqnarray*}
  f & = & \{ z \in \PP^1_{\Berk} : u_{\omega}(z,\mu) < V_(\mu) \} \ , \\
f_n & = & \{ z \in \PP^1_{\Berk} : u_{\omega}(z,\mu) < V_(\mu) - 1/n \}
       \ , \quad \text{for $n = 1, 2, 3, \ldots$.}
\end{eqnarray*}
Since $u_{\omega}(z,\mu)$ is lower semicontinuous, each $f_n$ is closed, hence
compact, so $f$ is an $F_{\sigma}$ set.  By Corollary~\ref{CorD10}, $f$ has
capacity $0$ if and only if each $f_n$ has capacity $0$.

Suppose $f_n$ has positive capacity for some $n$;  then there is
a probability measure $\sigma$ supported on $f_n$ such that
$I_{\zeta_0}(\sigma)
:= \iint j_{\zeta_0}(x,y) \, d\sigma(x) d\sigma(y) < \infty$.
Since $g_{\mu}(x,y)$ differs from $j_{\zeta_0}(x,y)$ by a bounded function,
$I_{\mu}(\sigma) < \infty$ as well.

Furthermore, since
\begin{equation*}
V(\mu)
\ = \ \iint g_{\mu}(x,y) \, d\omega(x) d\omega(y)
\ = \ \int u_{\omega}(x,\mu) \, d\omega(x) \ ,
\end{equation*}
there is a point $q \in \supp(\mu)$ with $u_{\omega}(q,\mu) \ge V(\mu)$.
Since $u_{\omega}(z,\mu)$ is lower semicontinuous, there is a neighborhood
$U$ of $q$ on which $u_{\omega}(z,\mu) > V(\mu) - 1/(2n)$.
After shrinking $U$ if necessary, we can assume its closure
$\Ubar$ is disjoint from $f_n$.
Since $q \in \supp(\omega)$,
it follows that $M := \omega(\Ubar) > 0$.
Define a measure $\sigma_1$ of total mass $0$ by
\begin{equation*}
\sigma_1 \ = \ \left\{ \begin{array}{ll}
                  M \cdot \sigma & \text{on $f_n$,} \\
                  -\mu & \text{on $\Ubar$,} \\
                  0 & \text{elsewhere.}
                       \end{array} \right.
\end{equation*}
We claim that $I_{\zeta}(\sigma_1)$ is finite.  Indeed
\begin{eqnarray*}
I_{\zeta}(\sigma_1) & = &
   M^2 \cdot \iint_{f_n \times f_n}
                   g_{\mu}(x,y) \, d\sigma(x) d\sigma(y)  \\
   & & \qquad
   - 2M \cdot \iint_{f_n \times \Ubar} g_{\mu}(x,y) \, d\sigma(x) d\omega(y) \\
   & & \qquad \qquad \qquad
      + \iint_{\Ubar \times \Ubar} g_{\mu}(x,y) \, d\omega(x) d\omega(y) \ .
\end{eqnarray*}
The first integral is finite by hypothesis.
The second is finite because $f_n$ and $\Ubar$ are disjoint,
so $g_{\mu}(z,w)$ is bounded on $f_n \times \Ubar$.
The third is finite because $I_{\mu}(\omega)$ is finite
and $g_{\mu}(z,w)$ is bounded below.

For each $0 \le t \le 1$, $\omega_t := \omega + t \sigma_1$
is a probability measure.  By an expansion like the one above,
\begin{eqnarray*}
I_{\mu}(\omega_t) - I_{\mu}(\omega) & = &
   2t \cdot \int_E u_{\omega}(z,\mu) \, d\sigma_1(z)
           + t^2 \cdot I_{\mu}(\sigma_1) \\
   & \le & 2t \cdot( (V(\mu) -1/n) - (V(\mu) - 1/(2n)) \cdot M
            + t^2 \cdot I_{\mu}(\sigma_1) \\
   & = & (-M/n) \cdot t + I_{\mu}(\sigma_1) \cdot t^2 \ .
\end{eqnarray*}
For sufficiently small $t > 0$, the right side is negative.
This contradicts the fact that $\omega$ minimizes the energy integral.
It follows that $f_n$ has capacity $0$, and hence $f = \cup_{n=1}^{\infty} f_n$
has capacity $0$ by Corollary~\ref{CorD10}.

\vskip .05 in
The second part requires showing that $u_{\omega}(z,\mu) \le V(\mu)$
for all $z \in \supp(\omega)$.  If $u_{\omega}(z,\mu) > V(\mu)$ for
some $q \in \supp(\omega)$, let $\varepsilon > 0$ be small enough that
$u_{\omega}(q,\mu) > V(\mu) + \varepsilon$.
The lower semicontinuity of $u_{\omega}(z,\mu)$ shows there is a
neighborhood $U$ of $q$ on which $u_{\omega}(z,\mu) > V(\mu) + \varepsilon$.
Then $T := \omega(U) > 0$, since $q \in \supp(\mu)$.
On the other hand, by Lemma~\ref{LemD9}, $\omega(f) = 0$
since $I_{\mu}(\omega) < \infty$ implies that $I_{\zeta_0}(\omega) < \infty$.
Since $u_{\omega}(z,\mu) \ge V(\mu)$ for all $z \notin f$,
\begin{eqnarray*}
V(\mu)
& = & \int_{U} u_{\omega}(z,\mu) \, d\omega(z)
       + \int_{\PP^1_{\Berk} \backslash U} u_{\omega}(z,\mu) \, d\omega(z) \\
& \ge & T \cdot (V(\mu) + \varepsilon) + (1-T) \cdot V(\mu)
\ = \ V(\mu) + T \varepsilon
\end{eqnarray*}
which is impossible.
Hence $u_{\omega}(z,\mu) \le V(\mu)$ on $\supp(\omega)$, and Maria's
theorem implies that $u_{\omega}(z,\mu) \le V(\mu)$ for all $z$.

As noted at the beginning of the proof,
the fact that $u_{\omega}(p,\mu) = V(\mu)$ for all $p \in \PP^1(\CC_v)$,
now follows from Proposition \ref{PropG18}.
\end{proof}

\vskip .1 in
We can now complete the proof of Theorem \ref{ThmG147}.
\vskip .1 in

\begin{proof} (of Theorem \ref{ThmG147}).

Suppose $\omega$ minimizes the energy integral $I_{\mu}(\nu)$.
By Proposition \ref{PropG151}, $u_{\omega}(z,\mu)$ is constant.
It follows that $\Delta(u_{\omega}(z,\mu)) \equiv 0$.  On the other hand, 
by Proposition \ref{PropG18}, $\Delta(u_{\omega}(z,\mu)) = \omega - \mu$.
Hence $\omega = \mu$.

However, $I_{\mu}(\mu) = 0$ by the normalization of $g_{\mu}(x,y)$, 
so $V(\mu) = 0$.  It follows that $I_{\mu}(\nu) \ge 0$ for all $\nu$,
and if $I_{\mu}(\nu) = 0$, then $\nu = \mu$.
\end{proof}

\vskip .1 in

Theorem~\ref{ThmG147} is one of the key results needed in \cite{B-R2} 
for proving a nonarchimedean equidistribution theorem for 
points of small dynamical height.

\end{document}